\documentclass[12pt,oneside]{book}
\usepackage{geometry}                		
\usepackage{graphicx}									
\usepackage{amssymb}\vspace{0.25cm}
\usepackage{mathtools}
\usepackage{setspace}
\usepackage{tikz-cd}
\usepackage[mathscr]{euscript}
\newcommand{\abs}[1]{\left\vert#1\right\vert}
\usepackage[normalem]{ulem}
\usepackage{manfnt}
\usepackage{pgfplots}
\usepackage{pifont}
\usepackage[safe]{tipa}
\usepackage{marvosym}
\usepackage{scalerel,stackengine}
\usepackage{titlesec}
\usepackage{makeidx}
\usepackage{centernot}
\usepackage{xspace}
\usepackage[noauto]{chappg}
\usepackage[OT2,T1]{fontenc}
\usepackage{hyperref}
\usepackage{amsmath}
\usepackage[amsmath,thmmarks]{ntheorem}
\usepackage{tabularx}
\usepackage{yfonts}
\usepackage{epstopdf}
\usepackage{pdfpages}
\newcolumntype{Y}{>{\raggedright\arraybackslash}X}

\newcommand\gphi{%
\text{$\Phi$}
}

\newcommand\E{%
\mathbb{E}
}

\newcommand\N{%
\mathbb{N}
}

\newcommand\Q{%
\mathbb{Q}
}
\newcommand\R{%
\mathbb{R}
}

\newcommand\X{%
\mathbb{X}
}
\newcommand\XX{%
\mathbb{X}
}
\newcommand\YY{%
\mathbb{Y}
}

\newcommand\bS{%
\textbf{S}\xspace
}

\newcommand\Hom{%
\text{Hom}
}

\newcommand\sA{%
\mathcal{A}
}
\newcommand\sB{%
\mathcal{B}
}
\newcommand\sC{%
\mathcal{C}
}
\newcommand\sD{%
\mathcal{D}
}
\newcommand\sE{%
\mathcal{E}
}
\newcommand\sF{%
\mathcal{F}
}
\newcommand\sH{%
\mathcal{H}
}

\newcommand\sK{%
\mathcal{K}
}
\newcommand\sO{%
\mathcal{O}
}
\newcommand\sM{%
\mathcal{M}
}

\newcommand\sP{%
\mathcal{P}
}

\newcommand\sV{%
\mathcal{V}
}

\newcommand\ap{%
\text{ap}
}

\newcommand\loc{%
\text{loc}
}
\newcommand\locx{%
\text{$\ell$oc}
}
\newcommand\elln{%
\text{$\ell$n}
}
\newcommand\elog{%
\text{$\ell$og\hsy}
}

\newcommand\Gr{%
\text{Gr}
}

\newcommand\nth{%
\text{th}
}

\newcommand\td{%
\text{d}
}
\newcommand\df{%
\text{d$f$}
}

\newcommand\tI{%
\text{I}
}

\newcommand\iI{
\text{I}\hspace{0.05cm}
}
\newcommand\tJ{%
\text{J}
}
\newcommand\tL{%
\text{L}
}

\newcommand\Lm{%
\text{L}
}
\newcommand\Lp{
\text{L}
}

\newcommand\subsetx{%
\ \subset \ 
}

\newcommand\tD{%
\text{D}
}
\newcommand\xD{%
\text{D}
}
\newcommand\dvg{%
\text{div}
}

\newcommand\ts{%
\text{s}
}

\newcommand\tu{%
\text{u}
}

\newcommand\tW{%
\text{W}
}

\newcommand\Tee{%
\mathsf{T}
}

\def\thereforex{\boldsymbol{\text{ }
\leavevmode
\lower0.4ex\hbox{\textbullet}
\kern-.9em\raise1.1ex\hbox{\textbullet}
\kern-0.9em\lower0.4ex\hbox{\textbullet}
\hspace{0.1cm}\thinspace\text{ }}}
\def\thereforez{\boldsymbol{\text{ }
\leavevmode
\lower0.4ex\hbox{$\circ$}
\kern-.9em\raise1.1ex\hbox{$\circ$}
\kern-0.9em\lower0.4ex\hbox{$\circ$}
\hspace{0.1cm}\thinspace\text{ }}}
\newcommand\ra{%
\rightarrow
}
\newcommand\lra{%
\longrightarrow
}

\newcommand\ds{%
\displaystyle
}

\newcommand\sptx{%
\text{spt}\hspace{0.05cm}
}

\newcommand\diam{%
\text{diam}\hspace{0.05cm}
}
\newcommand\dist{%
\text{dist}\hspace{0.05cm}
}
\newcommand\osc{%
\text{osc}\hspace{0.05cm}
}

\newcommand\Lip{%
\text{Lip}\hspace{0.05cm}
}

\newcommand\rank{%
\text{rank}
}

\newcommand\sign{%
\text{sign}\hspace{0.05cm}
}

\newcommand\un[1]{%
\underline{#1}\xspace
}

\newcommand\esssup{%
\text{ess\hspace{0.05cm}sup} \hspace{0.05cm}
}

\newcommand\restr[2]{%
{#1}|{#2}
}

\newcommand\hsx{%
\hspace{0.05cm}
}
\newcommand\hsy{%
\hspace{0.03cm}
}

\newcommand{\norm}[1]{\left\lVert #1 \right\rVert}

\newcommand{\normx}[1]{\big|\hspace{-.05cm}\big| #1 \big|\hspace{-.05cm}\big|}
\newcommand{\normxx}[1]{\bigg|\hspace{-.05cm}\bigg| #1 \bigg|\hspace{-.05cm}\bigg|}

\newcommand\BMO{%
\text{BMO}
}
\newcommand\BV{%
 \text{BV}
}

\newcommand\AC{%
\text{AC}
}
\newcommand\ACL{%
\text{ACL}
}
\newcommand\AT{%
\text{AT}
}

\newcommand\ov[1]{%
\overline{#1}
}

\newcommand\chisubE{\chi\raisebox{-.1cm}{$\scaleto{_{E}}{6.5pt}$}}
\newcommand\chisubOmega{\chi\raisebox{-.1cm}{$\scaleto{_{\Omega}}{6.5pt}$}}
\newcommand\chisubEcapOmega{\chi\raisebox{-.1cm}{$\scaleto{_{E \hsy \cap \hsy \Omega}}{6.5pt}$}}

\newcommand\chisubF{\chi\raisebox{-.1cm}{$\scaleto{_{F}}{6.5pt}$}}
\newcommand\chisubG{\chi\raisebox{-.1cm}{$\scaleto{_{G}}{6.5pt}$}}
\newcommand\chisubS{\chi\raisebox{-.1cm}{$\scaleto{_{S}}{6.5pt}$}}

\newcommand\chisubmQ{\chi\raisebox{-.1cm}{$\scaleto{_{\Q}}{6.5pt}$}}
\newcommand\chisubV{\chi\raisebox{-.1cm}{$\scaleto{_{V}}{6.5pt}$}}

\newcommand\chisubEsubt{\chi\raisebox{-.1cm}{$\scaleto{_{E_t}}{6.5pt}$}}

\newcommand\chisubfofEonem{\chi\raisebox{-.1cm}{$\scaleto{_{f(E_1^{(m)})}}{8.5pt}$}}
\newcommand\chisubfofEim{\chi\raisebox{-.1cm}{$\scaleto{_{f(E_i^{(m)})}}{8.5pt}$}}

\newcommand\chisubOmegasubtoff{\chi\raisebox{-.1cm}{$\scaleto{_{\Omega_t(f)}}{6.5pt}$}}
\newcommand\chisubOmegasubtofHr{\chi\raisebox{-.1cm}{$\scaleto{_{\Omega_t(H_r)}}{6.5pt}$}}
\newcommand\chisubOmegasubtoffsubk{\chi\raisebox{-.1cm}{$\scaleto{_{\Omega_t (f_k)}}{6.5pt}$}}

\newcommand\chisubU{\chi\raisebox{-.1cm}{$\scaleto{_U}{4.5pt}$}}
\newcommand\chisubK{\chi\raisebox{-.1cm}{$\scaleto{_K}{4.5pt}$}}

\newcommand\chisubalpha{\chi\raisebox{-.1cm}{$\scaleto{_\alpha}{4.0pt}$}}
\newcommand\chisubbeta{\chi\raisebox{-.1cm}{$\scaleto{_\beta}{7.5pt}$}}

\newcommand\chisubPhiInvofE{\chi\raisebox{-.1cm}{$\scaleto{_{\Phi^{-1}(E)}}{8.5pt}$}}
\newcommand\chisubTinvS{\chi\raisebox{-.1cm}{$\scaleto{_{T^{-1}(S)}}{8.5pt}$}}

\newcommand\chisubBofZeroOneInt{\chi\raisebox{-.1cm}{$\scaleto{_{B(0,1)^\circ}}{8.5pt}$}}

\newcommand\chisubRnbackslashEsubv{\chi\raisebox{-.1cm}{$\scaleto{_{\R^n \backslash E_{v}}}{6.5pt}$}}

\newcommand\xL{%
\ \mathsf{L} \ 
}

\newcommand\bdot{%
\text{\textbullet}
}

\definecolor{ultramarine}{RGB}{0, 32, 96}

\definecolor{darkcerulean}{rgb}{0.3, 0.27, 0.49}

\definecolor{forestgreen}{rgb}{0.0, 0.27, 0.13}
\definecolor{forestgreenweb}{rgb}{0.13, 0.55, 0.13}
\definecolor{deepjunglegreen}{rgb}{0.0, 0.29, 0.29}

\definecolor{midnightblue}{rgb}{0.1, 0.1, 0.44}
\definecolor{midnightgreen}{rgb}{0.0, 0.29, 0.33}
\definecolor{myrtle}{rgb}{0.13, 0.26, 0.12}
\definecolor{darkviolet}{rgb}{0.58, 0.0, 0.83}
\definecolor{darkgreen}{rgb}{0.0, 0.2, 0.13}
\definecolor{officegreen}{rgb}{0.0, 0.5, 0.0}

\definecolor{harvardcrimson}{rgb}{0.79, 0.0, 0.09}
\definecolor{hollywoodcerise}{rgb}{0.96, 0.0, 0.63}
\definecolor{debianred}{rgb}{0.84, 0.04, 0.33}
\definecolor{darkturquoise}{rgb}{0.0, 0.81, 0.82}

\definecolor{darktangernine}{rgb}{1.0, 0.66, 0.07}
\definecolor{aureolin}{rgb}{0.99, 0.93, 0.0}
\definecolor{canaryyellow}{rgb}{1.0, 0.94, 0.0}
\definecolor{amber}{rgb}{1.0, 0.75, 0.0}

\definecolor{urobilin}{rgb}{0.88, 0.68, 0.13}
\definecolor{uscgold}{rgb}{1.0, 0.8, 0.0}

\makeatletter
\newenvironment{sqcases}{%
  \matrix@check\sqcases\env@sqcases
}{%
  \endarray\right.%
}
\def\env@sqcases{%
  \let\@ifnextchar\new@ifnextchar
  \left\lbrack
  \def\arraystretch{1.2}%
  \array{@{}l@{\quad}l@{}}%
}
\makeatother

\usepackage{scalerel,stackengine}
\stackMath
\newcommand\reallywidehat[1]{%
\savestack{\tmpbox}{\stretchto{%
  \scaleto{%
    \scalerel*[\widthof{\ensuremath{#1}}]{\kern-.6pt\bigwedge\kern-.6pt}%
    {\rule[-\textheight/2]{1ex}{\textheight}}
  }{\textheight}%
}{0.5ex}}%
\stackon[1pt]{#1}{\tmpbox}%
}
\makeatletter
\DeclareRobustCommand\widecheck[1]{{\mathpalette\@widecheck{#1}}}
\def\@widecheck#1#2{%
    \setbox\z@\hbox{\m@th$#1#2$}%
    \setbox\tw@\hbox{\m@th$#1%
       \widehat{%
          \vrule\@width\z@\@height\ht\z@
          \vrule\@height\z@\@width\wd\z@}$}%
    \dp\tw@-\ht\z@
    \@tempdima\ht\z@ \advance\@tempdima2\ht\tw@ \divide\@tempdima\thr@@
    \setbox\tw@\hbox{%
       \raise\@tempdima\hbox{\scalebox{1}[-1]{\lower\@tempdima\box
\tw@}}}%
    {\ooalign{\box\tw@ \cr \box\z@}}}
\makeatother

\newtheoremstyle{xx}
  {4pt}
  {0pt}
  {\upshape}
  {\bfseries}
  {}
  { }
  {}
  
\makeatletter 
 \newtheoremstyle{myu}%
  {\upshape\item[ \indent\indent\bf\underline{\theorem@headerfont ##2:}]}%
\makeatother
\makeatletter 
 \newtheoremstyle{myn}%
  {\item[\hskip\labelsep \ \bf ##1 \theorem@headerfont ##2.]}%
\makeatother
\theoremstyle{myn}
\newtheorem{theoremn}{Theorem} 
\theoremstyle{myu}
{\upshape}
\newtheorem{x}[theoremn]{}

 \newtheoremstyle{mr}%
  {\upshape\item[ \indent{\theorem@headerfont ##2. \hspace{.2cm}}]}%
\makeatother
\theoremstyle{mr}
{\upshape}
\newtheorem{rf}[theoremn]{}

\title{\textbf{Analysis 101:\\
Functions of Several Variables}}
\author{Garth Warner\\
Department of Mathematics\\
University of Washington}
\date{}	

\titleformat{\chapter}[display]
{\normalfont\filcenter\huge\bfseries}{}{0pt}{\large}

\titleformat{\chapter}[display]
{\normalfont\filcenter\huge\bfseries}{}{0pt}{\large}

\usepackage[OT2,OT1]{fontenc} 
\newcommand\cyr
{
\renewcommand\rmdefault{wncyr} 
\renewcommand\sfdefault{wncyss} 
\renewcommand\encodingdefault{OT2} 
\normalfont
\selectfont
}

\DeclareTextFontCommand{\textcyr}{\cyr}

\makeindex 
\begin{document}

\maketitle                              

\titlespacing*{\chapter}{0pt}{-50pt}{40pt}
\setlength{\parskip}{0.1em}
\pagenumbering{bychapter}
\setcounter{chapter}{0}
\pagenumbering{bychapter}
\setcounter{chapter}{0}
\setcounter{section}{0}

\begingroup
\fontsize{11pt}{11pt}\selectfont

\[
\textbf{ABSTRACT}
\]
\\[-1cm]

Apart from an account of classical preliminaries, 
this volume contains a systematic introduction to Sobolev spaces and functions of bounded variation with selected applications.  
This is installment III of a four part discussion of certain aspects of Real Analysis: Functions of a Single Variable, Curves and Length, Functions of Several Variables, and Surfaces and Area.
\\[2cm]

\[
\textbf{ACKNOWLEDGEMENT}
\]
\\[-1cm]

Many thanks to David Clark for his rendering the original transcript into AMS-LaTeX.  
Both of us also thank Judith Clare for her meticulous proofreading.
\newpage

\[
\textbf{FUNCTIONS OF SEVERAL VARIABLES}
\]
\\[-.5cm]


\hspace{.65cm}
\text{SECTION 1: \quad MEASURE THEORY}
\\

\hspace{1.3cm} \ \S1.1. \quad FACTS%
\\[-.26cm]

\hspace{1.3cm} \ \S1.2. \quad BOREL MEASURES%
\\[-.26cm]

\hspace{1.3cm} \ \S1.3. \quad RADON MEASURES%
\\[-.26cm]

\hspace{1.3cm} \ \S1.4. \quad OUTER MEASURES%
\\[-.26cm]

\hspace{1.3cm} \ \S1.5. \quad LEBESGUE MEASURES%
\\[-.26cm]

\hspace{1.3cm} \ \S1.6. \quad HAUSDORFF MEASURES%
\\[.26cm]

\hspace{.65cm}
\text{SECTION 2: \quad DIFFERENTIATION THEORY}
\\

\hspace{1.3cm} \ \S2.1. \quad SCALAR FUNCTIONS%
\\[-.26cm]

\hspace{1.3cm} \ \S2.2. \quad VECTOR FUNCTIONS%
\\[-.26cm]

\hspace{1.3cm} \ \S2.3. \quad LIPSCHITZ FUNCTIONS%
\\[-.26cm]

\hspace{1.3cm} \ \S2.4. \quad RADEMACHER%
\\[-.26cm]

\hspace{1.3cm} \ \S2.5. \quad STEPANOFF%
\\[-.26cm]

\hspace{1.3cm} \ \S2.6. \quad LUSIN%
\\[.26cm]

\hspace{.65cm}
\text{SECTION 3: \quad DENSITY THEORY}
\\

\hspace{1.3cm} \ \S3.1. \quad LEBESGUE POINTS%
\\[-.26cm]

\hspace{1.3cm} \ \S3.2. \quad APPROXIMATE LIMITS%
\\[-.26cm]

\hspace{1.3cm} \ \S3.3. \quad APPROXIMATE DERIVATIVES%
\\[.26cm]

\hspace{.65cm}
\text{SECTION 4: \quad WEAK PARTIAL DERIVATIVES}
\\

\hspace{.65cm}
\text{SECTION 5: \quad MOLLIFIERS}
\\

\hspace{.65cm}
\text{SECTION 6: \quad $\text{W}^{1, \hsy \infty} (\R^n)$}
\\

\hspace{.65cm}
\text{SECTION 7: \quad SOBOLEV SPACES}
\\

\hspace{1.3cm} \ \S7.1. \quad FORMALITIES%
\\[-.26cm]

\hspace{1.3cm} \ \S7.2. \quad EMBEDDINGS: \ GNS%
\\[-.26cm]

\hspace{1.3cm} \ \S7.3. \quad EMBEDDINGS: \ BMO%
\\[-.26cm]

\hspace{1.3cm} \ \S7.4. \quad EMBEDDINGS: \ MOR%
\\[.26cm]

\hspace{.65cm}
\text{SECTION 8: \quad ACL}
\\

\hspace{.65cm}
\text{SECTION 9: \quad BV SPACES}
\\

\hspace{1.3cm} \ \S9.1. \quad PROPERTIES%
\\[-.26cm]

\hspace{1.3cm} \ \S9.2. \quad DECOMPOSITION THEORY%
\\[-.26cm]

\hspace{1.3cm} \ \S9.3. \quad DIFFERENTIATION%
\\[-.26cm]

\hspace{1.3cm} \ \S9.4. \quad BVL%
\\[.26cm]

\hspace{.65cm}
\text{SECTION 10: \quad ABSOLUTE CONTINUITY}
\\

\hspace{.65cm}
\text{SECTION 11: \quad MISCELLANEA}
\\

\hspace{1.3cm} \ \S11.1. \quad PROPERTY (N)%
\\[-.26cm]

\hspace{1.3cm} \ \S11.2. \quad THE MULTIPLICITY FUNCTION%
\\[-.26cm]

\hspace{1.3cm} \ \S11.3. \quad JACOBIANS%
\\[.26cm]

\hspace{.65cm}
\text{SECTION 12: \quad AREA FORMULAS}
\\

\hspace{1.3cm} \ \S12.1. \quad THE LINEAR CASE%
\\[-.26cm]

\hspace{1.3cm} \ \S12.2. \quad THE $C^1$-CASE%
\\[-.26cm]

\hspace{1.3cm} \ \S12.3. \quad PROOF%
\\[-.26cm]

\hspace{1.3cm} \ \S12.4. \quad THE DIFFERENTIABLE CASE%
\\[-.26cm]

\hspace{1.3cm} \ \S12.5. \quad THE LIPSCHITZ CASE%
\\[-.26cm]

\hspace{1.3cm} \ \S12.6. \quad THE SOBOLEV CASE%
\\[-.26cm]

\hspace{1.3cm} \ \S12.7. \quad THE APPROXIMATE CASE%
\\[-.26cm]

\[
\]



\endgroup 

\chapter{
\text{SECTION 1: \quad MEASURE THEORY}
\\[.5cm]
$\boldsymbol{\S}$\textbf{1.1}.\quad  FACTS}
\setlength\parindent{2em}
\renewcommand{\thepage}{1-\S1-\arabic{page}}

\qquad
Let $\XX$ be a nonempty set and let $\sE \subset \sP(\XX)$ be a collection of subsets of $\XX$.
\\

\qquad{\bf 1.1.1:}\quad  
{\small\bf DEFINITION} \ 
The pair $(\XX, \sE)$ is called a \un{measurable space} if $\sE$ is a $\sigma$-algebra.
\\

\qquad{\bf 1.1.2:}\quad  
{\small\bf EXAMPLE} \ 
If $(\XX, \tau)$ is a topological space, then $(\XX, \sB(\XX))$ is a measurable space, $\sB(\XX)$ the $\sigma$-algebra 
of Borel subsets of $\XX$, i.e., the $\sigma$-algebra generated by the open subsets of $\XX$.
\\

\qquad{\bf 1.1.3:}\quad  
{\small\bf DEFINITION} \ 
Let $(\XX, \sE)$ be a measurable space.
\\[-.5cm]

\qquad \textbullet \quad A function $\mu\hsx : \hsx \sE \ra [0, +\infty]$ is a \un{positive measure} provided
$\mu(\emptyset) = 0$ and $\mu$ is $\sigma$-additive on $\sE$.
\\

\qquad{\bf 1.1.4:}\quad  
{\small\bf LEMMA} \ 
Let $(\XX, \sE)$ be a measurable space and suppose that $\mu\hsx : \hsx \sE \ra [0, +\infty]$ is 
$\sigma$-subadditive and additive $-$then $\mu$ is $\sigma$-additive hence $\mu$ is a positive measure.
\\[-.5cm]

PROOF \ 
Let $E_1$, $E_2, \ldots$ be a sequence of pairwise disjoint elements of $\sE$ $-$then 
\allowdisplaybreaks
\begin{align*}
\mu \bigg(\bigcup\limits_{n = 1}^\infty \ E_n \bigg) \ 
&\leq \ 
\sum\limits_{n = 1}^\infty \ \mu(E_n)
\\[9pt]
&\leq \ 
\lim\limits_{N \ra \infty} \ \sum\limits_{n = 1}^N  \ \mu(E_n )
\\[9pt]
&= \ 
\lim\limits_{N \ra \infty} \ \mu \bigg(\bigcup\limits_{n = 1}^N  \ E_n \bigg)
\\[9pt]
&\leq \ 
\mu \bigg(\bigcup\limits_{n = 1}^\infty \ E_n \bigg).
\end{align*}
\\[-1cm]


\qquad{\bf 1.1.5:}\quad  
{\small\bf DEFINITION} \ 
Let $(\XX, \sE)$ be a measurable space. 
\\[-.5cm]

\qquad \textbullet \ 
A function $\mu \hsx : \hsx \sE \ra \R^m$ $(m \geq 1)$ is a \un{vector measure} provided
$\mu(\emptyset) = 0$ and $\mu$ is $\sigma$-additive. 
\\[-.5cm]

[Note: \ 
If $\mu$ is a vector measure and if $m = 1$, then $\mu$ is a \un{real measure}.  
Since $+\infty$ is admitted, a positve measure is not necessarily a real measure.]
\\

\qquad{\bf 1.1.6:}\quad  
{\small\bf REMARK} \ 
If $\mu$ is a vector measure and if $E_1$, $E_2, \ldots$ is a sequence of pairwise disjoiont elements of $\sE$, then
\[
\mu \bigg(\bigcup\limits_{n = 1}^\infty \ E_n \bigg)  
\ = \ 
\sum\limits_{n = 1}^\infty \ \mu(E_n).
\]
Here the series on the right is absolutely convergent since its sum does not depend on the order of its terms 
(this being the case of the union on the left).
\\

\qquad{\bf 1.1.7:}\quad  
{\small\bf DEFINITION} \ 
Suppose that $\mu \hsx : \hsx \sE \ra \R^m$ $(m \geq 1)$ is a vector measure $-$then its \un{total variation} $\norm{\mu}$ 
is the arrow $\sE \ra [0, +\infty]$ defined by the prescription
\[
\norm{\mu} (E) 
\ = \ 
\sup \ 
\left[ 
\sum\limits_{n = 1}^\infty \ \norm{\mu(E_n)} \hsx : \hsx \{E_n\} \ 
\text{pairwise disjoint}, \  E \hsx = \hsx \bigcup\limits_{n = 1}^\infty \ E_n 
\right].
\]
\\[-.75cm]

\qquad{\bf 1.1.8:}\quad  
{\small\bf THEOREM} \ 
$\norm{\mu}$ is a positive finite measure (hence $\norm{\mu} (\XX) < +\infty)$.
\\

\qquad{\bf 1.1.9:}\quad  
{\small\bf REMARK} \ 
Denote by $\sM (\XX; \R^m)$ the set of $\R^m$-valued vector measures 
$\mu : \X \ra \R^m$ $(m \geq 1)$ $-$then $\sM (\XX; \R^m)$ is a real vector space and the total variation 
is a norm on $\sM (\XX; \R^m)$ under which it is a Banach space.
\\

\qquad{\bf 1.1.10:}\quad  
{\small\bf NOTATION} \ 
Given a real measure $\mu$, let

\[
\begin{cases}
&\ds\mu^+ \ = \ \frac{\abs{\mu} + \mu}{2} \hspace{.7cm} \text{positive part}
\\[8pt]
&\ds\mu^- \ = \ \frac{\abs{\mu} - \mu}{2} \hspace{.7cm} \text{negative part}
\end{cases}
.
\]
\\[-1cm]

\qquad{\bf 1.1.11:}\quad  
{\small\bf \un{N.B.}} \ 
Therefore $\mu^+$ and $\mu^-$ are positive finite measures and 
\[
\mu 
\ = \ 
\mu^+ - \mu^-,
\]
the \un{Jordan decomposition} of $\mu$.
\\

\qquad{\bf 1.1.12:}\quad  
{\small\bf SCHOLIUM} \ 
If $\mu$ is an $\R^m$-valued vector measure, say
\[
\mu
\ = \ 
(\mu_1, \ldots, \mu_m),
\]
put
\[
\int\limits_\XX \ f \ \td \mu 
\ = \ 
\bigg(
\int\limits_\XX \ f \td \mu_1, \ldots, \int\limits_\XX \ f \td \mu_m
\bigg).
\]

\noindent
Then
\[
\normxx{\int\limits_\XX \ f \ \td \mu } 
\ \leq \ 
\int\limits_\XX \ \abs{f} \ \td \norm{\mu}. 
\]
\\[-.25cm]

\qquad{\bf 1.1.13:}\quad  
{\small\bf NOTATION} \ 
Let $\mu$ be a positive measure on $(\XX, \sE)$.  
Given an $f \in \Lp^1 (\XX, \mu)^m$, say
\[
f 
\ = \ 
(f_1, \ldots, f_m),
\]
and an $E \in \sE$, put
\[
\int\limits_E \ f \ \td \mu 
\ = \ 
\bigg(
\int\limits_E \ f_1 \ \td \mu, \ldots, \int\limits_E \ f_m  \ \td \mu
\bigg).
\]
\\[-.75cm]

\qquad{\bf 1.1.14:}\quad  
{\small\bf SUBLEMMA} \ 
The assignment
\[
E \lra 
\int\limits_E \ f \ \td \mu 
\]
is an $\R^m$-valued vector measure, call it $f  \mu$.
\\

\qquad{\bf 1.1.15:}\quad  
{\small\bf LEMMA} \ 
The total variation $\norm{f  \mu}$ of $f  \mu$ is the assignment
\[
E \lra 
\int\limits_E \ \norm{f} \ \td \mu.
\]
\\[-.75cm]

PROOF \ 
First
\[
\norm{f  \mu}
\ \leq \ 
\norm{f} \mu.
\]
This said, fix a countable dense set $\{\mu_k\} \subset \bS^{m - 1}$ ($\subset \R^m$) and let $E \in \sE$.  
Given $\varepsilon > 0$, put
\[
\sigma(x) 
\ = \ 
\min\{k \in \N \hsx : \hsx \langle f(x), \mu_k \rangle \geq (1 - \varepsilon) \norm{f(x)} \} 
\]
and write
\[
E_k 
\ = \ 
\sigma^{-1} (\{k\}) \cap E \in \sE.
\]
Then
\begin{align*}
(1 - \varepsilon) \hsy \norm{f} \hsy \mu(E) \ 
&=\ 
(1 - \varepsilon) \ \int\limits_E \ \norm{f} \ \td \mu
\\[15pt]
&=\ 
\sum\limits_{k = 1}^\infty \ (1 - \varepsilon) \ \int\limits_{E_k} \ \norm{f} \ \td \mu
\\[15pt]
&\leq\ 
\sum\limits_{k = 1}^\infty \ \langle f  \mu (E_k), \mu_k \rangle
\\[15pt]
&\leq\ 
\sum\limits_{k = 1}^\infty \ \norm{f  \mu (E_k)} 
\\[15pt]
&<\ 
\norm{f  \mu} \hsy (E).
\end{align*}
\\[-.75cm]

\qquad{\bf 1.1.16:}\quad  
{\small\bf DEFINITION} \ 
Let $(\XX, \sE)$ be a measurable space.
\\[-.5cm]

\qquad \textbullet \quad 
Let $\nu$ be a positive measure and let $\mu$ be a vector measure $-$then $\mu$ is 
\un{absolutely continuous} w.r.t $\nu$, denoted $\mu \ll \nu$, if for every $E \in \sE$, the implication 

\[
\nu(E) 
\ = \ 
0
\implies
\norm{\mu} \hsy (E) \ = \ 0
\]
obtains.
\\

\qquad{\bf 1.1.17:}\quad  
{\small\bf EXAMPLE} \ 
If $f \in \Lp^1 (\XX, \nu)^m$, then $f \hsy \nu \ll \nu$.
\\

\qquad{\bf 1.1.18:}\quad  
{\small\bf CRITERION} \ 
$\mu$ is absolutely continuous w.r.t $\nu$ iff for every $\varepsilon > 0$ there exists $\delta > 0$ 
such that for every $E \in \sE$,
\[
\nu (E) < \delta 
\implies 
\norm{\mu} \hsy (E) < \varepsilon.
\]
\\[-1.25cm]

\qquad{\bf 1.1.19:}\quad  
{\small\bf DEFINITION} \ 
Let $(\XX, \sE)$ be a measurable space.
\\[-.5cm]

\qquad \textbullet \quad 
Let $\mu_1$, $\mu_2$ be positive measures $-$then $\mu_1$, $\mu_2$ are 
\un{mutually singular}, 
denoted $\mu_1 \perp \mu_2$, if there exists $E \in \sE$ such that
\[
\mu_1(E) \ = \ 0
\quad \text{and} \quad 
\mu_2(\XX - E) \ = \ 0.
\]
\\[-1cm]

\qquad{\bf 1.1.20:}\quad  
{\small\bf \un{N.B.}} \ 
Vector measures $\mu_1$, $\mu_2$ are mutually singular provided this is the case of $\norm{\mu_1}$, $\norm{\mu_2}$.
\\[-.5cm]

[Note: \ 
If $\nu$, $\mu$ are as above, write $\nu \perp \mu$ when $\nu \perp \norm{\mu}$.]
\\

\qquad{\bf 1.1.21:}\quad  
{\small\bf RADON-NIKODYM} \ 
Let $(\XX, \sE)$ be a measurable space.
\\[-.5cm]

\qquad \textbullet \quad 
Let $\nu$ be a positive measure and let $\mu$ be a vector measure, 
say $\mu : \sE \ra \R^m$ $(m \geq 1)$.  
Assume: \ 
$\nu$ is $\sigma$-finite $-$then there is a unique pair $\mu^a$, $\mu^s$ of $\R^m$-valued vector measures such that 
\[
\begin{cases}
\ \mu^a
\ll
\nu
\\[4pt]
\ \mu^s
\perp
\nu
\end{cases}
\quad
\text{and} \quad 
\mu 
\ = \ 
\mu^a + \mu^s.
\]
\\[-.75cm]

\qquad{\bf 1.1.22:}\quad  
{\small\bf \un{N.B.}} \ 
In addition, there is a unique $f \in \Lp^1 (\XX, \nu)^m$ such that 
$\mu^a = f \hsy \nu$, the so-called \un{density} of $\mu$ w.r.t. $\nu$, denoted $\ds \frac{\td \mu}{\td \nu}$.
\\[-.5cm]

[Note: \ 
Uniqueness is taken to mean in the sense of equivalence classes of functions which agree $\nu$-a.e.]
\\

\qquad{\bf 1.1.23:}\quad  
{\small\bf LEMMA} \ 
Let $\mu$ be an $\R^m$-valued vector measure $-$then there is a unique $\bS^{m - 1}$-valued function 
$f \in \Lp^1(\XX, \norm{\mu})^m$ such that $\mu = f \hsy \norm{\mu}$.
\\[-.5cm]

PROOF \ 
Trivially, $\mu \ll \norm{\mu}$, so $\mu = f \hsy \norm{\mu}$  ($f \in \Lp^1(\XX, \norm{\mu})^m$).  
Therefore
\begin{align*}
\norm{\mu} \ 
&=\ 
\norm{f \hsy \norm{\mu}}
\\[11pt]
&=\ 
\norm{f} \hsx \norm{\mu}
\qquad (\text{cf. 1.1.15.})
\end{align*}
thus $\norm{f} = 1$ ($\norm{\mu}$-a.e.).
\\

\qquad{\bf 1.1.24:}\quad  
{\small\bf DEFINITION} \ 
Let $(\XX, \sE)$ be a measurable space.
\\[-.5cm]

\qquad \textbullet \quad 
A function $\mu:\sE \ra [-\infty, +\infty]$ is a \un{signed measure} provided $\mu(\emptyset) = 0$, 
$\mu$ takes at most one of the two values $\pm\infty$, i.e., either 
$\mu:\sE \ra ]-\infty, +\infty]$ 
or 
$\mu:\sE \ra [-\infty, +\infty[\hsy$, 
and $\mu$ is $\sigma$-additive on $\sE$.
\\

\qquad{\bf 1.1.25:}\quad  
{\small\bf \un{N.B.}} \ 
A positive measure is a signed measure.
\\

\qquad{\bf 1.1.26:}\quad  
{\small\bf \un{N.B.}} \ 
A real measure is a signed measure.
\\

\qquad{\bf 1.1.27:}\quad  
{\small\bf DEFINITION} \ 
Suppose that $\mu$ is a signed measure on $(\XX, \sE)$.  
\\[-.5cm]

\qquad \textbullet \quad 
A set $E \in \sE$ is a \un{positive set} for $\mu$ if $\mu(E_0) \geq 0$ for every $E_0 \in \sE$ such that $E_0 \subset E$.
\\[-.5cm]

\qquad \textbullet \quad 
A set $E \in \sE$ is a \un{negative set} for $\mu$ if $\mu(E_0) \leq 0$ for every $E_0 \in \sE$ such that $E_0 \subset E$.
\\

\qquad{\bf 1.1.28:}\quad  
{\small\bf DEFINITION} \ 
Suppose that $\mu$ is a signed measure on $(\XX, \sE)$ $-$then a set $E \in \sE$ is a \un{null set} for $\mu$ 
if $\mu(E_0) = 0$ for every $E_0 \in \sE$ such that $E_0 \subset E$.
\\[-.25cm]

\qquad{\bf 1.1.29:}\quad  
{\small\bf DEFINITION} \ 
Let $(\XX, \sE)$ be a measurable space.
\\[-.5cm]

\qquad \textbullet \quad 
Given a signed measure $\mu$, sets $E_+$, $E_-$ are said to constitute a \un{Hahn} \un{decomposition} for $\mu$
provided $E_+ \cap E_- = \emptyset$, $E_+ \cup E_- = \XX$, and 
\[
\begin{cases}
&E_+ \ \text{is a positive set for $\mu$}\\
&E_- \ \text{is a negative set for $\mu$}
\end{cases}
.
\]
\\[-.75cm]

\qquad{\bf 1.1.30:}\quad  
{\small\bf THEOREM} \ 
Hahn decompositions exist.  
Moreover, if $(E_+, E_-)$ and 
$(E_+^\prime, E_-^\prime,)$ are two such, then 
$E_+ \hsx \Delta \hsx E_+^\prime$ and 
$E_- \hsx \Delta \hsx E_-^\prime$
are null sets for $\mu$.
\\

\qquad{\bf 1.1.31:}\quad  
{\small\bf LEMMA} \ 
If $\mu$, $\nu$ are signed measures on $(\XX, \sE)$, at least one of which is finite, 
then the set function $\mu - \nu$ is  a signed measure on $(\XX, \sE)$.
\\

\qquad{\bf 1.1.32:}\quad  
{\small\bf LEMMA} \ 
Suppose that $\mu$ is a signed measure on $(\XX, \sE)$ and let $E_1$, $E_2 \in \sE$ 
with $E_1 \subset E_2$ $-$then $\mu(E_2) \in \R$ $\implies \mu(E_1) \in \R$.  
And 
\[
\begin{cases}
&\mu(E_1) \ = \ +\infty \implies \mu(E_2)  \ = \ +\infty
\\[4pt]
&\mu(E_1) \ = \ -\infty \implies \mu(E_2)  \ = \ -\infty
\end{cases}
.
\]
\\[-.75cm]

\qquad{\bf 1.1.33:}\quad  
{\small\bf LEMMA} \ 
Suppose that $\mu$ is a signed measure $-$then there exist unique positive measures 
$\mu^+$, $\mu^-$ such that $\mu = \mu^+ - \mu^-$ and $\mu^+ \perp  \mu^-$.
\\[-.5cm]

[Let $\XX = E_+ \cup E_-$ be a Hahn decomposition for $\mu$ and put
\[
\begin{cases}
&\mu^+(E) \ = \  \mu(E \cap E_+)
\\[4pt]
&\mu^-(E) \ = \  -\mu(E \cap E_-)
\end{cases}
\qquad (E \in \sE).]
\]
\\[-.75cm]

\qquad{\bf 1.1.34:}\quad  
{\small\bf REMARK} \ 
If $\mu$ omits the value $+\infty$ ($-\infty$), then $\mu^+$ ($\mu^-$) is a finite positive measure.  
So if the range of $\mu$ is contained in $\R$, then $\mu$ is bounded, i.e., $\mu$ is a real measure. 
\\[-.25cm]

\qquad{\bf 1.1.35:}\quad  
{\small\bf DEFINITION} \ 
Let $(\XX, \sE)$ be a measurable space.
\\[-.5cm]

\qquad \textbullet \quad 
A signed measure $\mu$ is \un{finite} if $\abs{\mu}$ is a finite positive measure.
\\

\qquad{\bf 1.1.36:}\quad  
{\small\bf LEMMA} \ 
$\mu$ is finite iff $\mu (\XX) \in \R$.
\\

\qquad{\bf 1.1.37:}\quad  
{\small\bf RESTRICTION} \ 
Let $(\XX, \sE)$ be a measurable space and suppose that
$\mu$ is a positive, real, signed, or vector measure in $(\XX, \sE)$.  
Given $E \in \sE$, put
\[
(\mu \xL E) \hsy (S)
\ = \ 
\mu(E \cap S) 
\qquad (S \in \sE).
\]
Then
\[
\mu \xL E 
\ = \ 
\chisubE \hsy \mu.
\]
In fact,
\[
(\chisubE \hsy \mu) \hsy (S) 
\ = \ 
\int\limits_S \ \chisubE \ \td \mu 
\ = \ 
\mu(E \cap S).
\]
\\[-1cm]

\qquad{\bf 1.1.38:}\quad  
{\small\bf EXAMPLE} \ 
Per Radon-Nikodym, consider $\nu$ and $\mu$ $-$then there exists a set $E \in \sE$ such that 
$\nu(E) = 0$ and $\mu^s = \mu \xL E$.

\chapter{
$\boldsymbol{\S}$\textbf{1.2}.\quad  BOREL MEASURES}
\setlength\parindent{2em}
\setcounter{paragraph}{1}
\renewcommand{\thepage}{1-\S2-\arabic{page}}

\qquad
Let $\XX$ be a locally compact Hausdorff space (LCH space).
\\

{\bf 1.2.1:}\quad  
{\small\bf NOTATION} \ 
\\[-.5cm]

\qquad \textbullet \quad 
$\sO(\XX)$ is the collection of open subsets of $\XX$.
\\[-.5cm]

\qquad \textbullet \quad 
$\sK(\XX)$ is the collection of compact subsets of $\XX$.
\\[-.5cm]

\qquad \textbullet \quad 
$\sB(\XX)$ is the collection of Borel subsets of $\XX$.
\\

{\bf 1.2.2:}\quad  
{\small\bf DEFINITION} \ 
A positive measure on $(\XX, \sB(\XX))$ is referred to simply as a \un{Borel measure on $\XX$}.
\\

{\bf 1.2.3:}\quad  
{\small\bf DEFINITION} \ 
Let $\mu$ be a Borel measure on $\XX$ and let $E \in \sB(\XX)$.
\\[-.5cm]

\qquad \textbullet \quad 
$\mu$ is \un{outer regular} on $E$ if
\[
\mu(E) 
\ = \ 
\inf \{\mu(U) \hsx : \hsx U \supset E, \ U \in \sO(\XX)\}.
\]
\\[-1.5cm]

\qquad \textbullet \quad 
$\mu$ is \un{inner regular} on $E$ if
\[
\mu(E) 
\ = \ 
\sup \{\mu(K) \hsx : \hsx K \subset E, \ K \in \sK(\XX)\}.
\]
\\[-.75cm]

{\bf 1.2.4:}\quad  
{\small\bf DEFINITION} \ 
Let $\mu$ be a Borel measure on $\XX$ and let $E \in \sB(\XX)$ $-$then $\mu$ is \un{regular} on $E$ if $\mu$ is both 
inner regular and outer regular on $E$.
\\

{\bf 1.2.5:}\quad  
{\small\bf DEFINITION} \ 
Let $\mu$ be a Borel measure on $\XX$ and let $\sC$ be a subset of $\sB(\XX)$ $-$then $\mu$ is 
outer regular for $\sC$, 
inner regular for $\sC$,
or regular for $\sC$
according to whether $\mu$ is outer regular, inner regular, or regular for every $E \in \sC$.
\\

{\bf 1.2.6:}\quad  
{\small\bf TERMINOLOGY} \ 
A Borel measure $\mu$ on $\XX$ is 
outer regular, inner regular, or regular if $\mu$ is 
outer regular, inner regular, or regular for $\sB(\XX)$.
\\

{\bf 1.2.7:}\quad  
{\small\bf LEMMA} \ 
If every open subset of $\XX$ is $\sigma$-compact, then every Borel measure on $\XX$ is inner regular for $\sO(\XX)$ 
($\subset \sB(\XX)$).
\\

{\bf 1.2.8:}\quad  
{\small\bf EXAMPLE} \ %
Every open subset of $\R^n$ is $\sigma$-compact.
\\

{\bf 1.2.9:}\quad  
{\small\bf REMARK} \ 
If $\XX$ is  $\sigma$-compact, then the $\sigma$-algebra generated by $\sK(\XX)$ is $\sB(\XX)$.
\\

{\bf 1.2.10:}\quad  
{\small\bf LEMMA} \ 
Let $\mu_1$ and $\mu_2$ be two Borel measures on $\XX$.
\\[-.5cm]

\qquad \textbullet \quad 
If $\mu_1$ and $\mu_2$ are outer regular for $\sB(\XX)$ and $\mu_1 = \mu_2$ on $\sO(\XX)$, then $\mu_1 = \mu_2$ on $\sB(\XX)$.
\\[-.5cm]

\qquad \textbullet \quad 
If $\mu_1$ and $\mu_2$ are inner regular for $\sB(\XX)$ and $\mu_1 = \mu_2$ on $\sK(\XX)$, then $\mu_1 = \mu_2$ on $\sB(\XX)$.
\\[-.25cm]

[To establish the first point, let $E \in \sB(\XX)$ and write
\\[-1cm]

\begin{align*}
\mu_1(E) \ 
&=\ 
\inf \{\mu_1(U) \hsx : \hsx U \supset E, \ U \in \sO(\XX)\}
\\[8pt]
&=\ 
\inf \{\mu_2(U) \hsx : \hsx U \supset E, \ U \in \sO(\XX)\}
\\[8pt]
&=\ 
\mu_2(E).]
\end{align*}
\\[-1cm]

\[
* \hspace{0.25cm}
* \hspace{0.25cm}
* \hspace{0.25cm}
* \hspace{0.25cm}
* \hspace{0.25cm}
* \hspace{0.25cm}
* \hspace{0.25cm}
* \hspace{0.25cm}
* \hspace{0.25cm}
* 
\]

\[
\text{$\mathcal{APPENDIX}$}
\]
\\

A locally compact Hausdorff space $\XX$ is $\sigma$-compact if $\XX$ can be expressed as the union of at most 
countably many compact subspaces.
\\[-.5cm]

[Note: \ 
$\Q = \bigcup\limits_{q \in \Q} \ \{q\}$ and $\forall \ q$, $\{q\}$ is compact but $\Q$ is not locally compact.]
\\[-.25cm]

{\small\bf LEMMA} \ 
Every open subset of a second countable LCH space $\XX$ is $\sigma$-compact.
\\[-.25cm]

E.g.: \ This is the case of $\R^n$.
\\[-.25cm]

{\small\bf RAPPEL} \ 
A second countable topological space is separable (but, in general, not conversely).
\\

{\small\bf RAPPEL} \ 
Every separable metric space is second countable.
\\

So, if $\XX$ is a metrizable separable LCH space, then every open subset of $\XX$ is $\sigma$-compact.
\\[-.25cm]

{\small\bf \un{N.B.}} \ 
Every compact metric space is separable.

\chapter{
$\boldsymbol{\S}$\textbf{1.3}.\quad  RADON MEASURES}
\setlength\parindent{2em}
\setcounter{paragraph}{1}
\renewcommand{\thepage}{1-\S3-\arabic{page}}

\qquad
Let $\XX$ be a locally compact Hausdorff space (LCH space).
\\

\qquad{\bf 1.3.1:}\quad  
{\small\bf DEFINITION} \ 
A Borel measure $\mu$ on $\XX$ is said to be \un{locally finite} if $\forall \ K \in \sK(\XX)$, $\mu(K) < +\infty$.
\\

\qquad{\bf 1.3.2:}\quad  
{\small\bf DEFINITION} \ 
A locally finite Borel measure $\mu$ on $\XX$ is a \un{Radon} \un{measure} provided
\\[-.5cm]

\qquad \textbullet \quad
$\mu$ is outer regular for $\sB(\XX)$
\\[-.5cm]

\qquad \textbullet \quad
$\mu$ is inner regular for $\sO(\XX)$.
\\

\qquad{\bf 1.3.3:}\quad  
{\small\bf \un{N.B.}} \ 
A Radon measure is a positive measure.
\\

\qquad{\bf 1.3.4:}\quad  
{\small\bf REMARK} \ 
A finite Borel measure $\mu$ on a compact Hausdorff space $\XX$ is locally finite but it need not be Radon.
\\

\qquad{\bf 1.3.5:}\quad  
{\small\bf EXAMPLE} \ 
Take $\XX = \R^n$ $-$then the restriction of Lebesgue measure to $\sB(\XX)$ is a Radon measure.
\\[-.5cm]

[Note: \ 
Counting measure on $\R^n$ is not locally finite, hence is not Radon.]
\\

\qquad{\bf 1.3.6:}\quad  
{\small\bf LEMMA} \ 
Every $\sigma$-finite Radon measure is inner regular for $\sB(\XX)$, 
hence is regular for $\sB(\XX)$.
\\

In particular: \ 
Every finite Radon measure is regular for $\sB(\XX)$, thus every Radon measure on a compact Hausdorff space is regular for $\sB(\XX)$.
\\

\qquad{\bf 1.3.7:}\quad  
{\small\bf LEMMA} \ 
Suppose that $\XX$ is $\sigma$-compact $-$then every Radon measure on $\XX$ is inner regular for $\sB(\XX)$, 
hence is regular for $\sB(\XX)$.
\\[-.5cm]

[A Radon measure is locally finite, so here
\[
\text{$\sigma$-compact}
\implies
\text{$\sigma$-finite}.]
\]
\\[-1.25cm]

\qquad{\bf 1.3.8:}\quad  
{\small\bf LEMMA} \ 
Suppose that $\XX$ is $\sigma$-compact, let $\mu$ be a Radon measure on $\XX$, 
and let $\nu$ be a locally finite Borel measure on $\XX$.  
Assume: \ $\nu = \mu$ on $\sO(\XX)$ $-$then $\nu$ is regular for $\sB(\XX)$.
\\

\qquad{\bf 1.3.9:}\quad  
{\small\bf RIESZ REPRESENTATION THEOREM (RRT)} \ 
If $\iI$ is a positive linear functional on $C_c(\XX)$, then there exists a unique Radon measure $\mu$ on $\XX$ such that 
\[
\iI (f) 
\ = \ 
\int\limits_\XX \ f \ \td \mu
\]
for all $f \in C_c(\XX)$.
\\

\qquad{\bf 1.3.10:}\quad  
{\small\bf SUBLEMMA} \ 
Let $U$ be an open $\sigma$-compact subset of $\XX$ $-$then there is an increasing sequence
$f_1, f_2, \ldots$ on $C_c(\XX)$ 
such that $\lim\limits_{n \ra \infty} \ f_n = \chisubU$.
\\[-.5cm]

[Note: \ 
An open subset of a compact Hausdorff space need not be $\sigma$-compact.]
\\

\qquad{\bf 1.3.11:}\quad  
{\small\bf THEOREM} \ 
If every open subset of $\XX$ is $\sigma$-compact, then every locally finite Borel measure $\nu$ on $\XX$ 
is a regular Radon measure.
\\[-.5cm]

PROOF \ 
The issue is outer and inner regularity for $\sB(\XX)$ per $\nu$.  
Define a positive linear functional in $\iI$ in 
$C_c(\XX)$ by the prescription
\[
\iI \hsy (f) 
\ = \ 
\int\limits_\XX \ f \ \td \nu.
\]
Then by the RRT, there exists a unique Radon masure $\mu$ on $\XX$ such that $\forall \ f \in C_c(\XX)$, 
\[
\iI \hsy (f) 
\ = \ 
\int\limits_\XX \ f \ \td \mu.
\]
The claim now is that $\mu = \nu$ on $\sO(\XX)$.  
So let $U \in \sO(\XX)$ and choose as above $\{f_n\}$, 
$\lim\limits_{n \ra \infty} \ f_n = \chisubU$  
$-$then by monotone convergence,

\allowdisplaybreaks
\begin{align*}
\nu(U) \ 
&=\ 
\int\limits_\XX \ \chisubU \ \td \nu
\\[15pt]
&=\ 
\lim\limits_{n \ra \infty} \ \int\limits_\XX \ f_n  \ \td \nu
\\[15pt]
&=\ 
\lim\limits_{n \ra \infty} \ \int\limits_\XX \ f_n  \ \td \mu
\\[15pt]
&=\ 
\int\limits_\XX \ \chisubU \ \td \mu
\\[15pt]
&=\ 
\mu(U).
\end{align*}
Therefore $\nu = \mu$ on $\sO(\XX)$, thus $\nu$ is regular for $\sB(\XX)$.
\\[-.5cm]

[Note: \ 
Consequently, $\nu = \mu$ (cf. 1.2.10).]
\\

\qquad{\bf 1.3.12:}\quad  
{\small\bf RAPPEL} \ 
Let $(\XX, \sE)$ be a measurable space $-$then a \un{simple} \un{function} 
is a finite linear combination with real coefficients of characteristic functions of sets in $\sE$.
\\

\qquad{\bf 1.3.13:}\quad  
{\small\bf LEMMA} \ 
For any positive measure $\mu: \sE \ra [0, +\infty]$, the simple functions are dense in $\Lp^p(\XX, \mu)$ $(1 \leq p < +\infty)$.
\\

\qquad{\bf 1.3.14:}\quad  
{\small\bf THEOREM} \ 
If $\mu$ is a Radon measure on $\XX$, then $C_c(\XX)$ is dense in $\Lp^p(\XX, \mu)$ $(1 \leq p < +\infty)$.
\\[-.25cm]

PROOF \ 
It is enough to show that for any Borel set $E$ with $\mu(E) < +\infty$, $\chisubE$ can be approximated in the $\Lp^p$-norm 
by elements of $C_c(\XX)$.  
Given $\varepsilon > 0$, choose a compact $K \subset E$ and an open $U \supset E$ such that $\mu(U - K) < \varepsilon$ 
and using Urysohn, choose an $f \in C_c(\XX)$ such that 
$\chisubK \leq f \leq \chisubU$ $-$then
\[
\norm{\chisubE - f}_p
\ \leq \ 
\mu(U - K)^{1/p} 
\ < \ 
\varepsilon^{1/p}.
\]

\chapter{
$\boldsymbol{\S}$\textbf{1.4}.\quad  OUTER MEASURES}
\setlength\parindent{2em}
\setcounter{paragraph}{1}
\renewcommand{\thepage}{1-\S4-\arabic{page}}

\qquad
Let $\XX$ be a nonempty set $-$then the pair $(\XX, \sP(\XX))$ is a measurable space.
\\

\qquad{\bf 1.4.1:}\quad  
{\small\bf DEFINITION} \ 
A monotone function $\mu^* \hsx : \hsx \sP(\XX) \ra [0, +\infty]$ is an \un{outer measure} 
provided $\mu^* (\emptyset) = 0$ and $\mu^*$ is $\sigma$-subadditive on $\sP(\XX)$.
\\

\qquad{\bf 1.4.2:}\quad  
{\small\bf DEFINITION} \ 
Let $\mu^*$ be an outer measure $-$then a set $E \in \sP(\XX)$ is \un{$\mu^*$-measurable} if for every $A \in \sP(\XX)$,
\[
\mu^*(A) 
\ = \ 
\mu^* (A \cap E) + \mu^* (A \cap E^\text{c}).
\]
\\[-1.cm]

\qquad{\bf 1.4.3:}\quad  
{\small\bf NOTATION} \ 
$\sM(\mu^*)$ is the collection of all $\mu^*$-measurable sets $E \in \sP(\XX)$.
\\

\qquad{\bf 1.4.4:}\quad  
{\small\bf THEOREM} \ 
$\sM(\mu^*)$ is a $\sigma$-algebra.
\\

\qquad{\bf 1.4.5:}\quad  
{\small\bf NOTATION} \ 
Let $\mu$  be the restriction of $\mu^*$ to $\sM(\mu^*)$.  
\\

\qquad{\bf 1.4.6:}\quad  
{\small\bf THEOREM} \ 
$\mu$ is a positive measure.
\\

\qquad{\bf 1.4.7:}\quad  
{\small\bf THEOREM} \ 
$\mu$ is a complete measure.
\\

\qquad{\bf 1.4.8:}\quad  
{\small\bf DEFINITION} \ 
An outer measure $\mu^*$  is said to be \un{regular} if every $E \in \sP(\XX)$ is contained in a $\mu^*$-measurable set $F$ 
of equal outer measure.  
\\[-.5cm]

[In symbols: \ 
$\forall \ E \in \sP(\XX)$ $\exists$ $F \in \sM(\mu^*)$ : $F \supset E$ \& $\mu^*(F) = \mu^*(E)$.]
\\

\qquad{\bf 1.4.9:}\quad  
{\small\bf DEFINITION} \ 
An outer measure $\mu^*$ on a topological space $(\XX, \tau)$ is 
\un{Borel} if $\sB(\XX) \subset \sM(\mu^*)$ and is 
\un{Borel regular} if in addition for every $E \in \sP(\XX)$ there exists an $F \in \sB(\XX)$ such that 
$F \supset E$ and $\mu^*(F) = \mu^* (E)$.
\\

\qquad{\bf 1.4.10:}\quad  
{\small\bf DEFINITION} \ 
An outer measure $\mu^*$ on a metric space $(\XX, d)$ is a 
\un{metric outer measure} if 
\[
\mu^* (E \cup F) 
\ = \ 
\mu^* (E) + \mu^* (F)
\]
for all sets $E$, $F \in \sP(\XX)$ such that $\dist (E, F) > 0$.
\\

\qquad{\bf 1.4.11:}\quad  
{\small\bf THEOREM} \ 
An outer measure on a metric space $(\XX, d)$ is Borel iff $\mu^*$ is a metric outer measure.

\chapter{
$\boldsymbol{\S}$\textbf{1.5}.\quad  LEBESGUE MEASURES}
\setlength\parindent{2em}
\setcounter{paragraph}{1}
\renewcommand{\thepage}{1-\S5-\arabic{page}}

\qquad
\qquad{\bf 1.5.1:}\  
{\small\bf NOTATION} \ 
$\tL^{n^*}$ is outer Lebesgue measure on $\R^n$.
\\

\qquad{\bf 1.5.2:}\ 
{\small\bf DEFINITION} \ 
$\sM_\tL^n$ $(= \sM(\tL^{n^*}))$ is the $\sigma$-algebra comprised of the $\tL^{n^*}$-measurable subsets of $\R^n$, 
the members of $\sM_\tL^n$ being referred to as the \un{Lebesgue measurable subsets} of $\R^n$.
\\

\qquad{\bf 1.5.3:}\  
{\small\bf THEOREM} \ 
$\tL^{n^*}$  is a metric outer measure, hence 
\[
\sB(\R^n)
\subset
\sM_\tL^n.
\]
\\[-1cm]

\qquad{\bf 1.5.4:}\ 
{\small\bf NOTATION} \ 
$\tL^n$ is the restriction of $\tL^{n^*}$ to $\sM_\tL^n$, the \un{Lebesgue} \un{measure} on $\R^n$.
\\

\qquad{\bf 1.5.5:}\ 
{\small\bf THEOREM} \ 
$\tL^n$ is a complete measure and is the completion of the restriction of $\tL^n$ to $\sB(\R^n)$.
\\

The restriction of $\tL^n$ to $\sB(\R^n)$ is locally finite and Borel regular, thus is Radon.
\\

\qquad{\bf 1.5.6:}\ 
{\small\bf NOTATION} \ 
Put
\[
\omega_n 
\ = \ 
\frac{\pi^{n/2}}{\Gamma (1 + n/2)},
\]
the Lebesgue measure of the unit ball in $\R^n$.
\\

\qquad{\bf 1.5.7:}\ 
{\small\bf ISODIAMETRIC INEQUALITY} \ 
For every bounded Borel set $E \subset \R^n$, 
\[
\tL^n (E) 
\ \leq \ 
\omega_n \hsx \left( \frac{\diam (E)}{2}\right)^n.
\]
\\[-1cm]

If 
\[
B(x, r) 
\ = \ 
\{y \in \R^n \hsx : \hsx \norm{y - x} \leq r\},
\]
then
\[
\tL^n (B(x, r))
\ = \ 
\omega_n \hsy r^n.
\]
\\[-1cm]

Thus the interpretation of the isodiametric inequality is that the Lebesgue measure of $E$ cannot exceed 
the Lebesgue measure of a ball with the same diameter as that of $E$, 
i.e., among all $E$ with a given diameter $d$, a ball $B(x, r)$ with diameter $d$ has Lebesgue measure 
$\ds \omega_n \hsy \left(\frac{d}{2}\right)^n$ 
and
\[
\tL^n (E) 
\ \leq \ 
\omega_n \hsx 
\left(\frac{d}{2}\right)^n.
\]
\\[-1cm]

\qquad{\bf 1.5.8:}\ 
{\small\bf NOTATION} \ 
Given a nonsingular linear transformation $T : \R^n \ra \R^n$, 
let $M_T$ be the matrix of $T$ per the standard basis of $\R^n$.
\\

\qquad{\bf 1.5.9:}\ 
{\small\bf LEMMA} \ 
\[
E \in \sB( \R^n)
\implies T(E) \in \sB( \R^n)
\]
and
\[
\tL^n (T(E)) 
\ = \ 
\abs{\det (M_T)} \hsx \tL^n (E).
\]
\\[-1cm]

\qquad{\bf 1.5.10:}\ 
{\small\bf LEMMA} \ 
\[
E \in \sM_\tL^n 
\implies T(E) \in \sM_\tL^n
\]
and
\[
\tL^n (T(E)) 
\ = \ 
\abs{\det (M_T)} \hsx \tL^n (E).
\]
\\[-1cm]

\qquad{\bf 1.5.11:}\ 
{\small\bf LEMMA} \ 
$\forall \ E \subset \R^n$, 
\[
\tL^{n^*} (T(E)) 
\ = \ 
\abs{\det (M_T)} \hsx \tL^{n^*} (E).
\]
\\[-1cm]

On general grounds, $C_c(\R^n)$ is dense in $\Lp^p(\R^n)$ $(1 \leq p < +\infty)$ 
(cf. 1.3.14).  
But more is true:
\\

\qquad{\bf 1.5.12:}\ 
{\small\bf THEOREM} \ 
$C_c^\infty (\R^n)$ is dense in $\Lp^p(\R^n)$ $(1 \leq p < +\infty)$.

\chapter{
$\boldsymbol{\S}$\textbf{1.6}.\quad  HAUSDORFF MEASURES}
\setlength\parindent{2em}
\setcounter{paragraph}{1}
\renewcommand{\thepage}{1-\S6-\arabic{page}}

\qquad
In what follows, take $\XX = \R^n$.
\\

\qquad{\bf 1.6.1:}\ 
{\small\bf NOTATION} \ 
Given $s \in [0, +\infty[$, put
\[
\omega_s 
\ = \ 
\frac{\pi^{s/2}}{\Gamma (1 + s/2)}.
\]
\\[-1cm]

\qquad{\bf 1.6.2:}\ 
{\small\bf NOTATION} \ 
Given $0 < \delta \leq +\infty$ and a subset $E \subset \XX$, put
\[
\sH_\delta^s (E) 
\ = \ 
\frac{\omega_s}{2^s} \ 
\inf 
\left\{
\sum\limits_{k = 1}^\infty \ (\diam (E_k))^s \hsx : \hsx E \subset \bigcup\limits_{k = 1}^\infty \ E_k, \ 
\diam(E_k) \leq \delta
\right\}.
\]
\\[-1cm]

\qquad{\bf 1.6.3:}\ 
{\small\bf SUBLEMMA} \ 
\[
\delta_1 
\ \leq \ 
\delta_2 
\implies 
\sH_{\delta_1}^s (E) 
\ \geq \ 
\sH_{\delta_2}^s (E).
\]
\\[-1cm]

\qquad{\bf 1.6.4:}\ 
{\small\bf LEMMA} \ 
$\forall \ E \subset \XX$, 
\begin{align*}
\sH^s (E) \ 
&\equiv\ 
\lim\limits_{\delta \downarrow 0} \ \sH_\delta^s (E) 
\\[8pt]
&=\ 
\sup\limits_{\delta > 0} \ \sH_\delta^s (E) 
\end{align*}
exists.
\\

\qquad{\bf 1.6.5:}\ 
{\small\bf THEOREM} \ 
\[
\sH^s (E) \hsx : \hsx \sP(\XX) \ra [0, +\infty]
\]
is a metric outer measure, 
the 
\un{$s$-dimensional Hausdorff outer measure} 
on $\XX$, hence 
$\sH^s$ is Borel, hence
\[
\sB(\XX) 
\subset 
\sM(\sH^s).
\]
\\[-.75cm]

\qquad{\bf 1.6.6:}\ 
{\small\bf LEMMA} \ 
$\sH^s $ is Borel regular.
\\[-.5cm]

[In fact, if $E \in \sP(\XX)$, then there exists a $G_\delta$ set $G \supset E$ such that 
$\sH^s (G) = \sH^s (E)$ and, of course , $G \in \sB(\XX)$.]
\\

\qquad{\bf 1.6.7:}\ 
{\small\bf \un{N.B.}} \ 
The restriction of $\sH^s$ to $\sM(\sH^s)$ is a complete measure.
\\

\qquad{\bf 1.6.8:}\ 
{\small\bf LEMMA} \ 
$\forall \ x \in \XX$, 
\[
\sH^s (x + E)  \ = \ \sH^s (E) 
\]
and $\forall \ t > 0$, 
\[
\sH^s (t \hsy E)  \ = \ t^s \hsy \sH^s (E).
\]
\\[-1cm]

\qquad{\bf 1.6.9:}\ 
{\small\bf LEMMA} \ 
\[
\sH^s
\ \equiv \ 
0
\]
if $s > n$.
\\

Therefore matters reduce to the  consideration of $\sH^s$ in the range $0 \leq s \leq n$.
\\

\qquad{\bf 1.6.10:}\ 
{\small\bf LEMMA} \ 
$\sH^0$ is counting measure.  
Moreover, $\sM(\sH^0) = \sP(\XX)$.
\\

Therefore matters reduce to the  consideration of $\sH^s$ in the range $0 < s \leq n$.
\\[-.5cm]

Recall that $\tL^{n^*}$ is outer Lebesgue measure on $\R^n$.
\\

\qquad{\bf 1.6.11:}\ 
{\small\bf THEOREM} \ 
\[
\tL^{n^*} 
\ = \ 
\sH^n.
\]
\\[-1cm]

Therefore matters reduce to the  consideration of $\sH^s$ in the range $0 < s < n$.
\\

\qquad{\bf 1.6.12:}\ 
{\small\bf LEMMA} \ 
Let $E \subset \R^n$ and let $0 \leq s < t < +\infty$.
\\[-.25cm]

\qquad \textbullet \quad 
If $\sH^s (E) < +\infty$, then $\sH^t (E) = 0$.
\\[-.5cm]

\qquad \textbullet \quad 
If $\sH^t (E) > 0$, then $\sH^s (E) = +\infty$.
\\[-.25cm]

PROOF \ 
The second point is implied by the first point.  
To arrive at the latter, choose sets $E_k$ such that 
$\diam (E_k) \leq \delta$, $E \subset \ds\bigcup\limits_{k = 1}^\infty \ E_k$, and 
\[
\frac{\omega_s}{2^s} \ \sum\limits_{k = 1}^\infty \ (\diam (E_k) )^s
\ \leq \ 
\sH_\delta^s (E) + 1 
\ \leq \ 
\sH^s (E) + 1 .
\]
Then
\begin{align*}
\sH_\delta^t (E)  \ 
&\leq\ 
\frac{\omega_t}{2^t} \ \sum\limits_{k = 1}^\infty \ (\diam (E_k))^t
\\[15pt]
&=\ 
\frac{\omega_t}{\omega_s} \hsy 2^{s - t} \hsy \frac{\omega_s}{2^s}  \ 
\sum\limits_{k = 1}^\infty \ (\diam (E_k))^s \hsx (\diam (E_k))^{t - s}
\\[15pt]
&\leq\ 
\frac{\omega_t}{\omega_s} \hsy 2^{s - t} \hsx 
(\sH^s (E) + 1) \delta^{t - s}.
\end{align*}
Noting that $t - s > 0$, send $\delta \downarrow 0$ to conclude that $\sH^t (E) = 0$.
\\

\qquad{\bf 1.6.13:}\ 
{\small\bf LEMMA} \ 
Let $E \subset \R^n$ $-$then there exists at most one point $s^* \in [0, +\infty[$ such that
$\sH^{s^*} (E) \in \hsx ]0, +\infty[$. 
\\[-.25cm]

PROOF \ 
Take two distinct points $s$, $t \in [0, +\infty[$ with $s < t$.  
If $\sH^s(E) \in \hsx ]0, +\infty[$, then $\sH^t (E) = 0$ 
while if $\sH^t (E)\in \hsx ]0, +\infty[$, then $\sH^s (E) = +\infty$.
\\

\qquad{\bf 1.6.14:}\ 
{\small\bf NOTATION} \ 
Given $E \in \sP(\XX)$, denote by $\sH^\bdot (E)$ the function 
\[
\begin{cases}
\ [0, +\infty[ \ &\ra \ [0, +\infty]
\\
\hspace{1cm} s  &\ra \ \sH^s (E)
\end{cases}
.
\]
\\[-1cm]

\qquad{\bf 1.6.15:}\ 
{\small\bf LEMMA} \ 
$\sH^\bdot (E)$ is a decreasing function on $[0, +\infty[$ which vanishes on $]n, +\infty[$.
\\

\qquad{\bf 1.6.16:}\ 
{\small\bf THEOREM} \ 
There are three possibilities for the range of $\sH^\bdot (E)$.
\\

\qquad (i) \hspace{0.45cm}
$\sH^\bdot (E)$ assumes one value, viz. 0.
\\[-.5cm]

\qquad (ii) \hspace{0.35cm}
 $\sH^\bdot (E)$ assumes two values, viz. $+\infty$, and 0.
\\[-.5cm]

\qquad (iii) \hspace{0.25cm}
 $\sH^\bdot (E)$ assumes three values, viz. $+\infty$,  0 and one finite positive value $s^*$.
\\

\qquad{\bf 1.6.17:}\ 
{\small\bf EXAMPLE} \ 

\[
\begin{cases}
&\sH^s (\R^n) \ = \ +\infty \quad (s \in [0,n])
\\[8pt]
&\sH^s (\R^n) \ = \ 0  \quad (s \in \hsx ]n,+\infty[)
\end{cases}
.
\]
\\[-.5cm]

\qquad{\bf 1.6.18:}\ 
{\small\bf LEMMA} \ 
If $\sH^\bdot (E)$ assumes a finite positive value at some point $s^* \in [0,+\infty[$, then 

\[
\begin{cases}
&\sH^s (E) \ = \ +\infty \quad (s \in [0,s^*[)
\\[8pt]
&\sH^s (E) \ = \ 0  \quad (s \in \hsx ]s^*,+\infty[)
\end{cases}
.
\]
\\[-.25cm]

\qquad{\bf 1.6.19:}\ 
{\small\bf \un{N.B.}}  \ 
$\sH^\bdot (E)$ is identically zero on $[0, +\infty[$ iff $E = \emptyset$.
\\

\qquad{\bf 1.6.20:}\ 
{\small\bf \un{N.B.}} \ 
If $E \neq \emptyset$, then $\sH^\bdot (E)$ has exactly one point of discontinuity in $[0, +\infty[$ and it belongs to $[0, n]$.
\\

\qquad{\bf 1.6.21:}\ 
{\small\bf DEFINITION} \ 
The \un{Hausdorff dimension} of a nonempty subset $E$ of $\R^n$, denoted $\dim_H (E)$, 
is the unique number $s^* \in [0, +\infty[$ at which $\sH^\bdot$ is discontinuous.
\\

\qquad{\bf 1.6.22:}\ 
{\small\bf \un{N.B.}} \ 

\[
\dim_H (\emptyset) \ = \ 0.
\]
\\[-1cm]

\qquad{\bf 1.6.23:}\ 
{\small\bf LEMMA} \ 
\begin{align*}
\dim_H (E)
&=\ 
\sup \ \{s \in [0, +\infty[ \ : \sH^s(E) > 0\}
\\[8pt]
&=\ 
\sup \ \{s \in [0, +\infty[ \ : \sH^s(E) = +\infty\}
\end{align*}
and
\begin{align*}
\dim_H (E)
&=\ 
\inf \ \{s \in [0, +\infty[ \ : \sH^s(E) = 0\}
\\[8pt]
&=\ 
\inf \ \{s \in [0, +\infty[ \ : \sH^s(E) < +\infty\}.
\end{align*}
\\[-1cm]

\qquad{\bf 1.6.24:}\ 
{\small\bf LEMMA} \ 
If $\sH^s(E) \in \hsx ]0, +\infty[$, then $s = \dim_H (E)$.
\\

\qquad{\bf 1.6.25:}\ 
{\small\bf EXAMPLE}\ 
\[
\dim_H (\R^n)
\ = \ 
n 
\quad 
\text{but} \quad 
\sH^n (\R^n) = +\infty.
\]
\\[-1cm]

\qquad{\bf 1.6.26:}\ 
{\small\bf LEMMA} \ 
If $E \in \sP(\XX)$ is countable, then $\dim_H (E) = 0$.
\\

\qquad{\bf 1.6.27:}\ 
{\small\bf LEMMA} \ 
If $E \in \sP(\XX)$ has a nonempty interior, then $\dim_H (E) = n$.
\\[-.5cm]

In particular: \ If $U \in \sP(\XX)$ is open and nonempty, then $\dim_H (U) = n$, so 
$\sH^s(U) = +\infty$ $(0 < s < n)$.
\\

\qquad{\bf 1.6.28:}\ 
{\small\bf THEOREM} \ 
For every $s \in ]0, n[$, there exists a compact $K \subset \R^n$ such that $\dim_H (K) = s$.
\\

\qquad{\bf 1.6.29:}\ 
{\small\bf EXAMPLE} \ 
Take $n = 1$ and let $C \subset \R^1$ be the Cantor set $-$then
\[
\dim_H (C) 
\ = \ 
\frac{\log \hsy 2}{\log \hsy 3}.
\]
\\[-1cm]

\qquad{\bf 1.6.30:}\ 
{\small\bf DEFINITION} \ 
A metric outer measure $\mu^*$ on $\XX$ is \un{locally finite} if 
$\mu^* (K) < +\infty$ for every $K \in \sK(\XX)$.
\\

\qquad{\bf 1.6.31:}\ 
{\small\bf THEOREM} \ 
Suppose that $\mu^*$ is locally finite $-$then for every Borel set $E \in \sB(\XX)$, 
\[
\mu^*(E) 
\ = \ 
\{\inf \ \mu^*(U) \hsx : \hsx U \supset E, \ U \in \sO(\XX)\}
\]
and
\[
\mu^*(E) 
\ = \ 
\{\sup \ \mu^*(K) \hsx : \hsx K \subset E, \ K \in \sK(\XX)\}.
\]
\\[-1cm]

If $U \subset \XX$ is open and nonempty, then
\[
\sH^s (U) 
\ = \ 
+\infty 
\qquad (0 < s < \dim_H (U) = n).  
\]
\\[-1.cm]

\qquad{\bf 1.6.32:}\ 
{\small\bf SCHOLIUM} \ 
$\sH^s$ is not locally finite if $0 < s < n$.
\\[-.5cm]

[Pretend it was $-$then for a generic $K \in \sK(\XX)$, 
\begin{align*}
\sH^s (K) \ 
&=\ 
\inf \{\sH^s (U) \hsx : \hsx U \supset K, \ U \in \sO(\XX)\}
\\[8pt]
&=\ 
+\infty \hsx \ldots \hsx .]
\end{align*}
\\[-1.cm]

\qquad{\bf 1.6.33:}\ 
{\small\bf \un{N.B.}}\ 
Bearing in mind that $\sB(\XX) \subset \sM(\sH^s)$, it follows that the restriction of $\sH^s$ to $\sB(\XX)$ is not Radon.
\\

\qquad{\bf 1.6.34:}\ 
{\small\bf THEOREM} \ 
Let $\Phi : \R^n \ra \R^n$ be an isometry (a distance preserving bijection) and suppose that $E \in \sP(\XX)$ $-$then
\[
\Phi(E) \in \sM(\sH^s) 
\iff 
E \in \sM(\sH^s) 
\quad (s \in [0, +\infty[).
\]

[Note: \ 
The assumption that $\Phi$ is an isometry implies that
\[
\sH^s (\Phi(E))
\ = \ 
\sH^s (E).]
\]

\chapter{
\text{SECTION 2: \quad DIFFERENTIATION THEORY}
\vspace{1.25cm}\\
$\boldsymbol{\S}$\textbf{2.1}.\quad  SCALAR FUNCTIONS}
\setlength\parindent{2em}
\renewcommand{\thepage}{2-\S1-\arabic{page}}


\qquad
Let $\Omega$ be a nonempty open subset of $\R^n$ and let 
$f: \Omega \ra \R$ be a function.
\\

\qquad{\bf 2.1.1:}\  
{\small\bf DEFINITION} \ 
$f$ is \un{differentiable} at a point $x_0 \in \Omega$ if there exists a linear function  
$T : \R^n \ra \R$ (depending on $x_0$) such that 

\[
\lim\limits_{h \ra 0} \ \frac{f(x_0 + h) - f(x_0) - T(h)}{\norm{h}} 
\ = \ 0.
\]
\\[-.75cm]

\qquad{\bf 2.1.2:}\  
{\small\bf \un{N.B.}} \ 
Consider the situation when $n = 1$, $\Omega = \R$ and suppose that $f: \R \ra \R$ is differentiable at $x_0$ 
in the traditional sense, i.e., 
\[
f^\prime (x_0) 
\ = \ 
\lim\limits_{h \ra 0} \ \frac{f(x_0 + h) - f(x_0)}{h} .
\]
Then $f$ is differentiable at $x_0$.
Thus view the number $f^\prime (x_0)$ as the linear map $\R \ra \R$ that sends $h$ to $f^\prime (x_0) (h)$, hence

\allowdisplaybreaks
\begin{align*}
\frac{f(x_0 + h) - f(x_0) - f^\prime (x_0) (h)}{h}\ 
&=\ 
\frac{f(x_0 + h) - f(x_0)}{h} - \frac{f^\prime (x_0) (h)}{h}
\\[11pt]
&=\ 
\frac{f(x_0 + h) - f(x_0)}{h} - f^\prime (x_0)  (1)
\\[11pt]
&\ra  
f^\prime (x_0)  - f^\prime (x_0)  \quad (h \ra 0)
\\[11pt]
&=\ 
0.
\end{align*}
\\[-.75cm]

$T$ is called the \un{differential} of $f$ at $x_0$ and is denoted by $\df (x_0)$.
\\[-.5cm]

[Note: \ 
The differential is unique, if it exists at all.  
Proof: \ 
Per the definition, suppose that $T = T_1$ and $T = T_2$ $-$then $\forall \ h \neq 0$, 

\[
\abs{(T_1 - T_2) (h)} 
\ < \ 
\abs{f(x_0 + h) - f(x_0) - T_1 (h)} + \abs{f(x_0 + h) - f(x_0) - T_2 (h)}
\]

\hspace{1.5cm} $\implies$

\[
\frac{\abs{(T_1 - T_2) (h)}}{\norm{h}} 
\ra 
0 
\quad \text{as $h \ra 0$}
\]

\hspace{1.5cm} $\implies$

\[
\frac{\abs{(T_1 - T_2) (t \hsy h)}}{\norm{t \hsy h}} 
\ra 
0 
\quad \text{as $t \ra 0$}
\]

\hspace{1.5cm} $\implies$

\[
(T_1 - T_2) (h)
\ = \ 
0 
\implies 
T_1 \ = \ T_2.]
\]
\\[-.75cm]

\qquad{\bf 2.1.3:}\  
{\small\bf \un{N.B.}} \ 
$f$ is \un{differentiable in $\Omega$} if $f$ is differentiable at every point of $\Omega$.
\\

\qquad{\bf 2.1.4:}\  
{\small\bf EXAMPLE} \ 
Take $\Omega = \R^n$ $-$then polynomials in several variables are everywhere differentiable.
\\

\qquad{\bf 2.1.5:}\  
{\small\bf EXAMPLE} \ 
Take $\Omega = \R^n$ and let $T : \R^n \ra \R$ be linear $-$then $\td T (x_0) = T$.
\\

\qquad{\bf 2.1.6:}\  
{\small\bf LEMMA} \ 
If $f$ is differentiable at $x_0 \in \Omega$, then $f$ is continuous at $x_0$.
\\[-.5cm]

[Given $h \neq 0$, write
\[
\abs{f(x_0 + h) - f(x_0) - T (h)} 
\ \leq \ 
\norm{h} \hsx 
\abs{\frac{f(x_0 + h) - f(x_0) - T (h)}{h}}
\]
to conclude that
\[
f(x_0 + h) - f(x_0) 
\ra 0 
\quad \text{as $h \ra 0$}.]
\]
\\[-1cm]


[Note: \ 
Since $T$ is linear, 
$\lim\limits_{h \ra 0} \ T(h) = 0$.]
\\

Given $x_0 \in \Omega$, suppose that $B(x_0, r_0)$ is contained in $\Omega$ $-$then for each nonzero 
$v \in \R^n$, $x_0 + t \hsy v \in \Omega$ for $\abs{t} \leq \ds\frac{r_0}{\norm{v}}$.
\\[.25cm]

\qquad{\bf 2.1.7:}\  
{\small\bf DEFINITION} \ 
The \un{directional derivative} of $f$ at $x_0$ in the direction $v$ is 

\[
\lim\limits_{t \ra 0} \ \frac{f(x_0 + t \hsy v) - f(x_0)}{t},
\]
denoted

\[
\frac{\partial f}{\partial v} (x_0).
\]
\\[-1cm]

[Note: \ 
The underlying assumption is that the limit exists and is finite.]
\\

\qquad{\bf 2.1.8:}\  
{\small\bf \un{N.B.}} \ 
$\forall \ \lambda \neq 0$, 
\[
\frac{\partial f}{\partial (\lambda \hsy v)} (x_0) 
\ = \ 
\lambda \hsx \frac{\partial f}{\partial v} (x_0).
\]
\\

\qquad{\bf 2.1.9:}\  
{\small\bf LEMMA} \ 
If $f$ is differentiable at $x_0$, then 
$\ds\frac{\partial f}{\partial v} (x_0)$ exists for all $v \neq 0$ and 

\[
\frac{\partial f}{\partial v} (x_0) 
\ = \ 
\df (x_0) (v).
\]

[Observe that 

\[
\abs{ \frac{f(x_0 + t \hsy v) - f(x_0)}{t} - T(v)}
\ = \ 
\abs{ \frac{f(x_0 + t \hsy v) - f(x_0) - T( t \hsy v)}{t}}.\big]
\]
\\


\qquad{\bf 2.1.10:}\  
{\small\bf EXAMPLE} \ 
The function 
\[
f(x, y) \ = \ 
\begin{cases}
\ \ds\frac{x^2 \hsy y}{x^2 + y^2} \hspace{0.7cm} \text{if $x \neq 0$}
\\[11pt]
\hspace{0.5cm}  0 \hspace{1.3cm} \text{if $x = 0$}
\end{cases}
\]
is continuous at $(0,0)$ and all its directional derivatives exist at $(0,0)$.  
Still, the differential $\td f(0,0)$ does not exist.
\\[-.25cm]

[To see this last point, suppose instead that $\td f(0,0)$ does exist, thus being linear, 
\[
\df (0, 0) (1,0) + \df (0, 0) (0,1) 
\ = \ 
\df (0, 0) (1,1).
\]
On the other hand, 
\[
\begin{cases}
\ \ds\df (0, 0) (1,0)  \ = \ \frac{\partial f}{\partial (1,0)} (0,0) \ = \ 0
\\[18pt]
\ \ds\df (0, 0) (0,1)  \ = \ \frac{\partial f}{\partial (0,1)} (0,0) \ = \ 0
\end{cases}
.
\]
Meanwhile
\[
\frac{\partial f}{\partial (1,1)} (0,0) \ = \ \frac{1}{2}
\]

\hspace{1.5cm} $\implies$
\[
\df (0, 0) (1,1)
\ = \ 
\frac{1}{2}.
\]
I.e.: 
\[
0 + 0 
\ = \ 
\frac{1}{2}.
\]
Contradiction.]
\\

\qquad{\bf 2.1.11:}\  
{\small\bf REMARK} \ 
If $\Omega$ is convex and if 
$f: \Omega \ra \R$ is convex, then $f$ is differentiable at $x_0$ iff $f$ has ordinary partial derivatives at $x_0$.
\\

Suppose that $\langle \hsx, \hsx \rangle$ is the standard inner product in $\R^n$.
Since the differential of $f$ at $x_0$ is a linear function from $\R^n$ to $\R$, there is a unique vector 
$\nabla f(x_0) \in \R^n$ such that for all $h \in \R^n$,
\[
\df (x_0) (h)
\ = \ 
\langle h, \nabla f(x_0) \rangle.
\]
\\

\qquad{\bf 2.1.12:}\  
{\small\bf DEFINITION} \ 
$\nabla f(x_0)$ is called the \un{gradient} of $f$ at $x_0$.
\\

\qquad{\bf 2.1.13:}\  
{\small\bf NOTATION} \ 
Let $(e_1, e_2, \ldots, e_n)$ be the standard basis for $\R^n$ and let 
$(x_1, x_2, \ldots, x_n)$ be the associated system of coordinates.
\\

\qquad{\bf 2.1.14:}\  
{\small\bf DEFINITION} \ 
The derivative of $f$ at $x_0$ in the direction $e_i$ is called the \un{partial derivative} of $f$ w.r.t. $x_i$, denoted
\[
\frac{\partial f}{\partial x_i}  (x_0).
\]
\\[-.75cm]

\qquad{\bf 2.1.15:}\  
{\small\bf LEMMA} \ 

\[
\nabla f(x_0)
\ = \ 
\left(
\frac{\partial f}{\partial x_1} (x_0), 
\frac{\partial f}{\partial x_2} (x_0), 
\ldots, 
\frac{\partial f}{\partial x_n} (x_0)
\right).
\]
\\[-.75cm]

\qquad{\bf 2.1.16:}\  
{\small\bf DEFINITION} \ 
The \un{Jacobian matrix} of $f$ at $x_0$ is the $1 \times n$ matrix

\[
\tD f(x_0) 
\ = \ 
\left[
\frac{\partial f}{\partial x_1} (x_0), 
\frac{\partial f}{\partial x_2} (x_0), 
\ldots, 
\frac{\partial f}{\partial x_n} (x_0)
\right].
\]
\\[-.75cm]

\qquad{\bf 2.1.17:}\  
{\small\bf LEMMA} \ 
For all $h \in \R^n$,
\begin{align*}
\df (x_0) (h) \ 
&=\ 
\sum\limits_{i = 1}^n \ \df (x_0) \hsy (e_i) \hsy h_i
\\[15pt]
&=\ 
\sum\limits_{i = 1}^n \ \frac{\partial f}{\partial x_i} \hsy (x_0)  \hsy h_i
\\[15pt]
&=\ 
\tD f(x_0) 
\begin{bmatrix}
h_1\\
h_2\\
\vdots\\
h_n
\end{bmatrix}
\\[15pt]
&=\ 
\langle h, \nabla f(x_0) \rangle .
\end{align*}
Consider two points $x_0$, $x_0 + h$ $-$then the line segment $\ell$ joining $x_0$ and $x_0 + h$ 
is the curve $x_0 + t \hsy h$ $(0 \leq t \leq 1)$.
\\

\qquad{\bf 2.1.18:}\  
{\small\bf MEAN VALUE THEOREM} \ 
Suppose that $f$ is continuous at the points of $\ell$ and differentiable at the points of $\ell$ except perhaps the endpoints $-$then 
there exists an $s \in ]0,1[$ such that 
\[
f(x_0 + h) - f(x_0) 
\ = \ 
\df (x_0 + s \hsy h) (h).
\]

PROOF \ 
Introduce
\[
\phi(t) 
\ = \ 
f(x_0 + t \hsy h) 
\qquad (0 \leq t \leq 1).
\]
Then $\phi$ is continuous in $[0,1]$ and 
\[
\begin{cases}
&\phi(0) \ = \ f(x_0)
\\[15pt]
&\phi(1) = f(x_0 + h)
\end{cases}
, \qquad \phi^\prime (t) = \df (x_0 + t \hsy h) (h) 
\qquad (0 < t < 1).
\]
By the mean value theorem for functions of one variable, there exists an $s \in ]0,1[$ such that
\[
\phi(1) - \phi(0) 
\ = \ 
\phi^\prime (s).
\]
\\[-1cm]

\qquad{\bf 2.1.19:}\  
{\small\bf APPLICATION} \ 
Suppose that $\Omega$ is not only open but is connected as well.  
Assume: $f$ is differentiable in $\Omega$ and that $\df (x) = 0$ for every $x \in \Omega$ $-$then $f$ is a constant function.
\\[-.5cm]

\[
* \ * \ * \ * \ * \ * \ * \ * \ * \ * \ 
\]
\[
\text{APPENDIX}
\]

What has been said in 2.1.10 can be substantially generalized.  
Indeed, there are continuous functions of 2 variables which have partial derivatives almost everywhere but for which the differential fails to exists anywhere.
\\

{\small\bf EXAMPLE}
(Bruckner and Goffman\footnote[2]{\vspace{.11 cm} \textit{Real Anal. Exchange}, \textbf{6} (1980/81) no. 1, 9-65.})
 \ 

Let $Q$ be the unit closed square and let 
$D = \{r_1, r_2, \ldots\} \subset Q$ be a countable dense subset in the interior of $Q$.  
(To be concrete one could take $D = \Q \times \Q \cap Q^\circ$.)

\noindent
Let $\{\varepsilon_\ell\}$ be a sequence of postive numbers such that 
\[
\sum\limits_\ell \ \varepsilon_\ell \ < \ \infty
\quad \text{and} \quad 
B(r_\ell, \varepsilon_\ell)^\circ \subset Q.
\]
(Again, being concrete, this may be easily achieved by taking
\[
\varepsilon_\ell
\ \equiv \ 
\min\big\{\frac{1}{2^\ell}, \dist(r_\ell, \partial Q)\big\}.)
\]
Let $k_\ell$ be a sequence of postive numbers such that 
\[
\sum\limits_\ell \ k_\ell \ = \ \infty.
\]
Let $F_\ell : Q \ra \R$ be defined: 
\[
F_\ell (x) \ = \ 
\begin{cases}
\ \text{Cone(height $k_\ell$, base $B(r_\ell, \varepsilon_\ell))$} \hspace{.5cm} r \in B(r_\ell, \varepsilon_\ell)
\\[4pt]
\ \hspace{0.5cm} 0 \hspace{5.31cm} r \in Q \backslash B(r_\ell, \varepsilon_\ell)
\end{cases}
\]
Define 
\[
\raisebox{1cm}{$F \ = \ \sum\limits_\ell \ F_\ell$}
\hspace{1.5cm} 
\begin{tikzpicture}[scale=0.35]



       \begin{scope}[shift = {(1,1)}][scale=2]
          \draw[clip] ({-sqrt(91)/10},9/100) coordinate (A) 
             arc[x radius=1, y radius=0.3, start angle=180-17.46, end angle=360+17.46]
             -- (0,3) -- cycle;
          \fill[cyan] (0,5) -- (-1.3628,-0.67) -- (1.3628,-0.67) -- cycle;
          \foreach \i in {0,1,...,60} {
             \fill[white!\i!cyan] (0,5) -- (-1.3628+\i*0.006814,-0.65) -- +(0.008,0) -- cycle;
             \fill[cyan!\i!white] (0,5) -- (-0.6814+\i*0.006814,-0.65) -- +(0.008,0) -- cycle;
             \fill[cyan!\i!black!50!cyan] (0,5) -- (0.6814-\i*0.006814,-0.65) -- +(-0.008,0) -- cycle;
             \fill[black!\i!cyan!50!black] (0,5) -- (0.6814+\i*0.006814,-0.65) -- +(-0.008,0) -- cycle;
          }
          \draw[dotted] (A) arc[x radius=2, y radius=0.3, start angle=180-17.46, end angle=17.46];
       \end{scope}

       
       \begin{scope}[shift = {(2.5,1)}][scale=2]
          \draw[clip] ({-sqrt(91)/10},9/100) coordinate (A) 
             arc[x radius=1, y radius=0.3, start angle=180-17.46, end angle=360+17.46]
             -- (0,3) -- cycle;
          \fill[cyan] (0,7) -- (-1.3628,-0.67) -- (1.3628,-0.67) -- cycle;
          \foreach \i in {0,1,...,60} {
             \fill[white!\i!cyan] (0,7) -- (-1.3628+\i*0.006814,-0.65) -- +(0.008,0) -- cycle;
             \fill[cyan!\i!white] (0,7) -- (-0.6814+\i*0.006814,-0.65) -- +(0.008,0) -- cycle;
             \fill[cyan!\i!black!50!cyan] (0,7) -- (0.6814-\i*0.006814,-0.65) -- +(-0.008,0) -- cycle;
             \fill[black!\i!cyan!50!black] (0,7) -- (0.6814+\i*0.006814,-0.65) -- +(-0.008,0) -- cycle;
          }
          \draw[dotted] (A) arc[x radius=2, y radius=0.3, start angle=180-17.46, end angle=17.46];
       \end{scope}
       
       \begin{scope}[shift = {(3.0,2)}][scale=2]
          \draw[clip] ({-sqrt(91)/10},9/100) coordinate (A) 
             arc[x radius=1, y radius=0.3, start angle=180-17.46, end angle=360+17.46]
             -- (0,4) -- cycle;
          \fill[cyan] (0,4) -- (-1.3628,-0.67) -- (1.3628,-0.67) -- cycle;
          \foreach \i in {0,1,...,60} {
             \fill[white!\i!cyan] (0,4) -- (-1.3628+\i*0.006814,-0.65) -- +(0.008,0) -- cycle;
             \fill[cyan!\i!white] (0,4) -- (-0.6814+\i*0.006814,-0.65) -- +(0.008,0) -- cycle;
             \fill[cyan!\i!black!50!cyan] (0,4) -- (0.6814-\i*0.006814,-0.65) -- +(-0.008,0) -- cycle;
             \fill[black!\i!cyan!50!black] (0,4) -- (0.6814+\i*0.006814,-0.65) -- +(-0.008,0) -- cycle;
          }
          \draw[dotted] (A) arc[x radius=2, y radius=0.3, start angle=180-17.46, end angle=17.46];
       \end{scope}


       \begin{scope}[scale=1]
          \draw[clip] ({-sqrt(91)/10},9/100) coordinate (A) 
             arc[x radius=1, y radius=0.3, start angle=180-17.46, end angle=360+17.46]
             -- (0,2.5) -- cycle;
          \fill[cyan] (0,2.5) -- (-1.3628,-0.3) -- (1.3628,-0.3) -- cycle;
          \foreach \i in {0,1,...,60} {
             \fill[white!\i!cyan] (0,2.5) -- (-1.3628+\i*0.006814,-0.3) -- +(0.008,0) -- cycle;
             \fill[cyan!\i!white] (0,2.5) -- (-0.6814+\i*0.006814,-0.3) -- +(0.008,0) -- cycle;
             \fill[cyan!\i!black!50!cyan] (0,2.5) -- (0.6814-\i*0.006814,-0.3) -- +(-0.008,0) -- cycle;
             \fill[black!\i!red!50!black] (0,2.5) -- (0.6814+\i*0.006814,-0.3) -- +(-0.008,0) -- cycle;
          }
          \draw[dotted] (A) arc[x radius=1, y radius=0.3, start angle=180-17.46, end angle=17.46];
       \end{scope}

       \begin{scope}[shift = {(2,0)}][scale=2]
          \draw[clip] ({-sqrt(91)/10},9/100) coordinate (A) 
             arc[x radius=1, y radius=0.3, start angle=180-17.46, end angle=360+17.46]
             -- (0,5) -- cycle;

          \fill[cyan] (-1,10) -- (-1.3628,-0.3) -- (1.3628,-0.3) -- cycle;
          \foreach \i in {0,1,...,60} {
             \fill[white!\i!cyan] (0,10) -- (-1.3628+\i*0.006814,-0.3) -- +(0.008,0) -- cycle;
             \fill[cyan!\i!white] (0,10) -- (-0.6814+\i*0.006814,-0.3) -- +(0.008,0) -- cycle;
             \fill[cyan!\i!black!50!cyan] (0,10) -- (0.6814-\i*0.006814,-0.3) -- +(-0.008,0) -- cycle;
             \fill[black!\i!cyan!50!black] (0,10) -- (0.6814+\i*0.006814,-0.3) -- +(-0.008,0) -- cycle;
          }
          \draw[dotted] (A) arc[x radius=1, y radius=0.3, start angle=180-17.46, end angle=17.46];
       \end{scope}
       
       \begin{scope}[shift = {(4,0)}][scale=.5]
          \draw[clip] ({-sqrt(91)/10},9/100) coordinate (A) 
             arc[x radius=1, y radius=0.3, start angle=180-17.46, end angle=360+17.46]
             -- (0,5) -- cycle;

          \fill[cyan] (-1,10) -- (-1.3628,-0.3) -- (1.3628,-0.3) -- cycle;
          \foreach \i in {0,1,...,60} {
             \fill[white!\i!cyan] (0,10) -- (-1.3628+\i*0.006814,-0.3) -- +(0.008,0) -- cycle;
             \fill[cyan!\i!white] (0,10) -- (-0.6814+\i*0.006814,-0.3) -- +(0.008,0) -- cycle;
             \fill[cyan!\i!black!50!cyan] (0,10) -- (0.6814-\i*0.006814,-0.3) -- +(-0.008,0) -- cycle;
             \fill[black!\i!cyan!50!black] (0,10) -- (0.6814+\i*0.006814,-0.3) -- +(-0.008,0) -- cycle;
          }
          \draw[dotted] (A) arc[x radius=1, y radius=0.3, start angle=180-17.46, end angle=17.46];
       \end{scope}


       \begin{scope}[shift = {(1,1)}][scale=2]
          \draw[clip] ({-sqrt(91)/10},9/100) coordinate (A) 
             arc[x radius=1, y radius=0.3, start angle=180-17.46, end angle=360+17.46]
             -- (0,3) -- cycle;
          \fill[cyan] (0,5) -- (-1.3628,-0.67) -- (1.3628,-0.67) -- cycle;
          \foreach \i in {0,1,...,60} {
             \fill[white!\i!cyan] (0,5) -- (-1.3628+\i*0.006814,-0.65) -- +(0.008,0) -- cycle;
             \fill[cyan!\i!white] (0,5) -- (-0.6814+\i*0.006814,-0.65) -- +(0.008,0) -- cycle;
             \fill[cyan!\i!black!50!cyan] (0,5) -- (0.6814-\i*0.006814,-0.65) -- +(-0.008,0) -- cycle;
             \fill[black!\i!cyan!50!black] (0,5) -- (0.6814+\i*0.006814,-0.65) -- +(-0.008,0) -- cycle;
          }
          \draw[dotted] (A) arc[x radius=2, y radius=0.3, start angle=180-17.46, end angle=17.46];
       \end{scope}

       
       \begin{scope}[shift = {(2.5,1)}][scale=2]
          \draw[clip] ({-sqrt(91)/10},9/100) coordinate (A) 
             arc[x radius=1, y radius=0.3, start angle=180-17.46, end angle=360+17.46]
             -- (0,3) -- cycle;
          \fill[cyan] (0,7) -- (-1.3628,-0.67) -- (1.3628,-0.67) -- cycle;
          \foreach \i in {0,1,...,60} {
             \fill[white!\i!cyan] (0,7) -- (-1.3628+\i*0.006814,-0.65) -- +(0.008,0) -- cycle;
             \fill[cyan!\i!white] (0,7) -- (-0.6814+\i*0.006814,-0.65) -- +(0.008,0) -- cycle;
             \fill[cyan!\i!black!50!cyan] (0,7) -- (0.6814-\i*0.006814,-0.65) -- +(-0.008,0) -- cycle;
             \fill[black!\i!cyan!50!black] (0,7) -- (0.6814+\i*0.006814,-0.65) -- +(-0.008,0) -- cycle;
          }
          \draw[dotted] (A) arc[x radius=2, y radius=0.3, start angle=180-17.46, end angle=17.46];
       \end{scope}


       \begin{scope}[shift = {(1,-1)}][scale=2]
          \draw[clip] ({-sqrt(91)/10},9/100) coordinate (A) 
             arc[x radius=1, y radius=0.3, start angle=180-17.46, end angle=360+17.46]
             -- (0,3) -- cycle;
          \fill[cyan] (0,5) -- (-1.3628,-0.67) -- (1.3628,-0.67) -- cycle;
          \foreach \i in {0,1,...,60} {
             \fill[white!\i!cyan] (0,5) -- (-1.3628+\i*0.006814,-0.65) -- +(0.008,0) -- cycle;
             \fill[cyan!\i!white] (0,5) -- (-0.6814+\i*0.006814,-0.65) -- +(0.008,0) -- cycle;
             \fill[cyan!\i!black!50!cyan] (0,5) -- (0.6814-\i*0.006814,-0.65) -- +(-0.008,0) -- cycle;
             \fill[black!\i!cyan!50!black] (0,5) -- (0.6814+\i*0.006814,-0.65) -- +(-0.008,0) -- cycle;
          }
          \draw[dotted] (A) arc[x radius=2, y radius=0.3, start angle=180-17.46, end angle=17.46];
       \end{scope}

       
       \begin{scope}[shift = {(2.75,-1)}][scale=2]
          \draw[clip] ({-sqrt(91)/10},9/100) coordinate (A) 
             arc[x radius=1, y radius=0.3, start angle=180-17.46, end angle=360+17.46]
             -- (0,2) -- cycle;
          \fill[cyan] (0,2) -- (-1.3628,-0.67) -- (1.3628,-0.67) -- cycle;
          \foreach \i in {0,1,...,60} {
             \fill[white!\i!cyan] (0,2) -- (-1.3628+\i*0.006814,-0.65) -- +(0.008,0) -- cycle;
             \fill[cyan!\i!white] (0,2) -- (-0.6814+\i*0.006814,-0.65) -- +(0.008,0) -- cycle;
             \fill[cyan!\i!black!50!cyan] (0,2) -- (0.6814-\i*0.006814,-0.65) -- +(-0.008,0) -- cycle;
             \fill[black!\i!cyan!50!black] (0,2) -- (0.6814+\i*0.006814,-0.65) -- +(-0.008,0) -- cycle;
          }
          \draw[dotted] (A) arc[x radius=2, y radius=0.3, start angle=180-17.46, end angle=17.46];
       \end{scope}
       
       \begin{scope}[shift = {(3.75,-1)}][scale=2]
          \draw[clip] ({-sqrt(91)/10},9/100) coordinate (A) 
             arc[x radius=1, y radius=0.3, start angle=180-17.46, end angle=360+17.46]
             -- (0,3.5) -- cycle;
          \fill[cyan] (0,3.5) -- (-1.3628,-0.67) -- (1.3628,-0.67) -- cycle;
          \foreach \i in {0,1,...,60} {
             \fill[white!\i!cyan] (0,3.5) -- (-1.3628+\i*0.006814,-0.65) -- +(0.008,0) -- cycle;
             \fill[cyan!\i!white] (0,3.5) -- (-0.6814+\i*0.006814,-0.65) -- +(0.008,0) -- cycle;
             \fill[cyan!\i!black!50!cyan] (0,3.5) -- (0.6814-\i*0.006814,-0.65) -- +(-0.008,0) -- cycle;
             \fill[black!\i!cyan!50!black] (0,3.5) -- (0.6814+\i*0.006814,-0.65) -- +(-0.008,0) -- cycle;
          }
          \draw[dotted] (A) arc[x radius=2, y radius=0.3, start angle=180-17.46, end angle=17.46];
       \end{scope}
 
    \end{tikzpicture}
    \quad
\raisebox{1cm}{.}
\]
\\[-1cm]

\un{Claim:} 
\\[-.5cm]

\quad 
\textbullet \quad 
The partial derivatives of $F$ exist almost everywhere.
\\[-.5cm]

\quad 
\textbullet \quad 
$F$ is nowhere totally differentiable.
\\

[$F$ is defined on almost every segment 
\quad
$
\begin{cases}
\ x = a \hspace{0.5cm} 0 \leq a \leq 1
\\
\quad \text{and} 
\\
\ y = b \hspace{0.5cm} 0 \leq b \leq 1
\end{cases}
, 
$
\\

\noindent
i.e.,  $F$ is defined on the intersection of $Q$ and almost every line parallel to a coordinate axis.   

\noindent
Let
\[
\tI_\ell \ = \ \Pi_x (B(r_\ell, \varepsilon_\ell)), 
\qquad
\text{($\Pi_x$ projection onto the $x$-axis).}
\]
For almost every $x \in [0,1]$, $x$ is in at most a finite number of the $\tI_\ell$.

\noindent
Therefore, $F(x,-)$ a sum of a finite number of the $F_\ell(x,-)$ for almost all $x \in [0,1]$
and the first part of the claim obtains.
\\[-.25cm]

The second part of the claim is immediate since $F$ is nowhere continuous.]
\\[-.25cm]

For an example of the more delicate situation cited above, where $F$ is constrained to be continuous, 
cf. Rademacher\footnote[2]{\vspace{.11 cm} \textit{Math. Ann.}, \textbf{79} (1919) 340-359.}

\chapter{
$\boldsymbol{\S}$\textbf{2.2}.\quad  VECTOR FUNCTIONS}
\setlength\parindent{2em}
\renewcommand{\thepage}{2-\S2-\arabic{page}}

\vspace{-.25cm}
\qquad
Let $\Omega$ be a nonempty open subset of $\R^n$ and let $f : \Omega \ra \R^m$ be a function. 
\\

\qquad{\bf 2.2.1:}\  
{\small\bf DEFINITION} \ 
$f$ is \un{differentiable} at a point $x_0 \in \Omega$ if there exists a linear function 
$T : \R^n \ra \R^m$ (depending on $x_0$) such that 

\[
\lim\limits_{h \ra 0} \ 
\frac{f(x_0 + h) - f(x_0) - T(h)}{\norm{h}}
\ = \ 
0.
\]
\\[-1cm]

$T$ is called the \un{differential} of $f$ at $x_0$ and is denoted by $\td f (x_0)$.
\\[-.5cm]

[Note: \ 
As in the scalar case, the differential is unique, if it exists at all.]
\\

\qquad{\bf 2.2.2:}\  
{\small\bf \un{N.B.}}\ 
$f$ is differentiable in $\Omega$ if $f$ is differentiable at every point of $\Omega$.
\\

Let $f^1 (x), \hsx f^2 (x), \hsx \ldots, \hsx f^m (x)$ be the components of $f$ 
and let $T^1, T^2, \ldots, T^m$ be the components  of $T$ $-$then 
the defining relation for the differential of $f$ at $x_0$ amounts to the relations
\\

\[
\begin{matrix*}[l]
&
\ds 
\lim\limits_{h \ra 0} \ \frac{f^1 (x_0 + h) - f^1 (x_0) - T^1 (h)}{\norm{h}} 
\ = \ 
0
\\[26pt]
&
\ds 
\lim\limits_{h \ra 0} \ \frac{f^2 (x_0 + h) - f^2 (x_0) - T^2 (h)}{\norm{h}} 
\ = \ 
0
\\[26pt]
&
\ds 
\quad \vdots
\\[26pt]
&
\ds 
\lim\limits_{h \ra 0} \ \frac{f^m (x_0 + h) - f^m (x_0) - T^m (h)}{\norm{h}} 
\ = \ 
0.
\end{matrix*}
\]
\\[-.75cm]

\noindent
Therefore $f$ is differentiable at $x_0$ iff all the components of $f$ are differentiable at $x_0$  
and when this is so, $f$ is continuous at $x_0$.
\\

\qquad{\bf 2.2.3:}\  
{\small\bf SCHOLIUM} \ 
For all $h = (h_1, h_2, \ldots, h_n)$, 

\[
\begin{matrix*}[l]
&
\ds 
T^1 (h) \ = \ \td f^1 (x_0) (h) \ = \ \frac{\partial f^1}{\partial h} (x_0)\ = \ \sum\limits_{i = 1}^n \ \frac{\partial f^1}{\partial x_i} (x_0) h_i
\\[26pt]
&
\ds 
T^2 (h) \ = \ \td f^2 (x_0) (h) \ = \ \frac{\partial f^2}{\partial h} (x_0)\ = \ \sum\limits_{i = 1}^n \ \frac{\partial f^2}{\partial x_i} (x_0) h_i
\\[26pt]
&
\ds
\quad \vdots
\\[26pt]
&
\ds
T^m (h) \ = \ \td f^m (x_0) (h) \ = \ \frac{\partial f^m}{\partial h} (x_0)\ = \ \sum\limits_{i = 1}^n \ \frac{\partial f^m}{\partial x_i} (x_0) h_i
\end{matrix*}
.
\]
\\[-.5cm]

\qquad{\bf 2.2.4:}\  
{\small\bf DEFINITION} \ 
The \un{Jacobian matrix} of $f$ at $x_0$ is the $m \times n$ matrix 
\\

\[
\begin{bmatrix*}[l]
&
\ds 
\frac{\partial f^1}{\partial x_1} (x_0) \ 
\frac{\partial f^1}{\partial x_2} (x_0)  \ 
\ \ldots \ \ 
\frac{\partial f^1}{\partial x_n} (x_0) \ 
\\[26pt]
&
\ds
\frac{\partial f^2}{\partial x_1} (x_0)  \
\frac{\partial f^2}{\partial x_2} (x_0)  \
\ \ldots \ \ 
\frac{\partial f^2}{\partial x_n} (x_0)  \
\\[26pt]
&
\ds
\quad \vdots
\\[26pt]
&
\ds
\frac{\partial f^m}{\partial x_1} (x_0)  \
\frac{\partial f^m}{\partial x_2} (x_0)  \
\ \ldots \ \ 
\frac{\partial f^m}{\partial x_n} (x_0)  \
\end{bmatrix*}
,
\]
denoted by $\tD f(x_0)$.
\\[-.5cm]

[Note: \ 
The partial derivatives of $f$ are the partial derivatives of its components, i.e., the

\[
\frac{\partial f^i}{\partial x_j} 
\ \equiv \ 
\tD_j f^i .]
\]
\\[-.5cm]


\qquad{\bf 2.2.5:}\  
{\small\bf DEFINITION} \ 
Suppose that $n = m$ $-$then the determinant of the Jacobian matrix $\tD f(x_0)$ is called the \un{Jacobian} of $f$ at $x_0$, 
denoted by

\[
\tJ_f (x_0) 
\quad \text{or} \quad 
\frac{\partial (f^1, f^2, \ldots, f^n)}{\partial (x_1, x_2, \ldots, x_n)}.
\]
\\[-.75cm]

\qquad{\bf 2.2.6:}\  
{\small\bf OPEN MAPPING THEOREM} \ 
Suppose that $n = m$ and suppose that $\tJ_f (x) \neq 0$ for all $x \in \Omega$ $-$then 
the image $f(U)$ of any open set $U \subset \Omega$ is open.
\\

\qquad{\bf 2.2.7:}\  
{\small\bf CHAIN RULE} \ 
Let $U \subset \R^n$ and $V \subset \R^m$ be nonempty open sets and let 
$f: U \ra \R^m$ and $g : V  \ra \R^p$ subject to $f(U) \subset V$.  
Assume: \ $f$ is differentiable at $x_0 \in U$ and $g$ is differentiable at $f(x_0)$ $-$then $g \circ f$ is differentiable at $x_0$ and 

\[
\td (g \circ f) (x_0) 
\ = \ 
\td g(f(x_0)) \circ \td f(x_0)
\]
or, in terms of Jacobian matrices, 

\[
\tD (g \circ f) (x_0) 
\ = \ 
\tD g (f(x_0)) \hsy \tD f (x_0).
\]
\\[-1cm]

\qquad{\bf 2.2.8:}\  
{\small\bf RAPPEL} \ 
The set $\Hom (\R^n, \R^m)$ of linear transformations from $\R^n$ to $\R^m$ is a vector space of dimension $n \hsy m$.  
Moreover, it is a Banach space under the norm 

\[
\norm{A} 
\ = \ 
\max \{\norm{A x} \hsy : \hsy \norm{x} \leq 1\}.
\]
And $\forall \ x$, 
\[
\norm{A x}
\ \leq \ 
\norm{A} \hsy \norm{x}.
\]
\\[-.5cm]

\qquad{\bf 2.2.9:}\  
{\small\bf EXAMPLE} \ 
Given $f : \Omega \ra \R^m$, 

\[
\td f(x_0) \in \Hom (\R^n, \R^m).
\]
\\[-1.5cm]

[Note: \ 
If $f = A \in \Hom (\R^n, \R^m)$, then $\td f(x_0) = A$.]
\\[-.5cm]

\qquad{\bf 2.2.10:}\  
{\small\bf DEFINITION} \ 
A differentiable function $f : \Omega \ra \R^n$ is \un{continuously} \un{differentiable} if 
\[
\td f : \Omega \ra \Hom (\R^n, \R^m) 
\]
is continuous.
\\[-.5cm]

[Spelled out, given $x_0 \in \Omega$ and $\varepsilon > 0$, there is a $\delta > 0$ such that 

\[
\norm{\td f (x) - \td f (x_0)} 
\ < \ 
\varepsilon
\]
if $\norm{x - x_0} < \delta$.]
\\

\qquad{\bf 2.2.11:}\  
{\small\bf NOTATION} \ 
$C^1(\Omega; \R^m)$ is the set of continuously differentiable functions from $\Omega$ to $\R^m$, 
often referred to as the $C^\prime$-functions (a vector space over $\R$).
\\

\qquad{\bf 2.2.12:}\  
{\small\bf THEOREM} \ 
$f : \Omega \ra \R^m$ is $C^\prime$ iff the partial derivatives of $f$ exist and are continuous throughout $\Omega$.
\\[-.5cm]

PROOF \ 
That the differentiability of $f$ implies the continuity of the partials can be seen by noting that 
\[
\abs{(\tD_j f^i) (y) - (\tD_j f^i) (x)}
\ \leq \ 
\norm{\td f(y) -\td f(x)}.
\]
In the other direction, take $m = 1$, fix $x_0 \in \Omega$, let $\varepsilon > 0$, and choose $r_0 > 0$: $B(x_0, r_0)^\circ \subset \Omega$ and 

\[
\abs{(\tD_j f) (x) - (\tD_j f) (x_0)}
\ < \ 
\frac{\varepsilon}{n}
\qquad (x \in B(x_0, r_0)^\circ, \ 1 \leq j \leq n).
\]
Write

\[
h 
\ = \ 
\sum\limits_{j = 1}^n \ 
h_j \hsy e_j, 
\qquad \norm{h} < r_0, 
\]
and put

\[
\begin{cases}
\ v_0 \ = \ 0,
\\[4pt] 
v_k \ = \ 
h_1 e_1 + \cdots + h_k e_k 
\quad (1 \leq k \leq n)
\end{cases}
.
\]
Then
\[
f(x_0 + h) - f(x_0) 
\ = \ 
\sum\limits_{j = 1}^n \ 
[f(x_0 + v_j) - f(x_0 + v_{j-1})].
\]
Since $\norm{v_k} < r_0$, the line segments with endpoints $x_0 + v_{j-1}$ and $x_0 + v_j$ lie in $B(x_0, r_0)^\circ$.  
Taking into account that
\[
v_j 
\ = \ 
v_{j-1} + h_j e_j, 
\]
the MVT implies that 
\[
f(x_0 + v_j) - f(x_0 + v_{j-1})
\ = \ 
h_j (\tD_j f) (x_0 + v_{j-1} + \theta_j h_j e_j)
\]
for some $\theta_j \in \hsx ]0,1[$.  
Next

\[
\abs{h_j (\tD_j f) (x_0 + v_{j-1} + \theta_j h_j e_j) 
-h_j (\tD_j f) (x_0)}
\ < \ 
\frac{\abs{h_j} \varepsilon}{n}.
\]
Consequently

\allowdisplaybreaks
\begin{align*}
\bigg | \hsx
f (x_0 + h) - f(x_0) 
&- 
\sum\limits_{j = 1}^n \ 
h_j (\tD_j f) (x_0)
\hsx \bigg | 
\
\\[15pt]
&=\ 
\bigg | \hsx
\sum\limits_{j = 1}^n \ 
\left[
f(x_0 + v_j) - f(x_0 + v_{j-1})
\right]
- 
\sum\limits_{j = 1}^n \ 
h_j (\tD_j f) (x_0)
\hsx \bigg | 
\\[15pt]
&=\ 
\bigg | \hsx
\sum\limits_{j = 1}^n \ 
h_j (\tD_j f) (x_0 + v_{j - 1} + \theta_j h_j e_j)
- 
\sum\limits_{j = 1}^n \ 
h_j (\tD_j f) (x_0)
\hsx \bigg | 
\\[15pt]
&\leq \ 
\sum\limits_{j = 1}^n \ 
\bigg | \hsx
h_j (\tD_j f) (x_0 + v_{j - 1} + \theta_j h_j e_j)
- 
h_j (\tD_j f) (x_0)
\hsx \bigg | 
\\[15pt]
&\leq \ 
\sum\limits_{j = 1}^n \ 
\frac{\abs{h_j} \varepsilon}{n}
\\[15pt]
&\leq \ 
\frac{1}{n}
\left(
\sum\limits_{j = 1}^n \ 
\abs{h_j} 
\right)
\varepsilon
\\[15pt]
&\leq \ 
\frac{1}{n} (\sqrt{n} \hsx \norm{h}) \hsx \varepsilon
\\[15pt]
&\leq \ 
\norm{h} \hsx \varepsilon.
\end{align*}

Therefore $f$ is differentiable at $x_0$: 
\begin{align*}
\td f(x_0) (h) \ 
&=\ 
\td f(x_0) 
\left(
\sum\limits_{j = 1}^n \ 
h_j e_j
\right)
\\[15pt]
&=\ 
\sum\limits_{j = 1}^n \ 
h_j (\tD_j f) (x_0).
\end{align*}

Since $m = 1$, the Jacobian matrix is a row:

\[
\tD f(x_0) 
\ = \ 
[
\tD_1 f(x_0), \tD_2 f(x_0), \ldots, \tD_n f(x_0) 
]
\]
or still, 

\[
\tD f(x_0) 
\ = \ 
\left[
\frac{\partial f}{\partial x_1} (x_0), \frac{\partial f}{\partial x_2} (x_0), \ldots, \frac{\partial f}{\partial x_n} (x_0)
\right].
\]

\noindent
Its entries are continuous functions of $x_0$, thus $f$ is a $C^\prime$-function.
\\

\qquad{\bf 2.2.13:}\  
{\small\bf DEFINITION} \ 
Suppose that $n = m$ and suppose that $f: \Omega \ra \R^n$ is a $C^\prime$-function $-$then 
a point $x_0 \in \Omega$ is a \un{critical point} for $f$ if the rank of $\tD f(x_0)$ is not maximal, i.e., 
if the rank of $\tD f(x_0)$ is $< n$ or still, if $\tJ_f (x_0) = 0$.
\\

\qquad{\bf 2.2.14:}\  
{\small\bf NOTATION} \ 
Write $Z_f$ for the set of critical points of $f$.
\\

\qquad{\bf 2.2.15:}\  
{\small\bf SARD} \ 
$f(Z_f)$ is a set of Lebesgue measure 0.
\\

There are numerous variants on this theme which need not be considered at this juncture.  
However: \ 
\\

\qquad{\bf 2.2.16:}\  
{\small\bf LEMMA} \ 
Under the above assumptions, for any Lebesgue measurable set $E \subset \Omega$, the set $f(E)$ is Lebesgue measurable and 

\[
\Lm^n (f(E)) 
\ \leq \ 
\int\limits_E \ \abs{\tJ_f} \ \td \Lm^n 
\qquad \text{(cf. 12.3.3).}
\]
\\[-.75cm]

\qquad{\bf 2.2.17:}\  
{\small\bf \un{N.B.}} \ 
SARD is an immediate consequence of this result.
\\

The mean value theorem does not hold in general for a vector valued function $f: \Omega \ra \R^m$ $(m > 1)$ 
(but it does hold if the number of auxiliary points is increased (details omitted)).  
However: 
\\

\qquad{\bf 2.2.18:}\  
{\small\bf THEOREM} \ 
Suppose that $f:[a,b] \ra \R^k$ is continuous and that its restriction to $]a, b[$ is differentiable $-$then 
there exists an $x \in ]a,b[$ such that 

\[
\norm{f(b) - f(a)}
\ \leq \ 
(b - a) \hsy \norm{f^\prime (x)}.
\]

PROOF \ 
Let
\[
\phi(t) 
\ = \ 
\langle f(b) - f(a), f(t) \rangle 
\qquad (a \leq t \leq b).
\]
Then $\phi$ satisfies the assumptions of the MVT, hence
\begin{align*}
\phi(b) - \phi(a) \ 
&=\ 
(b - a) \phi^\prime (x) 
\\[8pt]
&=\ 
(b - a) 
\hsx 
\langle f(b) - f(a), f^\prime(x) \rangle 
\end{align*}
for some $x \in ]a,b[$.  
On the other hand, 
\begin{align*}
\phi(b) - \phi(a) \ 
&=\ 
\langle f(b) - f(a), f(b) \rangle 
- 
\langle f(b) - f(a), f(a) \rangle 
\\[8pt]
&=\ 
\langle f(b) - f(a), f(b) -  f(a) \rangle
\\[8pt]
&=\ 
\norm{f(b) - f(a)}^2.
\end{align*}
Then
\begin{align*}
\norm{f(b) - f(a)}^2 \ 
&=\ 
(b - a) 
\hsx 
\langle f(b) - f(a), f^\prime(x) \rangle 
\\[8pt]
&\leq \ 
(b - a) 
\hsx 
\norm{f(b) - f(a)}
\norm{f^\prime(x)},
\end{align*}
hence
\[
\norm{f(b) - f(a)}
\ \leq \ 
(b - a) 
\hsx 
\norm{f^\prime(x)}.
\]

\chapter{
$\boldsymbol{\S}$\textbf{2.3}.\quad  LIPSCHITZ FUNCTIONS}
\setlength\parindent{2em}
\renewcommand{\thepage}{2-\S3-\arabic{page}}

\qquad
Let $E$ be a nonempty subset of $\R^n$.
\\

\qquad{\bf 2.3.1:}\  
{\small\bf DEFINITION} \ 
A function $f : E \ra \R^m$ is said to be \un{$L$-Lipschitz} ($L \geq 0$) if for all $x$, $y \in E$,
\[
\norm{f(x) - f(y)}
\ \leq \ 
L \hsy \norm{x - y}.
\]
\\[-1cm]

\qquad{\bf 2.3.2:}\  
{\small\bf EXAMPLE} \ 
A constant function $x \ra C \hsx (\in \R^m)$ is 0-Lipschitz.
\\

\qquad{\bf 2.3.3:}\  
{\small\bf EXAMPLE} \ 
$\norm{\hsx \cdot \hsx } : \R^n \ra \R$ is 1-Lipschitz.
\\[-.5cm]

[In fact,
\[
\abs{\hsx \norm{x} - \norm{y}\hsy}
\ \leq \ 
\norm{x - y}.]
\]
\\[-1cm]

\qquad{\bf 2.3.4:}\  
{\small\bf LEMMA} \ 
Let $\{f_i \hsy : \hsy i \in I\}$ be a collection of $L$-Lipschitz functions $f_i : E \ra \R$ $-$then the functions

\[
\begin{cases}
\ x \ra \sup\limits_{i \in I} \ f_i(x) \equiv F(x)
\\[11pt]
\ x \ra \inf\limits_{i \in I} \ f_i(x) \equiv f(x)
\end{cases}
\qquad (x \in E)
\]
are $L$-Lipschitz if finite at one point.
\\[-.25cm]

PROOF \ 
To establish the first assertion, note that for all $x$, $y \in E$, 
\[
f_i(y) 
\ \leq \ 
f_i(x) +  L \hsy \norm{x - y}.
\]
Take now the supremum on the RHS and then on the LHS to get: 
\[
F(y) 
\ \leq \ 
F(x) + L \hsy \norm{(x - y)}.
\]
If $F(x) < +\infty$, then $F(y) < +\infty$ for all $y \in E$, hence 
$F(y) - F(x) < L \hsy \norm{(x - y)}$, 
hence
$F(x) - F(y) \leq L \hsy \norm{(x - y)}$, hence $\abs{F(x) - F(y)} \leq L \hsy \norm{(x - y)}$.
\\

\qquad{\bf 2.3.5:}\  
{\small\bf APPLICATION} \ 
The function from $\R^n$ to $\R$ defined by the rule
\[
y \ra \dist (y, E) 
\ = \ \inf \{\norm{x - y} \hsy : \hsy x \in E\}
\]
is 1-Lipschitz.
\\

\qquad{\bf 2.3.6:}\  
{\small\bf THEOREM} \ 
If $f: E \ra \R$ is $L$-Lipschitz, then there is an $L$-Lipschitz function $F : \R^n \ra \R$ such that $\restr{F}{E} = f$.
\\[-.5cm]

[Consider
\[
F(y) 
\ = \ 
\inf\limits_{x \in E} \ \{f(x) + L \hsy \norm{x - y}\}.]
\]
\\[-.75cm]

\qquad{\bf 2.3.7:}\  
{\small\bf NOTATION} \ 
Given $f: E \ra \R^m$, put
\[
\Lip (f; E)
\ = \ 
\sup\limits_{\substack{x, y \in E\\ x \neq y}} \ \frac{\norm{f(x) - f(y)}}{\norm{x - y}} 
\quad (\equiv \inf \{L\}).
\]
\\[-1cm]

[Note: \ 
Omit the ``E'' if it is $\R^n \hsy : \hsy \Lip (f)$.]
\\

\qquad{\bf 2.3.8:}\  
{\small\bf DEFINITION} \ 
A function $f: E \ra \R^m$ is \un{Lipschitz} if it is $L$-Lipschitz for some $L \geq 0$. 
\\

\qquad{\bf 2.3.9:}\  
{\small\bf \un{N.B.}} \ 
If $f: E \ra \R$ is Lipschitz, then $f$ is uniformly continuous.
\\[-.5cm]

[Conversely, it can be shown that if $f: E \ra \R$ is bounded and uniformly continuous, then $f$ is the uniform limit of a sequence of 
Lipschitz functions.]
\\[-.5cm]

[Note: \ 
The function $f (x) = \sqrt{x}$ $(0 \leq x \leq 1)$ is not Lipschitz but it is uniformly continuous.]
\\

\qquad{\bf 2.3.10:}\  
{\small\bf THEOREM} \ 
Let $f : \R^n \ra \R^m$ be Lipschitz $-$then for any nonempty $E \subset \R^n$, 
\[
\sH^\ts (f(E)) 
\ \leq \ 
(\Lip (f))^\ts \hsx \sH^\ts(E)
\qquad (s \in [0, +\infty[).
\]
\\[-1cm]

PROOF \ 
Fix $\delta > 0$ and choose sets $\{E_k\} \subset \R^n$ such that
\[
E \subset 
\bigcup\limits_{k = 1} ^\infty \ E_k, 
\quad 
\diam (E_k) \leq \delta.
\]
Then
\[
\diam (f(E_k)) \leq \Lip (f) \hsx \diam (E_k).
\]
Since
\[
f(E) \subset 
\bigcup\limits_{k = 1} ^\infty \ f(E_k), 
\quad 
\diam (f(E_k)) \leq \Lip (f) \hsx \delta 
\]
it therefore follows that 

\begin{align*}
\sH_{\Lip (f) \hsx \delta}^\ts \hsx (f(E)) \ 
&\leq \ 
\frac{\omega_\ts}{2^\ts} \ 
\sum\limits_{k = 1}^\infty \ (\diam (f(E_k)))^\ts
\\[18pt]
&\leq \ 
\frac{\omega_\ts}{2^\ts} \hsx
\Lip (f)^\ts \ 
\sum\limits_{k = 1}^\infty \ (\diam (f(E_k)))^\ts.
\end{align*}
Now take the infimum over this data to arrive at

\[
\sH_{\Lip (f) \hsx \delta}^\ts \hsx (f(E)) 
\ \leq \ 
\Lip (f)^\ts \hsx \sH_\delta^\ts (E),
\]
from which the assertion upon sending $\delta \downarrow 0$.
\\

\qquad{\bf 2.3.11:}\  
{\small\bf EXAMPLE} \ 
If $n > m$ and if $P: \R^n \ra \R^m$ is the usual projection, then for all $E \subset \R^n$, 
\[
\sH^s (P(E)) 
\ \leq \ 
\sH^s (E).
\]
\\[-1cm]

[In fact $\Lip (P) = 1$.]
\\

\qquad{\bf 2.3.12:}\  
{\small\bf SUBLEMMA} \ 
Let
\[
\begin{cases}
\ \{x_1, \ldots, x_k\} \subset \R^n
\\[8pt]
\ \{y_1, \ldots, y_k\} \subset \R^m
\end{cases}
\]
\\[-1cm]

\noindent
subject to the condition that
\[
\norm{y_i - y_j} 
\ \leq \ 
\norm{x_i - x_j} 
\]
for all $i, \hsx j \in \{1, \ldots, k\}$.  
Suppose that $r_1, \ldots, r_k$ are positive numbers such that 
\[
\bigcap\limits_{i = 1}^k \ B(x_i, r_i) 
\ \neq \ 
\emptyset.
\]
Then
\[
\bigcap\limits_{j = 1}^k \ B(y_j, r_j) 
\ \neq \ 
\emptyset.
\]
\\[-.5cm]

\qquad{\bf 2.3.13:}\  
{\small\bf LEMMA} \ 
Let $E \subset \R^n$ be a nonempty finite set and let $f : E \ra \R^m$ be a 1-Lipschitz function $-$then 
for any $x \in \R^n$, there is an extension of $f$ to a 1-Lipschitz function on $E \cup \{x\}$.  
\\[-.5cm]

PROOF \ 
Let $E = \{x_1, \ldots, x_k\}$ and assume that $\forall \ i, x \neq x_i$.  
Put $r_i = \norm{x - x_i} > 0$ and let $y_i = f(x_i)$ $-$then there exists a point $y \in \R^m$ such that
\[
\norm{y - f(x_i)}
\ \leq \ 
\norm{x - x_i}
\]
for each $i$, so it remains only to let $f(x) = y$.
\\

\qquad{\bf 2.3.14:}\  
{\small\bf EXTENSION PRINCIPLE} \ 
Let $f : E \ra \R^m$ be an $L$-Lipschitz function $-$then there exists an $L$-Lipschitz function
$F : \R^n \ra \R^m$ such that $\restr{F}{E} = f$.  
\\[-.5cm]

PROOF \ 
Upon dividing $f$ by $L$, it can be assumed that $f$ is 1-Lipschitz and it will be enough to deal explicitly with the situation when 
$E$ and $\R^n \backslash E$ are both infinite.  
Accordingly, choose 
a countable dense set $\{x_1, x_2, \ldots\}$ in $E$ and 
a countable dense set $\{y_1, y_2, \ldots\}$ in $\R^n \backslash E$.  
This done, for each $k = 1, 2, \ldots,$ use the previous lemma repeatedly to obtain a 1-Lipschitz function
\[
f_k \hsy : \hsy \{x_1, \ldots, x_k, y_1, \ldots, y_k\} \ra \R^m
\]
such that $f_k(x_i) = f(x_i)$ $(i = 1, \ldots, k)$.  
Claim: \ 
The sequence $\{f_k (y_1)\} \subset \R^m$ is bounded.  
Proof: \ 
\[
f_k (y_1) 
\ = \ 
f_k(y_1) - f_k(x_1) + f_k(x_1) 
\]

$\implies$
\begin{align*}
f_k (y_1) \ 
&\leq \ 
\norm{f_k(y_1) - f_k(x_1)} + \norm{f_k(x_1)}
\\[11pt]
&\leq \ 
\norm{y_1 - x_1} + \norm{f_k(x_1)}
\\[11pt]
&< \ 
+\infty
\end{align*}
independently of $k$.  
Proceeding, extract a convergent subsequence, say $\{f_{k_j^1}(y_1)\}$, 
and then extract from it yet another convergent subsequence $\{f_{k_j^2}(y_2)\}$.  
ETC.  
Pass to the diagonal sequence $\{g_j\} \hsy : \hsy g_j = f_{k_j^j}$, hence for every
\[
c \in C = \{x_1, x_2, \ldots \} \cup \{y_1, y_2, \ldots\},
\]
there follows
\[
g(c) 
\ = \ 
\lim\limits_{j \ra \infty} \ g_j(c) \in \R^m.
\]
In addition, $g: C \ra \R^m$ is 1-Lipschitz and $g(x_i) = f(x_i)$ $(i = 1,2, \ldots)$.  
And finally, in view of the density of $C$ in $\R^m$ and the density of $\{x_1, x_2, \ldots\}$ in $E$, $g$ extends to a 1-Lipschitz 
function $F: \R^n \ra \R^m$ such that $\restr{F}{E} = f$.
\\

\qquad{\bf 2.3.15:}\  
{\small\bf THEOREM} \ 
Suppose that $f : \Omega \ra \R^m$ is differentiable (thus, by definition, $\Omega$ is open).  
Assume: \ $\Omega$ is convex and that there is an $L \geq 0$ such that 
\[
\norm {\td f(x)} 
\ \leq \ 
L
\]
for all $x \in \Omega$ $-$then $f$ is $L$-Lipschitz.
\\[-.5cm]

PROOF \ 
Given $x$, $y \in \Omega$, the convexity of $\Omega$ implies that 
\[
tx + (1 - t) y \in \Omega 
\quad (0 \leq t \leq 1).
\]
Let
\[
g(t) 
\ = \ 
f (tx + (1 - t) \hsy y ) 
\quad (0 \leq t \leq 1).
\]
Then
\[
\frac{\td}{\td t} \hsx g(t) 
\ = \ 
\td f (tx + (1 - t) \hsy y ) \hsx (x - y)
\]

$\implies$
\begin{align*}
\norm{\frac{\td}{\td t} \hsx g(t)} \ 
&\leq \ 
\norm{\td f (tx + (1 - t) \hsy y )}\hsx \norm{(x - y)}
\\[11pt]
&\leq \ 
L \hsy \norm{(x - y)}.
\end{align*}
Take now in 2.2.13, $[a,b] = [0,1]$ and apply it to $g$, 
thus for some $t_0 \in \hsx ]0, 1[$, 
\begin{align*}
\norm{g(1)  - g(0)} \ 
&\leq \ 
(1 - 0) \hsx \norm{g^\prime (t_0)}
\\[11pt]
&=\ 
\norm{g^\prime (t_0)}
\\[11pt]
&\leq\ 
L \hsy \norm{(x - y)}.
\end{align*}
But
\[
g(1) \ = \ f(x), 
\quad 
g(0) \ = \ f(y)
\]

$\implies$
\[
\norm{f(x) - f(y)}
\ \leq \ 
L \hsy \norm{(x - y)}.
\]
\\[-.75cm]

\qquad{\bf 2.3.16:}\  
{\small\bf EXAMPLE} \ 
The sine function is 1-Lipschitz (since its derivative is the cosine which is bounded by 1).
\\

\qquad{\bf 2.3.17:}\  
{\small\bf DEFINITION} \ 
A function $f : \Omega \ra \R^m$  is said to be \un{locally Lipschitz} if  for each compact set $K \subset \Omega$, 
there exists a constant $C_K \geq 0$ such that for all $x$, $y \in K$, 

\[
\norm{f(x) - f(y)}
\ \leq \ 
C_K \hsy \norm{(x - y)}.
\]
\\[-1cm]


[Note: \ 
If there exists $L$ such that $C_K = L$ for all $K$, then $f$ is said to be \un{locally $L$-Lipschitz}.]
\\

\qquad{\bf 2.3.18:}\  
{\small\bf EXAMPLE} \ 
In $\R^2$, let
\[
E \ = \ 
\{(r, \theta) \hsy : \hsy 0 < r < +\infty, -\pi < \theta < \pi\}.
\]
Then the function $E \ra \R^2$ given by
\[
(r, \theta) \ra (r, \theta/2)
\]
is locally 1-Lipschitz but not Lipschitz.
\\

\qquad{\bf 2.3.19:}\  
{\small\bf THEOREM} \ 
If $\Omega \subset \R^n$ is convex and if $f : \Omega \ra \R$ is convex, then $f$ is locally Lipschitz.
\\

\qquad{\bf 2.3.20:}\  
{\small\bf RAPPEL} \ 
Suppose that $E \subset \R^n$ is Lebesgue measurable $-$then there exists an increasing sequence $\{F_k\}$ of closed sets $F_k$ 
contained in $E$ and a set $N$ of Lebesgue measure 0 such that
\[
E 
\ = \ 
\bigg(
\bigcup\limits_k \ F_k
\bigg)
\hsx \cup \hsx 
N.
\]
\\[-1cm]

\qquad{\bf 2.3.21:}\  
{\small\bf \un{N.B.}}  \ 
A closed set is the union of a countable family of compact sets and a continuous function $f: \R^n \ra \R^n$ sends a countable union of 
compact sets to a countable union of 
compact sets.
\\

\qquad{\bf 2.3.22:}\  
{\small\bf DEFINITION}\ 
Let $\Omega \subset \R^n$ be nonempty and open $-$then a continuous function 
$f: \Omega \ra \R^n$ is said to have 
\un{property (N)} 
if $f$ sends Lebesgue sets of measure 0 to Lebesgue sets of measure 0.
\\

If $f: \R^n \ra \R^n$  is continuous and if $E \subset \R^n$ is closed, then $f(E)$ is Lebesgue measurable.  
Consequently, in the presence of property (N), it follows that $f$ sends 
Lebesgue measurable sets to Lebesgue measurable sets.
\\

\qquad{\bf 2.3.23:}\  
{\small\bf THEOREM} \ 
If $f: \R^n \ra \R^n$ is locally Lipschitz, then $f$ has property (N).
\\[-.5cm]

PROOF \ 
The claim is that 
\[
\Lm^n (N)  \ = \ 0 
\ \implies \  
\Lm^{n^*} (f (N)) \ = \  0.
\]
To this end, fix a closed cube $K$ in $\R^n$, write
\[
\norm{f(x) - f(y)} 
\ \leq \ 
C_K \hsy \norm{x - y} 
\quad (x, y \in K),
\]
and note that a cube $I$ of side $r$ in $K$ has diameter of length $r \sqrt{n}$.  
Since $f$ is Lipschitz, $f(I)$ has diameter at most $r \sqrt{n} \hsx C_K$, 
thus is contained in a cube of side length $r \sqrt{n} \hsx C_K$, and so
\[
\Lm^{n^*} (f (I))
\ \leq \ 
n^{n/2} \hsy C_K^n \hsy r^n 
\ = \ 
n^{n/2} \hsy C_K^n \hsy \Lm^{n^*} (I),
\]
or still, 
\[
\Lm^{n^*} (f (N \cap K))
\ \leq \ 
n^{n/2} \hsy C_K^n \hsy \Lm^{n^*} (I \cap K)
\]

$\implies$
\[
\Lm^{n^*} (f (N \cap K))
\ = \ 
0.
\]
Choose now an increasing sequence $\{K_j\}$ of closed cubes $K_j$ such that $\R^n = \bigcup\limits_j \ K_j$, hence

\[
f(N) 
\ = \ 
\bigcup\limits_j \ f(N \cap K_j), 
\]
and therefore

\[
\Lm^n (N) 
\ \leq \ 
\sum\limits_j \ 
\Lm^{n^*} (f (N \cap K_j))
\ = \ 
0.
\]
\\[-.75cm]

\qquad{\bf 2.3.24:}\  
{\small\bf LEMMA} \ 
Suppose that $f$ is a $C^\prime$-function, i.e., $f \in C^1 (\Omega, \R^m)$ $-$then $f$ is locally Lipschitz.

\chapter{
$\boldsymbol{\S}$\textbf{2.4}.\quad  RADEMACHER}
\setlength\parindent{2em}
\renewcommand{\thepage}{2-\S4-\arabic{page}}

\qquad
If $\Omega$ is a nonempty open subset of $\R$ and if $f : \Omega \ra \R$ is $L$-Lipschitz, then $f$ is absolutely continuous, 
hence is differentiable almost everywhere (per $\Lp^1$).
\\

\qquad{\bf 2.4.1:}\  
{\small\bf THEOREM} \ 
If $\Omega$ is a nonempty open subset of $\R^n$ and if $f : \Omega \ra \R^m$ is $L$-Lipschitz, then $f$ is differentiable at 
$\Lm^n$ almost all points in $\Omega$.
\\[-.25cm]

The proof will be given in the lines below.
\\[-.25cm]

\un{First Step:} \ 
It can be assumed that $m = 1$.
\\[-.5cm]

[For $f$ is Lipschitz (or differentiable) iff every component of $f$ is Lipschitz (or differentiable).]
\\[-.25cm]

\un{Second Step:} \ 
It can be assumed that $\Omega = \R^n$, so $f : \R^n \ra \R$.  
\\[-.5cm]

[Invoke the Extension Principle.]
\\

\qquad{\bf 2.4.2:}\  
{\small\bf RAPPEL} \ 
A Lebesgue measurable function $f : \R^n \ra \R$ is \un{locally integrable} if

\[
\int\limits_K \ 
\abs{f} 
\ \td \Lm^n 
\ < \ 
+\infty
\]
for every compact $K \subset \R^n$.
\\[-.25cm]

Denote the space of such by
\[
\Lp_\locx^1 (\R^n).
\]
\\[-1.25cm]

\qquad{\bf 2.4.3:}\  
{\small\bf LEMMA} \ 
If $f \in \Lp_\locx^1 (\R^n)$ and if 

\[
\int\limits_{\R^n} \ 
f \hsy \phi 
\ \td \Lm^n 
\ = \ 
0
\]
for all $\phi \in C_c^\infty (\R^n)$, then $f = 0$ almost everywhere.
\\

\un{Third Step:} \ 
Given $x \in \R^n$, $v \in \bS^{n-1}$, form
\[
f_{x, v} (t) 
\ = \ 
f(x + t v) 
\quad (t \in \R).
\]
Then $f_{x, v}$ as a function of $t$ is Lipschitz, hence is differentiable almost everywhere (per $\Lm^1$).
\\

\un{Fourth Step:} \ 
Recall that by definition, 
\[
\frac{\partial f}{\partial v} (x)
\ = \ 
\lim\limits_{t \ra 0} \ 
\frac{f(x + t v) - f(x)}{t}
\]
whenever the limit exists.
\\

\un{Fifth Step:} \ 
Let 
\[
E_v 
\ = \ 
\big\{x \in \R^n \hsy : \hsy \frac{\partial f}{\partial v} (x) \ \text{exists}\big\}.
\]
Then $E_v$ is Borel and the function
\[
\begin{cases}
\ E_v \ra \R
\\[8pt]
\ \ds x \ra \frac{\partial f}{\partial v} (x)
\end{cases}
\]
is Lebesgue measurable.
\\

\un{Sixth Step:} \ 
Write
\[
\R^n 
\ = \ 
\R v \hsx \oplus \hsx v^\perp.
\]

\noindent
Then
\allowdisplaybreaks
\begin{align*}
\Lm^n (\R^n \backslash E_v) \ 
&=\ 
\int\limits_{\R^n} \ 
\chisubRnbackslashEsubv 
\ \td \Lm^n
\\[15pt]
&=\ 
\int\limits_{v^\perp} \ 
\int\limits_{\R v}\ 
\chisubRnbackslashEsubv (tv + w) \ 
\ \td t \hsx \td w
\\[15pt]
&=\ 
\int\limits_{v^\perp} \ 
\Lm^1 (S_w) \ \td w
\\[15pt]
&=\ 
0,
\end{align*}

where
\[
S_w 
\ = \ 
\{t \in \R \hsy : \hsy t v + w \in \R^n \backslash E_v\}.
\]
\\[-0.75cm]

\un{Seventh Step:} \ 
Therefore $E_v$ is of full measure in that $\ds\frac{\partial f}{\partial v} (x)$ exists for almost every $x \in \R^n$ per $\Lm^n$.
\\

\un{Eighth Step:} \ 
In particular, the partial derivatives $\ds\frac{\partial f}{\partial x_i} (x)$ $(i = 1, 2, \ldots, n)$ 
exist for almost all $x$ per $\Lm^n$, hence the same is so of the formal gradient
\[
\nabla f (x) 
\ = \ 
\left(
\frac{\partial f}{\partial x_1} (x), \frac{\partial f}{\partial x_2} (x), \ldots, \frac{\partial f}{\partial x_n} (x)).
\right)
\]
\\[-0.75cm]

\qquad{\bf 2.4.4:}\  
{\small\bf LEMMA} \ 
For each $v \in \bS^{n-1}$, 
\[
\frac{\partial f}{\partial v} (x) 
\ = \ 
\langle v, \nabla f(x) \rangle
\]
almost everywhere (per $\Lp^1$).
\\[-.5cm]

PROOF \ 
Both functions are in $\Lp_\locx^1 (\R^n)$.  
E.g.: 

\allowdisplaybreaks
\begin{align*}
\frac{\abs{f(x + t v) - f(x)}}{\abs{t}} \ 
&\leq \ 
L \hsx \frac{\norm{x + t v - x}}{\abs{t}}
\\[11pt]
&=\ 
L \hsx \frac{\norm{t v}}{\abs{t}} 
\\[11pt]
&=\ 
L \hsx \norm{v}.
\end{align*}
Bearing in mind 2.4.3, it suffices to show that
\[
\int\limits_{\R^n} \ 
\frac{\partial f}{\partial v} (x)  \hsx \phi(x) 
\ \td \Lm^n (x) 
\ = \ 
\int\limits_{\R^n} \ 
\langle v, \nabla f (x) \rangle
\ \td \Lm^n (x) 
\]
for all $\phi \in C_c^\infty (\R^n)$.  
Start with the left hand side and proceed: 
\allowdisplaybreaks
\begin{align*}
\int\limits_{\R^n} \ 
\frac{\partial f}{\partial v} (x)  \hsx \phi(x) 
\ \td \Lm^n (x) \ 
&=\ 
\int\limits_{\R^n} \ 
\lim\limits_{t \ra 0} \ 
\frac{f(x + t v) - f(x)}{t} \hsx \phi(x) 
\ \td \Lm^n (x) \ 
\\[15pt]
&=\ 
\lim\limits_{t \ra 0} \
\int\limits_{\R^n} \ 
\frac{f(x + t v) - f(x)}{t} \hsx \phi(x) 
\ \td \Lm^n (x) \ 
\\[15pt]
&=\ 
\lim\limits_{t \ra 0} \
\int\limits_{\R^n} \ 
-f(x) \hsx
\frac{\phi (x) - \phi (x - tv)}{t} 
\ \td \Lm^n (x) \ 
\\[15pt]
&=\ 
-
\int\limits_{\R^n} \ 
f(x) \hsx 
\lim\limits_{t \ra 0} \ 
\frac{\phi (x) - \phi (x - tv)}{t} 
\ \td \Lm^n (x) \ 
\\[15pt]
&=\ 
-
\int\limits_{\R^n} \ 
f(x) \hsx \frac{\partial \phi}{\partial v} (x)
\ \td \Lm^n (x)
\\[15pt]
&=\ 
-\sum\limits_{i = 1}^n \ 
v_i \ 
\int\limits_{\R^n} \ 
f(x) \hsx \frac{\partial \phi}{\partial x_i} (x) 
\ \td \Lm^n (x)
\\[15pt]
&=\ 
\sum\limits_{i = 1}^n \ 
v_i \ 
\int\limits_{\R^n} \ 
\frac{\partial f}{\partial x_i} (x)  \hsx \phi(x) 
\ \td \Lm^n (x) \ 
\\[15pt]
&=\ 
\int\limits_{\R^n} \ 
\langle v, \nabla f (x) \rangle
\hsx \phi(x) 
\ \td \Lm^n (x).
\end{align*}

[Note: \ 
The justification of the formalities is left to the reader.]
\\

\un{Ninth Step:} \ 
Let $D \subset \bS^{n-1}$ be a countable dense set $-$then there is a Lebesgue measurable set $E \subset \R^n$ 
such that $\Lm^n (\R^n \backslash E) = 0$ and $\forall \ v \in D$, 
\[
\frac{\partial f}{\partial v} (x) 
\ = \ 
\langle v, \nabla f (x) \rangle
\qquad (x \in E).
\]

\un{Tenth Step:} \ 
Fix $x_0 \in E$ $-$then the claim is that $f$ is differentiable at $x_0$: 

\[
\lim\limits_{h \ra 0} \ 
\frac{f(x_0 + h) - f(x_0) - \langle h, \nabla f (x_0) \rangle}{\norm{h}}
\ = \ 
0, 
\]
the ambient linear function $T : \R^n \ra \R$ being the arrow
\[
h \ra \langle h, \nabla f (x_0) \rangle.
\]
\\[-1.25cm]

\qquad{\bf 2.45.:}\  
{\small\bf LEMMA} \ 
$f$ is differentiable at $x_0$ if $\forall \ \varepsilon > 0$, $\exists \ \delta > 0$ such that
\[
\abs{f(x_0 + t v) - f(x_0) - t \langle v, \nabla f(x_0) \rangle}
\ 
\ \leq \ 
\varepsilon \abs{t}
\]
provided $\abs{t} \leq \delta$ and $v \in \bS^{n-1}$.
\\[-.25cm]

To verify that this condition is satisfied, fix $\varepsilon > 0$ and choose a finite set $D_0 \subset D$ with the property that for every 
$v \in \bS^{n-1}$ there is a $v_0 \in D_0$ such that $\norm{v - v_0} \leq \varepsilon$.  
Since directional derivatives indexed by the $v_0 \in D_0$ are finite in number, there is a $\delta > 0$ such that $\forall \ v_0 \in D_0$, 
\[
\abs{f(x_0 + t v_0) - f(x_0) - t \langle v_0, \nabla f(x_0) \rangle }
\ 
\ \leq \ 
\varepsilon \abs{t}
\]
if $\abs{t} \leq \delta$.  
Given now $v \in \bS^{n-1}$, determine $v_0 \in D_0$ for which $\norm{v - v_0} \leq \varepsilon$ $-$then
\begin{align*}
|f(x_0 + t v) - &f(x_0) - t \langle v, \nabla f(x_0) \rangle| \ 
\\[11pt]
&\leq \ 
\varepsilon \abs{t} + 
\abs{f(x_0 + t v) - f(x_0 + t v_0)} 
+
\abs{t} \hsy \abs{\langle v - v_0, \nabla f(x_0)\rangle}
\\[11pt]
&\leq \ 
(1 + \Lip (f) + \norm{\nabla f(x_0)}) \hsx  \varepsilon \abs{t}
\end{align*}
for all  $\abs{t} \leq \delta$.

\chapter{
$\boldsymbol{\S}$\textbf{2.5}.\quad  STEPANOFF}
\setlength\parindent{2em}
\renewcommand{\thepage}{2-\S5-\arabic{page}}

\qquad 
Let $\Omega$ be a nonempty open subset of $\R^n$ and let $f : \Omega \ra \R^m$ be a Lebesgue measurable function
\\

\qquad{\bf 2.5.1:}\  
{\small\bf DEFINITION} \ 
The  \un{pointwise Lipschitz constant} of $f$ is 
\[
\Lip f(x) 
\ = \ 
\limsup\limits_{\substack{y \ra x \\y \in \Omega}} \ 
\frac{\norm{f(x) - f(y)}}{\norm{x - y}} 
\quad (x \in \Omega).
\]
\\[-.75cm]

\qquad{\bf 2.5.2:}\  
{\small\bf THEOREM} \ 
$f$ is differentiable almost everywhere in the set
\[
L_f 
\ = \ 
\{x \in \Omega \hsy : \hsy \Lip f(x) < +\infty\}.
\]
\\[-.75cm]

\qquad{\bf 2.5.3:}\  
{\small\bf REMARK} \ 

\[
L_f 
\ = \ 
\bigcup\limits_{k, \ell} \ E_{k, \ell}, 
\]
where
\[
E_{k, \ell}
\ = \ 
\left\{
x \in L_f \hsy : \hsy \norm{f(x)} \leq k 
\quad \text{and} \quad 
\frac{\norm{f(x) - f(y)}}{\norm{x - y}} 
\leq k 
\ \text{if } \ 
\norm{x - y} \leq \frac{1}{\ell}
\right\}
.
\]
\\[-1cm]

\noindent
Moreover $\restr{f}{E_{k, \ell}}$ is Lipschitz:
\\

\qquad \textbullet \quad 
$\ds \norm{x - y} \hsx < \hsx \frac{1}{\ell}
\implies 
\norm{f(x) - f(y)} \hsx\leq \hsx k \hsy \norm{x - y}$
\\

\qquad \textbullet \quad 
$\ds \norm{x - y} \hsx \geq \hsx \frac{1}{\ell}
\implies 
\norm{f(x) - f(y)} 
\hsx \leq  \hsx
2 \hsy k
\hsx \leq  \hsx 
2 \hsy k \hsy \ell \hsy \norm{x - y},$
\\

and it turns out that $f$ is differentiable almost everywhere in each $E_{k, \ell}$ (details omitted).
\\[-.5cm]

[Note: \ 
The $E_{k, \ell}$ are Lebesgue measurable, hence the same is true of $L_f$.]
\\

\qquad{\bf 2.5.4:}\  
{\small\bf SUBLEMMA} \ 
Let $g, \hsx f, \hsx h$ be functions from $\Omega$ to $\R$.  
Suppose that $g \leq f \leq h$, $g(x_0) = f(x_0) = h(x_0)$, and $g, \hsx h$ are differentiable at $x_0$ $-$then $f$ is differentiable at $x_0$.
\\[-.5cm]


PROOF \ 
Since $h - g \geq 0$ and $(h - g)(x_0) = 0$, it follows that $\td (h - g)(x_0) = 0$, hence $\td h(x_0) = \td g(x_0)$, call it $T$ $-$then 

\begin{align*}
\frac{\norm{g(x) - g(x_0) - \td g(x_0) (x - x_0)}}{\norm{x - x_0}} \ 
&\leq \ 
\frac{\norm{f(x) - f(x_0) - T(x - x_0)}}{\norm{x - x_0}}
\\[15pt]
&\leq \ 
\frac{\norm{h(x) - h(x_0) - \td h(x_0) (x - x_0)}}{\norm{x - x_0}}.
\end{align*}
The first and third terms converge to 0 when $x \ra x_0$, thus so does the second term.
\\

Passing to the proof of the theorem, take $m = 1$ and assume that $L_f$ is nonempty.  
Consider the countable collection $\{B_1, B_2, \ldots\}$ of all open balls $B(x,r)^\circ$ 
contained in $\Omega$ with $x \in \Q^n$ and 
$r \in \Q \hsx \cap \ ]0, +\infty[$ such that $\restr{f}{B(x,r)^\circ}$ is bounded $-$then 
$\ds \L_f \subset \bigcup\limits_{n = 1}^\infty \ B_n$.  
Given $x \in B_n$, introduce
\[
\begin{cases}
\ u_n (x) = \sup \{u(x) \hsy : \hsy  u \leq f \ \text{on} \ B_n, \ \Lip(u; B_n) \leq n\}
\\[8pt]
\ v_n (x) = \inf \{v(x) \hsy : \hsy  v \geq f \ \text{on} \ B_n, \ \Lip(v; B_n) \leq n\}
\end{cases}
.
\]
Here the ``sup'' (``inf'') is over all the $u$ $(v)$ with the stated properties, thus
\[
u_n 
\ \leq \ 
\restr{f}{B_n}
\ \leq \ 
v_n
\]
and
\[
\begin{cases}
\ \Lip(u_n; B_n) \leq n
\\[8pt]
\ \Lip(v_n; B_n) \leq n
\end{cases}
.
\]

\noindent
Let $E_n$ stand for the set of $x \in B_n$ at which both $u_n$ and $v_n$ are differentiable at $x$, 
hence by Rademacher, the set
\\[-1.2cm]

\[
Z 
\ = \ 
\bigcup_{n = 1}^\infty \ B_n \backslash E_n
\]
has Lebesgue measure 0.  
The claim now is that $f$ is differentiable at all points of $L_f \backslash Z$.  
So let $x_0 \in L_f \backslash Z$ 
$-$then it need only be shown that there is an index $n$ such that 
$x_0 \in E_n$ and $u_n(x_0) = v_n(x_0)$.  
This said, choose an $r_0 > 0$ and an $M > 0$ such that

\[
\abs{f(x) - f(x_0)} 
\ \leq \ 
M \hsy \norm{x - x_0} 
\qquad (x \in B(x_0, r_0)^\circ).
\]
Next, choose $n \geq M$: 
\[
x_0 \in B_n \subset B(x_0, r_0)^\circ.
\]
Then $x_0 \notin Z \implies x_0 \in E_n$ 
(for $x_0 \notin E_n \implies x_0 \in B_n \backslash E_n \subset Z$).  
Proceeding, $\forall \ x \in B_n$
\[
f(x) 
\ \leq \ 
f(x_0) + M \hsy \norm{x - x_0}
\ \leq \ 
f(x_0) + n \hsy \norm{x - x_0}
\]

$\implies$
\[
f(x) 
\ \leq \ 
v_n(x) 
\ \leq \ 
f(x_0) + n \hsy \norm{x - x_0}
\]

$\implies$
\[
f(x_0) 
\ \leq \ 
v_n(x_0) 
\ \leq \
f(x_0) 
\]

$\implies$
\[
f(x_0) 
\ = \ 
v_n(x_0).
\]
Therefore $u_n(x_0) = v_n(x_0)$, completing the proof. 
\\

\qquad{\bf 2.5.5:}\  
{\small\bf APPLICATION} \ 
Suppose that $\Lip f(x) < +\infty$ almost everywhere $-$then $f$ is differentiable almost everywhere.
\\

\qquad{\bf 2.5.6:}\  
{\small\bf EXAMPLE} \ 
Quasiconformal maps are differentiable almost everywhere.
\\

\qquad{\bf 2.5.7:}\  
{\small\bf REMARK} \ 
It can be shown that the subset $E \subset \Omega$ consisting of those $x$ 
at which $f$ is differentiable is Lebesgue measurable, as are the partial derivatives
\[
\frac{\partial f}{\partial x_i} \hsy : \hsy E \ra \R 
\qquad (i = 1, \ldots, n).
\]
\\[-.75cm]

\qquad{\bf 2.5.8:}\  
{\small\bf \un{N.B.}} \ 
The set of points where a given first order partial derivative $f$ exists need not be Lebesgue measurable. 
\\

\qquad{\bf 2.5.9:}\  
{\small\bf EXAMPLE} \ 
Let $S \subset \R$ be a non Lebesgue measurable set and let
\[
f(x, y) 
\ = \ 
\chisubmQ (x) \hsy \chisubS (y) 
\qquad ((x, y) \in \R^2).
\]
Then $f$ is Lebesgue measurable but the set of points $(x, y)$ at which $\ds \frac{\partial f}{\partial x_i}$ exists is not Lebesgue measurable.
\\

\qquad{\bf 2.5.10:}\  
{\small\bf REMARK} \ 
It can be shown that if $f : \Omega \ra \R$ is continuous and 
\\[-.5cm]

\noindent
if $E_i$ is the set of all $x \in \Omega$ such that 
$\ds \frac{\partial f}{\partial x_i}$ exists, 
then $E_i$ is a Borel set and $\ds \frac{\partial f}{\partial x_i}$ is a Borel function in $E_i$.
\\[-.5cm]

[Note: \ 
If instead $f : \Omega \ra \R$ is merely Borel measurable, then the $\ds \frac{\partial f}{\partial x_i}$ are Lebesgue measurable.]

\chapter{
$\boldsymbol{\S}$\textbf{2.6}.\quad  LUSIN}
\setlength\parindent{2em}
\renewcommand{\thepage}{2-\S6-\arabic{page}}

\qquad
Convention: \ Be it a set or a function, measurable means Lebesgue measurable.
\\

\qquad{\bf 2.6.1:}\  
{\small\bf THEOREM} \ 
Suppose given $f: \R^n \ra \R$ $-$then the following conditions are equivalent:  
$f$ is measurable or
\\[-.25cm]

\qquad \textbullet \quad 
For every $\varepsilon > 0$ and any compact $K \subset \R^n$, there is an open set $G \subset \R^n$  such that $\Lm^n(G) < \varepsilon$ 
and $\restr{f}{K \backslash G}$ is continuous.
\\[-.25cm]

\qquad \textbullet \quad 
For every $\varepsilon > 0$ and any compact $K \subset \R^n$, there exists a continuous function $\phi : \R^n \ra \R$ such that
\[
\Lm^n (\{x \in K \hsy : \hsy f(x) \neq \phi (x)\}) < \varepsilon.
\]
\\[-1.25cm]

\qquad \textbullet \quad 
For every compact $K \subset \R^n$, there exists a sequence $\{\phi_n\}$ of continuous functions  $\phi_n : \R^n \ra \R$ such that 
$\phi_n \ra f$ almost everywhere on $K$.  
\\

\qquad{\bf 2.6.2:}\  
{\small\bf THEOREM} \ 
Suppose given a function $f: \R^n \ra \R$ $-$then the  following conditions are equivalent:  
$f$ is measurable or 
\\[-.25cm]

\qquad \textbullet \quad 
For every $\varepsilon > 0$, there is an open set $G \subset \R^n$  such that $\Lm^n(G) < \varepsilon$ 
and $\restr{f}{\R^n \backslash G}$ is continuous.
\\[-.25cm]

\qquad \textbullet \quad 
For every $\varepsilon > 0$, there exists a continuous function $\phi : \R^n \ra \R$ and an open set $G \subset \R^n$ such that 
$\Lm^n(G) < \varepsilon$ and $\phi = f$ on $\R^n \backslash G$.
\\[-.25cm]

\qquad \textbullet \quad 
There exists a sequence $\{\phi_n\}$ of continuous functions  $\phi_n : \R^n \ra \R$ such that 
$\phi_n \ra f$ almost everywhere on $\R^n$.
\\

\qquad{\bf 2.6.3:}\  
{\small\bf CHARACTERIZATION} \ 
Let $f : E \ra \R$ be a function defined on a measurable set $E \subset \R^n$ $-$then $f$ is measurable iff for every $\varepsilon > 0$ there 
exists a closed set
$F \subset E$ such that $\Lm^n (E \backslash F) < \varepsilon$ and the restriction of $f$ to $F$ is continuous.
\\

\qquad{\bf 2.6.4:}\  
{\small\bf SCHOLIUM} \ 
Suppose that  $f: \R^n \ra \R$ is measurable $-$then for any measurable set $E$ and any  $\varepsilon > 0$, there exists a continuous function 
$g : \R^n \ra \R$ such that 
\[
\Lm^n (\{x \in E \hsy : \hsy f(x) \neq g (x)\}) < \varepsilon 
\]
and
\[
\norm{g}_{\infty, \R^n} 
\ \leq \ 
\norm{f}_{\infty, E}.
\]
\\[-1.cm]

In particular: \ Take $E = \R^n$ $-$then the conclusion is that a measurable function coincides with a continuous function outside a set of 
arbitrarily small measure.  
\\[-.5cm]

There is also a $C^\prime$ version of this result, the proof of which depends on an extension theorem due to Whitney.
\\

\qquad{\bf 2.6.5:}\  
{\small\bf THEOREM} \ 
Let $K \subset \R^n$ be a compact set and let $f : K \ra \R$, $T : K \ra \R^n$ be continuous functions.  
Assume: \ For every $\varepsilon > 0$, there is a $\delta > 0$ such that 

\[
\frac{\abs{f(y) - f(x) - T(x) (y -x)}}{\norm{y - x}} 
\ \leq \ 
\varepsilon
\]
whenever $x$, $y \in K$, $x \neq y$, and $\norm{y - x} \leq \delta$ $-$then there exists a $C^\prime$ function $g : \R^n \ra \R$ such that 

\[
\restr{g}{K} \ = \ f
\quad \text{and} \quad
\restr{\nabla g}{K} \ = \ T.
\]
\\[-1cm]

\qquad{\bf 2.6.6:}\  
{\small\bf NOTATION} \ 
As usual, $\Omega$ is a nonempty open subset of $\R^n$.
\\

\qquad{\bf 2.6.7:}\  
{\small\bf APPLICATION} \ 
Suppose that $f : \Omega\ra \R$ is measurable and differentiable almost everywhere $-$then for any $\varepsilon > 0$, 
there is a function $g \in C^1 (\Omega; \R)$ such that 
\[
\Lm^n (\{x \in \Omega \hsy : \hsy f(x) \neq g (x)\}) < \varepsilon.
\]
\\[-1cm]

\qquad{\bf 2.6.8:}\  
{\small\bf \un{N.B.}} \ 
Thanks to Rademacher, this applies in the special case when $f$ is Lipschitz.

\chapter{
\text{SECTION 3: \quad DENSITY THEORY}
\vspace{1.25cm}\\
$\boldsymbol{\S}$\textbf{3.1}.\quad  LEBESGUE POINTS}
\setlength\parindent{2em}
\renewcommand{\thepage}{3-\S1-\arabic{page}}

\qquad
Let $f \in \Lp_\locx^1 (\R^n)$.
\\

\qquad{\bf 3.1.1:}\   
{\small\bf DEFINITION} \ 
A point $x \in \R^n$ such that 
\[
\lim\limits_{r \ra 0} \ 
\frac{1}{\omega_n r^n} \ 
\int\limits_{B(x,r)} \ \abs{f - f (x)} \ \td \Lm^n 
\ = \ 
0
\]
is called a \un{Lebesgue point} of $f$.
\\[-.5cm]

[Note: \ 
Recall that
\[
\Lm^n (B(x,r)) \ = \ \omega_n r^n.
\]
In particular, if $n = 1$, then 
\[
\omega_1
\ = \ 
\frac{\pi^{1/2}}{\Gamma (1 + 1/2)}
\ = \ 
\frac{\pi^{1/2}}{\Gamma (3/2)}
\ = \ 
\frac{\pi^{1/2}}{\frac{1}{2} \hsx \Gamma (1/2)}
\ = \ 
\frac{\pi^{1/2}}{\frac{1}{2} \hsx\pi^{1/2}}
\ = \ 
2.]
\]
\\[-.5cm]

\qquad{\bf 3.1.2:}\   
{\small\bf DEFINITION} \ 
The \un{Lebesgue set} of $f$ is the set of its Lebesgue points, denoted $\Lambda (f)$.
\\

\qquad{\bf 3.1.3:}\   
{\small\bf THEOREM} \ 
\[
\Lambda (f) \in \sM_\Lm^n
\]
and
\[
\Lm^n (\R^n \backslash \Lambda (f)) 
\ = \ 
0.
\]
\\[-1cm]

\qquad{\bf 3.1.4:}\   
{\small\bf \un{N.B.}}  \ 
Every point of continuity of $f$ is a Lebesgue point of $f$. 
\\[-.5cm]

[Supposing that $f$ is continuous at $x$, given $\varepsilon > 0$, there exists $\delta > 0$ such that 
$\abs{f(y) - f(x)} < \varepsilon$ if $y \in B(x,\delta)$, so

$
\hspace{3cm}
r \in \hsx ]0, \delta[ 
\ \implies \ 
B(x,r) \subset B(x,\delta)
$
\\[-.5cm]

\hspace{1cm} $\implies$
\begin{align*}
\frac{1}{\omega_n r^n} \ 
\int\limits_{B(x,r)} \ \abs{f - f (x)} \ \td \Lm^n \ 
&\leq \ 
\frac{1}{\omega_n r^n} \ 
\int\limits_{B(x,r)} \ \varepsilon \ \td \Lm^n 
\\[15pt]
&=\ 
\varepsilon.]
\end{align*}
\\[-1cm]

\qquad{\bf 3.1.5:}\   
{\small\bf DEFINITION} \ 
If $E \subset \R^n$ is Lebesgue measurable, then the \un{density} of $E$ at a point $x \in \R^n$ (not necessarily in $E$) is
\[
\tD_E (x)
\ = \ 
\lim\limits_{r \ra 0} \ 
\frac{\Lm^n (E \cap B(x,r))}{\omega_n r^n}
\quad 
\big(\in [0,1]\big)
\]
provided the limit exists.
\\

\qquad{\bf 3.1.6:}\   
{\small\bf EXAMPLE} \ 
Work in $\R$ and let 
\[
E \ = \ \bigcup\limits_{k = 0}^\infty \ I_k, 
\]
where
\[
I_k \ = \ 
\left[
\frac{1}{2^{2 k + 1}}, \frac{1}{2^{2 k}}
\right]
\]
hence

\[
\Lm^1 (I_k ) 
\ = \ 
\frac{1}{2^{2 k + 1}},
\]
but $\tD_E (0)$ does not exist.  
In fact, 

\[
\frac{\Lm^1 (E \cap B(0,2^{-2k}))}{\Lm^1 (B(0,2^{-2k}))} 
\ = \ 
\frac{1}{3}
\]
and

\[
\frac{\Lm^1 (E \cap B(0,2^{-2k - 1}))}{\Lm^1 (B(0,2^{-2k-1}))} 
\ = \ 
\frac{1}{6}.
\]
Therefore
\[
\frac{\Lm^1 (E \cap B(0,r))}{\Lm^1  (B(0,r))}
\]
assumes the value $\ds \frac{1}{3}$ for $r = 2^{-2k}$ and the value $\ds \frac{1}{6}$ for the value $r = 2^{-2k - 1}$, 
so $\tD_E(0)$ does not exist.
\\

\qquad{\bf 3.1.7:}\   
{\small\bf \un{N.B.}}  \ 
\[
E \in \sM_\Lm^n
\ \implies \  
\chisubE \in \Lp_\locx^1 (\R^n).
\]
\\[-.75cm]

\qquad{\bf 3.1.8:}\   
{\small\bf LEMMA} \ 
Let $E \in \sM_\Lm^n$ $-$then
\[
E^\circ \cup (\R^n \backslash E)^\circ  \subsetx \Lambda(\chisubE).
\]

PROOF \ 
If $x \in E^\circ$ (or if $x \in \R^n \backslash E)^\circ$), then $\chisubE$ is continuous at $x$, thus $x \in \Lambda(\chisubE)$.  
\\

\qquad{\bf 3.1.9:}\   
{\small\bf EXAMPLE} \ 
It can happen that $\tD_E(x)$ exists for some $x \notin \Lambda(\chisubE)$.
\\[-.5cm]

[Work in $\R$ and let $E = [0, +\infty[$ $-$then 
\[
\R \backslash \{0\} 
\ = \ 
E^\circ \cup (\R^n \backslash E)^\circ 
\subset \Lambda(\chisubE).
\]
On the other hand, 
\begin{align*}
\lim\limits_{r \ra 0} \ \frac{1}{\Lm^1 (B(0,r))} \ \int\limits_{B(0,r)} \ \abs{\chisubE - \chisubE (0)} \ \td \Lm^1 \ 
&=\ 
\lim\limits_{r \ra 0} \ 
\frac{1}{2 r} \ 
\int\limits_{-r}^r \ (1 - \chisubE (y))  \ \td \Lm^1
\\[15pt]
&=\ 
\lim\limits_{r \ra 0} \ 
\frac{1}{2 r} \ 
\bigg(
2 r - \int\limits_0^r  \ \td \Lm^1
\bigg)
\\[15pt]
&=\ 
\frac{1}{2}
\\[15pt]
&\neq \ 
0
\\[15pt]
&\implies
0 \notin  \Lambda(\chisubE).
\end{align*}
Nevertheless
\begin{align*}
\tD_E(0) \ 
&=\ 
\lim\limits_{r \ra 0} \ \frac{\Lm^1 (E \cap B(0,r))}{\Lm^1 (B(0,r))} 
\\[15pt]
&=\ 
\lim\limits_{r \ra 0} \ 
\frac{r}{2 r} 
\\[15pt]
&=\ 
\frac{1}{2}.]
\end{align*}
\\[-1cm]

\qquad{\bf 3.1.10:}\   
{\small\bf LEMMA} \ 
If $E \subset \R^n$ is a set of Lebesgue measure 0, then $E \subset \R^n \backslash \Lambda(\chisubE)$.
\\[-.5cm]

PROOF \ 
The assertion is trivial if $E = \emptyset$, so take an $x \in E$ $-$then $\chisubE (x) = 1$, while $\chisubE = 0$ 
almost everywhere in $\R^n$, hence

\begin{align*}
\lim\limits_{r \ra 0} \ 
\frac{1}{\omega_n r^n} \ 
\int\limits_{B(x,r)} \ \abs{\chisubE - \chisubE (x)} \ \td \Lm^n \ 
&=\ 
\lim\limits_{r \ra 0} \ 
\frac{1}{\omega_n r^n} \ 
\int\limits_{B(x,r)} \ 1 \ \td \Lm^n
\\[15pt]
&=\ 
1 
\\[15pt]
&\neq \ 
0
\\[15pt]
&\implies
0 \notin  \R^n \backslash \Lambda(\chisubE).
\end{align*}
\\[-1cm]

\qquad{\bf 3.1.11:}\   
{\small\bf EXAMPLE} \ 
Take for $E$ the Cantor set in $\R$ $-$then $\Lambda(\chisubE) = \R \backslash E$.  
Moreover, $\forall \ x \in \R$, $\tD_E(x)$ exists and is equal to zero.
\\[-.5cm]

[Recall that $E$ is a closed subset of $\R$, thus $\R \backslash E$ is an open subset of $\R$, 
thus $\R \backslash E = (\R \backslash E)^\circ \subset \Lambda(\chisubE)$.  
But $E$ is also a set of Lebesgue measure 0, hence 
\[
E \subset \R \backslash \Lambda(E)
\]

\hspace{1.5cm} $\implies$
\[
\R \backslash E \supset \Lambda(E)
\]

\hspace{1.5cm} $\implies$
\[
\Lambda(E) \supset \R \backslash E 
\]

\hspace{1.5cm} $\implies$
\[
\Lambda(E\chisubE) \ = \  \R \backslash E.
\]

As for the other contention, simply note that
\[\Lm^1(E) 
\ = \ 0 
\implies
\Lm^1(E \cap B(x,r)) \ = \ 0.]
\]
\\[-1cm]

\qquad{\bf 3.1.12:}\   
{\small\bf LEMMA} \ 
Let $E \in \sM_\Lm^n$ $-$then
\[
\tD_E (x) \ = \ 
\begin{cases}
\ 1 \quad \text{for $x \in E \cap \Lambda (\chisubE)$}\\[4pt]
\ 0 \quad \text{for $x \in (\R^n \backslash E) \cap \Lambda (\chisubE)$}
\end{cases}
.
\]
\\[-.75cm]

\qquad{\bf 3.1.13:}\   
{\small\bf SCHOLIUM} \ 
\[
\tD_E (x) \ = \ 
\begin{cases}
\ 1 \quad \text{for almost all $x \in E$}\\[4pt]
\ 0 \quad \text{or almost all $x \in \R^n \backslash E$}
\end{cases}
.
\]

[It is a question of establishing that 
\\[-.25cm]

\qquad \textbullet \quad 
$\Lm^n (E \backslash (E \cap \Lambda (\chisubE))) = 0$
\\[-.5cm]

\noindent
and
\\[-.5cm]

\qquad \textbullet \quad 
$\Lm^n ((\R^n \backslash E) \backslash ((\R^n \backslash E) \cap \Lambda (\chisubE))) = 0. $
\\

E.g.:
\begin{align*}
E \backslash (E \cap \Lambda (\chisubE)) \ 
&=\ 
\R^n \backslash E \cap (E \cap \Lambda (\chisubE))
\\[15pt]
&=\ 
E \cap  ((\R^n \backslash E) \cup \R^n \backslash \Lambda (\chisubE))
\\[15pt]
&=\ 
E \cap (\R^n \backslash \Lambda (\chisubE)).
\end{align*}
But
\[
\Lm^n (\R^n \backslash \Lambda (\chisubE))
\ = \ 
0.]
\]
\\[-1cm]

\qquad{\bf 3.1.14:}\   
{\small\bf NOTATION} \ 
Given $E \in \sM_\Lm^n$, let
\[
\sD_E 
\ = \ 
\{x \in \R^n : \tD_E(x) \ \text{exists}\}.
\]
\\[-1cm]

\qquad{\bf 3.1.15:}\   
{\small\bf \un{N.B.}}  \ 

\[
\sD_E \supset \Lambda (\chisubE).
\]
\\[-1cm]

\qquad{\bf 3.1.16:}\   
{\small\bf LEMMA} \ 

\[
\sD_E \in \sM_\Lm^n.
\]
\\[-1cm]

\qquad{\bf 3.1.17:}\   
{\small\bf LEMMA} \ 
The function
\[
x \ra \tD_E (x) 
\quad (x \in \sD_E)
\]
is Lebesgue measurable.
\\

\qquad{\bf 3.1.18:}\   
{\small\bf THEOREM} \ 
\[
\int\limits_{\R^n} \ \tD_E \ \td \Lm^n 
\ = \ 
\Lm^n (E).
\]

PROOF \ 
Write
\[
\R^n 
\ = \ 
(E \cap \Lambda (\chisubE)) 
\cup 
(\R^n \backslash \E) \cap \Lambda (\chisubE) 
\cup 
\R^n \backslash \Lambda (\chisubE).
\]
Then this is a disjoint union of Lebesgue measurable sets, the third of which, viz. $\R^n \backslash \Lambda (\chisubE)$, 
being of Lebesgue measure 0.  
Therefore
\begin{align*}
\int\limits_{\R^n} \ \tD_E \ \td \Lm^n  \ 
&=\ 
\int\limits_{E \hsx \cap \hsx  \Lambda (\chisubE)} \ \tD_E \ \td \Lm^n 
\ + 
\int\limits_{(\R^n \backslash \E) \hsx \cap \hsx  \Lambda (\chisubE)} \ \tD_E \ \td \Lm^n 
\\[15pt]
&=\ 
\int\limits_{E \hsx \cap \hsx  \Lambda (\chisubE)} \ 1 \ \td \Lm^n 
\ + 
\int\limits_{(\R^n \backslash \E) \hsx \cap \hsx  \Lambda (\chisubE)}  \ 0\ \td \Lm^n 
\\[15pt]
&=\
\Lm^n (E \hsx \cap \hsx \Lambda (\chisubE)).
\end{align*}

Write
\begin{align*}
E \cap \Lambda (\chisubE) \ 
&=\ 
E \backslash (\R^n \backslash \Lambda (\chisubE))
\\[15pt]
&=\ 
E \backslash (E \cap (\R^n \backslash \Lambda (\chisubE))),
\end{align*}
from which 
\begin{align*}
\Lm^n (E \cap \Lambda (\chisubE))  \ 
&=\ 
\Lm^n (E) - \Lm^n (E \cap  (\R^n \backslash \Lambda (\chisubE)))
\\[15pt]
&=\ 
\Lm^n (E),
\end{align*}
thereby completing the proof.
\\

\qquad{\bf 3.1.19:}\   
{\small\bf DEFINITION} \ 
If $E \subset \R^n$ is Lebesgue measurable, then a point $x \in \R^n$ (not necessarily in $E$) is a 
\un{point of density 1} for $E$, denoted $x \in E^1$, if $\tD_E(x) = 1$ and a  
\un{point of density 0} for $E$, denoted $x \in E^0$, if $\tD_E(x) = 0$.
\\

\qquad{\bf 3.1.20:}\   
{\small\bf DEFINITION} \ 
\\[-.25cm]

\qquad \textbullet \quad 
$E^1$ is the \un{measure theoretic interior} of $E$.
\\

\qquad \textbullet \quad 
$E^0$ is the \un{measure theoretic exterior} of $E$.
\\

\qquad{\bf 3.1.21:}\   
{\small\bf DEFINITION} \ 
The \un{measure theoretic boundary} of $E$, denoted $\partial_\sM E$, is the set of points where the density is neither 0 nor 1.
\\

\qquad{\bf 3.1.22:}\   
{\small\bf DEFINITION} \ 
A Lebesgue measurable set $E \subset \R^n$ is $d$-open if each point of $E$ is a point of density 1, 
i.e., if $\forall \ x \in E$, $\tD_E(x) = 1$.
\\

\qquad{\bf 3.1.23:}\   
{\small\bf EXAMPLE} \ 
Take $n = 1$ $-$then the set of irrational numbers is $d$-open.
\\

\qquad{\bf 3.1.24:}\   
{\small\bf LEMMA} \ 
Every open subset of $\R^n$ is $d$-open.
\\

\qquad{\bf 3.1.25:}\   
{\small\bf THEOREM}\ 
The collection of all $d$-open sets forms a topology, the \un{density topology}.
\\

\qquad{\bf 3.1.26:}\   
{\small\bf \un{N.B.}} \ 
The density topology is strictly finer than the Euclidian 
topology.
\\

Let $E \subset \R^n$ be Lebesgue measurable $-$then

\begin{align*}
\frac{\Lm^n (E \cap B(x,r))}{\omega_n r^n} 
\ + \  
\frac{\Lm^n (\R^n \backslash E) \cap B(x,r))}{\omega_n r^n} \ 
&=\ 
\frac{\Lm^n (\R^n \cap B(x,r))}{\omega_n r^n} 
\\[15pt]
&=\ 
\frac{\Lm^n (B(x,r))}{\omega_n r^n} 
\\[15pt]
&=\ 
\frac{\omega_n r^n}{\omega_n r^n} 
\\[15pt]
&=\ 
1.
\end{align*}

Then it follows that 
\[
\begin{cases}
\ \tD_E(x) = 1 \quad \text{iff} \quad \tD_{(\R^n \backslash E)} (x) = 0\\[8pt]
\ \tD_E(x) = 0 \quad \text{iff} \quad \tD_{(\R^n \backslash E)} (x) = 1\
\end{cases}
.
\]


\chapter{
$\boldsymbol{\S}$\textbf{3.2}.\quad  APPROXIMATE LIMITS}
\setlength\parindent{2em}
\renewcommand{\thepage}{3-\S2-\arabic{page}}


\qquad
Let $f : \R^n \ra \R$ be a Lebesgue measurable function. 
\\

\qquad{\bf 3.2.1:}\   
{\small\bf DEFINITION} \ 
An element $\ell \in \R$ is the \un{approximate limit} of $f$ as $y \ra x$, denoted
\[
\ap \lim\limits_{y \ra x} \ f(y) 
\ = \ 
\ell, 
\]
if for every $\varepsilon > 0$, the set
\[
\{y  \hsy : \hsy \abs{f(y) - \ell} \geq \varepsilon\}
\]
has density 0 at $x$, i.e., 
\[
\lim\limits_{r \ra 0} \ 
\frac{\Lm^n (\{ \abs{f - \ell} \geq \varepsilon\} \cap B(x, r))}{\omega_n^{r^n}}
\ = \ 
0.
\]
\\[-.75cm]

\qquad{\bf 3.2.2:}\   
{\small\bf LEMMA} \ 
Approximate limits are unique (thereby justifying the use of ``the'' in the definition).  
\\[-.5cm]

PROOF \ 
Let $\ell_1$ and $\ell_2$ be two candidates for the approximate limit per the definition.  
Assume that $\ell_1  \neq \ell_2$ and take $\varepsilon = \abs{\ell_1 - \ell_2} / 3$ $-$then for each $y \in \R^n$, 
\[
3 \hsy \varepsilon 
\ = \ 
\abs{\ell_1 - \ell_2}
\ \leq \ 
\abs{f(y) - \ell_1} + \abs{f(y) - \ell_2}
\]

$\implies$
\[
B(x,r) 
\subset 
\{\abs{f(y) - \ell_1} \geq \varepsilon\}
\cup 
\{\abs{f(y) - \ell_2} \geq \varepsilon\}.
\]

Proof: \ 
If there were a $y \in B(x, r)$ which was not in the union, then
\[
\begin{cases}
\ \abs{f(y) - \ell_1} < \varepsilon
\\[4pt]
\ \abs{f(y) - \ell_2} < \varepsilon
\end{cases}
\]

$\implies$

\[
\abs{f(y) - \ell_1} + \abs{f(y) - \ell_2}
\ < \ 
2 \hsy \varepsilon 
\]

$\implies$

\[
3 \hsy \varepsilon < 2 \hsy \varepsilon 
\implies 
3 < 2  \hsx 
\ldots \hsx .
\]
Therefore
\begin{align*}
\omega_n^{r^n} \ 
&=\ 
\Lm^n (B(x, r)) 
\\[11pt]
&\leq \ 
\Lm^n (\{ \abs{f - \ell_1} \geq \varepsilon\} \cap B(x, r)) 
 + 
\Lm^n (\{ \abs{f - \ell_2} \geq \varepsilon\} \cap B(x, r)) .
\end{align*}
\\[-1.25cm]

\noindent
Now divide through by $\omega_n^{r^n}$ and send $r$ to zero to get $1 \leq 0$.
\\

\qquad{\bf 3.2.3:}\   
{\small\bf THEOREM} \ 
\\[-.25cm]
\[
\ap \lim\limits_{y \ra x} \ f(y) 
\ = \ 
\ell
\]
iff there exists a Lebesgue measurable set $E \subset \R^n$ with $\xD_E (x) = 1$ such that 
\[
\lim\limits_{\substack{y \ra x \\y \in E}} \ f(y) 
\ = \ 
\ell.
\]
\\[-1cm]

[The discussion infra supplies the proof.]
\\

\qquad{\bf 3.2.4:}\   
{\small\bf \un{N.B.}} \ 
In view of established principles, $x$ may or may not belong to $E$.  
As for the symbol
\[
\lim\limits_{\substack{y \ra x \\y \in E}} \ f(y) 
\ = \ 
\ell, 
\]

\noindent
it means: \ $\forall \ \varepsilon > 0$, $\exists \ r > 0$ such that
\[
\abs{f(y) - \ell}
\ < \ 
\varepsilon
\]
if $y \in E \cap (B(x,r) \backslash \{x\})$.
\\

Start matters by assuming that the limit above is in force $-$then the claim is that for every $\varepsilon > 0$, the set
\[
\{y \hsy : \hsy \abs{f(y) = \ell} \geq \varepsilon\}
\]
has density 0 at $x$ or, equivalently, that the set

\[
\{y \hsy : \hsy \abs{f(y) - \ell} < \varepsilon\}
\]
has density 1 at $x$.  
This set, however, contains $E \cap (B(x,r) \backslash \{x\})$ for small $r$.  
Therefore

\[
\frac{\Lm^n (\{\abs{f(y) - \ell} < \varepsilon\} \cap B(x,r))}{\omega_n^{r^n}}
\ \geq \ 
\frac{\Lm^n (E \cap (B(x,r)\backslash \{x\}))}{\omega_n^{r^n}}.
\]
But
\begin{align*}
\Lm^n (E \cap B(x,r)) \ 
&=\ 
\Lm^n (E \cap (B(x,r)\backslash \{x\})) + \Lm^n (E \cap \{x\})
\\[11pt]
&=\ 
\Lm^n (E \cap (B(x,r)\backslash \{x\})) 
\end{align*}

$\implies$
\begin{align*}
\frac{\Lm^n (\{\abs{f(y) - \ell} < \varepsilon\} \cap B(x,r))}{\omega_n^{r^n}} \ 
&\geq \ 
\frac{\Lm^n (E \cap B(x,r))}{\omega_n^{r^n}}
\\[15pt]
&\ra 1 \quad (r \ra 0).
\end{align*}

In the other direction, assume that 
\[
\ap \lim\limits_{y \ra x} \ f(y) 
\ = \ 
\ell, 
\]
the objective being to construct an $E \in \sM_{\Lm}^n$ with the stated property.  
To this end, choose a strictly decreasing sequence $\{r_k\}$ such that
\[
\Lm^n \big(\big\{\abs{f - \ell} \geq \frac{1}{k} \big\} \cap B(x,r_k)\big)
\ \leq \ 
\frac{\Lm^n ( B(x,r_k))}{2^k}
\]
and put
\[
E 
\ = \ 
\R^n \hsx \backslash \ 
\bigcup\limits_{k = 1}^\infty \ 
\left(
B(x, r_k) \backslash B(x, r_{k+1})
\right)
\hsx \cap \hsx 
\big\{\abs{f(y) - \ell} \geq \frac{1}{k} \big\} .
\]

\noindent
Then $E$ is Lebesgue measurable and 
\[
\lim\limits_{\substack{y \ra x \\y \in E}} \ f(y) 
\ = \ 
\ell.
\]
There remains the contention that $\xD_E (x) = 1$ or still, that
\[
\xD_{(\R^n \backslash E)} (x) \ = \ 0.
\]
By definition, 
\[
\R^n \backslash E
\ = \ 
\bigcup\limits_{k = 1}^\infty \ 
\left(
B(x, r_k) \backslash B(x, r_{k+1})
\right)
\hsx \cap \hsx 
\big\{y : \abs{f(y) - \ell} \geq \frac{1}{k} \big\} .
\]
Given $r > 0$, denote by $K$ the integer for which $r_{K + 1} < r < r_K$ $-$then
\allowdisplaybreaks
\begin{align*}
\Lm^n ((\R^n \backslash E) \cap B(x,r)) \ 
&\leq  \ 
\sum\limits_{k = K}^\infty \
\Lm^n ((B(x, r_k) \backslash B(x, r_{k+1})) \hsx \cap \hsx \big\{\abs{f - \ell} \geq \frac{1}{k} \big\})
\\[15pt]
&\leq \ 
\sum\limits_{k = K}^\infty \
\Lm^n(B(x, r_k) \hsx \cap \hsx \big\{\abs{f - \ell} \geq \frac{1}{k} \big\})
\\[15pt]
&\leq \ 
\sum\limits_{k = K}^\infty \
\frac{\Lm^n ( B(x,r_k) )}{2^k}
\\[15pt]
&\leq \ 
\Lm^n ( B(x,r) ) \ 
\sum\limits_{k = K}^\infty \
\frac{1}{2^k}
\\[15pt]
&\ra 0 \quad (r \ra 0).
\end{align*}

\qquad{\bf 3.2.5:}\   
{\small\bf DEFINITION} \ 
A Lebesgue measurable function $f : \R^n \ra \R$ is \un{approximately continuous}
at $x \in \R^n$ if $f$ is defined at $x$ and
\[
\ap \lim\limits_{y \ra x} \ f(y) 
\ = \ 
f(x). 
\]
\\[-.75cm]

\qquad{\bf 3.2.6:}\   
{\small\bf SCHOLIUM} \ 
$f$ is approximately continuous at $x$ iff there exists a Lebesgue measurable set $E \subset \R^n$ with $\xD_E(x) = 1$ 
such that $\restr{f}{E}$ is continuous at $x$.
\\

\qquad{\bf 3.2.7:}\   
{\small\bf REMARK} \ 
In terms of the density topology, $f$ is approximately continuous at $x$ iff $f$ is $d$-continuous at $x$.
\\

\qquad{\bf 3.2.8:}\   
{\small\bf THEOREM} \ 
A Lebesgue measurable function $f : \R^n \ra \R$ is approximately continuous $\Lm^n$ almost everywhere.  
\\[-.5cm]

PROOF \ 
Given $\varepsilon > 0$, there is a continuous function $g : \R^n \ra \R$ and an open set $G \subset \R^n$ such that 
$\Lm^n (G) < \varepsilon$ and $f = g$ in $\R^n \backslash G$.   
On general grounds, almost every point of $\R^n \backslash G$ is a point of density 1, thus $f$ is approximately continuous at 
almost all points of $\R^n$.
\\

\qquad{\bf 3.2.9:}\   
{\small\bf LEMMA} \ 
Let $f \in \Lp_\locx^1 (\R^n)$ and let $\ell \in \R^n$.  
Assume: 
\[
\lim\limits_{r \ra 0} \ 
\frac{1}{\omega_n^{r^n}}\ 
\int\limits_{B(x, r)} \ 
\abs{f - \ell} 
\ \td \Lm^n
\ = \ 
0.
\]
Then
\[
\ap \lim\limits_{y \ra x} \ f(y) 
\ = \ 
\ell.
\]
\\[-1.25cm]

PROOF \ 
Thanks to Chebyshev, $\forall \ \varepsilon > 0$,  
\[
\varepsilon \  
\frac{\Lm^n (\{\abs{f- \ell} \geq \varepsilon\} \cap B(x,r))}{\omega_n^{r^n}}
\ \leq \ 
\frac{1}{\omega_n^{r^n}}\ 
\int\limits_{B(x, r)} \ 
\abs{f - \ell} 
\ \td \Lm^n .
\]
\\[-.75cm]

\qquad{\bf 3.2.10:}\   
{\small\bf APPLICATION} \ 
The approximate limit exists at each Lebesgue point $x$ of $f$ and coincides with the value $f(x)$.
\\

The literal converse to 3.2.9 is false in general, i.e., it can happen that $\forall \ \varepsilon > 0$, 
\[
\lim\limits_{r \ra 0} \
\frac{\Lm^n (\{\abs{f - \ell} \geq \varepsilon\} \cap B(x,r))}{\omega_n^{r^n}}
\ = \ 
0,
\]
yet the relation

\[
\lim\limits_{r \ra 0} \ 
\frac{1}{\omega_n^{r^n}}\ 
\int\limits_{B(x, r)} \ 
\abs{f - \ell} 
\ \td \Lm^n 
\ = \ 
0
\]
fails or, what amounts to the same, it can happen that at some $x$, there is no $\ell$ such that 

\[
\lim\limits_{r \ra 0} \ 
\frac{1}{\omega_n^{r^n}}\ 
\int\limits_{B(x, r)} \ 
\abs{f - \ell} 
\ \td \Lm^n 
\ = \ 
0
\]
but for some $\ell$ and $\forall \ \varepsilon > 0$, 

\[
\lim\limits_{r \ra 0} \
\frac{\Lm^n (\{\abs{f - \ell} \geq \varepsilon\} \cap B(x,r))}{\omega_n^{r^n}}
\ = \ 
0.
\]
\\[-.75cm]

\qquad{\bf 3.2.11:}\   
{\small\bf EXAMPLE} \ 
In $\R^2$, take $\alpha > 0$ and consider
\[
f(x, y) \ = \ 
\begin{cases}
\ 0 \hspace{1.75cm} \text{if \ $y \leq 0$ \ or \  $y \geq x^2$}
\\[4pt]
\ \abs{y}^{-\alpha} \hspace{1cm} \text{otherwise}
\end{cases}
\quad (x, y) \in \R^2.
\]
Then

\[
\frac{1}{\omega_2^{r^2}}\ 
\int\limits_{B((0, 0),r)} \ 
\abs{f} 
\ \td x \td y
\]
tends to $+\infty$ as $r \ra 0$ if $\ds \frac{1}{2} < \alpha < 1$ while choosing $\ell = 0$, the sets $\{\abs{f - 0} \geq \varepsilon\}$
have density 0 at (0, 0).
\\

\qquad{\bf 3.2.12:}\   
{\small\bf LEMMA} \ 
Suppose that the sets 
\[
E_\varepsilon 
\ = \ 
\{\abs{f - \ell} \geq \varepsilon\}
\]
have density 0 at $x$ and $f$ is bounded in a neighborhood of $x$, say $\abs{f} \leq M$ $-$then

\[
\limsup\limits_{r \ra 0} \ 
\frac{1}{\omega_n^{r^n}}\ 
\int\limits_{B(x, r)} \ 
\abs{f - \ell} 
\ \td \Lm^n 
\ = \ 
0.
\]
\\[-1cm]

PROOF \ 
Write
\\[-.25cm]
\[
\frac{1}{\omega_n^{r^n}} \ 
\int\limits_{B(x, r)} \ 
\abs{f - \ell} 
\ \td \Lm^n 
\ \leq \ 
(M + \abs{\ell}) 
\hsx
\frac{\Lm^n (E_\varepsilon \cap B(x, r))}{\omega_n^{r^n}}
+ 
\varepsilon \hsx 
\frac{\Lm^n (B(x, r) \backslash E_\varepsilon)}{\omega_n^{r^n}}
\]

\noindent
from which

\[
\limsup\limits_{r \ra 0} \ 
\frac{1}{\omega_n^{r^n}}\ 
\int\limits_{B(x, r)} \ 
\abs{f - \ell} 
\ \td \Lm^n 
\ \leq \ 
\varepsilon.
\]
Now let $\varepsilon \ra 0$.
\\

\qquad{\bf 3.2.13:}\   
{\small\bf DEFINITION} \ 
$f$ has an \un{AFP approximate limit $\ell$} at $x$ if
\[
\lim\limits_{r \ra 0} \ \frac{1}{\omega_n^{r^n}} \ 
\int\limits_{B(x, r)} \ 
\abs{f - \ell} 
\ \td \Lm^n 
\ = \ 
0.
\]
\\[-.5cm]

\qquad{\bf 3.2.14:}\   
{\small\bf NOTATION} \ 
$S_f$ is the set of points $x$ which do not possess an AFP approximate limit. 
\\

\qquad{\bf 3.2.15:}\   
{\small\bf \un{N.B.}} \ 
If $f$ has an AFP approximate limit $\ell$ at $x$, then 
\[
\ap \limsup\limits_{y \ra x} \ f(y) 
\ = \ 
\ell,
\]
the converse being false in general (cf. supra).
\\


\qquad{\bf 3.2.16:}\   
{\small\bf LEMMA} \ 
AFP approximate limits are unique.
\\

\qquad{\bf 3.2.17:}\   
{\small\bf NOTATION} \ 
Write $\bar{f} (x)$ in place of $\ell$. 
\\

\qquad{\bf 3.2.18:}\   
{\small\bf OBSERVATION} \ 
If $\bar{f}(x) = f(x)$, then $x$ is a Lebesgue point of $f$. 
\\

\qquad{\bf 3.2.19:}\   
{\small\bf LEMMA} \ 
The set of points where the AFP approximate limit exists does not depend on the representative in the equivalence class of $f$, 
i.e., if $f = g \hsy \Lm^n$ almost everywhere in $\Omega$, 
then $x \notin S_f$ iff $x \notin S_g$ and $\bar{f} (x) = \bar{g} (x)$.
\\

\qquad{\bf 3.2.20:}\   
{\small\bf LEMMA} \ 
$S_f$ is a Borel set of Lebesgue measure 0.
\\[-.5cm]

PROOF \ 
The complement of the Lebesgue set $\Lambda (f)$ of $f$ is a set of Lebesgue measure 0, hence $\Lm^n (S_f) = 0$.   
As for $S_f$ being Borel, write
\\[-.5cm]
\[
\R^n \backslash  S_f 
\ = \ 
\bigcup\limits_{k = 1}^\infty \ 
\bigcup\limits_{q \in\ Q} \ 
\big\{
x \hsy : \hsy \limsup\limits_{r \ra 0} \ 
\frac{1}{\omega_n^{r^n}}\ 
\int\limits_{B(x, r)} \ 
\abs{f - q} 
\ \td \Lm^n 
< \frac{1}{k}
\big\}.
\]

[The inclusion $\subset$ is trivial.  
On the other hand, if $x$ belongs to the set on the RHS, then for any integer $k \geq 1$, there is a $q_k \in \Q$ such that

\[
\limsup\limits_{r \ra 0} \ 
\frac{1}{\omega_n^{r^n}}\ 
\int\limits_{B(x, r)} \ 
\abs{f - q_k} 
\ \td \Lm^n 
< \frac{1}{k}.
\]
The sequence $\{q_k\}$ obtained in this way is Cauchy and its limit $\ell$ has the property that

\[
\lim\limits_{r \ra 0} \ 
\int\limits_{B(x, r)} \ 
\abs{f - \ell} 
\ \td \Lm^n 
\ = \ 
0,
\]
so $x \notin S_f$, i.e., $x \in \R^n \backslash S_f$.]
\\

\qquad{\bf 3.2.21:}\   
{\small\bf LEMMA} \ 
$\bar{f} \hsy : \hsy \R^n \backslash S_f \ra \R$ is a Borel function which coincides with $\Lm^n$ almost everywhere with 
$\restr{f}{\R^n \backslash S_f }$.
\\[-.5cm]


PROOF \ 
In fact, for any $x \in \R^n \backslash S_f$, 

\[
\lim\limits_{r \ra 0} \ 
\frac{1}{\omega_n^{r^n}} \ 
\int\limits_{B(x, r)} \ 
f 
\ \td \Lm^n 
\ = \ 
\ell 
\ = 
\bar{f} (x), 
\]
thus $\bar{f}$ is the pointwise limit as $r \ra 0$ of the continuous function

\[
x \ra 
\frac{1}{\omega_n^{r^n}} \ 
\int\limits_{B(x, r)} \ 
f 
\ \td \Lm^n.
\]
\\[-1cm]

\qquad{\bf 3.2.22:}\   
{\small\bf EXAMPLE}  \ 
Suppose that $f = \chisubE$ is a characteristic function ($E$ a Lebesgue measurable set) $-$then $S_f$ is the measure theoretic 
boundary $\partial_M E$ of $E$.
\\[-.25cm]

On occasion, it will be necessary to consider a generalization of ``\ap lim''.
\\[-.25cm]

\qquad{\bf 3.2.23:}\   
{\small\bf DEFINITION} \ 
Let $f : \R^n \ra \R$ be a Lebesgue measurable function. 
\\

\qquad \textbullet \quad
An element $\ell \in \R$ is the \un{approximate $\limsup$} of $f$ as $y \ra x$, denoted
\[
\ap \limsup\limits_{y \ra x} \ f(y), 
\]
if $\ell$ is the infimum of the real numbers $t$ such that 

\[
\lim\limits_{r \ra 0} \ 
\frac{\Lm^n (\{ f > t\} \cap B(x, r))}{\omega_n^{r^n}}
\ = \ 
0.
\]

\qquad \textbullet \quad
An element $\ell \in \R$ is the \un{approximate $\liminf$} of $f$ as $y \ra x$, denoted

\[
\ap \liminf\limits_{y \ra x} \ f(y), 
\]
if $\ell$ is the supremum of the real numbers $t$ such that 
\[
\lim\limits_{r \ra 0} \ 
\frac{\Lm^n (\{ f < t\}  \cap B(x, r))}{\omega_n^{r^n}}
\ = \ 
0.
\]


Obviously

\[
\ap \liminf\limits_{y \ra x} \ f(y) 
\ \leq \ 
\ap \limsup\limits_{y \ra x} \ f(y) 
\]
and if 

\[
\ap \liminf\limits_{y \ra x} \ f(y) 
\ = \ 
\ap \limsup\limits_{y \ra x} \ f(y) 
\]
and if their common value is $\ell$, then the approximate limit exists and

\[
\ap \liminf\limits_{y \ra x} \ f(y) 
\ = \ 
\ell.
\]
\\[-1.5cm]

\[
* \ * \ * \ * \ * \ * \ * \ * \ * \ * \ * 
\]
\\[-1.75cm]

\[
\text{APPENDIX}
\]
\\[-1.5cm]

The preceding considerations have been formulated under the assumption that $f : \R^n \ra \R$ is Lebesgue measurable.  
Matters can be generalized.  
Thus let $S \subset \R^n$ be Lebesgue measurable and suppose that $f : S \ra \R$ is Lebesgue measurable.  
Fix a point $x \in \R^n$ such that $\xD_S (x) = 1$.
\\[-.25cm]

{\small\bf DEFINITION} \ 
An element $\ell \in \R$ is the \un{approximate limit} of $f$ as $y \ra x$ in $S$, denoted
\[
\ap \lim\limits_{\substack{y \ra x \\y \in S}} \ f(y) 
\ = \ 
\ell, 
\]
if for every $\varepsilon > 0$, the set
\[
\{y \in S \hsy : \hsy \abs{f(y) - \ell} \geq \varepsilon\}
\]
has density 0 at $x$.
\\

{\small\bf \un{N.B.}} \ 
If $S = \R^n$, then the demand that $\xD_{\R^n} (x) = 1$ is automatic.  
\\[-.5cm]

Proof: \ 
\begin{align*}
\xD_{\R^n} (x)  \ 
&=\ 
\frac{\Lm^n (\R^n \cap B(x, r))}{\omega_n^{r^n}}
\\[15pt]
&=\ 
\frac{\Lm^n (B(x, r))}{\omega_n^{r^n}}
\\[15pt]
&=\ 
\frac{\omega_n^{r^n}}{\omega_n^{r^n}}
\\[15pt]
&=\ 
1.
\end{align*}
\\[-1.25cm]

The earlier developments carry over modulo minor changes here and there.  
In particular: \ Approximate limits are unique and the notion of approximate continuity is clear.
\\[-.25cm]

{\small\bf THEOREM} \ 
\[
\ap \lim\limits_{\substack{y \ra x \\y \in S}} \ f(y) 
\ = \ 
\ell, 
\]
iff there exists a Lebesgue measurable set $E \subset \R^n$ with $\xD_E(x) = 1$ such that 
\[
\lim\limits_{\substack{y \ra x \\y \in E}} \ f(y) 
\ = \ 
\ell.
\]

\chapter{
$\boldsymbol{\S}$\textbf{3.3}.\quad  APPROXIMATE DERIVATIVES}
\setlength\parindent{2em}
\renewcommand{\thepage}{3-\S3-\arabic{page}}

\qquad
Let $f : \R^n \ra \R$ be a Lebesgue measurable function.
\\

\qquad{\bf 3.3.1:}\     
{\small\bf DEFINITION} \ 
$f$ is \un{approximately differentiable} at a point $x \in \R^n$ if there exists a linear function 
$T : \R^n \ra \R$ (depending on $x$) such that 

\[
\underset{y \ra x}{\ap \hsx \lim} \ 
\frac{\abs{f(y) - f(x) - T(y - x)}}{\norm{y - x}} 
\ = \ 
0.
\]
\\[-1cm]

$T$ is called the \un{approximate differential} of $f$ at $x$ and is denoted by
\[
\ap \hsx \td f (x).
\]
\\[-1.25cm]

[Note: \ 
If $f$ is differentiable at $x$ in the ordinary sense, then $f$ is approximately differentiable at $x$ and

\[
\td f(x) 
\ = \ 
\ap \hsx \td f (x).]
\]
\\[-1cm]

\qquad{\bf 3.3.2:}\     
{\small\bf \un{N.B.}} \ 
Existence is implied by demanding that

\[
\ap \hsx \limsup_{y \ra x}
\frac{\abs{f(y) - f(x) - T(y - x)}}{\norm{y - x}} 
\ = \ 
0.
\]
\\[-1cm]

\qquad{\bf 3.3.3:}\     
{\small\bf LEMMA} \ 
An approximate differential is unique (if it exists at all).  
\\[-.25cm]

PROOF \ 
Let $T_1$ and $T_2$ be two candidates for the approximate differential $-$then $\forall \ \varepsilon > 0$, 

\[
\lim\limits_{r \ra 0} \ 
\frac
{
\Lm^n 
\left(
\big\{
y \hsy : \hsy
\ds 
\frac{\abs{f(y) - f(x) - T_1(y - x)}}{\norm{y - x}} 
\ \geq \
\varepsilon
\}
\cap B(x, r)
\right)
}
{\omega_n r^n}
\ = \ 
0
\]
and

\[
\lim\limits_{r \ra 0} \ 
\frac
{
\Lm^n 
\left(
\big\{
y \hsy : \hsy
\ds 
\frac{\abs{f(y) - f(x) - T_2(y - x)}}{\norm{y - x}} 
\ \geq \
\varepsilon
\}
\cap B(x, r)
\right)
}
{\omega_n r^n}
\ = \ 
0. 
\]
To get a contradiction, suppose that $T_1 \neq T_2$ and take 
$\varepsilon = \norm{T_1 - T_2} / 6$.
Let

\[
S 
\ = \ 
\left\{
y \hsy : \hsy \abs{T_1 - T_2) (y - x)}
\ \geq \ 
\frac{\norm{T_1 - T_2} \hsx \norm{y - x}}{2}
\right\}.
\]
Then

\[
\frac{\Lm^n (S \cap B(x, r))}{\omega_n \hsy r^n}
\ \equiv \ 
C 
\ > \ 
0
\]
for all $r > 0$.  
On the other hand, 
\\[-.25cm]

$y \in S \implies$
\begin{align*}
3 \hsy \varepsilon \norm{y - x}  \ 
&=\ 
\frac{\norm{T_1 - T_2} \hsx \norm{y - x}}{2}
\\[11pt]
&\leq \
\abs{(T_1 - T_2) (y - x)}
\\[11pt]
&\leq \
\abs{f(y) - f(x) - T_1 (y -x)}
\hsx + \hsx 
\abs{f(y) - f(x) - T_2 (y -x)}
\end{align*}

$\implies$

\[
S \subset 
\left\{
y \hsy : \hsy 
\frac{\abs{f(y) - f(x) - T_1 (y -x)}}{\norm{y - x}}
\ \geq \ 
\varepsilon
\right\}
\ \cup \ 
\left\{
y \hsy : \hsy 
\frac{\abs{f(y) - f(x) - T_2 (y -x)}}{\norm{y - x}}
\ \geq \ 
\varepsilon
\right\}
\]
\\[-1cm]

$\implies$
\[
\lim\limits_{r \ra 0} \ 
\frac{\Lm^n (S \cap B(x, r))}{\omega_n \hsy r^n}
\ = \ 
0,
\]
a contradiction \ldots \hsx .
\\

[Note: \ 
Here is a different proof.  
Suppose that $T_1 \neq T_2$ and put $T = T_1 - T_2$, hence

\[
\underset{y \ra x}{\ap \hsx \lim} \ 
\frac{\abs{T(y - x)}}{\norm{y - x}} 
\ = \ 
0
\]
or still, 

\[
\underset{v \ra 0}{\ap \hsx \lim} \ 
\frac{\abs{T(v}}{\norm{v}} 
\ = \ 
0.
\]
So, if $0 < \varepsilon < 1$, then there exists $r > 0$ such that 

\[
\frac{\Lm^n (B(0 r) \cap \{v \hsy : \hsy \abs{T (v)} \geq \varepsilon \norm{v}\})}{\omega_n \hsy r^n}
\ < \ 
\varepsilon
\]
and for every $u \in \R^n$ with $\norm{u} = r - r \varepsilon$, there exists

\[
v \in B(u, \varepsilon r) 
\subset 
B(0,r)
\]
with

\[
\abs{T(v)}
\ \leq \ 
\varepsilon \norm{v}
\]
which implies that
\begin{align*}
\abs{T (u)} \ 
&=\ 
\abs{T (u - v + v)} 
\\[11pt]
&\leq \ 
\abs{T (u -v)} + \abs{T (v)} 
\\[11pt]
&\leq \ 
\norm{T} \norm{u - v} + \varepsilon  \norm{v}
\\[11pt]
&\leq \ 
\norm{T} \varepsilon r + \varepsilon  \norm{v}
\\[11pt]
&\leq \ 
\norm{T} \varepsilon r + \varepsilon r
\\[11pt]
&= \ 
(\norm{T}  + 1) \hsx \varepsilon r.
\end{align*}

\noindent
And

\[
r 
\ = \ 
\frac{r (1 - \varepsilon )}{1 - \varepsilon } 
\ = \ 
\frac{r  - r \varepsilon}{1 - \varepsilon } 
\ = \ 
\frac{\norm{u}}{1 - \varepsilon } 
\]

$\implies$
\[
\varepsilon r
\ = \ 
\frac{\varepsilon}{1 - \varepsilon} \norm {u}
\]

$\implies$
\begin{align*}
\abs{T (u)} \ 
&\leq\ 
(\norm{T} + 1) \hsx  \varepsilon r
\\[11pt]
&\leq \ 
(\norm{T} + 1) \hsx \frac{\varepsilon}{1 - \varepsilon} \norm {u}
\end{align*}

$\implies$

\[
\norm{T} 
\ \leq \ 
(\norm{T} + 1) \hsx \frac{\varepsilon}{1 - \varepsilon}
\]

$\implies$
\[
\norm{T} 
\ = \ 
0
\implies 
T_1 \ = \ T_2.]
\]
\\[-.75cm]

\qquad{\bf 3.3.4:}\     
{\small\bf REMARK} \ 
If $f$ is approximately differentiable at $x$, then $f$ is approximately continuous at $x$.
\\

\qquad{\bf 3.3.5:}\     
{\small\bf THEOREM} \ 
Let $f, \hsx g :\R^n \ra \R$ be Lebesgue measurable functions.  
Assume: \ $f$ is approximately differentiable almost everywhere and $f = g$ almost everywhere $-$then 
$g$ is approximately differentiable almost everywhere and 
\[
\ap \hsx\hsx \td f 
\ = \ 
\ap \ \td g
\]
almost everywhere in $\R^n$.
\\

Therefore the notion of approximate differentiability does not depend on the particular choice of the representative 
in the equivalence class.
\\

\qquad{\bf 3.3.6:}\     
{\small\bf THEOREM} \ 
$f$ is approximately differentiable at $x$ iff there exists a Lebesgue measurable set $E \subset \R^n$  
and a linear function $T : \R^n \ra \R$ with $\tD_E (x) = 1$ such that 
\\

\[
\lim\limits_{\substack{y \ra x \\y \in E}} \ 
\frac{\abs{f(y) - f(x) - T (y -x)}}{\norm{y - x}}
\ = \ 
0.
\]
\\[-.75cm]

\qquad{\bf 3.3.7:}\     
{\small\bf DEFINITION} \ 
For $i = 1, \ldots, n$, the \un{approximate partial derivative} $\ap \hsx \tD_i f (x)$ of $f$ at a point $x \in \R^n$ 
is defined by the condition
\\

\[
\underset{t \ra 0}{\ap \hsx \lim} \ 
\frac{\abs{f( x + t e_i) - f(x) - \ap \hsx \tD_i f (x) \hsx t}}{\abs{t}} 
\ = \ 
0.
\]
\\[-.75cm]

\qquad{\bf 3.3.8:}\     
{\small\bf THEOREM} \ 
The following conditions are equivalent.
\\[-.25cm]

\qquad (a) \quad 
The function $f$ has approximate partial derivatives almost everywhere in $\R^n$.
\\[-.5cm]

\qquad (b) \quad 
The function $f$ is approximately differentiable almost everywhere in $\R^n$.
\\[-.5cm]

[Note: \ 
Work in $\R^n$ $(n > 1)$ $-$then it can happen that the partial derivatives of $f$ exist almost everywhere in $\R^n$, 
yet $f$ might be nowhere differentiable 
(but, of course, $f$ will be approximately differentiable almost everywhere in $\R^n$).]
\\

\qquad{\bf 3.3.9:}\     
{\small\bf \un{N.B.}} \ 
The equivalent conditions (a) and (b) are also equivalent to 
\\[-.25cm]

\qquad (c) \quad For every $\varepsilon > 0$ there is a locally Lipschitz function $g : \R^n \ra \R$ such that 
\[
\Lm^n (\{x \hsy : \hsy f(x) \neq g(x))\}) 
\ < \ 
\varepsilon
\]
or even 
\\[-.25cm]

\qquad (d) \quad For every $\varepsilon > 0$ there is a $C^\prime$-function $g : \R^n \ra \R$ such that 
\[
\Lm^n (\{x \hsy : \hsy f(x) \neq g(x)\}) 
\ < \ 
\varepsilon.
\]
\\[-1cm]

\qquad{\bf 3.3.10:}\     
{\small\bf NOTATION} \ 
$A_\tD (f)$ is the domain of existence of $\ap \hsx \td f$.
\\

\qquad{\bf 3.3.11:}\     
{\small\bf LEMMA} \ 
If $f$ is approximately differentiable at $\Lm^n$ almost all points in $\R^n$, 
i.e., if $\Lm^n (\R^n \backslash A_\tD (f)) = 0$, then there exist Lebesgue measurable sets 
$E_0$, $E_k$ $(k = 1, 2, \ldots)$ such that 

\[
A_\tD (f) 
\ = \ 
E_0 \hsx \cup \hsx 
\bigcup\limits_{k = 1}^\infty \ E_k, 
\]
where $\Lm^n (E_0) = 0$ and for every $k$, the restriction $\restr{f}{E_k}$ is Lipschitz.
\\

\qquad{\bf 3.3.12:}\     
{\small\bf \un{N.B.}} \ 
\[
\ap \hsx \td f : A_\tD (f) \ra \R
\]
is Legesgue measurable.
\\

Owing to 2.5.1, $f$ is differentiable at almost all points where

\[
\limsup\limits_{y \ra x} \ \frac{\abs{f(x) - f(y)}}{\norm{x - y}} 
\ < \ 
+\infty.
\]
\\[-.75cm]

\qquad{\bf 3.3.13:}\     
{\small\bf LEMMA} \ 
$f$ is approximately differentiable at almost all points where

\[
\ap \hsx \limsup\limits_{y \ra x} \ \frac{\abs{f(x) - f(y)}}{\norm{x - y}} 
\ < \ 
+\infty.
\]
\\[-.75cm]

\[
* \ * \ * \ * \ * \ * \ * \ * \ * \ * \ 
\]

\[
\text{APPENDIX}
\]

Suppose that $f$ has ordinary partial derivatives almost everywhere $-$then 
$f$ is approximately differentiable almost everywhere, 
thus if $f$ is approximately differentiable at $x$, 
there exists a Lebesgue measurable set $E \subset \R^n$ with $\tD_E (x) = 1$ such that 
$\restr{f}{E}$ is differentiable at $x$ in the ordinary sense.  
Moreover, 
\[
\td (\restr{f}{E}) (x) 
\ = \ 
\ap \hsx \td f (x).
\]

Assume now that $n = 2$, $-$then in this special case $f$ admits a \un{regular} approximate differential at $x$.  
Here ``regular'' means that the ubiquitous set $E$ is comprised of the boundaries of oriented squares centered at $x$.  
\\

\textbf{SUMMARY} \ 
If $f : \R^2 \ra \R$ has ordinary partial derivatives almost everywhere, 
then it has a regular approximate differential almost everywhere.

\chapter{
SECTION \textbf{4}:\quad  WEAK PARTIAL DERIVATIVES}
\setlength\parindent{2em}
\setcounter{chapter}{4}
\renewcommand{\thepage}{4-\arabic{page}}

\vspace{-.75cm}
\qquad
Let $\Omega$ be a nonempty open subset of $\R^n$.
\\

\qquad{\bf 4.1:}\   
{\small\bf DEFINITION} \ 
A Lebesgue measurable function 
$f : \Omega \ra \R$ 
is \un{locally} \un{integrable} if

\[
\int\limits_K \ \abs{f} \ \td \Lm^n \ < \ +\infty
\]
for every compact $K \subset \Omega$.
\\

Denote the space of such by 

\[
\Lp_\locx^1 (\Omega).
\]
\\[-1cm]

\qquad{\bf 4.2:}\  
{\small\bf EXAMPLE} \ 
Take $\Omega = \R$ $-$then $\elln \abs{x} \in \Lp_\locx^1 (\R)$ but 
$x^{-1} \notin \Lp_\locx^1 (\R)$.
\\

\qquad{\bf 4.3:}\  
{\small\bf DEFINITION} \ 
Let $1 \leq p < +\infty$ $-$then a Lebesgue measurable function 
$f : \Omega \ra \R$ 
is \un{locally $\Lp^p$} if
\[
\int\limits_K \ \abs{f}^p \ \td \Lm^n 
\ < \ +\infty
\]
for every compact $K \subset \Omega$.
\\

\qquad{\bf 4.4:}\  
{\small\bf LEMMA} \ 
Every locally $\Lp^p$ function $f$ is locally $\Lp^1$ (i.e., is locally integrable).
\\[-.5cm]

PROOF \ 
Given a compact $K \subset \Omega$, by H\"older's inequality

\[
\int\limits_K \ \abs{f} \ \td \Lm^n 
\ \leq \ 
\norm{\chisubK \hsy f}_{\Lp^p} \hsx \norm{\chisubK \hsy f}_{\Lp^{p^\prime}},
\]
where 
$p^\prime = +\infty$ if $p = 1$ and 
$p^\prime = p / (p - 1)$ if $1 < p < +\infty$.
\\

\qquad{\bf 4.5:}\  
{\small\bf \un{N.B.}} \ 
The product of two functions in $\Lp_\locx^1 (\Omega)$ need not be locally integrable.
\\

Every $f \in \Lp_\locx^1 (\Omega)$ determines a distribution via the arrow
\[
\phi 
\ra
\int\limits_\Omega\ \phi \hsx f \hsy \td \Lm^n.
\]
Moreover two locally integrable functions define the same distribution iff they are equal almost everywhere.
\\[-.5cm]

[Note: \ 
A distribution $T$ ``is a function'' if there exists an element $f \in \Lp_\loc^1 (\Omega)$ such that $T = f \hsx \td \Lm^n$.]
\\

\qquad{\bf 4.6:}\  
{\small\bf NOTATION} \ 
Let $T : C_c^\infty (\Omega) \ra \R$ be a distribution $-$then

\[
\frac{\phi T}{\partial x_i}
\quad (i = 1, \ldots, n)
\]
is the distributional derivative of $T$ : $\forall \ \phi \in C_c^\infty (\Omega)$, 

\[
\langle \phi, \frac{\partial T}{\partial x_i} \rangle 
\ = \ 
-
\langle \frac{\partial \phi}{\partial x_i}, T \rangle.
\]
\\[-.75cm]

\qquad{\bf 4.7:}\  
{\small\bf DEFINITION} \ 
Given an $f \in \Lp_\locx^1 (\Omega)$, denote by

\[
\frac{\partial f}{\partial x_i} 
\quad (i = 1, \ldots, n)
\]
its distributional derivative (per $T = f \hsy \td \Lm^n$) $-$then $\ds\frac{\partial f}{\partial x_i}$ is said to be a 
\un{weak partial} \un{derivative} of $f$ if 
$\ds\frac{\partial f}{\partial x_i} \in \Lp_\locx^1 (\Omega)$, thus $\forall \ \phi \in C_c^\infty (\Omega)$, 

\[
\langle \phi, \frac{\partial f}{\partial x_i} \rangle 
\ = \ 
-
\langle \frac{\partial \phi}{\partial x_i}, f \hsy \td \Lm^n \rangle
\]
or still, 

\[
\int\limits_\Omega \ \phi \hsx \frac{\partial f}{\partial x_i} \ \td \Lm^n 
\ = \ 
-
\int\limits_\Omega \ 
\frac{\partial \phi}{\partial x_i} \hsx f
 \hsx \td \Lm^n.
\]
\\[-.75cm]

\qquad{\bf 4.8:}\  
{\small\bf EXAMPLE} \ 
Take $\Omega = \R$ and consider the function 
\[
h(x) \ = \ 
\begin{cases}
\ 0 \hspace{0.5cm} \text{if} \ x \leq 0\\[4pt]
\ x \hspace{0.5cm} \text{if} \ x > 0
\end{cases}
.
\]
Then $h \in \Lp_\locx^1 (\R)$ and its distributional derivative $\ds\frac{\partial h}{\partial x}$ is the Heaviside function

\[
H(x) \ = \ 
\begin{cases}
\ 0 \hspace{0.5cm} \text{if} \ x \leq 0\\[4pt]
\ 1 \hspace{0.5cm} \text{if} \ x > 0
\end{cases}
,
\]
which is therefore the weak derivative of $h$.  
Since $H \in \Lp_\locx^1 (\R)$, one can form its distributional derivative $\ds\frac{\partial H}{\partial x}$, so

\begin{align*}
\langle \phi, \frac{\partial H}{\partial x} \rangle \ 
&=\ 
-
\langle \frac{\partial \phi}{\partial x}, H \hsy \td \Lm^1 \rangle
\\[15pt]
&=\ 
-
\int\limits_\R \ \frac{\partial \phi}{\partial x} \hsy H(x) \ \td \Lm^1
\\[15pt]
&=\ 
-
\int\limits_0^\infty \ \frac{\partial \phi}{\partial x} \ \td \Lm^1
\\[15pt]
&=\ 
-
\left[
\phi(\infty) - \phi(0)
\right]
\\[15pt]
&=\ 
\phi(0).
\end{align*}
Consequently $\ds\frac{\partial H}{\partial x} = \delta$, the Dirac measure concentrated at the origin.  
However there 
\\[-.5cm]

\noindent
is  no $f \in \Lp_\locx^1 (\R)$ such that

\[
\int\limits_\R \ \phi  \hsx f \hsy  \td \Lm^1
\ = \ 
\phi(0)
\]
for all $\phi \in C_c^\infty(\Omega)$, hence $H$ does not have a weak derivative.
\\

\begin{spacing}{1.75}
\qquad{\bf 4.9:}\  
{\small\bf DEFINITION} \ 
Suppose that $f \in \Lp_\locx^1 (\Omega)$ admits weak partial derivatives 
$\ds\frac{\partial f}{\partial x_1}, \ldots, \frac{\partial f}{\partial x_n}$ $-$then the 
\un{distributional gradient} attached to $f$ is the $n$-tuple
\end{spacing}

\[
\nabla \hsy f
\ = \ 
\left(
\frac{\partial f}{\partial x_1}, \ldots, \frac{\partial f}{\partial x_n}
\right).
\]
\\[-.75cm]

\qquad{\bf 4.10:}\  
{\small\bf EXAMPLE} \ 
Working in $\R^n$, take $\Omega = B(0,1)^\circ \backslash \{0\}$ and define 
$f \in C^\infty (B(0,1)^\circ \backslash \{0\})$ by the rule

\[
f(x) = \norm{x}^{-\alpha}
\quad (\alpha > 0) 
\quad (\norm{x} = (x_1^2 + \cdots + x_n^2)^{1/2}).
\]
Then $f$ is unbounded in every neighborhood of the origin $(0 < \norm{x} < 1)$ and 
\[
\frac{\partial f}{\partial x_i} (x) 
\ = \ 
-\alpha \hsx \frac{x_i}{\norm{x}^{\alpha + 2}}
\quad (i = 1, \ldots, n).
\]
Therefore

\[
\nabla \hsy f (x)
\ = \ 
-\alpha \hsx \frac{x}{\norm{x}^{\alpha + 2}}
\]

\hspace{1.5cm} $\implies$ 
\[
\norm{\nabla \hsy f (x)}
\ = \ 
\frac{\abs{\alpha}}{\norm{x}^{\alpha + 1}}\hsx ,
\]
where
\[
\norm{\nabla \hsy f (x)}
\ = \ 
\left(
\sum\limits_{i = 1}^n \ 
\bigg|
\frac{\partial f}{\partial x_i} (x) 
\bigg|^2
\right)^{1/2}.
\]
\\[-.5cm]

\qquad{\bf 4.11:}\  
{\small\bf RAPPEL} \ 
Let $\bS^{n-1}$ $(= \partial B(0,1))$ be the unit sphere and let $\sigma^{n-1}$ be its surface measure, thus
\begin{align*}
\sigma^{n-1} (\bS^{n-1}) \ 
&=\ 
n \hsy \Lm^n (B(0,1))
\\[15pt]
&=\ 
n \hsx \frac{\pi^{n/2}}{\Gamma (1 + n/2)}
\\[15pt]
&=\
n \hsy \omega_n.
\end{align*}
\\[-.75cm]

\qquad{\bf 4.12:}\  
{\small\bf APPLICATION} \ 
Given $\alpha > 0$ subject to $n > \alpha$, put $f(x) = \norm{x}^{-\alpha}$ and write
\begin{align*}
\int\limits_{B(0,1)} \ 
\abs{f} \ 
\td \Lm^n 
&= \ 
\int\limits_{B(0,1)} \ 
\norm{x}^{-\alpha} \ 
\td \Lm^n 
\\[15pt]
&= \ 
\int\limits_0^1 \ 
\int\limits_{\bS^{n-1}} 
\norm{r x}^{-\alpha} \ 
\hsy r^{n-1} 
\td \sigma^{n-1} (x) \td r
\\[15pt]
&= \ 
\sigma^{n-1} (\bS^{n-1}) \ 
\int\limits_0^1 \ 
r^{-\alpha + n - 1} \ \td r
\\[15pt]
&= \ 
\sigma^{n-1} (\bS^{n-1}) \ 
\frac{r^{-\alpha+n}}{- \alpha +n} \bigg|_0^1
\\[15pt]
&< \ 
+\infty.
\end{align*}
\\[-.75cm]
Therefore $f \in \Lp^1 (B(0,1))$.
\\

\qquad{\bf 4.13:}\  
{\small\bf EXAMPLE} \ 
Consider again $f(x) = \norm{x}^{-\alpha}$ $(\alpha > 0)$ but replace 
$\Omega = B(0,1)^\circ \backslash \{0\}$ by $\Omega = B(0,1)^\circ$ $-$then
\allowdisplaybreaks
\begin{align*}
n > \alpha  
&\implies
f \in \Lp^1 (B(0,1))
\\[11pt]
&\implies
f \in \Lp^1 (\Omega)
\\[11pt]
&\implies
f \in \Lp_\locx^1 (\Omega).
\end{align*}
\\[-.75cm]

\noindent
Next
\allowdisplaybreaks
\begin{align*}
\int\limits_\Omega \ 
\abs{\frac{\partial f}{\partial x_i}} 
\ \td \Lm^n 
&\leq \ 
\int\limits_\Omega \ \norm{\nabla f} 
\td \Lm^n 
\\[15pt]
&= \ 
\abs{\alpha} \hsx 
\int\limits_\Omega \ 
\frac{1}{\norm{x}^{\alpha + 1}}
\ \td \Lm^n 
\\[15pt]
&= \ 
\abs{\alpha} \hsx \sigma^{n-1} (\bS^{n-1}) \  
\frac{r^n - \alpha -1}{n - \alpha - 1} \bigg|_0^1 
\\[15pt]
&< \ 
+\infty
\end{align*}
if $n > \alpha + 1$, so 

\[
\frac{\partial f}{\partial x_i}  \in \Lp^1 (\Omega) 
\implies 
\frac{\partial f}{\partial x_i}  \in \Lp_\locx^1 (\Omega).
\]
Let $T$ be the distribution corresponding to $f$, hence $\forall \ \phi \in C_c^\infty(\Omega)$, 

\begin{align*}
\langle \phi, \frac{\partial T}{\partial x_i} \rangle \ 
&=\ 
-
\langle \frac{\partial \phi}{\partial x_i}, T \rangle
\\[15pt]
&=\ 
-
\langle \frac{\partial \phi}{\partial x_i}, f \hsx \td \Lm^n\rangle
\\[15pt]
&=\ 
-
\langle \phi, \frac{\partial f}{\partial x_i} \hsx \td \Lm^n \rangle
\qquad \text{(dominated convergence).}
\end{align*}
Accordingly, as distributions, 

\[
\frac{\partial T}{\partial x_i}
\ = \ 
\frac{\partial f}{\partial x_i} \hsx \td \Lm^n.
\]
Therefore $f$ admits weak partial derivatives in $B(0,1)^\circ$ 
(and not just in $B(0,1) \backslash \{0\}$).
\\

\qquad{\bf 4.14:}\  
{\small\bf LEMMA} \ 
If $f \in C^1(\Omega)$, then the ordinary partial derivatives $\ds\frac{\partial f}{\partial x_i}$ of $f$ 
are also the corresponding weak partial derivatives of $f$.

\chapter{
SECTION \textbf{5}:\quad  MOLLIFIERS}
\setlength\parindent{2em}
\setcounter{chapter}{5}
\renewcommand{\thepage}{5-\arabic{page}}

\qquad
Let $\Omega$ be a nonempty open subset of $\R^n$.
\\

\qquad{\bf 5.1:}\  
{\small\bf NOTATION} \ 
Given $\varepsilon > 0$ and a nonnegative even bounded function $\phi \in \Lp^1(\R^n)$ with 

\[
\sptx \phi \subset B(0,1), \quad
\int\limits_{\R^n} \ \phi \ \td \Lm^n 
\ = \ 1,
\]
put

\[
\phi_\varepsilon (x) 
\ = \ 
\frac{1}{\varepsilon^n} \hsx \phi\bigg(\frac{x}{\varepsilon}\bigg) 
\qquad (x \in \R^n).
\]
\\[-.75cm]

\qquad{\bf 5.2:}\  
{\small\bf DEFINITION} \ 
The $\phi_\varepsilon$ are called \un{mollifiers}.
\\

\qquad{\bf 5.3:}\  
{\small\bf \un{N.B.}} \ 
Mollifiers exist \ldots \hsx .
\\[-.5cm]

[The \un{standard} choice for $\phi$ is 
\[
\phi(x) \ = \ 
\begin{cases}
&\ds C(n) \hsx \exp\left(\frac{1}{\norm{x}^2 - 1}\right) 
\hspace{1cm} \text{if} \ \norm{x} < 1\\
&0  \hspace{4.8cm} \text{if} \ \norm{x} \geq 1
\end{cases}
,
\]
where $C(n) > 0$ is chosen so that

\[
\int\limits_{\R^n} \ \phi \ \td \Lm^n 
\ = \ 1.
\]
Here $\phi \in C_c^\infty (\R^n)$.  
Another possibility is

\[
\phi (x) 
\ = \ 
\frac{1}{\omega_n} \hsx \chisubBofZeroOneInt \hsy (x).]
\]
\\[-.75cm]

\qquad{\bf 5.4:}\  
{\small\bf NOTATION} \ 
Put 

\[
\Omega_\varepsilon 
\ = \ 
\{x \in \Omega \hsx : \hsx \ \dist(x, \partial \Omega) > \varepsilon\}.
\]
\\[-1.5cm]

[Note: \ 
If $\Omega = \R^n$, then $\Omega_\varepsilon = \R^n$.]
\\

\qquad{\bf 5.5:}\  
{\small\bf DEFINITION} \ 
Given a function $f \in \Lp_\locx^1 (\Omega)$, write

\[
f_\varepsilon (x)
\ = \ 
(f * \phi_\varepsilon) (x) 
\ = \ 
\int\limits_\Omega \ 
\phi_\varepsilon (x - y)
\hsx f(y) \ \td \Lm^n,
\]
where $x \in \Omega_\varepsilon$.
\\[-.5cm]

[Note: \ 

\[
\int\limits_\Omega \ 
\phi_\varepsilon (x - y)
\hsx f(y) \ \td \Lm^n
\ = \ 
\int\limits_{B(x, \varepsilon)} \ 
\phi_\varepsilon (x - y)
\hsx f(y) \ \td \Lm^n.]
\]
\\[-.75cm]

The function $f_\varepsilon : \Omega_\varepsilon \ra \R$ is said to be a \un{mollification} of $f$, 
the \un{standard} \un{mollification} of $f$ being the $f_\varepsilon$ per the standard choice for $\phi$ per supra.
\\[-.5cm]

[Note: \ 
Given $x \in \Omega$, $f_\varepsilon (x)$ is well defined for all $0 < \varepsilon < \dist(x, \partial \Omega)$, 
thus it makes sense to consider
\[
\lim\limits_{\varepsilon > 0} \ f_\varepsilon (x).]
\]
\\[-.75cm]

\qquad{\bf 5.6:}\  
{\small\bf THEOREM} \ 
If $f_\varepsilon$ is the standard mollification of $f$, 
then $f_\varepsilon \in C^\infty (\Omega_\varepsilon)$ $(0 < \varepsilon < 1)$ 
and for every multi index $\alpha$ and for every $x \in \Omega_\varepsilon$, 
\[
\partial^\alpha f_\varepsilon (x) 
\ = \ 
(f * \partial^\alpha \phi_\varepsilon) (x)
\ = \ 
\int\limits_\Omega \ 
\frac
{\partial^{\abs{\alpha}} \hsy \phi_\varepsilon}
{\partial x^\alpha} (x - y) 
\hsx f(y) \ \td \Lm^n.
\]
\\[-.75cm]

\qquad{\bf 5.7:}\  
{\small\bf LEMMA} \ 
If the standard choice for $\phi$ is used and if $f \in \Lp_\locx^1 (\Omega)$ 
admits a weak partial derivative $\ds \frac{\partial f}{\partial x_i}$ 
(hence, by definition, 
$\ds \frac{\partial f}{\partial x_i} \in \Lp_\locx^1 (\Omega)$), 
then the derivative of the mollification 
coincides with the mollification of the weak partial derivative, i.e., 
\[
\frac{\partial f_\varepsilon}{\partial x_i}
\ = \ 
\frac{\partial f}{\partial x_i} \hsy * \hsy \phi_\varepsilon.
\]

In fact, 
\begin{align*}
\frac{\partial f_\varepsilon}{\partial x_i} (x) \ 
&=\ 
\frac{\partial }{\partial x_i}
\left(
\int\limits_\Omega \ \phi_\varepsilon(x - y) \hsy f(y) \ \td \Lm^n
\right)
\\[15pt]
&=\ 
\int\limits_\Omega \ 
\frac{\partial}{\partial x_i} \hsx
\phi_\varepsilon(x - y) \hsx f(y)
\ \td \Lm^n
\\[15pt]
&=\ 
(-1) \hsx 
\int\limits_\Omega \ 
\frac{\partial}{\partial y_i} \hsx
\phi_\varepsilon(x - y) \hsx f(y)
\ \td \Lm^n
\\[15pt]
&=\ 
(-1) \hsx (-1) \ 
\int\limits_\Omega \ \phi_\varepsilon(x - y) \hsx 
\frac{\partial f}{\partial y_i} (y)
\ \td \Lm^n
\\[15pt]
&=\ 
\left(
\frac{\partial f}{\partial x_i} \hsy * \hsy \phi_\varepsilon
\right)
(x).
\end{align*}
\\[-.75cm]

\begin{spacing}{1.45}
\qquad{\bf 5.8:}\  
{\small\bf APPLICATION} \ 
Work with the standard choice for $\phi$, 
suppose that $\Omega$ is connected, 
let $f \in \Lp_\locx^1 (\Omega)$, 
assume that the weak partial derivatives 
$\ds\frac{\partial f}{\partial x_i}$ $(i = 1, \ldots, n)$ 
exist and are equal to 0 almost everywhere 
$-$then $f$ coincides almost everywhere in $\Omega$ with a constant function.
\\[-.5cm]
\end{spacing}

[To begin with,

\[
\frac{\partial f_\varepsilon}{\partial x_i}
\ = \ 
\frac{\partial f}{\partial x_i} \hsy * \hsy \phi_\varepsilon
\ = \ 
0  \hsy * \hsy \phi_\varepsilon
\ = \ 
0,
\]
thus $f_\varepsilon$, being smooth, must be constant in each connected component of $\Omega_\varepsilon$ 
(in general, $\Omega_\varepsilon$ is not connected).  
Consider now a pair of points $x$, $y \in \Omega$ $-$then there exists a polygonal path 
$\gamma$ in $\Omega$ joining $x$ and $y$ and for small enough $\varepsilon$, $\gamma$ is in $\Omega_\varepsilon$, 
so 
$f_\varepsilon(x) = f_\varepsilon (y)$.  
But $f_\varepsilon \ra f$ $(\varepsilon \downarrow 0)$ almost everywhere (see below).  
Therefore 
$\lim f_\varepsilon$ $(\varepsilon \downarrow 0)$ is a constant function in $\Omega$.]
\\

\qquad{\bf 5.9:}\  
{\small\bf LEMMA} \ 
Suppose given $f \in C(\Omega)$ $(\subset \Lp_\locx^1 (\Omega)$ 
$-$then for every choice of $\phi$,

\[
f_\varepsilon \ra f
\qquad (\varepsilon \downarrow 0)
\]
uniformly on compact subsets of $\Omega$.
\\

\qquad{\bf 5.10:}\  
{\small\bf LEMMA} \ 
Suppose given $f \in \Lp_\locx^1 (\Omega)$ $-$then for every choice of $\phi$ and for every Lebesgue point 
$x \in \Omega$, 

\[
f_\varepsilon (x) 
\ra 
f(x) 
\qquad 
(\varepsilon \downarrow 0),
\]
hence $f_\varepsilon \ra f$ almost everywhere.
\\[-.5cm]

PROOF \ 
Write
\begin{align*}
\abs{f_\varepsilon (x)  - f(x)} \ 
&=\ 
\bigg|
\int\limits_\Omega \ \phi_\varepsilon (x - y \hsy f(y) \ \td \Lm^n \ - f(x)
\bigg|
\\[15pt]
&=\ 
\bigg|
\int\limits_{B(x,\varepsilon)} \ 
\phi_\varepsilon(x - y) \hsx f(y) \ \td \Lm^n
\hsx - \hsx 
f(x) \ 
\int\limits_{B(x,\varepsilon)} \
\phi_\varepsilon(x - y)
 \ \td \Lm^n
\bigg|
\\[15pt]
&=\ 
\bigg|
\int\limits_{B(x,\varepsilon)} \ 
\phi_\varepsilon(x - y) \hsx (f(y) - f(x)) 
\td \Lm^n
\bigg|
\\[15pt]
&\leq\ 
\frac{1}{\varepsilon^n} \hsx 
\int\limits_{B(x,\varepsilon)} \ 
\phi_\varepsilon\left(\frac{x - y}{\varepsilon}\right)
\abs{f(y) - f(x)} \td \Lm^n
\\[15pt]
&\leq\ 
\norm{\phi}_\infty \hsx \frac{\omega_n}{\omega_n} \hsx \frac{1}{\varepsilon^n} \hsy 
\int\limits_{B(x,\varepsilon)} \ \abs{f - f(x)} \td \Lm^n
\\[15pt]
&=\ 
\norm{\phi}_\infty \hsx \omega_n \hsx \frac{1}{ \omega_n \hsy \varepsilon^n} \  
\int\limits_{B(x,\varepsilon)} \ \abs{f - f(x)} \td \Lm^n
\\[15pt]
&\ra \ 
0 
\quad 
(\varepsilon \downarrow 0).
\end{align*}
\\[-.75cm]

\qquad{\bf 5.11:}\  
{\small\bf LEMMA} \ 
Suppose given $f \in \Lp_\locx^p (\Omega)$ $(1 \leq p < +\infty)$ 
$(\implies f \in \Lp_\locx^1 (\Omega))$ $-$then for every choice of $\phi$,
\[
\lim\limits_{\varepsilon \downarrow 0} \ \norm{f_\varepsilon - f}_{\Lp^p (\Omega)} 
\ = \ 0.
\]

\chapter{
SECTION \textbf{6}:\quad  $W^{1, \infty} (\R^n)$}
\setlength\parindent{2em}
\renewcommand{\thepage}{6-\arabic{page}}

\qquad
Let $(\XX, \sE)$ be a measurable space and let $\mu$ be a measure on $(\XX, \sE)$.
\\

\qquad{\bf 6.1:}\  
{\small\bf NOTATION} \ 
Given a measurable function $f : \XX \ra \R$, put

\[
\norm{f}_\infty
\ = \ 
\inf \{t \geq 0 : \mu (\{x : \abs{f(x)} > t\}) = 0\},
\]
with the convention that $\inf \emptyset = \infty$.
\\

\qquad{\bf 6.2:}\ 
{\small\bf DEFINITION} \ 
$\norm{f}_\infty$ is the \un{essential supremum} of $f$ and is written

\[
\norm{f}_\infty 
\ = \ 
\underset{x \in \XX}{\esssup} \abs{f(x)}.
\]

\qquad{\bf 6.3:}\ 
{\small\bf NOTATION} \ 
\[
\Lp^\infty (\XX) 
\ \equiv \ 
\Lp^\infty (\XX, \sE, \mu)
\]
is the set of measurable functions defined on $\XX$ for which $\norm{f}_\infty < \infty$.
\\[-.5cm]

[Note: \ 
Such functions are said to be \un{essentially bounded} and if $f$ is one such, then

\[
\abs{f(x)} 
\ \leq \ 
\norm{f}_\infty
\]
almost everywhere.]
\\

\qquad{\bf 6.4:}\ 
{\small\bf LEMMA} \ 
$f \in \Lp^\infty (\XX)$ iff there is a bounded measurable function $g$ such that $f = g$ almost everywhere.
\\

\qquad{\bf 6.5:}\ 
{\small\bf LEMMA} \ 
$\Lp^\infty (\XX)$ is a Banach space.
\\

Henceforth the pair $(\XX, \sE)$ will be the pair $(\R^n, \sM_\Lm^n)$, $\mu$ being $\Lm^n$.
\\

\qquad{\bf 6.6:}\ 
{\small\bf \un{N.B.}} \ 
The set of bounded  continuous functions $f:\R^n \ra \R$ carries the 
uniform norm
\[
\norm{f}_\tu 
\ = \ 
\sup\limits_{x \in \XX} \ \abs{f(x)}
\]
and
\[
\norm{f}_\tu 
\ = \ 
\norm{f}_\infty.
\]
\\[-.75cm]

\qquad{\bf 6.7:}\ 
{\small\bf NOTATION} \ 
$W^{1, \infty} (\R^n)$ is the space consisting of all essentially bounded functions $f : \R^n \ra \R$ whose distributional derivatives 
\[
\frac{\partial f}{\partial x_i}
\quad (i = 1, \ldots, n)
\]
are also essentially bounded functions in $\R^n$ as well.
\\

In $W^{1, \infty} (\R^n)$ introduce the norm 

\[
\norm{f}_{W^{1, \infty}} 
\ \equiv \ 
\norm{f}_\infty 
\hsx + \hsx 
\sum\limits_{i = 1}^n \ \norm{\frac{\partial f}{\partial x_i} }_\infty.
\]
\\[-.75cm]

\qquad{\bf 6.8:}\ 
{\small\bf THEOREM} \ 
$W^{1, \infty} (\R^n)$  is a Banach space.
\\

\qquad{\bf 6.9:}\ 
{\small\bf NOTATION} \ 
Given $f \in \Lp_\locx^1 (\R^n)$, the \un{$i^\nth$ difference quotient} is 

\[
\tD_i^h  f(x)
\ = \ 
\frac{f(x - h e_i) - f(x)}{h} 
\quad (i = 1, \ldots, n).
\]
\\[-.75cm]

\qquad{\bf 6.10:}\ 
{\small\bf LEMMA} \ 
$\forall \ \phi \in C_c^\infty(\R^n)$, 
\[
\int\limits_{\R^n} \ \frac{\phi(x + h e_i) - \phi(x)}{h} \hsx f(x) \ \td \Lm^n 
\ = \ 
- 
\int\limits_{\R^n} \ 
\phi(x) \hsx
\frac{f(x - h e_i) - f(x)}{-h} 
\ \td \Lm^n .
\]
Consequently
\[
\int\limits_{\R^n} \ \left(\tD_i^h \hsy \phi\right) \hsy f \ \td \Lm^n
\ = \ 
- 
\int\limits_{\R^n} \ \phi \left(\tD_i^{-h}  \hsy f\right) \ \td \Lm^n
\quad (i = 1, \ldots, n).
\]
\\[-.75cm]

\qquad{\bf 6.11:}\ 
{\small\bf THEOREM} \ 
Suppose that $f$ is a bounded Lipschitz continuous function, say
\[
\abs{f(x) - f(y)}
\ \leq \  
L \hsy \norm{x - y} 
\qquad (x, y \in \R^n).
\]
Then

\[
f \in W^{1, \infty} (\R^n).
\]
\\[-1cm]

Since by hypothesis $f \in \Lp^\infty (\R^n)$, the problem is to show that its distributional derivatives

\[
\frac{\partial f}{\partial x_i} 
\quad (i = 1, \ldots, n)
\]
are  (essentially) bounded functions as well.
\\[-.5cm]

To begin with, $\forall \ x \in \R^n$, 

\[
\abs
{
\frac{f(x - h e_i) - f(x)}{-h} 
}
\ \leq \ 
L 
\quad (i = 1, \ldots, n)
\]

\hspace{1.5cm} $\implies$

\[
\norm{\tD_i^{-h} f}_{\Lp^\infty (\R^n)}
\ \leq \ 
L 
\quad (i = 1, \ldots, n),
\]
so if $\Omega$ is open and bounded, 

\[
\norm{\tD_i^{-h} f}_{\Lp^2(\Omega)} 
\ \leq \ 
\norm{\tD_i^{-h} f}_{\Lp^\infty (\R^n)} 
\ \left(
\Lp^n(\Omega)
\right)^{1/2}
\leq 
L(\tL^n(\Omega))^{1/2}.
\]

Let $h_k = \ds\frac{1}{k}$ $(k = 1, 2, \ldots)$ $-$then 
\[
\{\tD_i^{-h_k} f \}
\]
is a bounded sequence in $\Lp^2 (\Omega)$, thus there is a subsequence

\[
\{\tD_i^{-h_{k_\ell}} f\}
\]
that converges weakly in $\Lp^2 (\Omega)$ as $\ell \ra \infty$ to $g_i \in \Lp^2 (\Omega)$.  
To simplify, put $h_j = h_{k_\ell}$, hence $h_j \ra 0$ $(j \ra \infty)$.
\\[-.5cm]

Accordingly, $\forall \ \phi \in C_c^\infty(\Omega)$,

\begin{align*}
- \langle \phi, \frac{\partial f}{\partial x_i} \rangle \ 
&= \ 
\int\limits_\Omega \  \frac{\partial \phi}{\partial x_i}  \ f \hsx \td \Lm^n
\\[15pt]
&= \ 
\int\limits_\Omega \  \left( \lim\limits_{j \ra \infty} \ \tD_i^{h_j} \phi \right) \hsy f \hsx \td \Lm^n
\\[15pt]
&= \ 
- \lim\limits_{j \ra \infty} \ \int\limits_\Omega \ \left(\tD_i^{h_j} \phi \right) \hsy f \hsx \td \Lm^n
\\[15pt]
&=\ 
- \lim\limits_{j \ra \infty} \ \int\limits_\Omega \ \phi \left(\tD_i^{-h_j} f \right)  \ \td \Lm^n
\\[15pt]
&= \ 
- \int\limits_\Omega \ \phi \hsy g_i \ \td \Lm^n,
\end{align*}
the last equality following from weak convergence.  
Therefore the weak partial derivative $\ds \frac{\partial f}{\partial x_i}$ exists and is represented by $g_i$.
\\

Because 

\[
f_j
\ \equiv \ 
\tD_i^{-h_j} f
\quad (i = 1, 2, \ldots)
\]
converges weakly in $\Lp^2 (\Omega)$ to $\ds\frac{\partial f}{\partial x_i}$, 
there exists a subsequence $\{f_{j_n}\}$ such that the convex combinations

\[
\sum\limits_{n = 1}^N \ a_n \hsx f_{j_n} 
\ra 
\frac{\partial f}{\partial x_i} 
\quad \text{in} \ \Lp^2 (\Omega)
\]
\\[-1cm]
\noindent
as $N \ra \infty$.  
Here

\allowdisplaybreaks
\begin{align*}
\normx{\sum\limits_{n = 1}^N \ a_n \hsx f_{j_n}}_{\Lp^\infty (\Omega)}  \ 
&\leq \ 
\sum\limits_{n = 1}^N \ a_n \hsx\normx{\tD_i^{-h_j} f}_{\Lp^\infty (\Omega)} 
\\[15pt]
&\leq \ 
\left(
\sum\limits_{n = 1}^N \ a_n
\right)
L
\\[15pt]
&= \ 
L.
\end{align*}
\\[-.75cm]

Summary: 

\[
\abs{\frac{\partial f}{\partial x_i}} 
\ \leq \ 
L 
\quad (i = 1, \ldots, n)
\]
for almost every $x \in \Omega$, hence

\[
\frac{\partial f}{\partial x_i} \in \Lp^\infty (\R^n).
\]
And then 

\[
f \in W^{1, \infty} (\R^n).
\]
\\[-.75cm]

\qquad{\bf 6.12:}\ 
{\small\bf REMARK} \ 
Let $f \in W^{1, \infty} (\R^n)$ $-$then on the basis of ``embedding theory'', 
it can be shown that $f$ has a bounded continuous representative $\bar{f}$.  
Moreover, $\bar{f}$ can be taken Lipschitz continuous.
\\[-.5cm]

[Working with the standard mollification $\bar{f}_\varepsilon$ and assuming that the support of $\bar{f}$ is compact, note that
\[
\norm{\nabla\hsy \bar{f}_\varepsilon}_{\Lp^\infty (\R^n)} 
\ \leq \ 
\norm{\nabla\hsy \bar{f}}_{\Lp^\infty (\R^n)} 
\qquad (\varepsilon > 0)
\]

\hspace{.5cm} $\implies$
\begin{align*}
\abs{\bar{f}_\varepsilon (x) - \bar{f}_\varepsilon (y)} \ 
&= \ 
\bigg|
\int\limits_0^1 \ 
\langle 
\nabla\hsy \bar{f}_\varepsilon (t x + (1 - t) y, x - y
\rangle
\ \td t
\bigg|
\\[15pt]
&\leq \ 
\norm{\nabla\hsy \bar{f}_\varepsilon}_{\Lp^\infty (\R^n)} \hsx \norm{x - y}
\\[15pt]
&\leq \ 
\norm{\nabla\hsy \bar{f}}_{\Lp^\infty (\R^n)}  \hsx \norm{x - y}
\end{align*}

\hspace{1.5cm} $\implies \ (\varepsilon \downarrow 0)$
\[
\abs{\bar{f} (x) - \bar{f} (y)} 
\ \leq \ 
\norm{\nabla\hsy \bar{f}}_{\Lp^\infty (\R^n)}  \hsx \norm{x - y}
\]
for all $x$, $y \in \R^n$.]

\chapter{
\text{SECTION 7: \quad SOBOLOV SPACES}
\vspace{.75cm}\\
$\boldsymbol{\S}$\textbf{7.1}.\quad  FORMALITIES}
\setlength\parindent{2em}
\renewcommand{\thepage}{7-\S1-\arabic{page}}
\vspace{-.5cm}

\qquad
Let $\Omega$ be a nonempty open subset of $\R^n$.
\\


\qquad{\bf 7.1.1:}
{\small\bf DEFINITION} \ 
Let $1 \leq p < +\infty$ $-$then the \un{Sobolov space} 
\[
W^{1, p} (\Omega)
\]
\begin{spacing}{1.55}
\noindent
consists of those $f \in \Lp_\locx^1 (\Omega)$ such that $f$ belongs to $\Lp^p (\Omega)$  and such that the 
distributional derivatives 
$\ds \frac{\partial f}{\partial x_i}$ are weak partial derivatives and also belong to 
$\Lp^p (\Omega)$ $(i = 1, \ldots, n)$.
\\
\end{spacing}

[Note: \ 
There is a local version of this definition, namely call 
\[
W_\locx^{1, p} (\Omega)
\]
the set comprised of all $f \in \Lp_\locx^1 (\Omega)$ with the property that the restriction 
$\restr{f}{\Omega^\prime} \in W^{1, p} (\Omega^\prime)$ for every nonempty open set 
$\Omega^\prime \subset \Omega$ whose closure is a compact subset of $\Omega$.]
\\

\begin{spacing}{1.65}
\qquad{\bf 7.1.2:}
{\small\bf \un{N.B.}} \ 
Spelled out, $W^{1, p} (\Omega)$ consists of those $f \in \Lp^p (\Omega)$ 
for which there exist functions 
$\ds \frac{\partial f}{\partial x_1}, \ldots, \frac{\partial f}{\partial x_n}$ in $\Lp^p (\Omega)$ such that 
$\forall \ \phi \in C_c^\infty (\Omega)$, 
\end{spacing}

\[
\int\limits_\Omega \ \phi \hsy \frac{\partial f}{\partial x_i} \ \td \Lm^n 
\ = \ 
-
\int\limits_\Omega \  \frac{\partial \phi}{\partial x_i} \hsx f \hsx \td \Lm^n 
\qquad (i = 1, \ldots, n).
\]

[Note: \ 
Another point is this: $W^{1, p} (\Omega)$ is closed under taking absolute values,
i.e., 

\[
f \in W^{1, p} (\Omega) 
\implies 
\abs{f} \in W^{1, p} (\Omega).]
\]
\\[-1cm]

Depending on the parameters, it can happen that there exists a function in $W^{1, p} (\Omega) $ which is nowhere continuous.
\\

\qquad{\bf 7.1.3:}
{\small\bf EXAMPLE} \ 
Take $\Omega = B(0,1)^\circ$, let $\{q_k\}$ be a countable dense subset of $\Omega$, and consider
\[
f(x) 
\ = \ 
\sum\limits_{k = 1}^\infty \ \frac{1}{2^k} \ \norm{x - q_k}^{-\alpha}
\qquad (\alpha > 0) \ 
\bigg(x \in \Omega \  \backslash \ \bigcup\limits_{k = 1}^\infty \ \{q_k\}\bigg).
\]
Then
\[
f \in W^{1, p} (\Omega)
\quad \text{if} \ \alpha < \frac{n - p}{p}
\]
but $f$ is unbounded in every nonempty open subset of $\Omega$.
\\

\qquad{\bf 7.1.4:}
{\small\bf \un{N.B.}} \ 
It will be seen later on that each function in $W^{1, p} (\Omega)$ $(p > n)$ coincides with a continuous function almost everywhere.
\\

\qquad{\bf 7.1.5:}
{\small\bf LEMMA} \ 
Let $f \in W^{1, p} (\Omega)$ $-$then there is a partition
\[
\Omega 
\ = \ 
\bigg(
\bigcup\limits_{k = 1}^\infty \ E_k \bigg)
\cup Z, 
\]
where the $E_k$ are Lebesgue measurable sets such that $\restr{f}{E_k}$ is Lipschitz and $Z$ has Lebesgue measure 0.
\\

\qquad{\bf 7.1.6:}
{\small\bf THEOREM} \ 
Let $f \in W^{1, p} (\Omega)$ $-$then $f$ is approximately differentiable almost everywhere.
\\[-.5cm]

[Extend $\restr{f}{E_k}$ to all of $\R^n$ and use Rademacher.]
\\

\qquad{\bf 7.1.7:}
{\small\bf THEOREM} \ 
The prescription

\[
\norm{f}_{W^{1, p}}
\ \equiv \ 
 \norm{f}_{\Lp^p} 
 \ + \ 
\sum\limits_{i = 1}^n \ \norm{\frac{\partial f}{\partial x_i}}_{\Lp^p} 
\]
endows $W^{1, p} (\Omega)$ with the structure of a Banach space.
\\[-.5cm]

[Note: \ 
An equivalent norm is the prescription

\[
f \ra \norm{f}_{\Lp^p} + \norm{\nabla f}_{\Lp^{p}},
\]
where $\nabla f$ is the weak gradient attached to $f$.]
\\

\qquad{\bf 7.1.8:}
{\small\bf LEMMA} \ 
$W^{1, p} (\Omega)$ is separable.
\\

\qquad{\bf 7.1.9:}
{\small\bf THEOREM} \ 
Let $f \in W^{1, p} (\Omega)$ $-$then there exists a sequence $\{f_k\} \subset W^{1, p} (\Omega) \cap C^\infty (\Omega)$ 
such that

\[
f_k \ra f 
\quad \text{in} \ W^{1, p} (\Omega).
\]
\\[-1cm]

\qquad{\bf 7.1.10:}
{\small\bf REMARK} \ 
It can be shown that $C_c^\infty (\R^n)$ is dense in $W^{1, p} (\R^n)$.
\\

\qquad{\bf 7.1.11:}
{\small\bf PRODUCT RULE} \ 
Let $f, \hsx g \in W^{1, p} (\Omega) \cap \Lm^\infty  (\Omega)$ $-$then
$f  g \in W^{1, p} (\Omega) \cap \Lm^\infty  (\Omega)$ 
and

\[
\frac{\partial ( f  g )}{\partial x_i} 
\ = \ 
\frac{\partial f}{\partial x_i} \hsx g + f \hsx \frac{\partial g}{\partial x_i} 
\qquad (i = 1, \ldots, n)
\]
\\
$\Lm^n$ almost everywhere in $\Omega$.
\\

\qquad{\bf 7.1.12:}
{\small\bf CHAIN RULE} \ 
Let $f \in W^{1, p} (\Omega)$ 
and let $g \in C^1(\R)$ subject to $g^\prime \in \Lp^\infty (\R)$, $g(0) = 0$ $-$then 
$g \circ f \in W^{1, p} (\Omega)$  and

\[
\frac{\partial (g \circ f}{\partial x_i} 
\ = \ 
(g^\prime \circ f) \hsx \frac{\partial f}{\partial x_i} 
\qquad (i = 1, \ldots, n)
\]
$\Lm^n$ almost everywhere in $\Omega$.
\\[-.5cm]

[Note: \ 
The assumption that $g(0) = 0$ is not needed if $\Lm^n (\Omega) < +\infty$.]
\\

To formulate the next result, let $\Omega^\prime$ be another nonempty open subset of $\R^n$.
\\

\qquad{\bf 7.1.13:}
{\small\bf CHANGE OF VARIABLE} \ 
Suppose that $\Psi : \Omega^\prime \ra \Omega$ is invertible, where $\Psi$ and $\Psi^\prime$ are Lipschitz continuous functions, 
and let $f \in W^{1, p} (\Omega)$ $-$then $f \circ \Psi \in W^{1, p} (\Omega^\prime)$ and 
\[
\frac{\partial (f \circ \Psi)}{\partial x_i^\prime} (x^\prime)
\ = \ 
\sum\limits_{k = 1}^n \ \frac{\partial f}{\partial x_k} 
\hsx (\Psi (x^\prime)) 
\hsx \frac{\partial \Psi_k}{\partial x_i^\prime} (x^\prime)
\qquad (i = 1, \ldots, n)
\]
$\Lm^n$ almost everywhere in $\Omega^\prime$.
\\

\qquad{\bf 7.1.14:}
{\small\bf TERMINOLOGY} \ 
Given normed spaces 
$(\XX, \norm{ \hsx \cdot \hsx}_\XX)$ and 
$(\YY, \norm{ \hsx \cdot \hsx}_\YY)$, 
one says that $\XX$ is \un{embedded} in $\YY$, denoted $\XX \hookrightarrow \YY$, if $\XX$ is a subspace of $\YY$ and there exists 
a constant $C > 0$ such that for all $x \in \XX$,
\[
\norm{x}_\YY 
\ \leq \ 
C \hsy \norm{x}_\XX. 
\]
\\[-1.25cm]

\qquad{\bf 7.1.15:}
{\small\bf EXAMPLE} \ 
Suppose that $\Lm^n (\Omega) < +\infty$ and $1 \leq p < q < +\infty$ $-$then $\forall \ f  \in \Lp^q(\Omega)$, 
\[
\norm{f}_{\Lp^p (\Omega)}
\ \leq \ 
\left(
\Lm^n (\Omega)
\right)^{(1/p) - (1/q)} \hsx \norm{f}_{\Lp^q (\Omega)} 
\]

\hspace{1.5cm} $\implies$
\[
\Lp^q(\Omega) 
\hookrightarrow
\Lp^p(\Omega).
\]

[Note: \ 
Here $\XX = \Lp^q(\Omega)$, $\YY = \Lp^p(\Omega)$, and 
\[
C 
\ = \ 
\left(
\Lm^n (\Omega)
\right)^{(1/p) - (1/q)}.]
\]

Looking ahead: 
\\[0.5cm]

\hspace{2cm}
\fbox
{
$
\begin{matrix*}[l]
&\\
&\ds W^{1,p} (\R^n) \hookrightarrow \Lp^{p^*} (\R^n) \hspace{0.7cm} \text{if} \hspace{0.7cm} p < n \left(\frac{1}{p^*} = \frac{1}{p} - \frac{1}{n}\right)
\\[11pt]
&W^{1,p} (\R^n) \hookrightarrow \Lp^q (\R^n)  \hspace{0.85cm} \text{if} \hspace{0.7cm} p = n (q \in [p, +\infty[)
\\[11pt]
&W^{1,p} (\R^n) \hookrightarrow \Lp^\infty (\R^n)  \hspace{0.72cm} \text{if} \hspace{0.7cm} p > n.\\
&
\end{matrix*}
$
}

\chapter{
$\boldsymbol{\S}$\textbf{7.2}.\quad  EMBEDDINGS: GNS}
\setlength\parindent{2em}
\renewcommand{\thepage}{7-\S2-\arabic{page}}

\qquad{\bf 7.2.1:}
{\small\bf DEFINITION} \ 
Let $1 \leq p < n$ $-$then the \un{conjugate exponent} of $p$ is 
\[
p^*
\ = \ 
\frac{np}{n - p}
\ = \ 
p + \frac{p^2}{n - p}.
\]
\\[-1.25cm]

[Note: \ 
$p^* > p$ and 
\[
\frac{1}{p^*} 
\ = \ 
\frac{1}{p} - \frac{1}{n}.]
\]
\\[-.75cm]

\qquad{\bf 7.2.2:}
{\small\bf THEOREM} \ 
Let $1 \leq p < n$ $-$then there exists a constant $C(n,p) > 0$ such that for all $f \in W^{1, p} (\R^n)$, 

\[
\bigg(\ 
\int\limits_{\R^n} \ \abs{f}^{p^*} \ \td \Lm^n
\bigg)^{1/p^*}
\ \leq \ 
C(n, p) \ 
\bigg( \ 
\int\limits_{\R^n} \ \norm{\nabla f}^p
 \ \td \Lm^n
\bigg)^{1/p}
\ < \ 
+\infty.
\]
\\[-.5cm]

\qquad{\bf 7.2.3:}
{\small\bf SCHOLIUM} \ 
When $1 \leq p < n$, 

\[
W^{1, p} (\R^n) \hookrightarrow \Lp^{p^*} (\R^n).
\]
\\[-1cm]

\qquad{\bf 7.2.4:}
{\small\bf RAPPEL} \ 
If $1 \leq p_1, \ldots, p_k < +\infty$ with 
\[
\frac{1}{p_1} + \cdots + \frac{1}{p_k}
\ = \ 1
\]
and if $f_j \in \Lp^{p_j} (\R^n)$ $(j = 1, \ldots, k)$, then
\[
\int\limits_{\R^n} \ \abs{f_1 \cdots f_k} \ \td \Lm^n 
\ \leq \ 
\prod\limits_{j = 1}^k \ 
\norm{f_j}_{\Lp^{p_j}}.
\]
\\[-.75cm]

The proof of the theorem can be divided into three parts.  
\\[-.25cm]

\qquad \un{Step 1:} \quad
$p = 1$, $f \in C_c^1 (\R^n)$.
\\[-.25cm]

[For each $i \in \{1, \ldots, n\}$ and each point $x = (x_1, \ldots, x_i, \ldots, x_n) \in \R^n$, write
\[
f (x_1, \ldots, x_i, \ldots, x_n) 
\ = \ 
\int\limits_{-\infty}^{x_i} \ \frac{\partial f}{\partial x_i} (x_1, \ldots, t_i, \ldots, x_n) \ \td t_i,
\]
hence
\[
\abs{f(x)}
\ \leq \ 
\int\limits_{\R} \ 
\bigg|
\frac{\partial f}{\partial x_i} (x_1, \ldots, t_i, \ldots, x_n)
\bigg|
\ \td t_i
\qquad (1 \leq i \leq n) 
\]
\\[-.75cm]

$\implies$
\[
\abs{f(x)}^n 
\ \leq \ 
\prod\limits_{i = 1}^n \ 
\int\limits_{\R} \ 
\bigg|
\frac{\partial f}{\partial x_i} (x_1, \ldots, t_i, \ldots, x_n)
\bigg|
\ \td t_i
\]
\\[-.75cm]

$\implies$
\[
\abs{f(x)}^{n / (n - 1)} 
\ \leq \ 
\prod\limits_{i = 1}^n \ 
\bigg( \ 
\int\limits_{\R} \ 
\bigg|
\frac{\partial f}{\partial x_i} (x_1, \ldots, t_i, \ldots, x_n)
\bigg|
\ \td t_i
\bigg)^{1/(n-1)}
\]
\\[-.75cm]

$\implies$
\allowdisplaybreaks
\begin{align*}
\int\limits_{\R} \ 
\abs{f}^{n / (n - 1)} \hsx 
\td x_1 \ 
&\leq \ 
\bigg( \ 
\int\limits_{\R} \ 
\bigg|
\frac{\partial f}{\partial x_i} 
\bigg|
\hsx \td t_1
\bigg)^{1/(n-1)}
\int\limits_{\R} \ 
\prod\limits_{i = 2}^n \ 
\bigg( \ 
\int\limits_{\R} \ \ 
\bigg|
\frac{\partial f}{\partial x_i} 
\bigg|
\ \td t_i
\bigg)^{1/(n-1)}
\hsx \td x_1
\\[15pt]
&\leq \ 
\bigg( \ 
\int\limits_{\R} \ 
\bigg|
\frac{\partial f}{\partial x_i} 
\bigg|
\td t_1
\bigg)^{1/(n-1)}
\prod\limits_{i = 2}^n \ 
\bigg( \ 
\int\limits_{\R} \  
\int\limits_{\R} \  
\bigg|
\frac{\partial f}{\partial x_i} 
\bigg|
\ \td t_i
\ \td x_1
\bigg)^{1/(n-1)}
\end{align*}
\\[-.75cm]

$\implies$ $\cdots$

\[
\int\limits_{\R} \  
\cdots
\int\limits_{\R} \  
\abs{f}^{n / (n - 1)} \ 
\td x_1 \ldots \ \td x_n 
\ \leq \ 
\prod\limits_{i = 1}^n \ 
\bigg( \ 
\int\limits_{\R} \  
\cdots
\int\limits_{\R} \
\bigg|
\frac{\partial f}{\partial x_i} 
\bigg|
\ 
\td x_1 \ldots \ \td x_n 
\bigg)^{1/(n-1)}
\]
\\[-.75cm]

$\implies$
\allowdisplaybreaks
\begin{align*}
\int\limits_{\R^n} \  
\abs{f}^{n / (n - 1)} \ \td \Lm^n \ 
&\leq \ 
\prod\limits_{i = 1}^n \ 
\bigg( \ 
\int\limits_{\R^n} \
\bigg|
\frac{\partial f}{\partial x_i} 
\bigg|  
\ \td \Lm^n 
\bigg)^{1/(n-1)}
\\[15pt]
&\leq \ 
\bigg( \ 
\int\limits_{\R^n} \  
\norm{ \nabla f}
\ \td \Lm^n 
\bigg)^{n/(n-1)}
\end{align*}

$\implies$

\[
\bigg(\ 
\int\limits_{\R^n} \  
\abs{f}^{n/(n-1)} 
\ \td \Lm^n 
\bigg)^{(n-1)/n} 
\ \leq \ 
\int\limits_{\R^n} \  
\norm{ \nabla f}
\ \td \Lm^n 
.]
\]
\\[-.5cm]

[Note: \  
\[
1^* 
\ = \ 
\frac{n}{n-1}
\implies 
1/p^* 
\ = \ 
\frac{n-1}{n}
\]
and
\[
C(n, 1) 
\ = \ 
1.]
\]
\\[-1cm]

\qquad \un{Step 2:} \quad
$1 < p < n$, $f \in C_c^1 (\R^n)$.
\\[-.25cm]

[Put
\[
\gamma 
\ = \ 
\frac{p (n - 1)}{n - p}.
\]
Then $\gamma > 1$ and 
\[
\frac{\gamma n}{n - 1} 
\ = \ 
p^*
\ = \ 
\frac{n p}{n - p} 
\ = \ 
\frac{(\gamma - 1) p}{p-1}.
\]
\\[-.75cm]

\noindent
Now apply Step 1 to $\abs{f^\gamma} = \abs{f}^\gamma$: 
\\[-.75cm]

\allowdisplaybreaks
\begin{align*}
\bigg(\ 
\int\limits_{\R^n} \ 
\abs{f}^{p^*} \ \td \Lm^n 
\bigg)^{(n-1)/n}\ 
&=\ 
\bigg(\ 
\int\limits_{\R^n} \ 
\abs{f^\gamma}^{n/(n-1)} \ \td \Lm^n 
\bigg)^{(n-1)/n}\ 
\\[15pt]
&\leq \ 
\int\limits_{\R^n} \ 
\norm{ \nabla \abs{f^\gamma}\hsy} \ 
\td \Lm^n
\\[15pt]
&=\ 
\gamma \ 
\int\limits_{\R^n} \ 
\abs{f}^{\gamma - 1} \norm{ \nabla f} \ 
\td \Lm^n
\\[15pt]
&\leq \ 
\gamma \ 
\bigg(\ 
\int\limits_{\R^n} \ 
\abs{f}^{(\gamma - 1) p / (p-1)}
\ \td \Lm^n
\bigg)^{(p-1)/p}
\ 
\bigg(\ 
\int\limits_{\R^n} \
\norm{ \nabla f}^p 
\ \td \Lm^n
\bigg)^{1/p}
\\[15pt]
&= \ 
\gamma \ 
\bigg(\ 
\int\limits_{\R^n} \ 
\abs{f}^{p^*}
\ \td \Lm^n
\bigg)^{(p-1)/p}
\ 
\bigg(\ 
\int\limits_{\R^n} \
\norm{ \nabla f}^p 
\ \td \Lm^n
\bigg)^{1/p}.
\end{align*}
And
\begin{align*}
\frac{n-1}{n} - \frac{p-1}{p} \ 
&=\ 
\frac{p (n - 1) - n (p - 1)}{p}
\\[11pt]
&=\ 
\frac{n - p}{n p}
\\[11pt]
&=\ 
\frac{1}{p^*}
\end{align*}

$\implies$
\[
\bigg(\ 
\int\limits_{\R^n} \ 
\abs{f}^{p^*} \ \td \Lm^n 
\bigg)^{1/p^*}\ 
\ \leq \ 
\gamma \ 
\bigg(\ 
\int\limits_{\R^n} \
\norm{ \nabla f}^p 
\ \td \Lm^n
\bigg)^{1/p},
\]
\\[-.5cm]

\noindent
where $C(n,p) = \gamma$.]
\\

\qquad \un{Step 3:} \quad
$1 \leq p < n$, $f \in W^{1, p} (\R^n)$.
\\[-.25cm]

[Given an $f \in W^{1, p} (\R^n)$, there exists a sequence $\{f_k\} \subset C_c^1(\R^n)$ such that
\[
\norm{f_k - f}_{W^{1, p}}
\ra 
0 \quad (k \ra \infty).
\]
So
\[
\norm{f_k - f}_{\Lp^p}
\ra 
0 \quad (k \ra \infty) 
\]
and, upon passing to a subsequence if necessary, it can be assumed that $f_k \ra f$ almost everywhere in $\R^n$.  
This said, we then claim that $\{f_k\}$ is a Cauchy sequence in $\Lp^{p^*} (\R^n)$.  
For $f_k - f_\ell \in C_c^1 (\R^n)$, thus it follows that
\allowdisplaybreaks
\begin{align*}
\norm{f_k - f_\ell}_{\Lp^{p^*}} \ 
&\leq \ 
C(n,p) \hsx \norm{\nabla f_k - \nabla f_\ell}_{\Lp^p}
\\[15pt]
&\leq \ 
C(n,p) \hsx 
\left(
\norm{\nabla f_k - \nabla f}_{\Lp^p}
+ 
\norm{\nabla f - \nabla f_\ell}_{\Lp^p}
\right)
\\[15pt]
&\ra 
0 \quad (k, \ell \ra \infty).
\end{align*}
Consequently there exists a $g \in \Lp^{p^*} (\R^n)$ such that $f_k \ra g$ in $\Lp^{p^*} (\R^n)$.  
Therefore $f = g$ almost everywhere in $\R^n$, which implies that $f_k \ra f$ in $\Lp^{p^*} (\R^n)$.  
Finally
\allowdisplaybreaks
\begin{align*}
\norm{f}_{\Lp^{p^*}} \ 
&\leq \ 
\norm{f - f_k}_{\Lp^{p^*}} + \norm{f_k}_{\Lp^{p^*}}
\\[15pt]
&\leq \ 
\norm{f - f_k}_{\Lp^{p^*}} + C(n,p) \hsx  \norm{\nabla f_k}_{\Lp^p}
\\[15pt]
&\leq \ 
\norm{f - f_k}_{\Lp^{p^*}} 
+ C(n,p) \hsx 
\left(
\norm{\nabla f_k - \nabla f}_{\Lp^p} + \norm{\nabla f}_{\Lp^p}
\right)
\\[15pt]
&\ra \ 
0 \hsy + \hsy C(n,p) \hsx 
\left(
0+ \norm{\nabla f}_{\Lp^p}
\right).
\end{align*}

I.e.: 
\[
\norm{f}_{\Lp^{p^*}}
\ \leq \ 
C(n,p) \hsx 
\left(
\norm{\nabla f}_{\Lp^p}
\right).]
\]
\\[-.75cm]

\qquad{\bf 7.2.5:}
{\small\bf APPLICATION} \ 
If $f \in W^{1, p} (\R^n)$ with $1 \leq p < n$ and if $\nabla f = 0$ almost everywhere in $\R^n$, 
then $f = 0$ almost everywhere in $\R^n$.
\\

\qquad{\bf 7.2.6:}
{\small\bf RAPPEL} \ 
If $1 \leq p < q < r < +\infty$, then

\[
\Lp^p \cap \Lp^r \subset \Lp^q
\]
and

\[
\norm{f}_q 
\ \leq \ 
\left(
\norm{f}_p 
\right)^\lambda 
\hsx 
\left(
\norm{f}_r 
\right)^{1 - \lambda}, 
\]
where
\[
\frac{1}{q}
\ = \ 
\frac{\lambda}{p} + \frac{1 - \lambda}{r} 
\qquad (0 < \lambda < 1).
\]
\\[-.75cm]

\qquad{\bf 7.2.7:}
{\small\bf RAPPEL} \ 
If $a \geq 0$, $b \geq 0$ and if $0 < \lambda < 1$, then
\[
a^\lambda \hsy b^{1 - \lambda}
\ \leq \ 
\lambda \hsy a + (1 - \lambda) \hsy b.
\]

Therefore
\begin{align*}
\norm{f}_q \ 
&\leq \ 
\lambda \hsy \norm{f}_p + (1 - \lambda) \hsy \norm{f}_r
\\[11pt]
&< \  
\norm{f}_p + \norm{f}_r.
\end{align*}
\\[-1cm]

Specialize now and take $r = p^*$ (recall that $p < p^*$) and let $p < q < p^*$ $-$then it follows that 
\[
W^{1, p} (\R^n) \hookrightarrow \Lp^q (\R^n).
\]

\chapter{
$\boldsymbol{\S}$\textbf{7.3}.\quad  EMBEDDINGS: BMO}
\setlength\parindent{2em}
\renewcommand{\thepage}{7-\S3-\arabic{page}}

\qquad 
Having dealt with the case when $1 \leq p < n$, the next item on the agenda is the case when $p = n$, which necessitates some preparation.
\\

\qquad{\bf 7.3.1:}
{\small\bf DEFINITION} \ 
The set 
\[
Q 
\ = \ 
[a_1, b_1] \times \ldots \times [a_n, b_n], 
\quad 
b_1 - a_1 = \cdots = b_n - a_n 
\]
is a \un{cube} in $\R^n$ (if $n = 1$, a cube is a bounded closed interval in $\R$, 
if $n = 2$, a cube is a square in $\R^2$ etc.).
The \un{side length} $\ell(Q)$ of $Q$ is the common value
\[
b_i - a_i 
\quad 
(i = 1, \ldots, n).
\]
\\[-1.25cm]

\qquad{\bf 7.3.2:}
{\small\bf NOTATION} \ 
\[
Q(x, \ell) 
\ = \ 
\{y \in \R^n \hsy : \hsy \abs{y_i - x_i} \leq \frac{\ell}{2} \ (i = 1, \ldots, n)\}
\]
is a cube with center $x$ and side length $\ell$.  
Here
\[
\Lm^n(Q(x, \ell))
\ = \ 
\ell^n 
\quad \text{and} \quad 
\diam  Q(x, \ell) 
\ = \ 
\sqrt{n} \hsx \ell.
\]
\\[-1.cm]

\qquad{\bf 7.3.3:}
{\small\bf DEFINITION} \ 
Given $f \in \Lp_\locx^1 (\R^n)$, its \un{integral average} over the cube $Q(x, \ell)$ is the entity
\[
f_{Q(x, \ell)} 
\ = \ 
\frac{1}{\ell^n} \ 
\int\limits_{Q(x, \ell)} \ 
f 
\ \td \Lm^n.
\]
\\[-1.cm]

\qquad{\bf 7.3.4:}
{\small\bf LEMMA} \ 
Let $1 \leq p < +\infty$ $-$then there exists a constant $C(n, p) > 0$ such that for all $f \in W^{1, p} (\R^n)$, 
\[
\bigg(
\int\limits_{Q(x, \ell)} \ 
\abs{f - f_{Q(x, \ell)}}^p 
\ \td \Lm^n
\bigg)^{1/p}
\ \leq \ 
C(n, p) \hsx \ell \hsx 
\bigg(
\int\limits_{Q(x, \ell)} \ 
\norm{\nabla f}^p
\ \td \Lm^n
\bigg)^{1/p}
\ < \ 
+\infty.
\]

PROOF \ 
Take $f \in C_c^1 (\R^n)$, put $Q = Q(x, \ell)$, let $z$, $y \in Q$, and write
\allowdisplaybreaks
\begin{align*}
\abs{f(z) - f(y)} \ 
&\leq \ 
\abs{f(z) - f(z_1, \ldots, z_{n-1}, y_n)} 
+ \cdots + 
\abs{f(z_1, y_2, \ldots, y_n) - f(y)}
\\[11pt]
&\leq \ 
\sum\limits_{i = 1}^n \ 
\int\limits_{a_i}^{b_i} \ 
\norm{\nabla f(z_1, \ldots, z_{i - 1}, t, y_{i+1}, \ldots, y_n)} 
\ \td t
\end{align*}

$\implies$
\\

\quad 
$\abs{f(z) - f(y)}^p$
\allowdisplaybreaks
\begin{align*}
\hspace{1.5cm}
&\leq \ 
\bigg(
\sum\limits_{i = 1}^n \ 
\int\limits_{a_i}^{b_i} \ 
\norm{\nabla f(z_1, \ldots, z_{i - 1}, t, y_{i+1}, \ldots, y_n)} 
\ \td t
\bigg)^p
\\[15pt]
&\leq \ 
\bigg(
\sum\limits_{i = 1}^n \ 
\bigg(
\int\limits_{a_i}^{b_i} \ 
\norm{\nabla f(z_1, \ldots, z_{i - 1}, t, y_{i+1}, \ldots, y_n)}^p 
\ \td t
\bigg)^{1/p}
(b_i - a_i)^{1 - 1/p}
\bigg)^p
\\[15pt]
&\leq \ 
n^p \hsy \ell^{p-1} \ 
\sum\limits_{i = 1}^n \ 
\int\limits_{a_i}^{b_i} \ 
\norm{\nabla f(z_1, \ldots, z_{i - 1}, t, y_{i+1}, \ldots, y_n)}^p 
\ \td t
\end{align*}

$\implies$

\allowdisplaybreaks
\begin{align*}
\int\limits_Q \ 
\abs{f - f_Q}^p
\td \Lm^n
&=\ 
\int\limits_Q \ 
\abs{f(z) - f_Q}^p
\td z
\\[15pt]
&=\ 
\int\limits_Q \ 
\bigg| \hsx 
\frac{1}{\ell^n} \ 
\int\limits_Q \ 
(f(z) - f(y)) 
\ \td y \hsx 
\bigg|^p
\td z
\\[15pt]
&\leq \ 
\int\limits_Q \ 
\bigg( \hsx 
\frac{1}{\ell^n} \ 
\int\limits_Q \ 
\abs{f(z) - f(y)} 
\ \td y
\bigg)^p
\td z
\\[15pt]
&\leq \ 
\int\limits_Q \ 
\frac{1}{\ell^n} \ 
\int\limits_Q \ 
\abs{f(z) - f(y)}^p 
\ \td z \hsx \td y
\\[15pt]
&\leq \ 
\frac{n^p \ell^{p-1}}{\ell^n} \ 
\sum\limits_{i = 1}^n \ 
\int\limits_Q \ 
\int\limits_Q \ 
\int\limits_{a_i}^{b_i} \ 
\norm{\nabla f(z_1, \ldots, z_{i - 1}, t, y_{i+1}, \ldots, y_n)}^p 
\ \td t \hsx \td y \hsx \td z 
\\[15pt]
&\leq \ 
\frac{n^p \ell^{p-1}}{\ell^n} \hsx
\sum\limits_{i = 1}^n \ 
(b_i - a_i) \ 
\int\limits_Q \ 
\int\limits_Q \ 
\norm{\nabla f (z)}^p 
\td z \hsx \td y
\\[15pt]
&\leq \ 
\frac{n^p \ell^{p-1}}{\ell^n} \hsx \cdot \hsx 
n \hsy \ell^{n+1} \ 
\int\limits_Q \ 
\norm{\nabla f}^p 
\ \td \Lm^n
\\[15pt]
&= \ 
n^{p+1} \hsx \ell^p \ 
\int\limits_Q \ 
\norm{\nabla f}^p 
\ \td \Lm^n
\end{align*}

$\implies$
\[
\bigg(
\int\limits_Q \ 
\abs{f - f_Q}^p
\ \td \Lm^n
\bigg)^{1/p}
\ \leq \ 
C(n, p) \hsx \ell \ 
\bigg(
\int\limits_Q \ 
\norm{\nabla f}^p 
\ \td \Lm^n
\bigg)^{1/n}
\]
if 
\[
C(n, p) 
\ = \ 
(n^{p+1})^{1/p}.
\]
\\[-.75cm]

\qquad{\bf 7.3.5:}
{\small\bf SCHOLIUM} \ 
If $f \in W^{1, n} (\R^n)$, then for every cube $Q$, 
\\[-.75cm]

\begin{align*}
\frac{1}{\Lm^n (Q)} \
\int\limits_Q \
\abs{f - f_Q}
\ \td \Lm^n \ 
&\leq \ 
\bigg(
\frac{1}{\Lm^n (Q)} \ 
\int\limits_Q \
\abs{f - f_Q}^n
\ \td \Lm^n
\bigg)^{1/n}
\\[15pt]
&\leq \ 
C(n) \hsx 
\ell (Q) \hsx 
\bigg(
\frac{1}{\Lm^n (Q)} \ 
\int\limits_Q \ 
\norm{\nabla f}^n
\ \td \Lm^n
\bigg)^{1/n}
\\[15pt]
&= \ 
C(n) \hsx 
\ell (Q) \hsx 
\bigg(
\frac{1}{\ell(Q)^n} \ 
\int\limits_Q \ 
\norm{\nabla f}^n
\ \td \Lm^n
\bigg)^{1/n}
\\[15pt]
&= \ 
C(n) \hsx 
\bigg(
\int\limits_Q \ 
\norm{\nabla f}^n
\ \td \Lm^n
\bigg)^{1/n}
\\[15pt]
&\leq \ 
C(n) \hsx \norm{\nabla f}_{\Lm^n (\R^n)}
\\[15pt]
&< +\infty.
\end{align*}
\\[-1cm]

\qquad{\bf 7.3.6:}
{\small\bf DEFINITION} \ 
A function $f \in \Lp_\locx^1 (\R^n)$ is of \un{bounded mean oscillation} provided
\[
\norm{f}_\BMO 
\ \equiv \ 
\sup\limits_Q \ 
\frac{1}{\Lm^n(Q)} \ 
\int\limits_Q \ 
\abs{f - f_Q} 
\ \td \Lm^n
\ < \ 
+\infty, 
\]
where the supremum is taken over all cubes $Q$ in $\R^n$.
\\

\qquad{\bf 7.3.7:}
{\small\bf NOTATION} \ 
$\BMO (\R^n)$ is the set of functions of bounded mean oscillation.
\\

\qquad{\bf 7.3.8:}
{\small\bf \un{N.B.}} \ 
$\norm{\hsx \cdot \hsx} _\BMO$ is a seminorm, not a norm (constant functions have vanishing bounded mean oscillation).
\\

\qquad{\bf 7.3.9:}
{\small\bf LEMMA} \ 
$\BMO (\R^n)$ is a vector space over $\R$.
\\[-.5cm]

[If $f$, $g \in \BMO (\R^n)$, then
\[
\norm{f + g}_\BMO 
\ \leq \ 
\norm{f}_\BMO + \norm{g}_\BMO
\]
and $\norm{\hsx \cdot \hsx} _\BMO$ is scale invariant, i.e., $\forall \ r \in \R$, 
\[
\norm{f (r \hsy \cdot)}_\BMO 
\ = \ 
\norm{f}_\BMO.]
\]
\\[-1.25cm]

\qquad{\bf 7.3.10:}
{\small\bf THEOREM} \ 
\[
\BMO (\R^n)  / \R
\]
is a Banach space.
\\

\qquad{\bf 7.3.11:}
{\small\bf LEMMA}  \ 
\[
f \in \BMO (\R^n) 
\implies 
\abs{f} \in \BMO (\R^n).
\]
\\[-1.25cm]

\qquad{\bf 7.3.12:}
{\small\bf THEOREM} \ 
$\Lp^\infty (\R^n)$ is contained in $\BMO (\R^n)$.
\\[-.5cm]

[If $f \in \Lp^\infty (\R^n)$, then
\[
\norm{f}_\BMO 
\ \leq \ 
2 \hsx \norm{f}_{\Lp^\infty}.]
\]
\\[-1.25cm]

[Note: \ 
Therefore
\[
\Lp^\infty (\R^n) \hookrightarrow \BMO (\R^n).]
\]
\\[-1.25cm]

\qquad{\bf 7.3.13:}
{\small\bf \un{N.B.}} \ 
The containment is strict.
\\[-.5cm]

[The unbounded function 
\[
\elog \norm{x}
\]
belongs to $\BMO (\R^n)$.]
\\

\qquad{\bf 7.3.14:}
{\small\bf EXAMPLE} \ 
Take $n = 1$ $-$then the function 
\[
f(x) \ = \ 
\begin{cases}
\ \elog x \hspace{0.85cm} (x > 0)
\\[4pt]
\ \ 0 \hspace{1.5cm} (x \leq 0)
\end{cases}
\]
is not of bounded mean oscillation.
\\

Let $f \in W^{1, n} (\R^n)$ $-$then, as has been seen above, for every cube $Q$, 
\[
\frac{1}{\Lm^n (Q)} \ 
\int\limits_Q \ 
\abs{f - f_Q} 
\ \td \Lm^n 
\ \leq \ 
C(n) \hsy
\norm{\nabla f}_{\Lm^n (\R^n)},
\]
so upon taking the supremum over $Q$, it follows that $f \in \BMO (\R^n)$, where
\[
\norm{f}_\BMO
\ \leq \ 
C(n) \hsy 
\norm{\nabla f}_{\Lm^n (\R^n)}.
\]
\\[-1cm]

\qquad{\bf 7.3.15:}
{\small\bf SCHOLIUM} \ 
\[
W^{1, n} (\R^n) \hookrightarrow \BMO (\R^n).
\]
\\[-1.25cm]

\qquad{\bf 7.3.16:}
{\small\bf APPLICATION} \ 
If $f \in W^{1, n} (\R^n)$ and if $\nabla f = 0$ almost everywhere in $\R^n$, then $f$ = some constant almost everywhere in $\R^n$.

\chapter{
$\boldsymbol{\S}$\textbf{7.4}.\quad  EMBEDDINGS: MOR}
\setlength\parindent{2em}
\setcounter{theoremn}{0}
\renewcommand{\thepage}{7-\S4-\arabic{page}}

\qquad 
It remains to consider the situation when $p > n$.
\\

\qquad{\bf 7.4.1:}
{\small\bf RAPPEL} \ 
Let $E$ be a nonempty subset of $\R^n$ $-$then a function $f:E \ra \R$ is 
\un{H\"older continuous} with exponent $0 < \alpha \leq 1$ if there is a constant $C > 0$ such that 

\[
\abs{f(x) - f(y)} 
\ \leq \ 
C \hsx \norm{x - y}^\alpha
\]
for all $x$, $y \in E$.
\\[-.5cm]

[Note: \ 
Of course if $\alpha = 1$, then it is a question of Lipschitz continuous.]
\\

A H\"older continuous function is continuous but it can be nowhere differentiable.
\\

\qquad{\bf 7.4.2:}
{\small\bf NOTATION} \ 
$C^{0, \alpha} (E)$ is the set of all bounded functions that are H\"older continuous with exponent $\alpha$ and norm

\[
\norm{f}_{C^{0, \alpha}} 
\ = \ 
\sup\limits_{x \in E} \ \abs{f(x)} 
+ 
\sup\limits_{\substack{x, y \in E \\x \neq y}} \ \frac{\abs{f(x) - f(y)}}{\norm{x - y}^\alpha}.
\]
\\[-1cm]

[Note: \ 
When so equipped, $C^{0, \alpha} (E)$ is a Banach space.]
\\

\qquad{\bf 7.4.3:}
{\small\bf LEMMA} \ 
Let $p > n$ $-$then there is a constant $C(n, p) > 0$ such that for all $f \in W^{1,p} (\R^n) \cap C^\infty (\R^n)$ and 
for all $z$, $y \in \R^n$, 
\[
\abs{f(z) - f(y)}
\ \leq \ 
C(n, p) \hsx 
\norm{z - y}^{1 - n/p} \hsx \norm{\nabla f}_{\Lp^p (\R^n)}.
\]
\\[-1cm]

PROOF \ 
To begin with, 
\allowdisplaybreaks
\begin{align*}
f(z) - f(y) \ 
&=\ 
\int\limits_0^1 \ 
\frac{\partial }{\partial t} (f (tz + (1-t)y)) 
\ \td t
\\[15pt]
&=\ 
\int\limits_0^1 \ 
\langle
\nabla f (t z + (1 - t) y), z - y
\rangle
\ \td t.
\end{align*}

Assume now that $z$, $y\in Q(x, \ell)$ $-$then
\allowdisplaybreaks
\begin{align*}
\abs{f_{Q(x, \ell)} - f(y)} \ 
&=\
\bigg| \hsx
\frac{1}{\ell^n} \hsx 
\int\limits_{Q(x, \ell)} \ 
(f(z) - f(y))
\ \td z
\hsx \bigg|
\\[15pt]
&=\
\bigg| \hsx
\frac{1}{\ell^n} \hsx 
\int\limits_{Q(x, \ell)} \ 
\int\limits_0^1 \ 
\langle
\nabla f (t z + (1 - t) y), z - y
\rangle
\ \td t \hsx \td z
\hsx \bigg|
\\[15pt]
&\leq\
\sum\limits_{i = 1}^n \ 
\frac{1}{\ell^n} \hsx 
\int\limits_{Q(x, \ell)} \ 
\int\limits_0^1 \ 
\abs{\frac{\partial f}{\partial x_i} (tz + (1 - t) y)}
\abs{z_i - y_i}
\ \td t \hsx \td z
\\[15pt]
&\leq\
\sum\limits_{i = 1}^n \ 
\frac{1}{\ell^{n-1}} \hsx 
\int\limits_0^1 \ 
\int\limits_{Q(x, \ell)} \ 
\abs{\frac{\partial f}{\partial x_i} (tz + (1 - t) y)}
\ \td t \hsx \td z
\\[15pt]
&=\
\sum\limits_{i = 1}^n \ 
\frac{1}{\ell^n} \hsx 
\int\limits_0^1 \ 
\frac{1}{t^n} \ 
\int\limits_{Q(tx + (1 - t)y, t \ell)} \
\abs{\frac{\partial f}{\partial x_i} (w)}
\ \td w \hsx \td t
\\[15pt]
&\leq\
\sum\limits_{i = 1}^n \ 
\frac{1}{\ell^{n-1}} \hsx 
\int\limits_0^1 \ 
\frac{1}{t^n} \ 
\bigg(\ 
\int\limits_{Q(tx + (1 - t)y, t \ell)} \
\abs{\frac{\partial f}{\partial x_i} (w)}^p 
\ \td w
\bigg)^{1/p}
\\[15pt]
&\hspace{3cm} \times
\bigg(
\Lm^n (Q(t x + (1 - t)y, t \ell)
\bigg)^{1 - 1/p} 
\ \td t
\\[15pt]
&\leq\
n \hsx 
\norm{\nabla f}_{\Lp^p(Q(x, \ell))} \ 
\frac{\ell^{n(1 - 1/p)}}{\ell^{n-1}} \ 
\int\limits_0^1 \ 
\frac{t^{n(1 - 1/p)}}{t^n} \ 
\ \td t
\\[15pt]
&=\
\frac{n p }{p - n} \ 
\ell^{1 - n/p} \ 
\norm{\nabla f}_{\Lp^p(Q(x, \ell))}.
\end{align*}

Since the same estimate obtains if the roles of $y$ and $z$ are interchanged, write

\allowdisplaybreaks
\begin{align*}
\abs{f(z) - f(y)} \ 
&\leq \ 
\abs{f(z) - f_{Q(x, \ell)}} + \abs{f_{Q(x, \ell)} - f(y)}
\\[11pt]
&\leq \ 
2 \hsx 
\frac{n p }{p - n} \ 
\ell^{1 - n/p} \ 
\norm{\nabla f}_{\Lp^p(Q(x, \ell))},
\end{align*}
where $z$, $y \in Q(x, \ell)$.  
Proceed finally to when $z$, $y \in \R^n$ are arbitrary $-$then there exists a cube $Q(x, \ell)$ such that $z$, $y \in Q(x, \ell)$ 
and $\ell = \norm{z - y}$ (e.g., take $\ds x = \frac{z + y}{2}$), hence

\allowdisplaybreaks
\begin{align*}
\abs{f(z) - f(y)} \ 
&\leq \ 
C(n, p) \hsx 
\norm{z - y}^{1 - n/p} \hsx
\norm{\nabla f}_{\Lp^p(Q(x, \ell))}
\\[11pt]
&\leq \ 
C(n, p) \hsx 
\norm{z - y}^{1 - n/p} \hsx
\norm{\nabla f}_{\Lp^p(\R^n)},
\end{align*}
where

\[
C(n, p) 
\ = \ 
2 \hsx 
\frac{n p }{p - n}.
\]
\\[-.75cm]

In the foregoing, replace $z$, $y$ by $x$, $y$.
\\

\qquad{\bf 7.4.4:}
{\small\bf THEOREM} \ 
Let $f \in W^{1,p} (\R^n)$  $(p > n)$ and let $x$, $y \in \Lambda (f)$ $-$then

\[
\abs{f(x) - f(y)} 
\ \leq \ 
C(n, p) \hsx \norm{x - y}^{1 - n/p} \hsx \norm{\nabla f}_{\Lp^p (\R^n)}.
\]
\\[-1cm]

[Utilize the standard mollification $f_\varepsilon$ of $f$ and apply it to Lebesgue points $x$, $y$ of $f$.]
\\

The restriction $\restr{f}{\Lambda (f)}$ can be extended uniquely to $\R^n$ as a H\"older continuous function $\bar{f}$ of exponent 
$1 - n/p$ in such a way that 
\[
\abs{\bar{f}(x) - \bar{f}(y)} 
\ \leq \ 
C(n, p) \hsx \norm{x - y}^{1 - n/p} \hsx \norm{\nabla \bar{f}}_{\Lp^p (\R^n)}
\]
for all $x$, $y \in \R^n$.
\\

[Bearing in mind that $\Lambda (f)$ is dense, given an $x \in \R^n$, choose a sequence $\{f_k\} \subset \Lambda (f)$ 
such that $x_k \ra x$ $(k \ra \infty)$.  
From what has been said above, $\{f(x_k)\}$ is 
a Cauchy sequence, thus the prescription
\[
\bar{f}(x)
\ = \ 
\lim\limits_{k \ra \infty} \ f(x_k)
\]
makes sense and has the desired property.]
\\

\qquad{\bf 7.4.5:}
{\small\bf THEOREM} \ 
Let $f \in W^{1,p} (\R^n)$  $(p > n)$ $-$then in the equivalence class of $f$ there is a unique function $\bar{f}$ 
which is H\"older continuous with exponent $1 - n/p$. 
\\

In particular: Every element of $W^{1,p} (\R^n)$  $(p > n)$ coincides with a continuous function almost everywhere.
\\

\qquad{\bf 7.4.6:}
{\small\bf LEMMA} \ 
If $f \in W^{1,p} (\R^n)$ with $p > n$, then $f$ is essentially bounded. 
\\

PROOF \ 
Let $y \in Q(x, 1)$, take $f \in W^{1,p} (\R^n) \cap C^\infty (\R^n)$, and write

\allowdisplaybreaks
\begin{align*}
\abs{f(x)} \ 
&\leq \ 
\abs{f(z) - f_{Q(x, 1)}} + \abs{f_{Q(x, 1)}}
\\[15pt]
&\leq \ 
\frac{n p}{p - n} \hsx 
1^{1 - n/p} \hsx 
\norm{\nabla f}_{\Lp^p (Q(x, 1))} 
+
\int\limits_{Q(x, 1)} \ 
\abs{f(y)}\ \td y
\\[15pt]
&\leq \ 
\frac{n p}{p - n} \hsx 
\norm{\nabla f}_{\Lp^p (\R^n)} 
+
\bigg( \ 
\int\limits_{Q(x, 1)} \ 
\abs{f(y)}^p \ \td y
\bigg)^{1/p}
\\[15pt]
&\leq \ 
\frac{n p}{p - n} \hsx 
\norm{\nabla f}_{\Lp^p (\R^n)} 
+ 
\norm{f}_{\Lp^p (\R^n)}.
\end{align*}
But one choice for the norm of $f$ in $W^{1,p} (\R^n)$ is
\[
\norm{f}_{W^{1,p} (\R^n)}
\ = \ 
\norm{f}_{\Lp^p (\R^n)}
+ 
\norm{\nabla f}_{\Lp^p (\R^n)} .
\]
This said, there are then two possibilities.
\\

\qquad \textbullet \quad 
$\ds \frac{n p}{p - n}  \geq 1$
\\

$\implies$
\allowdisplaybreaks
\begin{align*}
\abs{f(z)} \ 
&\leq \ 
\frac{n p}{p - n} \hsx 
\norm{\nabla f}_{\Lp^p (\R^n)} 
+ 
\norm{f}_{\Lp^p (\R^n)}
\\[15pt]
&\leq \
\frac{n p}{p - n} \hsx 
\norm{\nabla f}_{\Lp^p (\R^n)} 
+ 
\frac{n p}{p - n} \hsx 
\norm{f}_{\Lp^p (\R^n)}
\\[15pt]
&= \
\frac{n p}{p - n} \hsx 
\bigg(
\norm{f}_{W^{1,p} (\R^n)}
\bigg).
\end{align*}
\\

\qquad \textbullet \quad 
$\ds \frac{n p}{p - n}  < 1$
\\

$\implies$
\allowdisplaybreaks
\begin{align*}
\abs{f(z)} \ 
&\leq \ 
\norm{\nabla f}_{\Lp^p (\R^n)} 
+ 
\norm{f}_{\Lp^p (\R^n)}
\\[15pt]
&= \
\norm{f}_{W^{1,p} (\R^n)}.
\end{align*}
\\[-1cm]

\noindent
Therefore the $\Lp^\infty$-norm of $f$ is bounded by a constant
\[
\bar{C} (n, p) \ = \ 
\begin{cases}
\ \ds\frac{n p}{p - n} \hspace{0.75cm} \text{(if $\geq 1$)}
\\[4pt]
\ \ 1 \hspace{1.4cm} \text{(if $< 1$)}
\end{cases}
\]
depending on $n$ and $p$ times the $W^{1, p}$-norm of $f$.
\\

\qquad{\bf 7.4.7:}
{\small\bf SCHOLIUM} \ 
When $p > n$, 

\[
W^{1,p} (\R^n)
\hookrightarrow 
\Lp^\infty (\R^n).
\]

\qquad{\bf 7.4.8:}
{\small\bf THEOREM} \ 
\[
\bar{f} \in C^{0, 1 - n/p} (\R^n) 
\quad (p > n).
\]

PROOF \ 

\allowdisplaybreaks
\begin{align*}
\norm{\bar{f}}_{C^{0, 1 - n/p}} \ 
&=\ 
\sup\limits_{x \in \R^n} \ \abs{\bar{f} (x)} 
+ 
\sup\limits_{\substack{x, y \in \R^n \\x \neq y}} \ \frac{\abs{\bar{f}(x) - \bar{f}(y)}}{\norm{x - y}^{1 - n/p}}
\\[15pt]
&\leq \ 
\bar{C}(n, p) \hsx \norm{f}_{W^{1,p}} 
+ 
C(n, p) \norm{f}_{W^{1,p}}. 
\end{align*}
\\[-1cm]

\qquad{\bf 7.4.9:}
{\small\bf \un{N.B.}} \ 
If $p > n$, then this state of affairs is symbolized by writing

\[
W^{1,p}(\R^n) \hookrightarrow C^{0, 1 - n/p} (\R^n) 
\]
since

\[
\norm{\bar{f}}_{C^{0, 1 - n/p}}
\ \leq \ 
C \hsx \norm{f}_{W^{1,p}}.
\]
\\[-.75cm]

\qquad{\bf 7.4.10:}
{\small\bf THEOREM} \ 
$(p > n)$
\[
\bar{f} (x) \ra 0 
\ \text{as} \ 
\norm{x} \ra +\infty.
\]
\\[-1cm]

PROOF \ 
Given $f \in W^{1,p} (\R^n)$, choose a sequence $\{f_k\}$ in $C_c^\infty (\R^n)$ that converges to $f$ in $W^{1,p} (\R^n)$ 
$-$then
\[
\norm{f - f_k}_{\Lp^\infty} 
\ \leq \ 
C \hsx \norm{f - f_k}_{W^{1,p}}
\ra 0 
\quad (k \ra \infty).
\]
Fix $\varepsilon > 0$ and choose $\bar{k}$ such that 

\[
\norm{f - f_k}_{\Lp^\infty} 
\ \leq \ 
\varepsilon
\quad(k \geq \bar{k}).
\]
Next choose $R_{\bar{k}} > 0$ such that $f_{\bar{k}} (x) = 0$ for all $x$ : $\norm{x} \geq R_{\bar{k}}$.  
So, for $\Lm^n$ almost every $x \in \R^n$ with $\norm{x} \geq R_{\bar{k}}$, it follows that
\[
\abs{\bar{f} (x)} 
\ = \ 
\abs{\bar{f} (x) - f_{\bar{k}} (x) } 
\ \leq \ 
\norm{f - f_k}_{\Lp^\infty} 
\ \leq \ 
\varepsilon
\]
and since $\bar{f}$ is continuous, this inequality holds for all $\ds x : \norm{x} \geq R_{\bar{k}}$.
\\

\qquad{\bf 7.4.11:}
{\small\bf LEMMA} \ 
$(p > n)$ $\forall \ x, \hsx y \in \R^n$, 
\[
\abs{\bar{f} (x) - \bar{f} (y)} 
\ \leq \ 
C(n, p) \hsx 
\norm{x - y}^{1 - n/p} \ 
\bigg(\ 
\int\limits_{B(x, \norm{x - y})} \ 
\norm{\nabla \bar{f}}^p 
\ \td \Lm^n
\bigg)^{1/p}.
\]
\\[-.75cm]

\qquad{\bf 7.4.12:}
{\small\bf THEOREM} \ 
$(p > n)$ $\bar{f}$ is differentiable almost everywhere.
\\[-.5cm]

PROOF \ 
It suffices to show that $\bar{f}$ is differentiable at every $\Lp^p$-Lebesgue point $x_0$ of $\nabla f$, where by definition

\[
\lim\limits_{r \ra 0} \ 
\frac{1}{\omega_n^{r^n}} \ 
\int\limits_{B(x, r)} \ 
\norm{\nabla f - \nabla f(x_0)}^p 
\ \td \Lm^n  
\ = \ 
0.
\]
To this end, note that 
\begin{align*}
\big|\bar{f} (x) - \bar{f} (x_0) &- \langle \nabla f(x_0), x - x_0 \rangle \big| \ 
\\[15pt]
&\leq \ 
C(n, p) \hsx 
\norm{x - x_0}^{1 - n/p} \ 
\bigg(\ 
\int\limits_{B(x, \norm{x - x_0})} \ 
\norm{\nabla f - \nabla f(x_0)}^p 
\ \td \Lm^n
\bigg)^{1/p}.
\end{align*}
And
\allowdisplaybreaks
\begin{align*}
C(n, p) \hsx \norm{x - x_0}^{1 - n/p} \ 
&=\ 
C(n, p) \hsx 
\norm{x - x_0} \hsx 
\left(
\frac{1}{\norm{x - x_0}^n}
\right)^{1/p}
\\[15pt]
&=\ 
C(n, p) \hsx 
(\omega_n)^{1/p} \hsx 
\norm{x - x_0} \hsx
\left(
\frac{1}{\omega_n \norm{x - x_0}^n}
\right)^{1/p}
\end{align*}

$\implies$
\\

\hspace{1cm}
$
\ds
\frac{\abs{\bar{f} (x) - \bar{f} (x_0) - \langle \nabla f(x_0), x - x_0 \rangle}}{\norm{x - x_0}}
$
\\[-.5cm]

\allowdisplaybreaks
\begin{align*}
&\leq \ 
C(n, p) \hsx 
(\omega_n)^{1/p} \hsx 
\bigg(
\frac{1}{\omega_n \norm{x - x_0}^n} \ 
\int\limits_{B(x, \norm{x - x_0})} \ 
\norm{\nabla f - \nabla f(x_0)}^p 
\ \td \Lm^n
\bigg)^{1/p}
\\[11pt]
&\ra 0
\qquad (\text{as $x \ra x_0$}).
\end{align*}
\\[-1cm]


\qquad{\bf 7.4.13:}
{\small\bf \un{N.B.}} \ 
The weak derivatives of $f$ coincide with the ordinary partial derivatives of $\bar{f}$ almost everywhere in $\R^n$.

\chapter{
SECTION \textbf{6}:\quad  ACL}
\setlength\parindent{2em}
\setcounter{chapter}{8}
\renewcommand{\thepage}{8-\arabic{page}}

\qquad
Working in $\R^n$, let
\[
I_k 
\ = \ 
[a_k, b_k]
\quad (k = 1, \ldots, n)
\]
be closed intervals and put
\[
Q \ = \ 
[a_1, b_1] \times \ldots \times [a_n, b_n],
\]
a rectangular box.
\\

\qquad{\bf 8.1:}\  
{\small\bf DEFINITION} \ 
A function $f: Q \ra \R$ is said to be ACL (\un{absolutely} \un{continuous on lines}) if for each 
$k = 1, \ldots, n$, and almost every point
\[
(x_1, \ldots, x_{k - 1}, x_{k + 1}, \ldots, x_n)
\in 
I_1 \times \cdots \times I_{k - 1} \times I_{k + 1} \times \cdots \times I_n
\subset \R^{n - 1} 
\]
with respect to $\Lm^{n-1}$ measure, the function 
\[
x_k \ra 
f(x_1, \ldots, x_{k - 1}, x_{k + 1}, \ldots, x_n)
\qquad (a_k \leq x_k \leq b_k)
\]
is absolutely continuous.   
\\

\qquad{\bf 8.2:}\  
{\small\bf \un{N.B.}} \ 
In the literature, the foregoing situation is sometimes referred to as saying that $f$ is absolutely continuous on 
almost every line segment in $Q$ parallel to the coordinate axes.
\\

Let $\Omega$ be a nonempty open subset of $\R^n$.
\\

\qquad{\bf 8.3:}\  
{\small\bf DEFINITION} \ 
A function $f: \Omega \ra \R$ is ACL if the restriction $\restr{f}{Q}$ is ACL for every $Q \subset \Omega$.
\\

\qquad{\bf 8.4:}\  
{\small\bf NOTATION} \ 
$\ACL(\Omega)$ is the set of ACL functions in $\Omega$.
\\


\qquad{\bf 8.5:}\  
{\small\bf EXAMPLE} \ 
A quasiconformal map belongs to $\ACL(\Omega)$.
\\

\begin{spacing}{1.5}
\qquad{\bf 8.6:}\
{\small\bf NOTATION} \ 
Let $1 \leq p < +\infty$ $-$then $\Lp^{1,p} (\Omega)$ consists of those $f \in \Lp_\locx^1 (\Omega)$ 
such that the distributional derivatives $\ds\frac{\partial f}{\partial x_i}$ are weak partial derivatives 
and also belong to $\Lp^p (\Omega)$ $(i = 1, \ldots, n)$.
\\[-.5cm]
\end{spacing}

[Note: \ 
Evidently

\[
\Lp^{1,p} (\Omega)
\subset
\Lp_\locx^p (\Omega)
\subset
\Lp_\locx^1 (\Omega).]
\]
\\[-1cm]

\qquad{\bf 8.7:}\
{\small\bf \un{N.B.}} \ 
Obviously

\[
\tW^{1, p} (\Omega)
\subset
\Lp^{1,p} (\Omega).
\]
\\[-1cm]

\qquad{\bf 8.8:}\
{\small\bf THEOREM} \ 
Let $1 \leq p < +\infty$ $-$then a function $f \in \Lp^{1,p} (\Omega)$ admits a representative 
$\ov{f} : \Omega \ra \R$ in $\ACL (\Omega)$.
\\[-.5cm]

[Note: \ 
The ordinary partial derivatives of $\bar{f}$ exist almost everywhere.]
\\

The proof in general is notationally involved so to simplify the bookkeeping, take $n = 2$, 
assume that $f$ is continuous, suppose that 

\[
Q = [0,1] \times [0,1] 
\subset 
\Omega 
\ = \ 
]-\varepsilon, 1 + \varepsilon [ \ \times \ ]-\varepsilon, 1 + \varepsilon [\hsy ,
\]
let $(x_1, x_2) = (x, y)$, thus the distributional derivatives are 
$\ds\frac{\partial f}{\partial x}$
and 
$\ds\frac{\partial f}{\partial y}$, 
and the claim is that 

\[
\begin{cases}
\ x \ra f(x,y) \ \text{is absolutely continuous for almost every $y \in [0,1]$}
\\[4pt]
\ y \ra f(x,y) \ \text{is absolutely continuous for almost every $x \in [0,1]$}
\end{cases}
.
\]

\noindent
The discussion in either case is conceptually the same, hence it will suffice to deal with the second of these.
\\

For use below: 
\\

\qquad{\bf 8.9:}\
{\small\bf CRITERION} \ 
If $f \in \Lp_\locx^1 ]a,b[$ and if $\forall \ \phi \in C_c^\infty ]a,b[$, 

\[
\int\limits_{]a,b[} \ \phi \hsy f  \ \td \Lm^1 
\ = \ 
0,
\]
then $f = 0$ almost everywhere in $]a,b[$.
\\

\qquad{\bf 8.10:}\
{\small\bf APPLICATION} \ 
If $f \in \Lp^1 ]0,1[$ and if $\forall \ \phi \in C_c^\infty ]0,1[$, 

\[
\int\limits_0^1 \ f \hsx \phi^\prime \ \td \Lm^1 
\ = \ 
0,
\]
then there exists a constant $C$ such that $f = C$ almost everywhere in $]0,1[$.
\\[-.5cm]

[Let $\Phi$, $\psi$ be functions in $C_c^\infty ]0,1[$ with $\ds\int\limits_0^1 \ \psi = 1$.  
Put
\[
\Psi (x) 
\ = \ 
\int\limits_0^x \ \Phi - \bigg(\int\limits_0^1 \ \Phi \bigg) \hsx \int\limits_0^x \ \psi.
\]
Then $\Psi (0) = 0$ and
\vspace{-.5cm}
\begin{align*}
\Psi(1) \ 
&=\ 
\int\limits_0^1 \ \Phi - \bigg(\int\limits_0^1 \ \Phi \bigg) \hsx \int\limits_0^1 \ \psi
\\[15pt]
&=\ 
0.
\end{align*}
Therefore $\Psi \in C_c^\infty ]0,1[$ and 
\vspace{-.5cm}
\[
\Psi^\prime 
\ = \ 
\Phi -  \bigg(\int\limits_0^1 \ \Phi \bigg) \hsx \psi
\]
\\[-1.5cm]

\hspace{1.cm} $\implies$
\allowdisplaybreaks
\begin{align*}
0 \ 
&=\ 
\int\limits_0^1 \ f \hsy \Psi^\prime 
\qquad (\text{by assumption})
\\[15pt]
&=\ 
\int\limits_0^1 \ \bigg(f  \hsx - \hsx \int\limits_0^1 \ f \hsy \psi\bigg) \hsx \Phi 
\end{align*}

\hspace{1.cm} $\implies$
\[
f - \int\limits_0^1 \ f \hsy \psi 
\ = \ 
0 
\qquad \text{almost everywhere in $]0,1[$}
\]

\hspace{1.cm} $\implies$
\[
f 
\ = \ 
\int\limits_0^1 \ f \hsy \psi 
\qquad \text{almost everywhere in $]0,1[$.}
\]
\\[-.75cm]

By hypothesis, 

\[
f \in \Lp^{1,p} (\Omega) 
\implies 
\restr{\frac{\partial f}{\partial y}}{Q} \in \Lp^p (Q), 
\]
from which it follows that 

\[
\int\limits_0^1 \ \abs{\frac{\partial f}{\partial y} (x, y)} \td y 
\ < \ 
+\infty
\]
for almost all $x \in [0,1]$ (Fubini).  
Therefore the function 

\[
y \ra g(x,y) 
\ = \ 
\int\limits_0^y \ \frac{\partial f}{\partial t} (x,t) \ \td t
\]
\begin{spacing}{1.5}
\noindent
is absolutely continuous on the segment $[0,1]$ for almost all $x \in [0,1]$ and its ordinary derivative 
$\ds g^\prime = \frac{\partial g}{\partial y}$ coincides with the distributional derivative 
$\ds\frac{\partial f}{\partial g}$ for almost all $y \in ]0,1[$.  
Consider now a test function $\phi$ of the form
\end{spacing}
\[
\phi 
\ = \ 
\xi \hsy\eta  
\qquad (\xi \in C_c^\infty (]0,1[, \ \eta \in C_c^\infty (]0,1[), 
\]
i.e., 
\[
\phi(x, y) 
\ = \ 
\xi (x) \hsy \eta(y).
\]
Then
\allowdisplaybreaks
\begin{align*}
\int\limits_Q \ f(x,y) \hsy \xi (x) \hsy \eta^\prime (y) \ \td x \hsy \td y \ 
&=\ 
- \int\limits_Q \ \xi (x) \hsy \eta (y) \hsy \frac{\partial f}{\partial y} \ \td x \hsy \td y 
\\[15pt]
&=\ 
- \int\limits_Q \ \phi(x, y) \hsy \frac{\partial g}{\partial y}  \ \td x \hsy \td y.
\end{align*}
On the other hand,
\[
\int\limits_0^1 \ g(x, y) \hsy \eta^\prime (y) \td y 
\ = \ 
- \int\limits_0^1 \ \eta(y) \hsy \frac{\partial g}{\partial y} (x,y) \  \td y
\]

$\implies$
\[
\xi (x)  \hsx
\int\limits_0^1 \ g(x,y) \eta^\prime (y) \ \td y 
\ = \ 
\xi(x)  \hsx 
\bigg(
- \int\limits_0^1 \ \eta(y) \hsy \frac{\partial g}{\partial y} (x,y) \ \td y
\bigg)
\]

$\implies$

\allowdisplaybreaks
\begin{align*}
\int\limits_Q \ g(x,y) \hsy \xi (x) \hsy \eta^\prime (y) \ \td x \hsy \td y \ 
&=\ 
- \int\limits_Q \ \xi (x) \hsy \eta (y) \hsy \frac{\partial g}{\partial y} \ \td x \hsy \td y 
\\[15pt]
&=\ 
- \int\limits_Q \ \phi(x, y) \hsy \frac{\partial g}{\partial y}  \ \td x \hsy \td y
\end{align*}

$\implies$

\[
\int\limits_Q \ f(x,y) \hsy \xi (x) \hsy \eta^\prime (y) \ \td x \hsy \td y
\ = \ 
\int\limits_Q \ g(x,y) \hsy \xi (x) \hsy \eta^\prime (y) \ \td x \hsy \td y
\]

$\implies$

\[
\int\limits_Q \ \xi (x) \hsy f(x,y) \hsy  \eta^\prime (y) \ \td x \hsy \td y
\ = \ 
\int\limits_Q \ \xi (x) \hsy g(x,y) \hsy  \eta^\prime (y) \ \td x \hsy \td y
\]

$\implies$

\[
\int\limits_0^1 \ [f(x,y) - g(x, y)] \eta^\prime (y) \ \td y
\ = \ 
0
\]
for almost all $x \in [0,1]$.  
Denote by $E \subset [0,1]$ the set of $x$ for which equality obtains $-$then $\forall \ x \in E$, 

\[
y \ra f(x,y)
\]
is absolutely continuous. 
In fact, for any such $x$,

\[
\int\limits_0^1 \ [f(x,y) - g(x, y)] \eta^\prime (y) \ \td y
\ = \ 
0
\]

$\implies$

\[
f(x,y) - g(x,y)  
\ = \ 
C_x 
\ (\in \R)
\]
for almost all $y \in [0,1]$
\\[-.25cm]

$\implies$
\begin{align*}
f(x,y) \ 
&=\ 
g(x, y) + C_x
\\[15pt]
&=\ 
\int\limits_0^y \ \frac{\partial f}{\partial t} (x, t) \ \td t \hsx + \hsx C_x
\end{align*}
for almost all $y \in [0,1]$.  
The right hand side is an absolutely continuous function of $y$ and the left hand side is a continuous function of $y$.  
Since equality holds for a subset of $[0,1]$ of full measure and such a set is dense in $[0,1]$, 
the conclusion is that

\[
f(x,y) 
\ = \ 
\int\limits_0^y \ \frac{\partial f}{\partial t} (x, t) \ \td t \hsx + \hsx C_x
\]
for all $y \in [0,1]$.
\\

Summary: \ 
\\[-.5cm]

\qquad $y \ra f(x,y)$ is absolutely continuous for almost all $x \in [0,1]$ (viz., $\forall \ x \in E$), 
thereby completing the proof.  
\\

The preceding result also admits an easy converse (where, as above, $1 \leq p < +\infty$).
\\

\qquad{\bf 8.11:}\
{\small\bf THEOREM} \ 
If $f : \Omega \ra \R$ has an \ACL representative $\bar{f}$ whose ordinary partial derivatives belong to $\Lp^p(\Omega)$, 
then these derivatives coincide almost everywhere with the corresponding distributional derivatives of $f$, 
hence $f \in \Lp^{1,p} (\Omega)$.
\\

As for Sobolev spaces, there is a characterization.
\\

\qquad{\bf 8.12:}\
{\small\bf THEOREM} \ 
Let $1 \leq p < +\infty$ $-$then a fucntion $f \in \Lp^p (\Omega)$ belongs to $W^{1,p} (\Omega)$ 
iff it has a representative $\bar{f}$ that is \ACL and whose ordinary partial derivatives belong to $\Lp^p(\Omega)$.
\\

\qquad{\bf 8.13:}\
{\small\bf CRITERION} \ 
Suppose that $f : \Omega \ra \R$ is continuous and \ACL $-$then the ordinary partial derivatives of $f$ 
exist almost everywhere in $\Omega$ and they are Borel functions. 
\\[-.5cm]

PROOF \ 
Fix $i \in \{1, \ldots, n\}$, put

\[
\R_i^{n-1} 
\ = \ 
\{x \in \R^n : x_i = 0\},
\]
\begin{spacing}{1.5}
\noindent
and let $P_i$ be the orthogonal projection of $\R^n$ onto $\R_i^{n - 1}$ 
(so $P_i(x) = x - x_i \hsy e_i)$. 
Let $E_i$ be the set of all $x \in \Omega$ at which $\ds\frac{\partial f}{\partial x_i}$ does not exist, 
the claim being that $\Lm^n(E \cap Q) = 0$ for all $Q \subset \Omega$.  
Since $f$ is continuous, $E$ is Borel and by Fubini, 
\end{spacing}

\[
\Lm^n(E_i \cap Q)
\ = \  
\int\limits_{P_i \Omega} \ \Lm^1 (P_i^{-1} (x) \cap E_i \cap Q) 
\ \td \Lm^{n - 1} (x).
\]

\begin{spacing}{1.6}
\noindent
If $f$ is absolutely continuous on the segment $P_i^{-1} (x) \cap Q$, then $\ds\frac{\partial f}{\partial x_i}$ 
exists almost everywhere on this segment, hence 
$\Lm^1 (P_i^{-1} (x) \cap E_i \cap Q) = 0$, implying thereby that 
$\Lm^n(E_i \cap Q) = 0$, $f$ being \ACL.
\\
\end{spacing}

\begin{spacing}{1.65}
\qquad{\bf 8.14:}\
{\small\bf REMARK} \ 
Recall that without some assumption, the set of points 
in $\Omega$ where 
$\ds\frac{\partial f}{\partial x_i}$ exists need not be Lebesgue measurable (let alone Borel).
\\
\end{spacing}

\qquad{\bf 8.15:}\
{\small\bf CRITERION} \ 
Let 
\[
R 
\ = \ 
]a_1, b_1[ \ \times \ldots \times \ ]a_n, b_n[ 
\hsx \subset \hsx 
\R^n
\]

\begin{spacing}{1.65}
\noindent
be an open rectangle.  
Fix $i \in \{1, \ldots,n\}$ and let $f: R \ra \R$ be a Lebesgue measurable function that is monotone 
on almost every line in $R$ parallel to the $x_i$ axis 
$-$then the ordinary partial derivative $\ds\frac{\partial f}{\partial x_i}$ exists almost everywhere in $R$ 
(and is Lebesgue measurable).
\\
\end{spacing}

\qquad{\bf 8.16:}\
{\small\bf REMARK} \ 
The assumption ``monotone on almost every line in $R$'' cannot be replaced by 
``bounded variation  on almost every line in $R$''
\\[-.5cm]

\begin{spacing}{1.65}
[Note: \ 
But if $f$ is of bounded variation on almost every line in $R$ parallel to the $x_i$ axis, then there is an equivalent function 
$\bar{f}$ which does not have an ordinary partial derivative 
$\ds\frac{\partial \bar{f}}{\partial x_i}$
almost everywhere in $R$ (and is Lebesgue measurable).
\end{spacing}

\chapter{
\text{SECTION 9: \quad BV SPACES}
\vspace{1.cm}\\
$\boldsymbol{\S}$\textbf{9.1}.\quad  PROPERTIES}
\setlength\parindent{2em}
\renewcommand{\thepage}{9-\S1-\arabic{page}}

\vspace{-.75cm}
\qquad
Let $\Omega$ be a nonempty open subset of $\R^n$.
\\

\qquad{\bf 9.1.1:}\quad  
{\small\bf DEFINITION} \ 
Let $f \in \Lp^1 (\Omega)$ $-$then $f$ is said to be a \un{function} \un{of bounded variation} 
if its partial derivatives in the sense of distributions are finite signed Radon measures in $\Omega$ of finite total variation.
\\

\qquad{\bf 9.1.2:}\quad 
{\small\bf NOTATION} \ 
$\BV (\Omega)$ is the set of functions of bounded variation in $\Omega$.
\\[-.5cm]

[Note: \ 
There is a local version of this definition, namely call 
\[
\BV_\locx (\Omega)
\]
the set comprised of all $f \in \Lp_\locx^1 (\Omega)$ with the property that 
$\restr{f}{\Omega^\prime} \in \BV (\Omega^\prime)$ for every nonempty open set
$\Omega^\prime \subset \Omega$ whose closure is a compact subset of $\Omega$.]
\\

\qquad{\bf 9.1.3:}\quad 
{\small\bf \un{N.B.}} \ 
Let $f \in \BV (\Omega)$ $-$then there are finite signed Radon measures
\[
\tD_i f 
\quad (i = 1, \ldots, n)
\]
of finite total variation such that $\forall \ \phi \in C_c^\infty (\Omega)$,
\[
\int\limits_\Omega \ \phi \hsy \td \tD_i f 
\ = \ 
-
\int\limits_\Omega \ \frac{\partial \phi}{\partial x_i} \ f \td \Lm^n 
\quad (i = 1, \ldots, n).
\]

[Note: \ 
\[
\tD f 
\ \equiv \ 
(\tD_1 f, \ldots, \tD_n f)
\]
is an $\R^n$-valued vector measure and on general grounds, $\norm{\tD f}$ is a positive finite measure 
(hence $\norm{\tD f} (\Omega) < +\infty$).]
\\

\qquad{\bf 9.1.4:}\quad {\small\bf LEMMA} \ 
Let $f \in \BV (\Omega)$ $-$then $\norm{\tD f}$ is absolutely continuous w.r.t. Lebesgue measure iff
each of the $\tD_i f$ is, in which case the distributional partial derivatives can be represented by $\Lp^1$ functions.
\\


\qquad{\bf 9.1.5:}\quad 
{\small\bf LEMMA} \ 
Let $f \in \BV (\Omega)$ $-$then
$f \in W^{1, 1} (\Omega)$ iff $\norm{\tD f}$ is absolutely continuous w.r.t. Lebesgue measure, in which case
\[
\norm{\tD f} (\Omega) 
\ = \ 
\int\limits_\Omega \ \norm{\nabla f} \ \td \Lm^n.
\]
\\[-1cm]

\qquad{\bf 9.1.6:}\quad 
{\small\bf \un{N.B.}}  \ 
The containment
\[
W^{1, 1} (\Omega)
\subset 
\BV (\Omega)
\]
is strict and every $f \in C^\infty (\Omega) \cap \BV (\Omega)$ belongs to $W^{1, 1} (\Omega)$.
\\

\qquad{\bf 9.1.7:}\quad 
{\small\bf NOTATION} \ 
Given $\Phi \in C_c^\infty (\Omega; \R^n)$, put
\[
\dvg \hsy \Phi 
\ = \ 
\sum\limits_{i = 1}^n \  \frac{\partial \Phi_i}{\partial x_i},
\]
the \un{divergence} of $\Phi$.
\\

\qquad{\bf 9.1.8:}\quad 
{\small\bf DEFINITION} \ 
Let $f \in \Lp^1 (\Omega)$ $-$then the \un{variation} of $f$ in $\Omega$ is the entity
\[
\sV (f ; \Omega) 
\ = \ 
\sup \bigg\{
\int\limits_\Omega \ f \hsx \dvg \ \hsy  \Phi \ \td \Lm^n \hsy : \hsy 
\Phi \in C_c^\infty (\Omega; \R^n), \norm{\Phi}_\infty \leq 1
\bigg\}.
\]
\\[-.75cm]

\qquad{\bf 9.1.9:}\quad 
{\small\bf THEOREM} \ 
Let $f \in \Lp^1 (\Omega)$ $-$then
\[
\sV (f ; \Omega) 
\ < \ 
+\infty
\]
iff $f \in \BV (\Omega)$.  
And when this is so, 
\[
\sV (f ; \Omega) 
\ = \ 
\norm{\tD f} (\Omega).
\]
\\[-1cm]

\qquad{\bf 9.1.10:}\quad 
{\small\bf LSC PRINCIPLE} \ 
Suppose that $\{f_k\}$ is a sequence in $\BV (\Omega)$ which converges in $\Lp^1 (\Omega)$ to a function $f$ $-$then
\[
\norm{\tD f} (\Omega) 
\ \leq \ 
\liminf\limits_{k \ra \infty} \ \norm{\tD_k f} (\Omega).
\]

PROOF \ 
Choose a $\Phi \in C_c^\infty (\Omega; \R^n)$ with $\norm{\Phi}_\infty \leq 1$, thus 
\begin{align*}
\int\limits_\Omega \ f \hsx \dvg \ \hsy  \Phi \ \td \Lm^n \ 
&=\ 
\lim\limits_{k \ra \infty}
\int\limits_\Omega \ f_k \hsx \dvg \ \hsy  \Phi \ \td \Lm^n
\\[15pt]
&\leq \ 
\liminf\limits_{k \ra \infty} 
\sV(f_k; \Omega) 
\\[15pt]
&= \ 
\liminf\limits_{k \ra \infty}  
\norm{\tD f_k} (\Omega).
\end{align*}
Now take the supremum over $\Phi$.
\\

\qquad{\bf 9.1.11:}\quad 
{\small\bf REMARK}  \ 
To conclude that $f \in \BV (\Omega)$, it suffices to assume that the $f_k$ have equibounded total variations, say $\forall \ k$, 
\[
\norm{\tD f_k} (\Omega)
\ \leq \ 
M.
\]
For then
\\[-1.cm]
\begin{align*}
\sV(f; \Omega) \ 
&= \ 
\norm{\tD f} (\Omega)
\\[15pt]
&\leq \ 
\liminf\limits_{k \ra \infty} \ \norm{\tD_k f} (\Omega) 
\\[15pt]
&\leq \ 
M
\\[15pt]
&< \ 
+\infty.
\end{align*}

\qquad{\bf 9.1.12:}\quad 
{\small\bf NOTATION} \ 
Given $f \in \BV (\Omega)$, put
\[
\norm{f}_\BV \ 
\ = \ 
\norm{f}_{\Lp^1} + \norm{\tD f} (\Omega).
\]

\qquad{\bf 9.1.13:}\quad 
{\small\bf THEOREM} \ 
Under the norm $\norm{\hsx \cdot \hsx}_\BV$, $\BV (\Omega)$ is a Banach space.
\\[-.5cm]

PROOF \ 
Completeness is the issue so suppose that $\{f_k\}$ is a Cauchy sequence in $\BV (\Omega)$ $-$then
by the definition of $\norm{\hsx \cdot \hsx}_\BV$, it must also be a Cauchy sequence in $\Lp^1 (\Omega)$, 
hence by the completeness of $\Lp^1 (\Omega)$, there exists a function $f \in \Lp^1 (\Omega)$ such that 
$f_k \ra f$ in $\Lp^1 (\Omega)$.  
On the other hand, since $\{f_k\}$ is a Cauchy sequence in $\BV (\Omega)$, 
$\{\norm{f_k}_\BV\}$ is bounded: $\exists \ M > 0$ such that $\forall \ k$: 

\begin{align*}
\norm{f_k}_\BV \ 
&= \ 
\norm{f_k}_{\Lp^1} + \norm{\tD f_k} (\Omega) 
\\[15pt]
&\leq\ 
M
\end{align*}

\hspace{1.5cm} $\implies$
\[
\norm{\tD f_k} (\Omega) 
\ \leq \ 
M
\]

\hspace{1.5cm} $\implies$
\[
\norm{\tD f} (\Omega) 
\ \leq \ 
M
\]

\hspace{1.5cm} $\implies$
\[
f \in \BV (\Omega).
\]
The claim now is that $f_k \ra f$ in $\BV (\Omega)$.  
Because we already have convergence in $\Lp^1(\Omega)$, matters reduce to showing that 
\[
\norm{\tD (f_k - f)} (\Omega)  \ra 0 \quad (k \ra \infty).
\]
To this end, let $\varepsilon > 0$ $-$then there exists $N$: 
\begin{align*}
k, \hsy j \geq N 
&\implies 
\norm{f_k - f_j}_\BV < \varepsilon
\\[15pt]
&\implies 
\norm{\tD (f_k - f_j)} (\Omega)  < \varepsilon.
\end{align*}
By construction, 
\[
f_j \ra f 
\quad \text{in} \ \Lp^1(\Omega), 
\]
so
\[
f_k - f_j \ra f_k - f 
\quad \text{in} \ \Lp^1(\Omega), 
\]
thus
\[
\norm{\tD (f_k - f)} (\Omega) 
\ \leq \ 
\liminf\limits_{j \ra \infty} \ 
\norm{\tD (f_k - f_j)} (\Omega) 
\ \leq \ 
\varepsilon
\]
from which the conclusion, $\varepsilon > 0$ being arbitrary.
\\

\qquad{\bf 9.1.14:}\quad 
{\small\bf REMARK} \ 
$\BV (\Omega)$ is not separable.
\\[-.5cm]

[To illustrate, work in $\R$ and consider the family $\sF$ of characteristic functions $\chisubalpha$ 
of the interval $]\alpha,1[$ $(0 < \alpha < 1)$ $-$then $\sF \subset \BV]0,1[$ and for $\alpha \neq \beta$, 
\[
\norm{\chisubalpha - \chisubbeta}_\BV 
\ = \ 
2 + \abs{\alpha - \beta}.]
\]
\\[-.75cm]

\qquad{\bf 9.1.15:}\quad 
{\small\bf \un{N.B.}} \ 
The closure of $\BV (\Omega) \cap C^\infty (\Omega)$ in $\BV (\Omega)$  is $W^{1,1} (\Omega)$, hence is not 
dense in $\BV (\Omega)$.
\\[-.5cm]

[Note: \ 
By way of comparison, recall that $W^{1,1} (\Omega) \cap C^\infty (\Omega)$ is dense in $W^{1,1} (\Omega)$.]
\\

\qquad{\bf 9.1.16:}\quad 
{\small\bf THEOREM} \ 
Let $f \in \BV (\Omega)$ $-$then there exists a sequence $\{f_k\} \subset \BV (\Omega) \cap C^\infty (\Omega)$  such that
\[
f_k \ra f 
\quad \text{in} \ \Lp^1(\Omega) \quad (k \ra \infty)
\]
and
\[
\lim\limits_{k \ra \infty} \ \norm{\tD f_k} (\Omega) 
\ = \ 
\norm{\tD f} (\Omega).
\]

[Note: \ 
It is not claimed nor is it true in general that
\[
\norm{\tD (f_k - f)} (\Omega)  \ra 0 \quad (k \ra \infty).]
\]
\\[-1cm]

\qquad{\bf 9.1.17:}\quad 
{\small\bf APPLICATION} \ 
Take $\Omega = \R^n$ and in GNS, take $p = 1$, hence
\[
1^\star 
\ = \ 
\frac{n}{n-1} 
\quad (n \geq 2).
\]
Choose the $f_k$ above in $C_c^\infty (\R^n)$ $-$then

\allowdisplaybreaks
\begin{align*}
\norm{f_k}_{\Lm^{n/(n-1)}} \ 
&\leq \ 
C(n) \hsx \norm{\nabla f_k}_{\Lp^1}
\\[15pt]
&=\ 
C(n) \hsx \int\limits_{\R^n} \ \norm{\nabla f_k} \ \td \Lm^n
\\[15pt]
&=\ 
C(n) \hsx \norm{\tD f_k} (\R^n),
\end{align*}
so, upon passing to the limit, it follows that $\forall \ f \in \BV (\Omega)$,
\[
\norm{f}_{\Lm^{n/(n-1)}}
\ \leq \ 
C(n) \hsx \norm{\tD f} (\R^n), 
\]

\hspace{1.5cm} $\implies$
\[
\BV(\R^n) 
\hookrightarrow 
\Lm^{n/(n-1)} (\R^n).
\]
\\[-.75cm]

\qquad{\bf 9.1.18:}\quad 
{\small\bf HEURISTICS} \ 
Let $(\XX, \sE)$ be a measurable space, let $\mu$ be a $\sigma$-finite positive measure on $(\XX, \sE)$, 
and let $f : \XX \ra [0,+\infty]$ be a $\mu$-measurable function $-$then (Cavalieri)
\[
\int\limits_\XX \ f \td \mu 
\ = \ 
\int\limits_0^\infty \ \mu(\{x \in \XX : f(x) > t\}) \ \td \Lm^1.
\]

[Let
\[
E_t 
\ = \ 
\{x \in \XX : f(x) > t\}.
\]
Then $\chisubEsubt (x) = 1$ iff $x \in E_t$ iff $f(x) > t$
\\[-.25cm]

\hspace{.7cm} $\implies$
\[
\int\limits_0^\infty \ \chisubEsubt (x) \td t 
\ = \ 
\int\limits_0^{f(x)} \ 1 \ \td t 
\ = \ 
f(x)
\]
\\[-1.25cm]

\hspace{.7cm} $\implies$
\allowdisplaybreaks
\begin{align*}
\int\limits_\XX \ f \ \td \mu  \ 
&=\ 
\int\limits_\XX \ 
\bigg(\int\limits_0^\infty \ \chisubEsubt (x) \ \td t \bigg) 
\td \mu 
\\[15pt]
&=\ 
\int\limits_0^\infty \
\bigg(\int\limits_\XX \ \chisubEsubt (x) \ \td \mu \bigg)
\ \td t 
\qquad \text{(Fubini)}
\\[15pt]
&=\ 
\int\limits_0^\infty \ \mu(E_t) \ \td t .]
\end{align*}
\\[-.75cm]

\qquad{\bf 9.1.19:}\quad 
{\small\bf LEMMA} \ 
Let $E \subset \R^n$ be Lebesgue measurable $-$then $\restr{\chisubE}{\Omega} \in \Lp^1 (\Omega)$ iff 
$\Lm^n(E \cap \Omega) < +\infty$.
\\[-.5cm]

PROOF \ 

\allowdisplaybreaks
\begin{align*}
\int\limits_\Omega \ \restr{\chisubE}{\Omega} \ \td \Lm^n \ 
&=\ 
\int\limits_{\R^n} \ \chisubE \hsx \chisubOmega  \ \td \Lm^n \ 
\qquad \text{(by definition)}
\\[15pt]
&=\ 
\int\limits_{\R^n} \ \chisubEcapOmega  \ \td \Lm^n \ 
\\[15pt]
&=\ 
\Lm^n (E \cap \Omega).
\end{align*}

[Note: \ 
For the record, 
\[
\Lm^n (E \cap \Omega)
\ = \ 
(\Lm^n \Lm E) (\Omega)
\ = \ 
(\chisubE \Lm^n) (\Omega).]
\]
\\[-.75cm]

\qquad{\bf 9.1.20:}\quad 
{\small\bf DEFINITION}\ 
Let $E \subset \R^n$ be a Lebesgue measurable set and suppose that 
$\restr{\chisubE}{\Omega} \in \Lp^1(\Omega)$ 
$-$then the \un{perimeter} of $E$ in $\Omega$, denoted $P(E; \Omega)$, 
is the variation of $\restr{\chisubE}{\Omega}$ in $\Omega$, i.e., 
\[
P(E; \Omega)
\ = \ 
\sV(\restr{\chisubE}{\Omega}; \Omega).
\]
\\[-1cm]

[Note: \ 
The set $E$ is said to have \un{finite perimeter} in $\Omega$ if $P(E,\Omega) < +\infty$.]
\\

\qquad{\bf 9.1.21:}\quad 
{\small\bf NOTATION} \ 
Given $f \in \BV (\Omega)$, put
\[
\Omega_t (f) 
\ = \ 
\{x \in \Omega : f(x) > t\}.
\]
\\[-1cm]

\qquad{\bf 9.1.22:}\quad 
{\small\bf SUBLEMMA} \ 
The function
\[
(x, t) \ra \chisubOmegasubtoff (x)
\]
is $(\Lm^n \times \Lm^1)$-measurable, thus for each $\Phi \in C_c^\infty (\Omega; \R^n)$, the function
\[
t \ra 
\int\limits_{\Omega_t } \ \dvg \hsy \Phi \ \td \Lm^n 
\ = \ 
\int\limits_\Omega \ \chisubOmegasubtoff \hsx \dvg \hsy \Phi \ \td \Lm^n 
\]
is $\Lm^1$-measurable.
\\

\qquad{\bf 9.1.23:}\quad 
{\small\bf LEMMA} \ 
The function
\[
t \ra \norm{\tD \chisubOmegasubtoff} (\Omega)
\]
is Lebesgue measurable.
\\

\qquad{\bf 9.1.24:}\quad 
{\small\bf THEOREM (COAREA)} \ 
Let $f \in \BV (\Omega)$ $-$then
the set $\Omega_t(f)$ has finite perimeter in $\Omega$ for almost all $t$ and
\[
\norm{\tD f}  (\Omega)
\ = \ 
\int\limits_\R \ \norm{\tD \chisubOmegasubtoff} (\Omega) \ \td \Lm^1.
\]

The proof proceeds in two steps.
\\

\un{Step 1}: \ 
Consider
\[
\int\limits_\Omega \ f \hsx \dvg \hsx  \Phi \ \td \Lm^n,
\]
where
\[
\Phi \in C_c^\infty (\Omega; \R^n), 
\quad 
\norm{\Phi}_\infty \leq 1,
\]
and recall that
\[
\int\limits_\Omega \ \dvg \hsx  \Phi \ \td \Lm^n
\ = \ 0.
\]

\qquad \textbullet \quad
$f \geq 0 \implies$
\[
f(x) 
\ = \ 
\int\limits_0^\infty \ \chisubOmegasubtoff (x) \ \td t
\]

\hspace{1.5cm} $\implies$

\begin{align*}
\int\limits_\Omega \ f \hsx \dvg \hsx  \Phi \ \td \Lm^n\ 
&=\ 
\int\limits_\Omega \
\bigg(
\int\limits_0^\infty \ \chisubOmegasubtoff (x) \td t
\bigg)
\dvg \hsy \Phi \ \td \Lm^n
\\[15pt]
&=\ 
\int\limits_0^\infty  \ 
\bigg(
\int\limits_\Omega \
\chisubOmegasubtoff (x) \ \dvg \hsx \Phi \ \td \Lm^n
\bigg)
\ \td t
\\[15pt]
&=\ 
\int\limits_0^\infty  \ 
\bigg(
\int\limits_{\chisubOmegasubtoff} \ \dvg \hsx \Phi \ \td \Lm^n
\bigg)
\ \td t.
\end{align*}

\qquad \textbullet \quad
$f \leq 0 \implies$
\[
f(x) 
\ = \ 
\int\limits_{-\infty}^0 \ (\chisubOmegasubtoff (x) -1) \ \td t
\]

\qquad $\implies$

\begin{align*}
\int\limits_\Omega \ f \hsx \dvg \hsx  \Phi \ \td \Lm^n \ 
&=\ 
\int\limits_\Omega \
\bigg(
\int\limits_{-\infty}^0 \ (\chisubOmegasubtoff (x) -1) \td t
\bigg)
\dvg \hsy \Phi \ \td \Lm^n
\\[15pt]
&=\ 
\int\limits_{-\infty}^0 \ 
\bigg(
\int\limits_\Omega \
(\chisubOmegasubtoff (x) -1) \ \dvg \hsx \Phi \ \td \Lm^n
\bigg)
\ \td t
\\[15pt]
&=\ 
\int\limits_{-\infty}^0 \ 
\bigg(
\int\limits_{\chisubOmegasubtoff} \ \dvg \hsx \Phi \ \td \Lm^n
\bigg)
\ \td t.
\end{align*}
So, upon writing $f = f^+  + (- f^-)$, it follows that 
\[
\int\limits_\Omega \ f \dvg \hsx \Phi \ \td \Lm^n 
\ = \ 
\int\limits_\R 
\bigg(
\int\limits_{\chisubOmegasubtoff} \ \dvg \hsx \Phi \ \td \Lm^n
\bigg)
\td \Lm^1
\]
or still, by the definition of the variation of the perimeter of $\Omega_t(f)$ in $\Omega$,
\[
\int\limits_\Omega \ f \dvg \hsx \Phi \ \td \Lm^n 
\ \leq \ 
\int\limits_\R \ P(\Omega_t (f); \Omega) \ \td \Lm^1
\]
or still, upon taking the supremum over $\Phi$,

\begin{align*}
\norm{\tD f}  (\Omega) \ 
&\leq \ 
\int\limits_\R \ P(\Omega_t(f) ; \Omega) \td t
\\[15pt]
&=\ \int\limits_\R \
\norm{\tD \chisubOmegasubtoff} (\Omega) \ \td \Lm^1.
\end{align*}

It remains to reverse this inequality and for that, as an intermediary, one first shows that for $f \in \BV (\Omega) \cap C^\infty (\Omega)$,
\[
\norm{\tD f}  (\Omega)
\ \geq \ 
\int\limits_\R \ \norm{\tD \chisubOmegasubtoff} (\Omega) \ \td \Lm^1,
\]
a point of detail that will be admitted without proof. 
\\

\un{Step 2}: \ 
Choose a sequence $\{f_k\} \subset \BV (\Omega) \cap C^\infty (\Omega)$ such that
\[
f_k \ra f  
\quad \text{in} \ \Lp^1(\Omega) \ (k \ra \infty)
\]
and
\[
\lim\limits_{k \ra \infty} \ \norm{\tD f_k} (\Omega) 
\ = \ 
\norm{\tD f} (\Omega).
\]
Then $\forall \ k$, 
\[
\norm{\tD f_k} (\Omega) 
\ = \ 
\int\limits_\R \ 
\norm{\tD \chisubOmegasubtoffsubk} (\Omega) 
 \ \td \Lm^1.
\]
Next
\[
f_k(x) - f(x) 
\ = \ 
\int\limits_\R \ (\chisubOmegasubtoffsubk (x) - \chisubOmegasubtoff (x) )  \ \td \Lm^1
\]
and moreover
\[
\abs{f_k(x) - f(x)}
\ = \ 
\int\limits_\R \ \abs{\chisubOmegasubtoffsubk (x) - \chisubOmegasubtoff (x) }  \ \td \Lm^1
\]
since
\[
\sign(f_k(x) - f(x))
\ = \ 
\sign(\chisubOmegasubtoffsubk (x) - \chisubOmegasubtoff (x) )
\]
for all $t$.  
Therefore
\[
\int\limits_\Omega \ \abs{f_k(x) - f(x)} \ \td \Lm^1
\ = \ 
\int\limits_\R \ 
\bigg(
\int\limits_\Omega \ 
\abs{\chisubOmegasubtoffsubk (x) - \chisubOmegasubtoff (x)}
\ \td \Lm^n
\bigg)
 \ \td \Lm^1.
\]
Bearing in mind that $f_k \ra f$ in $\Lp^1 (\Omega)$, there exists a subsequence, not relabeled, with the property that
\[
\chisubOmegasubtoffsubk \ra \chisubOmegasubtoff (f) 
\quad \text{in} \ \Lp^1(\Omega) \ (k \ra \infty)
\]
for almost every $t$.  
Finally

\allowdisplaybreaks
\begin{align*}
\int\limits_\R \ \norm{\tD \chisubOmegasubtoff} (\Omega) \ \td \Lm^1 \ 
&\leq \ 
\int\limits_\R \ 
\liminf\limits_{k \ra \infty} \ 
\norm{\tD \chisubOmegasubtoffsubk} (\Omega) \ \td \Lm^1 \ 
\hspace{0.5cm} \text{(LSC)}
\\[15pt]
&\leq \ 
\liminf\limits_{k \ra \infty} \ 
\int\limits_\R \ 
\norm{\tD \chisubOmegasubtoffsubk} (\Omega) \ \td \Lm^1 \ 
\hspace{0.5cm} \text{(Fatou)}
\\[15pt]
&= \ 
\liminf\limits_{k \ra \infty} \ 
\norm{\tD f_k} (\Omega)
\hspace{2.5cm} \text{(cf. supra)}
\\[15pt]
&= \ 
\lim\limits_{k \ra \infty} \ 
\norm{\tD f_k} (\Omega)
\\[15pt]
&= \ 
\norm{\tD f} (\Omega).
\end{align*}
\\[-1cm]

\qquad{\bf 9.1.25:}\quad 
{\small\bf EXAMPLE} \ 
Given $f \in \BV (\Omega)$, let
\[
f_r (x) \ = \ 
\begin{cases}
\ r \hspace{0.5cm} \text{if} \ f(x) > r\\[4pt]
\ f(x) \hspace{0.5cm} \text{if} \ -r \leq f(x) \leq  r\\[4pt]
\ -r \hspace{0.5cm} \text{if} \ f(x) < -r
\end{cases}
\]
and put
\[
H_r (x) 
\ = \ 
f(x) - f_r(x).
\]
Then $H_r \in \BV (\Omega)$ and 
\begin{align*}
\norm{\tD H_r} (\Omega) \ 
&=\ 
\int\limits_\R \ \norm{\tD \chisubOmegasubtofHr} (\Omega) \ \td \Lm^1
\\[15pt]
&=\ 
\int\limits_{\abs{t} > r} \ \norm{\tD \chisubOmegasubtoff} (\Omega) \ \td \Lm^1.
\end{align*}
\\[-1cm]

\qquad{\bf 9.1.26:}\quad 
{\small\bf PRODUCT RULE} \ 
If $f, \hsx g \in \BV (\Omega)\cap \Lp^\infty (\Omega)$ $-$then $f  g \in \BV (\Omega)\cap \Lp^\infty (\Omega)$ and 
\[
\norm{\tD (f \hsy g)}  (\Omega) 
\ \leq \ 
\norm{f}_\infty \norm{\tD g}  (\Omega) 
+ 
\norm{g}_\infty \norm{\tD f}  (\Omega).
\]
\\[-1cm]

\qquad{\bf 9.1.27:}\quad 
{\small\bf REMARK} \ 
$\BV (\Omega) \cap \Lp^\infty (\Omega)$ is dense in $\BV (\Omega)$.
\\

\qquad{\bf 9.1.28:}\quad 
{\small\bf PRODUCT RULE} \ 
If $f \in \BV (\Omega)$ and if $\phi \in C_c^\infty (\Omega)$, then $\phi \hsy f \in \BV (\Omega)$ and 
\[
\norm{\tD (\phi \hsy f)}  (\Omega) 
\ = \ 
\ldots \hsx .
\]

\chapter{
$\boldsymbol{\S}$\textbf{9.2}.\quad  DECOMPOSITION THEORY}
\setlength\parindent{2em}
\renewcommand{\thepage}{9-\S2-\arabic{page}}

\qquad 
We shall first review matters in $\R$, with $\Omega = ]a,b[$.  
So fix an $f \in \BV (\Omega)$.
\\

\qquad \textbullet \quad 
$\tD f = \tD^a f + \tD^s f$
is the decomposition of $\tD f$ into its absolutely continuous part w.r.t Lebesgue measure $\Lm^1$ and its singular part $\tD^s f$.
\\

Recall next that $\AT_f$ is the set of atoms of the theory, i.e., the $x \in ]a,b[$ such that $\tD (\{x\}) \neq 0$.
\\

\qquad \textbullet \quad 
$\tD^s f = \tD^j f + \tD^c f$, 
\\[-.5cm]

\noindent
where

\[
\begin{cases}
\ \tD^j f \ = \ \tD^s f \lfloor \AT_f
\\[4pt]
\ \tD^c f \ = \ \tD^s f \lfloor (\Omega \backslash\AT_f)
\end{cases}
.
\]
\\[-.75cm]

\qquad{\bf 9.2.1:}\quad  
{\small\bf \un{N.B.}} \ 
The measures $\tD^a f$, $\tD^j f$, $\tD^c f$ are mutually singular and 

\[
\norm{\tD f} (\Omega)
\ = \ 
\norm{\tD^a f} (\Omega) + \norm{\tD^j f} (\Omega) + \norm{\tD^c f} (\Omega).
\]
\\[-.75cm]

\qquad{\bf 9.2.2:}\quad  
{\small\bf DEFINITION} \ 
$f$ if a \un{jump function} if $\tD f = \tD^j f$.
\\

\qquad{\bf 9.2.3:}\quad  
{\small\bf DEFINITION} \ 
$f$ if a \un{Cantor function} if $\tD f = \tD^c f$.
\\

\qquad{\bf 9.2.4:}\quad  
{\small\bf THEOREM} \ 
Each $f \in \BV (\Omega)$ can be represented as a sum

\[
f^a + f^j + f^s,
\]
where $f^a$ belongs to $W^{1, 1} (]a,b[)$, $f^j$ is a jump function, and $f^s$ is  a Cantor function,
\\


\qquad{\bf 9.2.5:}\quad  
{\small\bf \un{N.B.}}  \ 
These functions are uniquely determined up to additive constants and if $\bar{f}$ is an admissible representative of $f$, then
\[
\norm{\tD^a f} (\Omega) 
\ = \ 
\int\limits_a^b \ 
\abs{\bar{f}^\prime} 
\ \td \Lm^1
\]
and
\[
\norm{\tD^j f} (\Omega) 
\ = \ 
\sum\limits_{x \in \AT_f} \ 
\abs{\bar{f} (x+) - \bar{f} (x-)}.
\]
\\[-.5cm]

\qquad{\bf 9.2.6:}\quad  
{\small\bf EXAMPLE} \ 
Work in $\Omega = ]0,1[$ and let $\{r_n\} \subset ]0,1[$ be a sequence.  
Define $f \in \BV (\Omega)$ by the prescription
\[
f(x) 
\ = \ 
\sum\limits_{\{n : r_n < x\}} \ 
2^{-n}.
\]
Then $f$ is a jump function and its distributional derivative $\tD f$ is 
\[
\sum\limits_n \ 
2^{-n} \hsx \delta_{r_n}.
\]
\\[-.75cm]

\qquad{\bf 9.2.7:}\quad  
{\small\bf EXAMPLE} \ 
Work in $\Omega = ]0,1[$ and take for $f$ the Cantor function $-$then its distributional derivative has no absolutely continuous part and 
no jump part.  
\\[-.5cm]

[If $C$ is the Cantor set, then $\tD f$ is (a constant multiple of) $\sH^\gamma \lfloor C$, where 
\\

\noindent
$\ds \gamma = \frac{\log 2}{\log 3}$.]
\\

Assume henceforth that $n > 1$, where as usual $\Omega \subset \R^n$ is nonempty and open.
\\

\qquad{\bf 9.2.8:}\quad  
{\small\bf NOTATION} \ 
Given an $f \in \BV (\Omega)$, put
\[
\begin{cases}
\ n_-(x) \ = \ \ap \liminf\limits_{\substack{y \ra x\\ y \in \Omega}} \  f(y)
\\[15pt]
\ n_+(x) \ = \ \ap \limsup\limits_{\substack{y \ra x\\ y \in \Omega}} \  f(y)
\end{cases}
\quad (x \in \Omega).
\]


\qquad{\bf 9.2.9:}\quad  
{\small\bf LEMMA} \ 
The functions
\[
\begin{cases}
\ x \ra n_-(x)
\\[4pt]
\ x \ra n_+ (x)
\end{cases}
\quad (x \in \Omega)
\]
are Borel measurable functions in $\Omega$.
\\

\qquad{\bf 9.2.10:}\quad  
{\small\bf NOTATION} \ 
\[
\tJ_f 
\ = \ 
\{x \in \Omega \hsy : \hsy n_- (x) < n_+ (x)\}.
\]
\\[-1.25cm]

[Accordingly, $\tJ_f$ is the set of points at which the approximate limit of $f$ does not exist.]

\qquad{\bf 9.2.11:}\quad  
{\small\bf \un{N.B.}} \ 
\[
\Lm^n (\tJ_f) 
\ = \ 
0.
\]
\\[-1cm]

\qquad{\bf 9.2.12:}\quad  
{\small\bf THEOREM} \ 
$\tJ_f$ is $\sH^{n-1}$-measurable.
\\

\qquad{\bf 9.2.13:}\quad  
{\small\bf THEOREM} \ 
There exist countably many $C^1$-hypersurfaces $S_k$ such that 
\[
\sH^{n-1} \big(\tJ_f \hsx - \hsx \bigcup\limits_{k = 1}^\infty \ S_k\big) 
\ = \ 
0.
\]
\\[-.75cm]

\qquad{\bf 9.2.14:}\quad  
{\small\bf THEOREM} \ 
\[
\norm{\tD f} \lfloor \tJ_f 
\ = \ 
(n_+ - n_-) \hsy \sH^{n-1} \lfloor \tJ_f.
\]
\\[-.75cm]

Let 
\[
\tD f 
\ = \ 
\tD^a f + \tD^s f
\]
be the decomposition of $\tD f$ into its absolutely continuous part w.r.t Lebesgue
measure $\Lm^n$ and its singular part $\tD^s f$.  
So
\[
\tD^a f 
\ = \ 
f^a \hsy \Lm ^n,
\]
where $f^a : \Omega \ra \R^n$ is the density of $\tD^a f$ and 
\[
\tD^s f \perp \Lm^n.
\]
\\[-.75cm]

\qquad{\bf 9.2.15:}\quad  
{\small\bf DEFINITION} \ 
\\[-.5cm]

\qquad \textbullet \quad 
The \un{jump part} of $\tD f$ is 
\[
\tD^j f 
\ \equiv \ 
\tD^s \hsy f \hsy \lfloor   \tJ_f.
\]

\qquad \textbullet \quad 
The \un{Cantor part} of $\tD f$ is 
\[
\tD^c f 
\ \equiv \ 
\tD^s \hsy f \hsy \lfloor  (\Omega \backslash \tJ_f).
\]

Therefore
\[
\tD f 
\ = \ 
\tD^a f + \tD^j f + \tD^c f.
\]
\\[-.75cm]

\qquad{\bf 9.2.16:}\quad  
{\small\bf THEOREM} \ 
\[
\tD^j f 
\ = \  
(n_+ - n_-) \hsy \nu_f \hsy \sH^{n-1} \lfloor \tJ_f.
\]
\\[-1cm]

[Note: 
Here
\[
\nu_f (x) 
\ = \ 
\frac{\td \tD^j f}{\td \norm{\tD^j f}} (x)
\]
for $\norm{\tD f}$ almost every $x$ in $\tJ_f$.]
\\

Therefore
\[
\tD f
\ = \ 
\tD^a f + (n_+ - n_-) \hsy \nu_f \hsy \sH^{n-1} \lfloor \tJ_f + \tD^c f.
\]
\\[-.75cm]

\qquad{\bf 9.2.17:}\quad  
{\small\bf REMARK} \ 
Earlier, under the assumption that $n = 1$, we exhibited a decomposition of a BV function but a result of this type does not obtain for BV 
functions of two or more variables.

\chapter{
$\boldsymbol{\S}$\textbf{9.3}.\quad  DIFFERENTIATION}
\setlength\parindent{2em}
\renewcommand{\thepage}{9-\S3-\arabic{page}}

\qquad
Let $\Omega$ be a nonempty open subset of $\R^n$.
\\

\qquad{\bf 9.3.1:}\quad  
{\small\bf RAPPEL} \ 
Let $f \in W^{1, p} (\Omega)$ $-$then $f$ is approximately differentiable almost everywhere (cf. 7.1.6).
\\

\qquad{\bf 9.3.2:}\quad  
{\small\bf THEOREM} \ 
Let $f \in \BV (\Omega)$ $-$then $f$ is approximately differentiable almost everywhere.
\\

[Note: \ 
Let
\[
\tD f 
\ = \ 
\tD^a f + \tD^s f
\]
be the decomposition of $\tD f$ into its absolutely continuous part $\tD^a f$ w.r.t. Lebesgue measure $\Lm^n$ 
and its singular part $\tD^s f$ $-$then 
\[
\tD^a f 
\ = \ 
f^a \hsy \Lm^n, 
\]
where $f^a : \Omega \ra \R^n$ is the density of $\tD^a f$, and 
\[
\ap \hsx \td f 
\ = \ 
f^a
\]
almost everywhere.]
\\

For the moment, take $n = 1$ and let $\Omega = ]a,b[$.  
Suppose that $f \in \BV (\Omega)$ $-$then there is a $g \in \BV (\Omega)$ such that $g = f$ 
almost everywhere and $g$ has an ordinary derivative almost everywhere.
\\[-.5cm]

[To see this, choose $g$ admissible, thus
\[
T_g ]a,b[ 
\ = \ 
e - T_f ]a,b[ 
\ < \ 
+\infty,
\]
so $g$ is of bounded variation in the traditional sense, thus has an ordinary derivative almost everywhere.]
\\

These considerations can be extended to arbitrary $n > 1$.
\\


\qquad{\bf 9.3.3:}\quad  
{\small\bf THEOREM} \ 
Let $f \in \BV (\Omega)$ $-$then there is a $g \in \BV (\Omega)$ which is equivalent to $f$ with the property that 
its ordinary partial derivative $\ds\frac{\partial g}{\partial x_i}$ $(i = 1, \ldots, n)$ exists almost everywhere.
\\

\qquad{\bf 9.3.4:}\quad  
{\small\bf REMARK} \ 
It follows that $f$ has approximate partial derivatives almost everywhere, hence has an approximate differential almost everywhere.
\\[-.5cm]

[Note: \ 
Neither $f$ nor any equivalent function need have ordinary partial derivatives at any point.]
\\

\qquad{\bf 9.3.5:}\quad  
{\small\bf \un{N.B.}}\ 
If in addition $f$ is continuous, then $f$ does have ordinary partial derivatives almost everywhere.
\\

\qquad{\bf 9.3.6:}\quad  
{\small\bf NOTATION} \ 
Generically, 
\[
(x_1, \ldots, x_n) 
\ = \ 
(x_i^\prime, x_i),
\]
where
\[
x_i^\prime 
\ = \ 
(x_1, \ldots, x_{i-1}, x_{i+1}, \ldots, x_n)
\quad (i = 1, \ldots, n).
\]
\\[-1cm]

\qquad{\bf 9.3.7:}\quad  
{\small\bf NOTATION} \ 
An open rectangle 
\[
R 
\ = \ 
]a_1, b_1[ \times \cdots \times ]a_n, b_n[ 
\subsetx \R^n
\]
can be viewed as the product of a rectangle
\[
\begin{cases}
\ R_i^\prime \subsetx \R^{n-1} \hspace{0.5cm} \text{(variable $x_i^\prime$)}
\\[4pt]
\ R_i \subsetx \R \hspace{0.5cm} \text{(variable $x_i$)}
\end{cases}
\]
and we write
\[
R 
\ = \ 
R_i^\prime \times R_i.
\]

Let $f \in \BV (\Omega)$ $-$then
\[
\tD f 
\ = \ 
\tD^a f + \tD^s f,
\]
where
\[
\tD^a f 
\ = \ 
f^a \hsx \Lm^n.
\]
\\[-1cm]

Here it is a question of $\R^n$-valued vector measures: 
\[
f^a \in \Lp^1 (\Omega ; \R^n)
\]
and for every Borel set $E \subset \Omega$, 
\[
\tD f (E) 
\ = \ 
\int\limits_E \ 
f^a 
\ \td \Lm^n 
\hsx + \hsx 
\tD^s f (E).
\]
\\[-.75cm]

\qquad{\bf 9.3.8:}\quad  
{\small\bf THEOREM} \ 
Let $f \in \BV (\Omega)$, let $1 \leq i \leq n$, and let $g$ be any function equivalent to $f$ for which the ordinary partial 
derivative $\ds \frac{\partial g}{\partial x_i}$ exists almost everywhere $-$then
\[
\frac{\partial g}{\partial x_i} 
\ = \ 
f_i^a
\]
almost everywhere.
\\

\begin{spacing}{1.5}
The proof is on the lengthy side and will be broken up into 3 steps.  
Write for brevity $\partial_i$ in place of $\ds\frac{\partial}{\partial x_i}$.
\\[-.25cm]
\end{spacing}

\qquad  \un{Step 1:} \quad
Consider a convex function $\Phi : \R \ra [0, +\infty[$ and let $R$ be an open rectangle whose closure $\bar{R}$ is 
contained in $\Omega$ $-$then the claim is that
\[
\int\limits_R \ 
\Phi (\partial_i g) 
\ \td \Lm^n 
\ \leq \ 
\liminf\limits_{\varepsilon \downarrow 0} \ 
\int\limits_{R_\varepsilon} \ 
\Phi (\partial_i f_\varepsilon) 
\ \td \Lm^n , 
\]
where

\[
R_\varepsilon 
\ = \ \{x \in R \hsy : \hsy \dist(x, \partial R) > \varepsilon\}.
\]
\\[-.75cm]

\noindent
To begin with, if $\phi$ is sufficiently smooth and if $h > 0$ is sufficiently small, then 

\allowdisplaybreaks
\begin{align*}
\Phi
\bigg(
\frac{\phi (x_i^\prime, x_i + h) - \phi (x_i^\prime, x_i)}{h}
\bigg)
 \ 
&=\ 
\Phi
\bigg(
\frac{1}{h} \ 
\int\limits_{x_i}^{x_i + h} \ 
\partial_i \phi (x_i^\prime, t) 
\ \td t
\bigg)
\\[15pt]
&\leq \ 
\frac{1}{h} \ 
\int\limits_{x_i}^{x_i + h} \ 
\Phi (\partial_i \phi (x_i^\prime, t)) 
\ \td t
\qquad \text{(Jensen).}
\end{align*}

\noindent
Now integrate this along $(R_i)_h$, hence 

\allowdisplaybreaks
\begin{align*}
\int\limits_{{(R_i)}_h} \ 
\Phi
\bigg(
\frac{\phi (x_i^\prime, x_i + h) - \phi (x_i^\prime, x_i)}{h}
\bigg)
\ \td x_i \ 
&\leq \ 
\frac{1}{h} \ 
\int\limits_{{(R_i)}_h} \ 
\int\limits_{x_i}^{x_i + h} \ 
\Phi (\partial_i \phi (x_i^\prime, t)) 
\ \td t \hsy \td x_i
\\[15pt]
&\leq \ 
\frac{1}{h} \ 
\int\limits_{R_i} \ 
\int\limits_{t - h}^t \ 
\Phi (\partial_i \phi (x_i^\prime, t)) 
\ \td x_i \hsy \td t
\\[15pt]
&= \ \int\limits_{R_i} \ 
\Phi (\partial_i \phi (x_i^\prime, t)) 
\ \td t
\end{align*}

$\implies$

\[
\int\limits_{R_h} \ 
\Phi
\bigg(
\frac{\phi (x_i^\prime, x_i + h) - \phi (x_i^\prime, x_i)}{h}
\bigg)
\ \td x
\ \leq \ 
\int\limits_R \ 
\Phi (\partial_i \phi (x)) 
\ \td x.
\]
\\[-.75cm]

\noindent
Specialize and take $\gphi = f_\varepsilon$, the mollification of $f$ $-$then almost everywhere

\[
 f_\varepsilon (x_i^\prime, x_i + h) -  f_\varepsilon (x_i^\prime, x_i) 
\ra
 f (x_i^\prime, x_i + h) -  f (x_i^\prime, x_i) 
 \quad (\varepsilon \downarrow 0), 
\]
thus by Fatou, 

\[
\int\limits_{R_h} \ 
\Phi
\left(
\frac{f(x_i^\prime, x_i + h) - f(x_i^\prime, x_i)}{h}
\right)
\ \td x
\ \leq \ 
\liminf\limits_{\varepsilon \downarrow 0} \ 
\int\limits_{R_\varepsilon} \ 
\Phi (\partial_i f_\varepsilon (x)) 
\ \td x
\]
or still, 
\[
\int\limits_{R_h} \ 
\Phi
\left(
\frac{g(x_i^\prime, x_i + h) - g(x_i^\prime, x_i)}{h}
\right)
\td x
\ \leq \ 
\liminf\limits_{\varepsilon \downarrow 0} \ 
\int\limits_{R_\varepsilon} \ 
\Phi (\partial_i f_\varepsilon (x)) 
\ \td x.
\]
To finish, it remains only to send $h$ to 0.
\\

\qquad  \un{Step 2:} \quad
The next claim is that
\[
\int\limits_{R_\varepsilon} \ 
\Phi (\partial_i f_\varepsilon) 
\ \td \Lm^n  \
\ \leq \ 
\int\limits_{R} \ 
\Phi (f_i^a) 
\ \td \Lm^n  \
\hsx + \hsx 
\norm{\tD_i^s f} (R),
\]
where the convex function $\Phi$ is subject to the condition
\[
\Phi (s + t) 
\ \leq \ 
\Phi (s) + \abs{t}
\]
for all $s, \hsx t \in \R$.  
First, for every $x \in \Omega_\varepsilon$ and any $i \in \{1, \ldots, n\}$, 

\allowdisplaybreaks
\begin{align*}
\partial_i f_\varepsilon (x) \ 
&=\ 
\frac{\partial f_\varepsilon}{\partial x_i} (x)
\\[15pt]
&=\ 
\int\limits_\Omega \ 
\frac{\partial \phi_\varepsilon}{\partial x_i} (x - y) \hsx f(y)
\ \td y
\\[15pt]
&=\ 
-
\int\limits_\Omega \ 
\frac{\partial \phi_\varepsilon}{\partial y_i} (x - y) \hsx f(y)
\ \td y
\\[15pt]
&=\ 
\int\limits_\Omega \ 
\phi_\varepsilon (x - y) 
\ \td \tD_i f
\\[15pt]
&=\ 
\int\limits_\Omega \ 
\phi_\varepsilon (x - y) 
f_i^a (y) 
\ \td y 
\hsx + \hsx 
\int\limits_\Omega \ 
\phi_\varepsilon (x - y) 
\ \td \tD_i^s f(y)
\end{align*}

$\implies$
\[
\Phi(\partial_i f_\varepsilon (x)) 
\ \leq \ 
\Phi
\left(
\int\limits_\Omega \ 
\phi_\varepsilon (x - y) 
f_i^a (y) 
\ \td y
\right)
\hsx + \hsx 
\abs{
\int\limits_\Omega \ 
\phi_\varepsilon (x - y) 
\ \td \tD_i^s f(y)
}
.
\]
Since
\[
\int\limits_\Omega \ 
\phi_\varepsilon (x - y) 
\ \td y
\ = \ 1,
\]
the first term can be estimated by Jensen, so
\[
\Phi(\partial_i f_\varepsilon (x)) 
\ \leq \ 
\int\limits_\Omega \ 
\phi_\varepsilon (x - y) 
\Phi (f_i^a (y))
\ \td y
\hsx + \hsx 
\int\limits_\Omega \ 
\phi_\varepsilon (x - y) 
\ \td \norm{\tD_i^s f} (y).
\]
Therefore
\[
\int\limits_{R_\varepsilon} \ 
\Phi (\partial_i f_\varepsilon) 
\ \td \Lm^n  \
\ \leq \ 
\int\limits_{R} \ 
\Phi (f_i^a) 
\ \td \Lm^n  \
\hsx + \hsx 
\int\limits_{R} \ 
\ \td \norm{\tD_i^s f}.
\]
And
\[
\int\limits_{R} \ 
\ \td \norm{\tD_i^s f}
\ = \ 
\norm{\tD_i^s} (R).
\]
\\

\qquad{\bf 9.3.9:}\quad  
{\small\bf \un{N.B.}} \  
Step 1 and Step 2 
\\[-.25cm]

$\implies$
\allowdisplaybreaks
\begin{align*}
\int\limits_R \ 
\Phi (\partial_i g) 
\ \td \Lm^n  \ 
&\leq \ 
\liminf\limits_{\varepsilon \downarrow 0} \ 
\int\limits_{R_\varepsilon} \ 
\Phi (\partial_i f_\varepsilon) 
\ \td \Lm^n  \
\\[15pt]
&\leq \ 
\int\limits_{R} \ 
\Phi (f_i^a) 
\ \td \Lm^n  \
\hsx + \hsx 
\norm{\tD_i^s f} (R).
\end{align*}
\\[-1cm]

\qquad  \un{Step 3:} \quad
Work with 
\[
R(x_0, r)
\ \equiv \ 
x_0 \hsx + \hsx 
\left]
-\frac{r}{2}, \frac{r}{2} \hsy
\right[^{\hsx n}.
\]
Then
\[
\lim\limits_{r \ra 0} \ 
\frac{\norm{\tD_i^s f}(R(x_0, r))}{r^n} 
\ = \ 
0
\]
\\[-1cm]

\noindent
for almost all $x_0\in \Omega$ (differentiation principle).  
Fix such an $x_0$ and take it to be a Lebesgue point for $\Phi (\partial_i g)$ and $\Phi (f_i^a)$.  
Since
\[
\Lm^n (R(x_0, r)) 
\ = \ 
r^n, 
\]
when $r \ra 0$, 

\[
\begin{cases}
\ds 
\frac{1}{r^n} \ 
\ \int\limits_{R(x_0, r)} \ 
\Phi (\partial_i g) 
\ \td \Lm^n 
\ra 
\Phi (\partial_i g) (x_0)
\\[26pt]
\ds 
\frac{1}{r^n} \ 
\ \int\limits_{R(x_0, r)} \ 
\Phi (f_i^a) 
\ \td \Lm^n 
\ra 
\Phi (f_i^a) (x_0)
\end{cases}
\quad \ra 0,
\]
\\[-.75cm]

$\implies$
\[
\Phi (\partial_i g) (x_0)
\ \leq \ 
\Phi (f_i^a (x_0)).
\]
\\[-1cm]

\noindent
Choose now for $\Phi$ the function

\[
\Psi (t) \ = \ 
\begin{cases}
\ e^t \hspace{1.05cm} (t < 0)
\\[4pt]
\ t + 1 \hspace{0.5cm} (t \geq 0)
\end{cases}
.
\]
Because $\Psi$ is monotone increasing, it follows that $\partial_i g(x) \leq f_i^a(x)$ almost everywhere in $\Omega$ 
and consideration of $\Psi (-t)$ implies that $\partial_i g(x) \geq f_i^a (x)$ almost everywhere in $\Omega$.
\\

\qquad{\bf 9.3.10:}\quad  
{\small\bf SCHOLIUM} \ 
Start with an $f \in \BV (\Omega)$, replace it by an equivalent $g \in \BV (\Omega)$ with the property that 
the ordinary partial derivatives $\ds \frac{\partial g}{\partial x_i}$ $(i = 1, \ldots, n)$ exist almost everywhere $-$then 
\allowdisplaybreaks
\begin{align*}
\nabla g \ 
&=\ 
\left(
\frac{\partial g}{\partial x_1}, \ldots, \frac{\partial g}{\partial x_n}
\right)
\\[11pt]
&=\ 
(f_1^a , \ldots, f_n^a)
\\[11pt]
&=\ 
f^a
\end{align*}
almost everywhere.
\\[-.5cm]

[Note: \ 
\[
f^a \in \Lp^1 (\Omega; \R^n)
\implies 
\nabla g \in \Lp^1 (\Omega; \R^n).] 
\]

\chapter{
$\boldsymbol{\S}$\textbf{9.4}.\quad  BVL}
\setlength\parindent{2em}
\renewcommand{\thepage}{9-\S4-\arabic{page}}

\qquad
Let us first review the situation when $n = 1$.
\\

\qquad{\bf 9.4.1:}\quad
{\small\bf RAPPEL} \ 
If $\Omega \subset \R$ is open and nonempty and if $f \in \Lp^1 (\Omega)$, then the essential variation of $f$, 
denoted $e - T_f \Omega$, is the set

\[
\inf  \ \{T_g \Omega \hsy : \hsy g = f \ \text{almost everywhere}\}.
\]
Moreover $f \in \BV (\Omega)$ iff 
\[
e - T_f \Omega \ < \ +\infty.
\]
And then

\[
e - T_f \Omega 
\ = \ 
\norm{\tD f} (\Omega).
\]
\\[-.75cm]

[Note: \ 
Recall that $\Omega$ is the union of its connected components, these being intervals (finite or infinite).]
\\

Let $\Omega$ be a nonempty open subset of $\R^n$.
\\

\qquad{\bf 9.4.2:}\quad
{\small\bf NOTATION} \ 
Per $x_i^\prime \in \R^{n-1}$, put

\[
\Omega_{x_i^\prime} 
\ = \ 
\{x_i \in \R \hsy : \hsy (x_i^\prime, x_i) \in \Omega\}
\]
and if $\Omega_{x_i^\prime} \neq \emptyset$ and if $f : \Omega \ra \R$ is Lebesgue measurable, denote by 
\[
e - T_{f(x_i^\prime,-)} \Omega_{x_i^\prime} 
\]
the essential variation of the function 

\[
x_i \ra f(x_i^\prime, x_i).
\]
\\[-.75cm]


\qquad{\bf 9.4.3:}\quad
{\small\bf NOTATION} \ 
If $f \in \Lp^1 (\Omega)$, write

\[
\int\limits_\Omega \ 
f 
\ \td \Lm^n 
\ = \ 
\int\limits_{\R^{n-1}} \ 
\int\limits_{\Omega_{x_i^\prime}} \ 
f(x_i^\prime, x_i) 
\ \td x_i \hsy \td x_i^\prime.
\]
\\[-.75cm]

[Put

\[
\int\limits_{\Omega_{x_i^\prime}} \ 
f(x_i^\prime, x_i) 
\ \td x_i 
\ = \ 
0
\]
if $\Omega_{x_i^\prime} = \emptyset$.]
\\

\qquad{\bf 9.4.4:}\quad
{\small\bf CRITERION} \ 
Let $f \in \Lp^1 (\Omega)$.  
Suppose that there exists an equivalent function $g \in \Lp^1 (\Omega)$ and nonnegative functions 
$V_1, \ldots V_n$ in $\Lp^1 (\R^{n-1})$ 
such that 

\[
e - T_{g(x_i^\prime,-)} \Omega_{x_i^\prime} 
\ \leq \ 
V_i (x_i^\prime)
\]
for all $x_i^\prime \in \R^{n-1}$ such that $\Omega_{x_i^\prime}$ is nonempty $(i = 1, \ldots, n)$ $-$then $f \in \BV (\Omega)$.
\\

PROOF \ 
Fix $i \in \{1, \ldots, n\}$ and note that $V_i (x_i^\prime)$ is finite almost everywhere (being in $\Lp^1 (\R^{n-1})$), hence

\[
e - T_{g(x_i^\prime,-)} \Omega_{x_i^\prime} 
\]
is finite almost everywhere.  
But

\[
g(x_i^\prime,-) \in \Lp^1  (\Omega_{x_i^\prime})
\]
for almost all $x_i^\prime$ (Fubini), so the conclusion is that for almost all $x_i^\prime$, 

\[
g(x_i^\prime,-)\in \BV (\Omega_{x_i^\prime}),
\]
from which a finite signed Radon measure $\tD_{x_i^\prime}$ in $\Omega_{x_i^\prime}$ of finite total variation such that 
\[
\partial_i g(x_i^\prime,-) 
\ = \ 
\tD_{x_i^\prime}
\]
with

\[
\norm{\tD_{x_i^\prime}} (\Omega_{x_i^\prime}) 
\ = \ 
e - T_{g(x_i^\prime,-)} \Omega_{x_i^\prime} 
\]

\noindent
Proceeding, define a linear functional $\Lambda$ on $C_c^1 (\Omega)$ by the rule

\[
\Lambda (\phi) 
\ = \ 
\int\limits_\Omega \ 
f \hsy \partial_i \hsy \phi
\ \td \Lm^n.
\]

\noindent
Then $\Lambda$ is continuous w.r.t. the topology of $C_0(\Omega)$.  
Proof: \ 
\begin{align*}
\Lambda (\phi) \ 
&= \ 
\int\limits_{\R^{n-1}} \ 
\int\limits_{\Omega_{x_i^\prime}} \ 
f \hsy \partial_i \hsy \phi
\ \td x_i \hsy \td x_i^\prime
\\[15pt]
&=\ 
-
\int\limits_{\R^{n-1}} \ 
\int\limits_{\Omega_{x_i^\prime}} \ 
\phi 
\ \td \tD_{x_i^\prime} \hsy \td x_i
\end{align*}

$\implies$

\begin{align*}
\abs{\Lambda (\phi)} \ 
&\leq\ 
\max\limits_\Omega \ 
\abs{\phi} \ 
\int\limits_{\R^{n-1}} \ 
\norm{\tD_{x_i^\prime}} (\Omega_{x_i^\prime}) 
\ \td x_i
\\[15pt]
&=\ 
\max\limits_\Omega \ 
\abs{\phi} \ 
\int\limits_{\R^{n-1}} \
e - T_{g(x_i^\prime,-)} \Omega_{x_i^\prime} 
\\[15pt]
&\leq\ 
\max\limits_\Omega \ 
\abs{\phi} \ 
\int\limits_{\R^{n-1}} \
V_i(x_i^\prime)
\ \td x_i^\prime.
\end{align*}
Now extend it to a continuous linear functional on $C_0(\Omega)$ and use the ``$C_0$'' version of the RRT to get a finite 
signed Radon measure $\tD_i$ such that for all $\phi \in C_0 (\Omega)$, 
\[
\Lambda (\phi)
\ = \ 
-\int\limits_\Omega \ 
\phi 
\ \td \tD_i,
\]
or still, $\forall \ \phi \in C_c^1 (\Omega)$, 

\[
\int\limits_\Omega \ 
f \hsy \partial_i \hsy \phi 
\ \td \Lm^n 
\ = \ 
-
\int\limits_\Omega \ 
\phi 
\ \td \hsy \tD_i.
\]
Since 

\[
\norm{\tD_i}  (\Omega) 
\ = \ 
\norm{\Lambda}^* 
\ \leq \ 
\int\limits_{\R^{n-1}} \ 
V_i (x_i^\prime) 
\ \td x_i^\prime 
\ < \ 
+\infty, 
\]
it follows that $f \in \BV (\Omega)$.
\\

\qquad{\bf 9.4.5:}\quad
{\small\bf REMARK} \ 
Take $n = 1$ and suppose that $\Omega = ]a,b[$ $-$then $x_i^\prime$ is just an abstract point, 
call it $*$, and $\Omega_*$ can be identified with $]a,b[$.  
Starting with $f \in \Lp^1 (]a,b[)$, the assumption above amounts to saying that there exists a 
$g \in \Lp^1 (]a,b[)$ and a constant $C \geq 0$ such that 

\[
e - T_{g(*, -)} ]a,b[ 
\ \leq \ 
C 
\ < \ 
+\infty.
\]
But then
\[
g \in \BV (]a,b[)
\]
which implies that 

\[
f \in \BV (]a,b[).
\]
\\[-.75cm]

Take 
\[
Q 
\ = \ 
[0,1] \times [0,1]
\]
and let $f : Q \ra \R$ be a continuous function (hence $f \in \Lp^1 (Q))$.
\\

\qquad{\bf 9.4.6:}\quad
{\small\bf NOTATION} \ 

\[
\begin{cases}
V_x (f; y) \ = \ T_{f(-, y)} [0,1] \hspace{0.25cm} (0 \leq y \leq 1)
\\[8pt]
V_y (f; x) \ = \ T_{f(x,-)} [0,1] \hspace{0.25cm} (0 \leq x \leq 1)
\end{cases}
.
\]

[Note: 
Since $f$ is continuous, 
\[
\begin{cases}
T_{f(-,y)} [0,1] \ = \ T_{f(-,y)} ]0,1[
\\[8pt]
T_{f(x,-)} [0,1] \ = \ T_{f(x,-)} ]0,1[
\end{cases}
.
\]
\\[-.5cm]

\qquad{\bf 9.4.7:}\quad
{\small\bf LEMMA} \ 
$V_x (f; y)$ and $V_y (f; x)$ are Lebesgue measurable.
\\

\qquad{\bf 9.4.8:}\quad
{\small\bf DEFINITION} \ 
$f$ is said to be of \un{bounded variation in the sense} \un{of Tonelli} if 

\[
\begin{cases}
\ \ds\ \int\limits_0^1 \ V_x (f; y) \ \td y \ < \ +\infty 
\\[21pt]
\ \ds\ \int\limits_0^1 \ V_y (f; x) \ \td x \ < \ +\infty 
\end{cases}
.
\]
\\[-.25cm]

\qquad{\bf 9.4.9:}\quad
{\small\bf \un{N.B.}} \ 
When dealing with essential variation on open subsets of the line, if the function in question is continuous, 
one can work instead with the usual variation, the reason being that the approximation via approximate points of continuity 
amounts to approximation via points of continuity.
\\

\qquad{\bf 9.4.10:}\quad
{\small\bf SCHOLIUM} \ 
If $f$ is of bounded variation in the sense of Tonelli, then $\restr{f}{Q^\circ} \in \BV (Q^\circ)$ 
and the ordinary partial derivatives 

\[
\frac{\partial f}{\partial x}, 
\quad 
\frac{\partial f}{\partial y}
\]
exist almost everywhere in $Q^\circ$ ($\restr{f}{Q^\circ}$ being continuous).
\\[0.25cm]

Relax the assumption that $f : Q \ra \R$ is continuous to merely that $f \in \Lp^1 (Q^\circ)$.
\\

\qquad{\bf 9.4.11:}\quad
{\small\bf NOTATION}\ 
\[
\begin{cases}
\ e - V_x (f; y) \ = \ e - T_{f(-,y)} ]0,1[ \hspace{0.25cm} (0 < y < 1)
\\[8pt]
\ e - V_y (f; x) \ = \ e - T_{f(x,-)} ]0,1[ \hspace{0.25cm} (0 < x < 1)
\end{cases}
.
\]
\\[-.75cm]

[Note: \ The essential variations here are per $]0, 1[$ \hsy :

\[
\begin{cases}
\ f(-, y) \ \text{is the function} \  x \ra f(x, y) \hspace{0.25cm} (0< x < 1)
\\[8pt]
\ f(x, -) \ \text{is the function} \  y \ra f(x, y) \hspace{0.25cm} (0< y < 1)
\end{cases}
.
\]
\\[-.5cm]

\qquad{\bf 9.4.12:}\quad
{\small\bf LEMMA} \ 
$e - V_x (f; y)$ and $e - V_y (f; x)$ are Lebesgue measurable. 
\\

\qquad{\bf 9.4.13:}\quad
{\small\bf DEFINITION} \ 
$f$ is said to be of \un{bounded variation in the sense} \un{of Cesari} if

\[
\begin{cases} 
\ 
\ds
\ \int\limits_0^1 \ 
e - V_x (f; y) 
\ \td y 
\ < \ 
+\infty
\\[21pt]
\ 
\ds
\ \int\limits_0^1 \ 
e - V_y (f; x)
\ \td x 
\ < \ 
+\infty
\end{cases}
.
\]
\\[-.5cm]

\qquad{\bf 9.4.14:}\quad
{\small\bf REMARK} \ 
Under the preceding circumstances, it can be shown that there exists a function $g$ equivalent to $f$ such that 

\[
\begin{cases}
\ 
\ \int\limits_0^1 \ 
V_x (g; y) 
\ \td y 
\ < \ 
+\infty
\\[21pt]
\ 
\ds
\ \int\limits_0^1 \ 
V_y (g; x) 
\ \td x 
\ < \ 
+\infty
\end{cases}
.
\]
\\[-.75cm]

\qquad{\bf 9.4.15:}\quad
{\small\bf SCHOLIUM} \ 
If $f$ is of bounded variation in the sense of Cesari, then $f \in \BV(Q^\circ)$ 
and there is a $g \in \BV(Q^\circ)$ 
which is equivalent to $f$ with the property that 
the ordinary partial derivatives 
\[
\frac{\partial g}{\partial x}, 
\quad 
\frac{\partial g}{\partial y}
\]
exist almost everywhere in $Q^\circ$.

\chapter{
SECTION \textbf{10}:\quad  ABSOLUTE CONTINUITY}
\setlength\parindent{2em}
\setcounter{chapter}{10}
\renewcommand{\thepage}{10-\arabic{page}}

\vspace{-.25cm}
\qquad
\qquad{\bf 10.1:}\ 
{\small\bf RAPPEL} \ 
Let $\Omega$  be a nonempty open subset of $\R$ $-$then a function 
$f : \Omega \ra \R$ is absolutely continuous if $\forall \ \epsilon > 0$, $\exists \ \delta > 0$ such that for every finite collection 
$[a_1, b_1], \hsx , \ldots, \hsx [a_k, b_k]$ of pairwise disjoint closed intervals in $\Omega$, 

\[
\sum\limits_{j = 1}^k \ \Lm^1 ([a_j, b_j]) 
\ < \ 
\delta 
\hsx \implies \hsx 
\sum\limits_{j = 1}^k \ \abs{f(a_j) - f(b_j)} 
\ < \ 
\epsilon. 
\]
\\[-.75cm]

Here is one extrapolation from $\R$ to $\R^n$ ($n > 1)$, where now $\Omega$ is a nonempty open subset of $\R^n$.
\\[-.25cm]

\qquad{\bf 10.2:}\ 
{\small\bf DEFINITION} \ 
A function $f : \Omega \ra \R^n$ is \un{$n$-absolutely continuous} if 
$\forall \ \epsilon > 0$, $\exists \ \delta > 0$ such that for every finite collection 
$B(x_1, r_1), \hsx \ldots, \hsx B(x_k, r_k)$ of pairwise disjoint closed balls in $\Omega$, 
\[
\sum\limits_{j = 1}^k \ \Lm^n (B(x_j, r_j))
\ < \ 
\delta 
\hsx \implies \hsx 
\sum\limits_{j = 1}^k \ \osc (f, B(x_j, r_j))^n
\ < \ 
\epsilon. 
\]

[Note: \ 
If $E$ is a subset of $\R^n$, then

\[
\osc (f, E) 
\ = \ 
\diam (f(E)).]
\]

Obviously, 
\[
\text{$n$-absolute continuity} 
\implies
\text{continuity}.
\]
\\[-1cm]

\qquad{\bf 10.3:}\ 
{\small\bf NOTATION} \ 
Put 
\[
V_n (f, \Omega) 
\ = \ 
\sup \ 
\left\{
\sum\limits_{j = 1}^k \ (\osc (f, B(x_j, r_j)))^n  
\right\}.
\]
\\[-.75cm]

\qquad{\bf 10.4:}\ 
{\small\bf DEFINITION} \ 
$f$ is of \un{bounded $n$-variation} in $\Omega$ if
\[
V_n (f, \Omega) 
\ < \ 
+\infty.
\]

\qquad{\bf 10.5:}\ 
{\small\bf NOTATION} \ 
$\BV^n(\Omega)$ is the set of all functions $f : \Omega \ra \R^n$ of bounded $n$-variation in $\Omega$.
\\

\qquad{\bf 10.6:}\ 
{\small\bf NOTATION} \ 
$\AC^n(\Omega)$ is the set of all functions $f \in \BV^n(\Omega)$ that are $n$-absolutely continuous.
\\

\qquad{\bf 10.7:}\ 
{\small\bf REMARK} \ 
The notion of $n$-absolutely continuous uses closed balls.  
One could also work with closed cubes.  
Here, however, one has to be careful: 
Examples have been constructed which show that working with closed balls is not the same as working with closed cubes, 
thus that these two concepts are incomparable.  
\\

\qquad{\bf 10.8:}\ 
{\small\bf DEFINITION} \ 
A function $f : \Omega \ra \R^n$ satisifies the \un{condition $RR$} if there is a nonnegative function $\phi \in \Lp^1(\Omega)$, 
a so-called \un{weight}, such that 
\[
(\osc (f, B(x,r)))^n 
\ \leq \ 
\int\limits_{B(x, r) } \ \phi \ \td \Lm^n
\]
for every $B(x, r) \subset \Omega$.
\\

\qquad{\bf 10.9:}\ 
{\small\bf NOTATION} \ 
Denote by $RR^n(\Omega)$ the set of all functions $f : \Omega \ra \R^n$ which satisfy condition $RR$.
\\

\qquad{\bf 10.10:}\ 
{\small\bf THEOREM} \ 
\[
RR^n(\Omega)
\ = \ 
\AC^n(\Omega).
\]

\qquad{\bf 10.11:}\ 
{\small\bf THEOREM}  \ 
Let $f \in \BV^n(\Omega)$ $-$then $f$ is differentiable almost everywhere.
\\

Matters can be generalized, thus suppose that $0 < \lambda \leq 1$. 
\\

\qquad{\bf 10.12:}\ 
{\small\bf DEFINITION} \ 
A function $f : \Omega \ra \R^n$  is \un{$\lambda$, $n$-absolutely continuous} if 
$\forall \ \epsilon > 0$, $\exists \ \delta > 0$ such that for any finite collection 
$B(x_1, r_1), \hsx \ldots, \hsx B(x_k, r_k)$ of pairwise disjoint closed balls in $\Omega$, 
\[
\sum\limits_{j = 1}^k \ \Lm^n (B(x_j, r_j))
\ < \ 
\delta 
\hsx \implies \hsx 
\sum\limits_{j = 1}^k \ \osc (f, B(x_j, \lambda \hsy r_j))^n
\ < \ 
\epsilon. 
\]
\\[-.75cm]

\qquad{\bf 10.13:}\ 
{\small\bf \un{N.B.}}  \ 
1, $n$-absolute continuity coincides with $n$-absolute continuity.
\\

\qquad{\bf 10.14:}\ 
{\small\bf NOTATION} \ 
Put
\[
V_{\lambda, n} (f, \Omega) 
\ = \ 
\sup \ 
\left\{
\sum\limits_{j = 1}^k \ (\osc (f, B(x_j, \lambda \hsy r_j)))^n
\right\}.
\]
\\[-.75cm]

\qquad{\bf 10.15:}\ 
{\small\bf DEFINITION} \ 
$f$ is of \un{bounded $\lambda$, $n$-variation} in $\Omega$ if 
\[
V_{\lambda, n} (f, \Omega) 
\ < \ 
+\infty.
\]
\\[-.75cm]

\qquad{\bf 10.16:}\ 
{\small\bf NOTATION} \ 
$\BV^{\lambda, n} (\Omega)$ is the set of all functions $f : \Omega \ra \R^n$ of bounded $\lambda$, $n$-variation in $\Omega$.
\\

\qquad{\bf 10.17:}\ 
{\small\bf NOTATION} \ 
$\AC^{\lambda, n} (\Omega)$ is the set of all functions $f \in \BV^{\lambda, n} (\Omega)$ which are $\lambda$, $n$-absolutely continuous.
\\

\qquad{\bf 10.18:}\ 
{\small\bf \un{N.B.}} \ 
\begin{align*}
\AC^n (\Omega) \ 
&=\ 
\AC^{1, n} (\Omega)
\\[11pt]
&\subset
\AC^{\lambda, n} (\Omega)
\qquad (0 < \lambda < 1)
\end{align*}
and it can be shown that the containment is proper.
\\

\qquad{\bf 10.19:}\ 
{\small\bf LEMMA} \ 
Let $f : \Omega \ra \R^n$  and suppose that $0 < \lambda_1 < \lambda_2 < 1$ 
$-$then $f$ is $\lambda_1$, $n$-absolutely continuous iff 
$f$ is $\lambda_2$, $n$-absolutely continuous.
\\

\qquad{\bf 10.20:}\ 
{\small\bf THEOREM} \ 
Suppose that $0 < \lambda_1 < \lambda_2 < 1$ $-$then
\[
\BV^{\lambda_1, n} (\Omega)
\ = \ 
\BV^{\lambda_2, n} (\Omega).
\]
\\[-.75cm]

\qquad{\bf 10.21:}\ 
{\small\bf SCHOLIUM} \ 
There are but two classes of $\lambda$, $n$-absolutely continuous functions, viz.  
those corresponding to $\lambda = 1$ and to $0 < \lambda < 1$.  
\\

\qquad{\bf 10.22:}\ 
{\small\bf DEFINITION} \ 
Let $f : \Omega \ra \R^n$  and suppose that $0 < \lambda< 1$ $-$then $f$ satisfies the 
\un{condition $RR_\lambda$} if there is a nonnegative function $\phi \in \Lp^1 (\Omega)$, a so-called \un{weight}, such that
\[
(\osc (f, B(x, \lambda \hsy r)))^n 
\ \leq \ 
\int\limits_{B(x,r)} \ \phi \ \td \Lm^n
\]
for every $B(x,r) \subset \Omega$. 
\\

[Note:  \ 
Formally, $RR_1 = RR$.]
\\

\qquad{\bf 10.23:}\ 
{\small\bf NOTATION} \ 
Denote by $RR^{\lambda, n} (\Omega)$ the set of all functions $f : \Omega \ra \R^n$ which satisfy the condition $RR_\lambda$.
\\

\qquad{\bf 10.24:}\ 
{\small\bf THEOREM} \ 
\[
RR^{\lambda, n} (\Omega)
\ = \ 
\AC^{\lambda, n} (\Omega).
\]
\\[-.75cm]

\qquad{\bf 10.25:}\ 
{\small\bf THEOREM} \ 
Let $f \in B^{\lambda, n} (\Omega)$ $-$then $f$ is differentiable almost everywhere.
\\

Return now to the beginning.
\\

\qquad{\bf 10.26:}\ 
{\small\bf LEMMA}  \ 
Suppose that $\Omega$ is bounded and $f : \Omega \ra \R^n$ is Lipschitz, say
\[
\norm{f(x_1) - f(x_2)} 
\ \leq \ 
M \hsx \ \norm{x_1 - x_2} 
\]
for all $x_1$, $x_2 \in \Omega$ $-$then $f \in R R^n (\Omega)$.

[Define $\phi : \Omega \ra \R$ by the rule
\[
\phi(x) 
\ = \ 
\frac{M^n \hsy 2^n}{\omega_n} 
\qquad (\implies \phi \in \Lp^1(\Omega)).
\]
Then for any $B(x, r) \subset \Omega$,

\begin{align*}
\int\limits_{B(x,r)} \ \phi \ \td \Lm^n \ 
&=\ 
\frac{M^n \hsy 2^n}{\omega_n} \ 
\int\limits_{B(x,r)} \ 1\ \td \Lm^n
\\[15pt]
&=\ 
\frac{M^n \hsy 2^n}{\omega_n} \ 
\omega_n \hsy r^n
\\[15pt]
&=\
M^n \hsy 2^n \hsy r^n
\end{align*}

\hspace{.7cm} $\implies$

\[
\bigg(
\int\limits_{B(x,r)} \ \phi \ \td \Lm^n
\bigg)^{1/n}
\ = \ 
M \hsx (2 \hsy r).
\]
But

\[
x_1, \hsx x_2 \in B(x, r) 
\implies 
\norm{x_1 - x_2} 
\ \leq \ 
2 \hsy r
\]

\hspace{.7cm} $\implies$
\[
\norm{f(x_1) - f(x_2)} 
\ \leq \ 
M \hsx \ \norm{x_1 - x_2} 
\ \leq \ 
M \hsx (2 \hsy r)
\]

\hspace{.7cm} $\implies$
\begin{align*}
\osc (f, B(x,r)) \ 
&\leq \ 
M( 2\hsy r) 
\\[11pt]
&=\ 
\bigg(
\int\limits_{B(x,r)} \ \phi \ \td \Lm^n
\bigg)^{1/n}
\end{align*}

\hspace{.7cm} $\implies$

\[
(\osc (f, B(x,r)))^n 
\ \leq \ 
\int\limits_{B(x,r)} \ \phi \ \td \Lm^n.]
\]
\\[-.75cm]

\qquad{\bf 10.27:}\ 
{\small\bf LEMMA} \ 
Suppose that $\Omega$ is bounded and $f \in W^{1, p} (\Omega; \R^n)$ $(p > n)$ is continuous 
$-$then $f \in RR^n (\Omega)$.  
\\[-.5cm]

[Upon consideration of components, one can take $n = 1$.  
This said, for any $B(x, r) \subset \Omega$, (cf. 7.4.11), 
\[
\osc (f, B(x,r)) 
\ \leq \ 
C \hsy r^{1 - n/p} \ 
\bigg(
\int\limits_{B(x,r)} \ \norm{\nabla f}^p \ \td \Lm^n
\bigg)^{1/p}
\]

\hspace{.25cm} $\implies$
\begin{align*}
(\osc (f, B(x,r)))^n  \ 
&\leq \ 
C \hsy r^{n(1 - n/p)} \ 
\bigg(
\int\limits_{B(x,r)} \ \norm{\nabla f}^p \ \td \Lm^n
\bigg)^{n/p}
\\[18pt]
&\leq \ 
C \hsx
\bigg(
r^n + \ 
\int\limits_{B(x,r)} \ \norm{\nabla f}^p \ \td \Lm^n
\bigg)
\\[18pt]
&\leq \ 
C \hsx
\int\limits_{B(x,r)} \ 
\left(
1 + \norm{\nabla f}^p
\right)
\td \Lm^n.
\end{align*}
So, for the weight, one can take
\[
\phi 
\ = \ 
C (1 + \norm{\nabla f}^p).]
\]

[Note: \ 
The usual convention on the constant ``$C$'' is in force, i.e., it may change from line to line.]

\chapter{
\text{SECTION 11: \quad MISCELLANEA}
\vspace{.75cm}\\
$\boldsymbol{\S}$\textbf{11.1}.\quad  PROPERTY (N)}
\setlength\parindent{2em}
\renewcommand{\thepage}{11-\S1-\arabic{page}}

\vspace{-.5cm}
\qquad
Let $\Omega$ be a nonempty open subset of $\R^n$.
\\

\qquad{\bf 11.1.1:}\quad 
{\small\bf DEFINITION} \ 
A function $f : \Omega \ra \R^n$ is said to have \un{property (N)} if $f$ sends sets of Lebesgue 
measure 0 to sets of Lebesgue measure 0:
\[
\Lm^n (E) 
\ = \ 
0
\quad (E \subset \Omega) 
\ \implies \ 
\Lm^n (f(E))
\ = \ 0.
\]
\\[-1.25cm]

\qquad{\bf 11.1.2:}\quad
{\small\bf SUBLEMMA} \ 
If $E \in \sM_\Lm^n$, then there exists an $F_\sigma$-set $F \subset E$ such that $\Lm^n (E \backslash F) = 0$.  
Choose next a countable collection of compact subsets $C_j$ for which 
$F = \bigcup\limits_j \ C_j$ and put 
$K_j = \bigcup\limits_{k = 1}^j \ C_k$, thus $\{K_j\}$ is an increasing sequence of compact sets with 
$\bigcup\limits_j \ K_j = \bigcup\limits_j \ C_j = F$.  
Finally, since $E$ is the disjoint union of $F$ and $E \backslash F$, we have 

\vspace{-.75cm}
\allowdisplaybreaks
\begin{align*}
\Lm^n (E) 
&= \ 
\Lm^n (F) + \Lm^n (E \backslash F) 
\\[11pt]
&= \ 
\Lm^n (F) 
\\[11pt]
&= \ 
\lim\limits_{j \ra \infty} \ \Lm^n (K_j ).
\end{align*}
\\[-1cm]

\qquad{\bf 11.1.3:}\quad
{\small\bf LEMMA} \ 
Suppose that $f : \Omega \ra \R^n$ has property (N) $-$then the implication

\[
E \in \sM_\Lm^n \quad (E \subset \Omega) 
\implies
f(E) \in \sM_\Lm^n
\]
obtains.
\\[-.5cm]

PROOF \ 
As above, write
\[
E 
\ = \ 
F \cup (E \backslash F) 
\quad (F = \bigcup\limits_j \ K_j).
\]
Then
\begin{align*}
f(E) \ 
&=\ 
f(F) \cup f(E \backslash F)
\\[11pt]
&=\ 
\bigcup\limits_j \ f(K_j) + f(E / F).
\end{align*}
Since $f$ is continuous, the $f(K_j)$ are compact, hence measurable, so the union $\bigcup\limits_j \ f(K_j)$ is measurable.  
On the other hand, 
\begin{align*}
\Lm^n (E \backslash F) = 0 
&\implies
\Lm^n (f(E \backslash F)) = 0 
\\[11pt]
&\implies
f(E \backslash F) \in \sM_\Lm^n.
\end{align*}
All told therefore, 
\[
f(E) \in \sM_\Lm^n.
\]
\\[-1cm]

\qquad{\bf 11.1.4:}\quad 
{\small\bf LEMMA} \ 
Take $n = 1$, $\Omega = \hsx ]a,b[$, and suppose that $f : \hsx ]a,b[ \ra \R$ is absolutely continuous $-$then 
$f$ has property (N).
\\[-.25cm]

\qquad{\bf 11.1.5:}\quad 
{\small\bf EXAMPLE} \ 
If  $f : \Omega \ra \R^n$ is locally Lipschitz, then $f$ has property (N) (cf. 2.3.23).
\\[-.5cm]

[Note: \ 
In particular: A $C^\prime$-function $f$ has property (N) (being locally Lipschitz),]
\\[-.25cm]

\qquad{\bf 11.1.6:}\quad 
{\small\bf \un{N.B.}}  \ 
The preceding consideration is false if $f$ is merely continuous or even H\"older continuous with exponent $0 < \alpha < 1$.
\\[-.5cm]

\begin{spacing}{1.75}
[The Cantor function $f$ sends the Cantor set 
$C$ $(\Lm^1(C) = 0)$ to $f(C)$ $(\Lm^1(f(C)) = 1)$.  
And $f$ is H\"older continuous with exponent 
$\ds \alpha = \frac{\elog 2}{\elog 3}$.]
\\
\end{spacing}

\qquad{\bf 11.1.7:}\quad 
{\small\bf RAPPEL} \ (Vitali) \ 
Let $\sB$ be a system of closed balls in $\R^n$ such that 
\[
\sup \{\diam (B) \hsy : \hsy B \in \sB\} 
\ < \ 
+\infty.
\]
Then there exists a pairwise disjoint, at most countable subsystem $\{B(x_i, r_i)\} \subset \sB$ such that 
\[
\bigcup\limits_{B \in \sB} 
\subsetx 
\bigcup\limits_i \ B(x_i, 5 \hsy r_i).
\]
\\[-.75cm]

\qquad{\bf 11.1.8:}\quad 
{\small\bf THEOREM} \ 
Suppose that  $f : \Omega \ra \R^n$ is $n$-absolutely continuous $-$then  $f$ has property (N).
\\[-.5cm]

PROOF \ 
Fix an $E \subset \Omega$ of Lebesgue measure 0.  
Given $\varepsilon > 0$, choose $\delta > 0$ per the definition of $n$-absolute continuity, subject to $\delta < \varepsilon$.  
Let $G \subset \Omega$  be an open set containing $E$ with $\Lm^n(G) < \delta$.  
Given an $x \in E$, choose $r(x) > 0$ such that
\[
B(x, r(x)) \subset G, 
\quad 
r(x) < \frac{\varepsilon}{10}, 
\quad \text{and} \ 
\rho(x) \equiv \osc(f, B(x, r(x))) < \frac{\varepsilon}{10}.
\]
Using Vitali, determine a disjoint system

\[
\{B(f(x_i), \rho(x_i))\} 
\subset 
\{B(f(x), \rho(x)) \hsy : \hsy x \in E\} 
\]
such that 

\[
f(E) 
\subset
\bigcup\limits_i \ 
B(f(x_i), 5 \hsy \rho(x_i)).
\]
Since

\[
f (B(x_i, r(x_i))) 
\subset
B(f(x_i), \rho(x_i)),
\]
the $B(x_i, r(x_i))$ are pairwise disjoint, hence

\allowdisplaybreaks
\begin{align*}
\sH_\varepsilon^n(f(E)) \ 
&\leq 
C \hsx \sum\limits_i \ \rho(x_i)^n
\\[15pt]
&\leq \ 
C \hsx \sum\limits_i \ \osc(f, B(x_i, r(x_i))^n
\\[15pt]
&\leq \ 
C \hsy \varepsilon.
\end{align*}
Now let $\varepsilon \ra 0$ to conclude that 
\[
\sH^n (f(E)) 
\ = \ 
0
\]
or still, that
\[
\Lm^n (f(E)) 
\ = \ 
0.
\]
\\[-.75cm]

\qquad{\bf 11.1.9:}\quad 
{\small\bf APPLICATION} \ 
Suppose that $f \in W^{1, p} (\Omega ; \R^n)$ $(p > n)$ is continuous $-$then $f$ has property (N).
\\

[In fact, 
\[
f \in R R^n (\Omega) 
\quad \text{(cf. 10.27).}
\]
But
\[
R R^n (\Omega) 
\ = \ 
\AC^n (\Omega) 
\quad \text{(cf. 10.10).]}
\]
\\[-.75cm]

\qquad{\bf 11.1.10:}\quad 
{\small\bf REMARK} \ 
There are continuous functions in $W^{1, n} (\Omega ; \R^n)$ $(n > 1)$ that do not have property (N).
\\[-.5cm]

[E.g., it is possible to construct a continuous $f \in W^{1, n} (\R^n ; \R^n)$ $(n > 1)$ which sends $[0,1]$ to $[0,1]^n$.  
Therefore $f$ does not have property (N).]
\\

\qquad{\bf 11.1.11:}\quad
{\small\bf THEOREM} \ 
If $f \in W^{1, n} (\Omega ; \R^n)$ is continuous and open, then $f$ has  property (N).
\\

\qquad{\bf 11.1.12:}\quad 
{\small\bf \un{N.B.}}  \ 
There exists a homeomorphism
\[
f \in W^p ((]-1, 1[)^n, (]-1, 1[)^n) 
\quad (p < n)
\]
which does not have property (N).
\\

\qquad{\bf 11.1.13:} 
{\small\bf THEOREM} \ 
If $f \in W^{1, n} (\Omega ; \R^n)$ is H\"older continuous, then $f$ has property (N).
\\[1cm]


\begin{center}
* \ * \ * \ * \ * \ * \ * \ * \ * \ * \ * \ * \ 
\end{center}

\[
\text{APPENDIX}
\]
\\[-.75cm]

{\small\bf LEMMA} \ 
Let $1 \leq k \leq n$, let $\Omega$ be a nonempty open subset of $\R^k$, and let 
$T : \Omega \ra \R^n$ be continuous and one-to-one $-$then
\[
E \in \sB(\R^k) 
\quad (E \subset \Omega) 
\ \implies \ 
T(E) \in \sB(\R^n).
\]

PROOF \ 
$\Omega$ is a $\sigma$-compact subset of $\R^k$, hence $T(\Omega)$ is a $\sigma$-compact subset of $\R^n$ 
($T$ being continuous), 
hence $T(\Omega) \in \sB (\R^n)$.  
Let now 
\[
\sA 
\ = \ 
\{E \subset \Omega : T(E) \in \sB(\R^n)\},
\]
a $\sigma$-algebra of subsets of $\Omega$ (as regards complementation, 
note that $T(\Omega \backslash E) = T(\Omega) \backslash T(E)$, 
$T$ being one-to-one).
It is clear that $\sA$ contains the open subsets of $\Omega$, (per the initial observation), 
so $\sA$ contains the Borel $\sigma$-algebra $\sB(\Omega)$.  
But $\sB(\sA) = \sB(\R^k) \cap \Omega$, thus $\forall \ E \in \sB(\R^k) \cap \Omega$, $T(E) \in \sB(\R^n)$.

\chapter{
$\boldsymbol{\S}$\textbf{11.2}.\quad  THE MULTIPLICITY FUNCTION}
\setlength\parindent{2em}
\renewcommand{\thepage}{11-\S2-\arabic{page}}

\vspace{-.25cm}
\qquad
Let $\Omega$ be a nonempty open subset of $\R^n$ and let $f : \Omega \ra \R^n$ be a continuous function.
\\[-.25cm]

\qquad{\bf 11.2.1:}\quad
{\small\bf DEFINITION} \ 
If $E \subset \Omega$, then 
\[
N(f, E, y) \ = \ 
\#\{x \in E \hsy : \hsy f(x) = y\}
\]
or still, 
\[
N(f, E, y) \ = \ \sH^0 (E \cap f^{-1} (y) )
\]
is the \un{multiplicity function} of $f$ at $y \in \R^n$ w.r.t. $E$.
\\[-.5cm]

[Note: \ 
$N(f, E, y)$ is the cardinality of $E \cap f^{-1} (y) $ and if this set is infinite, then we put
\[
N(f, E, y) \ = \ +\infty.]
\]
\\[-1cm]

\qquad{\bf 11.2.2:}\quad
{\small\bf LEMMA} \ 
\[
E_1 \subset E_2 
\implies 
N(f, E_1, y) 
\ \leq \ 
N(f, E_2, y).
\]
\\[-1cm]

\qquad{\bf 11.2.3:}\quad
{\small\bf LEMMA} \ 
If $\{E_k\}$ is an increasing sequence of subsets of $\Omega$, then 
\[
N(f, E, y)
\ = \ 
\lim\limits_{k \ra \infty} \ N(f, E_k, y), 
\]
where $\ds E = \bigcup\limits_{k = 1}^\infty \ E_k$.
\\[0.25cm]

\qquad{\bf 11.2.4:}\quad
{\small\bf THEOREM} \ 
Suppose that $f : \Omega \ra \R^n$ has property $(N)$ $-$then for any Lebesgue measurable set $E \subset \Omega$, 
the multiplicity function
\[
y \ra N(f, E, y)
\]
is Lebesgue measurable in $\R^n$.
\\[-,5cm]

PROOF \ 
Take $E$ bounded and for every $m \in \N$ construct a partition of $E$ into pairwise nonintersecting 
measurable sets
\[
E_1^{(m)}, \ldots, E_{k_m}^{(m)}
\]
such that 

\[
\diam (E_i^{(m)}) 
\ \leq \ 
\frac{1}{m}
\quad (i = 1, \ldots, k(m)).
\]
Put

\[
N(f, m, -)
\ = \ 
\chisubfofEonem + \cdots + \chisubfofEim
\]
and note that each of the sets $f(E_i^{(m)})$ is measurable 
(since $f$ has property (N)), 
hence $N(f, m, -)$ is measurable.  
Accordingly it need only be shown that

\[
\limsup\limits_{m \ra \infty} \ 
N(f, m, y) 
\ = \ 
N(f, E, y)
\]
to establish the contention.  
Given $y \in \R^n$, there are two possibilities for $E \cap f^{-1} (y)$: 
It is either finite or infinite.  
To treat the first of these, say

\[
E \cap f^{-1} (y) 
\ = \ 
\{x_1, \ldots, x_k\},
\]
take
\[
m_0 
\ > \ 
\min\limits_{i \neq j} \ \frac{1}{\norm{x_i - x_j}}.
\]
If $m > m_0$, then none of the $E_i^{(m)}$ contain two distinct $x_r$, $x_s$ , so it can be assumed that
\[
x_1 \in E_1^{(m)}, \ldots, x_k \in E_k^{(m)}.
\]
Next, $\forall \ m > m_0$, 

\[
N(f, m, y) 
\ \geq \ 
N(f, E, y)
\]

\qquad
$\implies$ 
\[
\liminf\limits_{m \ra \infty} \ 
N(f, m, y) 
\ \geq \ 
N(f, E, y).
\]
On the other hand, 
\[
N(f, m, y) 
\ \leq \ 
N(f, E, y)
\]

\qquad
$\implies$ 
\[
\limsup\limits_{m \ra \infty} \ 
N(f, m, y) 
\ \leq \ 
N(f, E, y).
\]
Therefore
\[
\lim\limits_{m \ra \infty} \ 
N(f, m, y) 
\ = \ 
N(f, E, y).
\]
\\[-.75cm]

\qquad{\bf 11.2.5:}\quad
{\small\bf LEMMA} \ 
If $f : \Omega \ra \R^n$ is continuous and open, then for every open set $G \subset \Omega$, the function 
$y \ra N(f, G, y)$ is lower semicontinuous in $\R^n$.

\chapter{
$\boldsymbol{\S}$\textbf{11.3}.\quad  JACOBIANS}
\setlength\parindent{2em}
\setcounter{theoremn}{0}
\renewcommand{\thepage}{11-\S3-\arabic{page}}

\qquad
Let $\Omega$ be a nonempty open subset of $\R^n$.
\\

\qquad{\bf 11.3.1:}\quad
{\small\bf DEFINITION} \ 
Let $f = (f^1, \ldots, f^n) \in W^{1, n} (\Omega; \R^n)$ $-$then the \un{Jacobian} of $f$, denoted 
$\tJ_f$, is the determinant
\[
\det (\nabla f^1, \ldots, \nabla f^n).
\]
\\[-1.cm]

\qquad{\bf 11.3.2:}\quad
{\small\bf \un{N.B.}} \ 
The coordinate functions $f^i$ $(1 \leq i \leq n)$ of $f$ and their first order distributional derivatives belong to 
$\Lp^n (\Omega)$.
\\[-.5cm]

[Note: \ 
Nevertheless, an element of $W^{1, n} (\Omega; \R^n)$ may be nowhere continuous, hence nowhere differentiable.]
\\

\qquad{\bf 11.3.3:}\quad
{\small\bf THEOREM} \ 
If $f \in W^{1, n} (\Omega; \R^n)$, then $\tJ_f \in \Lp^1 (\Omega)$.
\\[-.5cm]

PROOF \ 
$\tJ_f$ is a sum of terms, each of which is (plus or minus) the product of $n$ weak partial derivatives of the components of $f$ 
and, 
as noted above, each of these is in $\Lp^n (\Omega)$.  
The product of $n$ $\Lp^n (\Omega)$ functions is in $\Lp^1 (\Omega)$ (apply H\"older), hence $f \in \Lp^1 (\Omega)$.
\\

\qquad{\bf 11.3.4:}\quad
{\small\bf FACT} \ 
If $f : \R^n \ra \R^n$ is Lipschitz, then for any $E \in \sM_\Lm^n$, 
\[
\int\limits_E \ 
\abs{\tJ_f} 
\ \td \Lm^n 
\ = \ 
\int\limits_{\R^n} \ 
N(f, E, y) 
\ \td y.
\]
\\[-1cm]

[This will be established in \S12.5.]
\\

\qquad{\bf 11.3.5:}\quad
{\small\bf RAPPEL} \ 
Let $f \in W^{1, p} (\Omega)$ $-$then there is a partition
\[
\Omega 
\ = \ 
\big(
\bigcup\limits_{k = 1}^\infty \ E_k
\big)
\cup Z,
\]
where the $E_k$ are Lebesgue measurable sets such that $\restr{f}{E_k}$ is Lipschitz and $Z$ has Lebesgue measure 0 
(cf. 7.1.5).
\\

\qquad{\bf 11.3.6:}\quad
{\small\bf THEOREM} \ 
Suppose that $f \in W^{1, n} (\Omega; \R^n)$ has property $(N)$ $-$then 
\[
\int\limits_\Omega \ 
\abs{\tJ_f} 
\ \td \Lm^n 
\ = \ 
\int\limits_{\R^n} \ 
N(f, \Omega, y)
\ \td y.
\]
\\[-1cm]

PROOF \ 
The foregoing decomposition principle extends from $W^{1, n} (\Omega)$ to $W^{1, n} (\Omega; \R^n)$ 
and the Lipschitz function 
$\restr{f}{E_k}$ extends to a Lipschitz function $f_k : \R^n \ra \R^n$, hence per supra

\[
\int\limits_{E_k} \ 
\abs{\tJ_{f_k}} 
\ \td \Lm^n 
\ = \ 
\int\limits_{\R^n} \ 
N(f_k, E_k, y)
\ \td y.
\]
Put $E_0 = \bigcup\limits_{k = 1}^\infty \ E_k$ $-$then $\Omega = E_0 \cup Z$ $(\Lm^n (Z) = 0)$, so

\allowdisplaybreaks
\begin{align*}
\int\limits_\Omega \ 
\abs{\tJ_f} 
\ \td \Lm^n \ 
&=\ 
\int\limits_E \ 
\abs{\tJ_f} 
\ \td \Lm^n \ 
\\[15pt]
&=\ 
\sum\limits_{k = 1}^\infty \ 
\int\limits_{E_k} \ 
\abs{\tJ_f}
\ \td \Lm^n \ 
\\[15pt]
&=\ 
\sum\limits_{k = 1}^\infty \ 
\int\limits_{E_k} \ 
\abs{\tJ_{f_k}}
\ \td \Lm^n \ 
\\[15pt]
&=\ 
\sum\limits_{k = 1}^\infty \ 
\int\limits_{\R^n} \ 
N(f_k, E_k, y) 
\ \td y 
\\[15pt]
&=\ 
\int\limits_{\R^n} \ 
N(f, E, y) 
\ \td y 
\\[15pt]
&\leq\ 
\int\limits_{\R^n} \ 
N(f, \Omega, y) 
\ \td y 
\\[15pt]
&=\ 
\int\limits_{\R^n} \ 
N(f, E, y) 
\ \td y 
\hsx + \hsx 
\int\limits_{\R^n} \ 
N(f, Z, y) 
\ \td y 
\\[15pt]
&=\ 
\int\limits_{\R^n} \ 
N(f, E, y) 
\ \td y 
\hsx + \hsx 
\int\limits_{f(Z)} \ 
N(f, Z, y) 
\ \td y 
\\[15pt]
&=\ 
\int\limits_{\R^n} \ 
N(f, E, y) 
\ \td y 
\qquad (\Lm^n (f(Z) = 0)
\end{align*}

$\implies$
\[
\int\limits_\Omega \ 
\abs{\tJ_f} 
\ \td \Lm^n 
\ = \ 
\int\limits_{\R^n} \ 
N(f, \Omega, y)
\ \td y.
\]

[Note: \ 
\[
f(Z) 
\ = \ 
\{y \hsy : \hsy N(f, Z, y) \neq 0\}.]
\]
\\[-.75cm]

\qquad{\bf 11.3.7:}\quad
{\small\bf \un{N.B.}} \ 
The assumption that $f$ has property $(N)$ implies that the relevant multiplicity functions are Lebesgue measurable.
\\

\qquad{\bf 11.3.8:}\quad
{\small\bf THEOREM} \ 
If $f \in W^{1, n} (\Omega; \R^n)$  is continuous and if $\tJ_f > 0$ almost everywhere in $\Omega$, then 
$f$ has property $(N)$.
\\

\qquad{\bf 11.3.9:}\quad
{\small\bf REMARK} \ 
Examples have been constructed of continuous functions $f \in W^{1, n} (\Omega; \R^n)$ such that 
$\tJ_f = 0$ almost everywhere in $\Omega$ but such that $f$ fails to have property $(N)$.
\\

On general grounds, an $f \in W^{1, p} (\Omega; \R^n)$ $(1 \leq p < +\infty)$ is approximately differentiable almost 
everywhere in $\Omega$.  
More is true if $p = n$, namely
\[
\ap \hsx - \td f
\]
is ``regular'' (i.e., ``$E$'' can be written as a union of concentric spheres centered at $x$).

\chapter{
\text{SECTION 12: \quad AREA FORMULAS}
\vspace{.75cm}\\
$\boldsymbol{\S}$\textbf{1}.\quad  THE LINEAR CASE}
\setlength\parindent{2em}
\renewcommand{\thepage}{12-\S1-\arabic{page}}

\qquad
\qquad{\bf 12.1.1:}\quad
{\small\bf RAPPEL} \ 
Let $T: \R^n \ra \R^n$ be a nonsingular linear transformation $-$then
\[
E \in \sM_\Lm^n 
\implies
T (E) \in \sM_\Lm^n 
\]
and
\[
\Lm^n (T (E)) 
\ = \ 
\abs{\det (M_T)} \hsx \Lm^n (E).
\]
\\[-.75cm]

\qquad{\bf 12.1.2:}\quad
{\small\bf \un{N.B.}} \ 
This is the simplest instance of what is known as an ``area formula''.  
As will be shown below, it leads to a ``change of variable formula''.
\\

Retaining $T$ and $E$, suppose given a function $f: E \ra [-\infty, +\infty]$.
\\

\qquad{\bf 12.1.3:}\quad
{\small\bf LEMMA} \ 
If $f$ is Lebesgue measurable on $E$, then $f \circ T$ is Lebesgue measurable on $T^{-1} (E)$.
\\

\qquad{\bf 12.1.4:}\quad
{\small\bf THEOREM} \ 
\[
\int\limits_E \ f \ \td \Lm^n 
\ = \ 
\abs{\det (M_T)} \ 
\int\limits_{T^{-1} (E)} \ f \circ T \ \td \Lm^n \ 
\]
in the sense that if one of the two integrals exists then so does the other and the two are equal.
\\[-.5cm]

PROOF \ 
Let $S \in \sM_\Lm^n$ $-$then
\[
T^{-1} (E \cap S) 
\ = \ 
T^{-1} (E) \cap T^{-1} (S)
\]

\hspace{1.cm} $\implies$
\\[-1.5cm]

\allowdisplaybreaks
\begin{align*}
\Lm^n (E \cap S)  \ 
&=\ 
\Lm^n ((T \circ T^{-1}) \quad (E \cap S))
\\[11pt]
&=\ 
\Lm^n (T (T^{-1} (E \cap S)))
\\[11pt]
&=\ 
\Lm^n (T (T^{-1} (E) \cap  T^{-1}(S)))
\\[11pt]
&=\ 
\abs{\det (M_T)} \hsx \Lm^n (T^{-1} (E) \cap  T^{-1}(S)).
\end{align*}

Take for $f$ the characteristic function $\chisubS$ of $S$, hence
\\[-1.25cm]

\allowdisplaybreaks
\begin{align*}
\int\limits_E \ f \ \td \Lm^n \ 
&=\ 
\int\limits_E \ \chisubS \ \td \Lm^n
\\[11pt]
&=\ 
\Lm^n (E \cap S)
\\[11pt]
&=\ 
\abs{\det (M_T)} \hsx
\Lm^n (T^{-1} (E) \cap T^{-1} (S))
\\[11pt]
&=\ 
\abs{\det (M_T)} \hsx
\int\limits_{T^{-1} (E)} \ \chisubTinvS \ \td \Lm^n
\\[11pt]
&=\ 
\abs{\det (M_T)} \hsx  
\int\limits_{T^{-1} (E)} \ \chisubS \circ T \ \td \Lm^n
\\[11pt]
&=\ 
\abs{\det (M_T)} \hsx 
\int\limits_{T^{-1} (E)} \ f \circ T \ \td \Lm^n.
\end{align*}
\\[-.75cm]

One can then proceed from here to a nonnegative simple function on $E$ and then to a nonnegative extended real valued 
Lebesgue measurable function on $E$ and finally to the general case 
(write $f = f^+ - f^-$ and work separately with $f^+$ and $f^-$).  
\\[-.5cm]

[Note: \ 
By way of a justification, monotone convergence is used when coupled with the fact that there exists an increasing sequence 
$\{f_j\}$ of nonnegative simple functions such that $f_j \uparrow f$.]
\\

\qquad{\bf 12.1.5:}\quad
{\small\bf REMARK} \ 
Matters can be restated, viz.

\[
\int\limits_{T (E)} \ f \ \td \Lm^n 
\ = \ 
\abs{\det (M_T)} \ 
\int\limits_E \ f \circ T \  \td \Lm^n,
\]
the underlying supposition being that in this context, $f: T(E) \ra [-\infty, +\infty]$ is Lebesgue measurable.
\\

Assume: \ $k$, $n \in \N$, $1 \leq k \leq n$.
\\

\qquad{\bf 12.1.6:}\quad
{\small\bf DEFINITION} \ 
Suppose that  $T: \R^k \ra \R^n$ is a linear transformation $-$then the \un{adjoint} of $T$ is the linear tranformation 
$T^*: \R^n \ra \R^k$ characterized by the condition
\[
\langle x, T^* y \rangle
\ = \ \langle T x, y \rangle
\]
for all $x \in \R^k$ and for all $y \in \R^n$.
\\

[Note: \ 
In terms of matrices,
\[
M_{T^*}
\ = \ 
M_T^\Tee, 
\]
the transpose of $M_T$.]
\\

\qquad{\bf 12.1.7:}\quad
{\small\bf NOTATION} \ 
Given a linear transformation $T: \R^k \ra \R^n$, put
\[
J (T) 
\ = \ 
\sqrt{ \det (T^* \hsy T)}.
\]

[Note: \ 
$J (T)$ is nonzero iff $T$ is nonsingular.]
\\

\qquad{\bf 12.1.8:}\quad
{\small\bf \un{N.B.}}  \ 
If $k = n$, then 

\begin{align*}
\det (T^* \hsy T) \ 
&=\ 
\det (M_{T^* \hsy T})
\\[11pt]
&=\ 
\det (M_{T^*} \hsy M_T)
\\[11pt]
&=\ 
\det (M_T^\Tee \hsy M_T)
\\[11pt]
&=\ 
\det (M_T^\Tee) \hsx \det (M_T)
\\[11pt]
&=\ 
\det (M_T)  \det (M_T)
\\[11pt]
&=\ 
\det (M_T)^2
\end{align*}

\hspace{1.5cm} $\implies$
\begin{align*}
J (T) \ 
&=\ 
\sqrt{\det (T^* \hsy T)}
\\[11pt]
&=\ 
\sqrt{\det (M_T)^2}
\\[11pt]
&=\ 
\abs{\det (M_T)}
\\[11pt]
&=\ 
\abs{\det (T)}.
\end{align*}
\\[-.75cm]

\qquad{\bf 12.1.9:}\quad
{\small\bf DEFINITION} \ 

\qquad \textbullet \quad
A linear map $U: \R^k \ra \R^n$  is said to be \un{orthogonal} if 
$\langle Ux, Uy \rangle = \langle x, y \rangle$ for all $x$, $y \in \R^k$.
\\[-.5cm]

\qquad \textbullet \quad
A linear map $S: \R^k \ra \R^k$  is said to be \un{symmetric} if 
$\langle x, Sy \rangle = \langle Sx, y \rangle$ for all $x$, $y \in \R^k$.
\\

\qquad{\bf 12.1.10:}\quad
{\small\bf POLAR DECOMPOSITION} \ 
Let $T: \R^k \ra \R^n$ be an injective linear transformation $-$then there exists a symmetric map 
$S: \R^k \ra \R^k$ and an orthogonal map U$: \R^k \ra \R^n$ such that $T = U \hsy S$.
\\

\qquad{\bf 12.1.11:}\quad
{\small\bf \un{N.B.}} \ 
If $T: \R^k \ra \R^n$ is an injective linear transformation and $E \in \sM_\Lm^k$, then $T(E) \in \sM(\sH^k)$ and 
\[
\sH^k (T(E)) 
\ = \ 
J (T) \Lm^k (E).
\]

PROOF \ 
To establish the purported equality, consider first the case when $k = n$, thus
\\[-1.25cm]

\allowdisplaybreaks
\begin{align*}
\sH^n (T(E)) \ 
&=\ 
\Lm^n (T(E))
\\[11pt]
&=\ 
\abs{ \det (M_T)} \hsy \Lm^n(E)
\\[11pt]
&=\ 
\abs{ \det (T)} \hsy \Lm^n(E)
\\[11pt]
&=\ 
J (T)\hsy \Lm^n(E).
\end{align*}
Suppose now that $k < n$, write
\begin{align*}
\sH^k (T(E))  \ 
&=\ 
\sH^k (U \hsy S(E))  
\\[11pt]
&=\ 
\sH^k (S(E)),
\end{align*}
$U$ being an isometry.  
But
\begin{align*}
\sH^k (S(E)) \ 
&=\ 
\Lm^k (S(E))
\\[11pt]
&=\ 
\abs{ \det (S)} \hsy \Lm^k(E).
\end{align*}
And
\begin{align*}
T ^* \hsy T \ 
&=\ 
S^* \hsy U^* \hsy U \hsy S
\\[11pt]
&=\ 
S^* \hsy S \qquad (U^* \hsy U = \text{id})
\\[11pt]
&=\ 
S^2 \qquad (S^* = S)
\end{align*}

\hspace{.5cm} $\implies$

\[
\det (T^* \hsy T) 
\ = \ 
\det (S)^2
\]

\hspace{.5cm} $\implies$

\[
J (T)
\ = \ 
\sqrt{ \det (T^* \hsy T)}
\ = \ 
\sqrt{\det (S)^2}
\ = \ 
\abs{\det (S)}.
\]
\\[-.75cm]

\qquad{\bf 12.1.12:}\quad
{\small\bf REMARK} \ 
If $T: \R^k \ra \R^n$ is Lipschitz, then (cf. 12.3.1)
\[
E \in \sM_\Lm^k 
\implies 
T (E) \in \sM (\sH^k).
\]

[Note: \ 
\[
\text{$T$ linear $\implies$ $T$ Lipschitz}.
\]

Proof: \ 
\[
\norm{T x - T y} 
\ \leq \ 
\norm{T} \hsy \norm{x - y} 
\qquad (x, \ y \in \R^k).]
\]
\\[-.75cm]

\qquad{\bf 12.1.13:}\quad
{\small\bf SCHOLIUM} \ 
\[
\int\limits_E \ J (T) \ \td \Lm^k
\ = \ 
\int\limits_{T(E)} \ \sH^0 (E \cap T^{-1} (y) ) \ \td \sH^k (y).
\]
To repeat: \ $k$, $n \in \N$, $1 \leq k \leq n$.
\\

\qquad{\bf 12.1.14:}\quad
{\small\bf NOTATION} \ 
\[
\Lambda_{k, n}
\ = \ 
\{\lambda \in \N^k \hsy : \hsy 1 \leq \lambda_1 \hsy < \hsy \cdots \hsy < \hsy \lambda_k \hsx \leq \hsx n\}.
\]
\\[-1.5cm]

The matrix $M_T$ associated with $T$ is $n \times k$.  
Given $\lambda \in \N^k$, let $M_T^\lambda$ be the $k \times k$ submatrix of $M_T$ made up of the rows 
$\lambda_1, \hsy \ldots, \hsy \lambda_k$ of $M_T$.
\\

\qquad{\bf 12.1.15:}\quad
{\small\bf CAUCHY-BINET FORMULA} \ 
\[
J (T)^2 
\ = \ 
\sum\limits_{\alpha \in \Lambda_{k, n}} \ (\det (M_T^\lambda))^2.
\]

Therefore $J(T)$ is the square root of the sum of the squares of the $k \times k$ subdeterminants of $\det (M_T)$.
\\

\qquad{\bf 12.1.16:}\quad
{\small\bf EXAMPLE} \ 
Suppose that $k = 2$, $n = 3$, and 
\\[-.5cm]

\[
M_T \quad = \quad
\begin{bmatrix}
a &b {\hsx}\\
c &d {\hsx}\\
e &f  {\hsx}
\end{bmatrix}
.
\]
Put

\[
u
\ = \
\begin{bmatrix}
\hsx a {\hsx} \\
\hsx c {\hsx}\\
\hsx e {\hsx}
\end{bmatrix}
, \quad
v
\ = \
\begin{bmatrix}
\hsx b {\hsx} \\
\hsx d {\hsx}\\
\hsx f {\hsx}
\end{bmatrix}
\]

\noindent
and set

\[
\begin{sqcases}
&E \hsx = \hsx \norm{u}^2 \hsx = \hsx a^2 + c^2 + e^2
\\[4pt]
&F \hsx = \hsx \langle u, v \rangle \hsx = \hsx ab + cd + ef
\\[4pt]
&G \hsx = \hsx \norm{v}^2 \hsx = \hsx b^2 + d^2 + f^2
\end{sqcases}
.
\]
Then

\[
\det (M_T^* \hsx M_T)
\ = \ 
\det
\begin{bmatrix}
\hsx E &F{\hsx} \\
\hsx F &G {\hsx}
\end{bmatrix}
\quad = \quad 
E \hsy G - F^2.
\]
On the other hand, 

\[
\Lambda_{2, 3} 
\ = \ 
\{(1, 2), (2, 3), (1,3)\},
\]
so

\[
M_T^{(1,2)} \quad = \quad
\begin{bmatrix}
\hsx a &b {\hsx}\\
\hsx c &d {\hsx}
\end{bmatrix}
\quad \equiv \quad 
A
\]
\[
M_T^{(2,3)} \quad = \quad
\begin{bmatrix}
\hsx c &d {\hsx}\\
\hsx e &f {\hsx}
\end{bmatrix}
\quad \equiv \quad 
B
\]
\[
M_T^{(1,3)} \quad = \quad
\begin{bmatrix}
\hsx a &b{\hsx}\\
\hsx e &f {\hsx}
\end{bmatrix}
\quad \equiv \quad 
C,
\]
and by Cauchy-Binet, 
\begin{align*}
\det (M_T^* \hsx M_T) \ 
&=\ 
\det (A)^2 + \det (B)^2 + \det (C)^2
\\[11pt]
&=\ 
(ad - bc)^2 + (cf - ed)^2 + (af - be)^2.
\end{align*}
Consider now $u \times v$, the vector cross product of $u$ and $v$: 
\[
\begin{bmatrix}
\hsx cf - ed {\hsx}\\
\hsx eb - af {\hsx}\\
\hsx ad - cb {\hsx}
\end{bmatrix}
\]

\hspace{.5cm} $\implies$

\[
\norm{u \times v}^2
\ = \ 
\det (M_T^* \hsx M_T)
\]

\hspace{.5cm} $\implies$

\[
\norm{u \times v} 
\ = \ 
J(T).
\]

\chapter{
$\boldsymbol{\S}$\textbf{12.2}.\quad  THE $C^\prime$ CASE}
\setlength\parindent{2em}

\renewcommand{\thepage}{12-\S2-\arabic{page}}

\qquad
It was shown in the previous \S that if 
$T \hsy : \hsy \R^k \ra \R^n$ $(k \leq n)$ is an injective linear transformation, then 

\[
\begin{cases}
&E \in \sM_L^k  \implies T(E) \in \sM(\sH^k)
\\[4pt]
&\sH^k (T(E)) \ =  \ J(T) \  \Lm^k (E)
\end{cases}
.
\]
\\[-.75cm]

This conclusion can be generalized: 
\\[-.25cm]

\qquad (1) \quad Replace $\R^k$ by a nonempty open subset $\Omega \subset \R^k$.
\\[-.5cm]

\qquad (2) \quad Replace $T$ by a one-to-one function $\Phi \hsy : \hsy  \Omega \ra \R^n$ of class $C^\prime$.
\\[-.5cm]

After a fair amount of effort, matters then will read
\[
\begin{cases}
&E \in \sM_\Lm^k \ (E \subset \Omega) \implies \Phi(E) \in \sM(\sH^k)
\\[15pt]
&\ds\sH^k (\Phi(E)) \ =  \ \int\limits_E \ J(\Phi) \ \td \Lm^k
\end{cases}
.
\]
\\[-.75cm]

Setting aside the proof until later, we shall first deal with some preliminaries and consider some examples.
\\

\qquad{\bf 12.2.1:}\quad
{\small\bf LEMMA} \ 
Let $T \hsy : \hsy \R^k \ra \R^n$ $(k \leq n)$ be a linear transformation $-$then
\[
\rank T 
\ \leq \ 
k, 
\quad 
\ker T 
\ = \ 
\ker T^* T,
\]
and the following are equivalent:
\\[-.25cm]

\qquad (a) \quad $J(T) \equiv \sqrt{\det(T^* \hsy T)} \ = \ 0$,
\\[-.5cm]

\qquad (b) \quad $\ker T \neq \{0\}$,
\\[-.5cm]

\qquad (c) \quad $\rank \hsx T < k$.
\\[-.5cm]

[Note: \ 
If $T$ is injective, then $\forall \ x \in \R^k$, 
\allowdisplaybreaks
\begin{align*}
\norm{x} \ 
&=\ 
\norm{T^{-1} \hsy T x}
\\[8pt]
&\leq \ 
\norm{T^{-1}} \hsx \norm{T x}
\end{align*}

\hspace{.5cm} $\implies$
\[
\norm{x} \hsx \frac{1}{\norm{T^{-1}} } 
\ \leq \ 
\norm{T x}.]
\]
\\[-.75cm]

\qquad{\bf 12.2.2:}\quad
{\small\bf NOTATION} \ 
Given $x_0 \in \Omega$, put
\begin{align*}
J(\Phi) (x_0) \ 
&=\ 
J( \td \Phi(x_0))
\\[11pt]
&=\ 
\sqrt{\det(\td \Phi(x_0)^* \  \td \Phi(x_0))},
\end{align*}
from which a function $J(\Phi) \hsy : \hsy \Omega \ra \R$.
\\[-.25cm]

[Note: \ 
$J(\Phi) (x_0) \neq 0$ iff $\td \Phi(x_0)$ is injective.]
\\

\qquad{\bf 12.2.3:}\quad
{\small\bf RAPPEL} \ 
If $\Phi \hsy : \hsy  \Omega \ra \R^n$ is differentiable at a point $x_0 \in \Omega$, then the $n \times k$ matrix

\[
\tD \hsy \Phi(x_0) 
\ = \ 
\nabla \hsy \Phi(x_0) 
\ = \ 
\begin{bmatrix}
\hsx \nabla \hsy \Phi_1 (x_0) {\hsx}\\
\\
\vdots\\
\\
\hsx \nabla \hsy \Phi_n (x_0) {\hsx}
\end{bmatrix}
\]

\noindent
is the \un{Jacobian matrix} of $\Phi$ at $x_0$.
\\[-.5cm]

[Note: \ 
$\tD \hsy \Phi(x_0)$ is the matrix that represents $\td \Phi(x_0)$.]
\\

\qquad{\bf 12.2.4:}\quad
{\small\bf \un{N.B.}}  \ 
Suppose that $k = n$ $-$then

\[
\det (\tD \hsy \Phi(x_0)) 
\ = \ 
J_\Phi \hsy (x_0).
\]
\\[-.75cm]

What follows are some particular cases of the relation
\[
\sH^k(\Phi(E)) 
\ = \ 
\int\limits_E \ J (\Phi) \ \td \Lm^k.
\]


\qquad{\bf 12.2.5:}\quad
{\small\bf EXAMPLE} \ 
Take $k = 1$, $n > 1$, take $\Omega = ]a,b[$\hsy, so $\Phi \hsx : \hsx ]a,b[ \hsx \ra \R^n$ is a curve:
\[
\Phi (x)
\ = \ 
(\Phi^1(x), \hsx \ldots , \hsx \Phi^n(x)) 
\quad (a < x < b).
\]
And $\tD \hsy \Phi (x)$ is an $n \times 1$ matrix or still, upon switching the column vector to a row vector, 
$\tD \hsy \Phi (x)$ becomes a $1 \times n$ matrix, viz.

\[
\left(
\frac{\td \Phi^1}{\td x} 
, \hsx \ldots , \hsx 
\frac{\td \Phi^n}{\td x} 
\right), 
\]
thus

\[
J (\Phi) 
\ = \ 
\sqrt
{
\left(
\frac{\td \Phi^1}{\td x} 
\right)^2
\hsx + \cdots + \hsx
\left(
\frac{\td \Phi^n}{\td x} 
\right)^2
}
\ = \ 
\normx{\dot{\phi}}.
\]

\noindent
If therefore $\Phi$ is one-to-one, then

\[
\sH^1(\Phi( ]a,b[)) 
\ = \ 
\int\limits_a^b \ \normx{\dot{\Phi}} \ \td t.
\]
E.g.: Let 

\[
\Phi(x) 
\ = \ 
(\cos x, \sin x, x) 
\quad (0 < x < 1),
\]
hence

\[
\dot{\Phi} (x) 
\ = \ 
(-\sin x, \cos x, 1) 
\implies 
\normx{\dot{\Phi}} 
\ = \ 
\sqrt{2},
\]
hence

\[
\sH^1(\Phi( ]0,1[)) 
\ = \ 
\int\limits_0^1 \ \normx{\dot{\Phi}} \ \td t
\ = \ 
\sqrt{2}.
\]
\\[-.75cm]

\qquad{\bf 12.2.6:}\quad
{\small\bf EXAMPLE} \ 
The \un{graph} of a $C^\prime$ function $f \hsy : \hsy \Omega \ra \R$ $(\Omega \subset \R^k)$ is the subset 
of $\R^{k + 1}$ defined by

\[
\Gr_f 
\ = \ 
\{ (x, f(x)) \in \R^k \times \R \hsx : \hsx x \in \Omega\},
\]
i.e., $\Gr_f$ is the image of the injective map $\Phi(x) = (x, f(x))$ from $\Omega$ to $\R^{k + 1}$.  
Here

\[
\tD \Phi (x) \ = \ 
\begin{bmatrix}
\fbox{\text{id}}\hsx
\\[4pt]
\tD f(x) \hsx
\end{bmatrix}
,
\]
thus by Cauchy-Binet, 

\[
J (\tD \Phi(x)) 
\ = \ 
\sqrt
{
1 + 
\left(\frac{\partial f}{\partial x_1}\right)^2 
\hsx + \hsx \cdots \hsx + \hsx 
\left(\frac{\partial f}{\partial x_k}\right)^2 
}
\]
or still, 
\[
J (\tD \Phi(x)) 
\ = \ 
\sqrt{1 + \norm{\tD f(x)}^2}
\]
\hspace{.5cm} $\implies$
\begin{align*}
\sH^k (\Gr_f) \ 
&=\ 
\sH^k (\Phi(\Omega)) 
\\[15pt]
&=\ 
\int\limits_\Omega \ 
\sqrt{1 + \norm{\tD f(x)}^2} \ \td \Lm^k.
\end{align*}
\\[-.75cm]

\qquad{\bf 12.2.7:}\quad
{\small\bf EXAMPLE} \ 
Let $\Phi \hsy : \hsy \Omega \ra \R^{k + 1}$ $(\Omega \subset \R^k)$ be a one-to-one map of class $C^\prime$ $-$then 
the Jacobian matrix of $\Phi$ has $k + 1$ rows and $k$ columns and its 
$k \times k$ submatrices can be indexed by the missing row.  
If 

\[
\frac{\partial ( \Phi_1, \hsx \ldots , \hsy \Phi_{i - 1}, \Phi_{i + 1},\hsx \ldots , \hsy, \Phi_k)}{\partial (x_1,\hsx \ldots , \hsy, x_k)}
\]
denotes the determinant of the submatrix obtained by removing the $i^\nth$ row, it therefore follows that

\[
\sH^k (\Phi(\Omega)) 
\ = \ 
\int\limits_\Omega \ 
\left( \ 
\sum\limits_{i = 1}^k \ 
\left(\ 
\frac{\partial ( \Phi_1, \hsx \ldots , \hsy \Phi_{i - 1}, \Phi_{i + 1},\hsx \ldots , \hsy, \Phi_k)}{\partial (x_1,\hsx \ldots , \hsy, x_k)}
\right)^2 
\right)^{1/2} \ 
\td \Lm^k.
\]

\qquad{\bf 12.2.8:}\quad
{\small\bf SCHOLIUM} \ 
Take $k = n$ and let $\Phi \hsy : \hsy \Omega \ra \R^n$ be a one-to-one function of 
class $C^\prime$ $-$then $\forall \ E \in \sM_\Lm^n$ $(E \subset \Omega)$,

\begin{align*}
\Lm (\Phi(E)) \ 
&=\ 
\int\limits_E \ J(\Phi) \ \td \Lm^n
\\[15pt]
&=\ 
\int\limits_\Omega \ \chisubE \hsy J(\Phi) \ \td \Lm^n
\end{align*}
from which 

\[
\int\limits_{\Phi(\Omega)} \ \chisubE  \ \td \Lm^n
\ = \ 
\int\limits_\Omega \ (\chisubE \circ \Phi) \hsy J(\Phi) \ \td \Lm^n.
\]
\\[-.75cm]

[The first point is a special case of the general theory and the second point follows from the first. 
To see this, assume to begin with that $E$ is Borel, hence that 
$\Phi^{-1}(E) = \{x \in \Omega \hsy : \hsy \Phi(x) \in E\}$ is Borel ($\Phi$ being continuous), so

\begin{align*}
\int\limits_{\Phi(\Omega)} \ \chisubE  \ \td \Lm^n
&=\ 
\Lm^n (E\cap \Omega))
\\[15pt]
&=\ 
\Lm^n (\Phi (\Phi^{-1}(E) \cap \Omega))
\\[15pt]
&=\ 
\int\limits_\Omega  \ \chisubPhiInvofE \hsy J(\Phi) \ \td \Lm^n
\\[15pt]
&=\ 
\int\limits_\Omega \ (\chisubE \circ \Phi) \hsy J(\Phi) \ \td \Lm^n.
\end{align*}
To proceed in general, let $E \in \sM_\Lm^n$ $(E \subset \Omega)$ and write $E = F \cup N$, 
where $F \cap N = \emptyset$, 
$F$ is an $F_\sigma$-set, and $N$ is a subset of a $G_\delta$-set $G$ with $\Lm^n(G) = 0$.  
Since $F$ and $G$ are Borel,
\[
\int\limits_{\Phi(\Omega)} \ \chisubF \ \td \Lm^n 
\ = \ 
\int\limits_\Omega \ (\chisubF \circ \Phi) \hsy J(\Phi) \ \td \Lm^n
\]
and
\[
0
\ = \ 
\int\limits_{\Phi(\Omega)} \ \chisubG \ \td \Lm^n 
\ = \ 
\int\limits_\Omega \ (\chisubG \circ \Phi) \hsy J(\Phi) \ \td \Lm^n.
\]

\noindent
From here, it remains only to incorporate $N$ \ldots \hsx .]
\\

\qquad{\bf 12.2.9:}\quad
{\small\bf THEOREM} \ 
If $\Phi \hsy : \hsy \Omega \ra \R^n$ is one-to-one and if 
$f \hsy : \hsy \Phi(\Omega) \ra [-\infty, +\infty]$ is Lebesgue measurable, then

\[
\int\limits_{\Phi(\Omega)} \ f \ \td \Lm^n 
\ = \ 
\int\limits_\Omega \ (f \circ \Phi) \hsy J(\Phi) \ \td \Lm^n .
\]
\\[-1cm]

[This is true when $f = \chisubE$ and the general case follows by a standard approximation argument.]
\\

\qquad{\bf 12.2.10:}\quad
{\small\bf \un{N.B.}} \ 
The relation
\[
\Lm^n (\Phi(E)) 
\ = \ 
\int\limits_E \ J(\Phi) \ \td \Lm^n 
\]
is an instance of a so-called ``area formula''.
\\

\qquad{\bf 12.2.11:}\quad
{\small\bf EXAMPLE}  \ 
Work in $\R^2$, take
\[
\Omega 
\ = \ 
]0, +\infty[ \ \times \ ]-\pi, \pi[\hsy,
\]
and for $(r, \theta) \in \Omega$ define $\Phi \hsy : \hsy \Omega \ra \R^2$ by the rule

\[
\begin{cases}
&\Phi^1 (r, \theta) \ = \ r \cos \theta \ = \ x
\\[4pt]
&\Phi^2 (r, \theta) \ = \ r \sin \theta \ = \ y
\end{cases}
.
\]
Then $\Phi$ is one-to-one, of class $C^\prime$, and its range $\Phi(\Omega)$ is $\R^2 \backslash \Lambda$, where

\[
\Lambda 
\ = \ 
]-\infty, 0] \times \{0\} \subsetx \R \times \R \quad (\implies \Lp^2 (\Lambda) = 0).
\]

\noindent
The Jacobian matrix $\tD \hsy \Phi(r, \theta)$ is given by 
\[
\begin{bmatrix}
\hsx
\cos \theta 
& - r \sin \theta
{\hsx}
\\
\hsx
\sin \theta 
& r \cos \theta 
{\hsx}
\end{bmatrix}
\]
and the Jacobian $\tJ_\Phi (r, \theta)$, i.e., $\det (\tD \hsy \Phi(r, \theta))$, equals $r$.  
So formally
\begin{align*}
\int\limits_{\R^2} \ f(x, y) \ \td \Lm^n (x, y)
&=\ 
\int\limits_{\R^2 \backslash \Lambda} \ f(x, y) \ \td \Lm^n (x, y)
\\[15pt]
&=\ 
\int\limits_{\Phi(\Omega)} \ f \ \td \Lm^n
\\[15pt]
&=\ 
\int\limits_\Omega \ (f \circ \Phi) \abs{\det J_\Phi}\ \td \Lm^n
\\[15pt]
&=\ 
\int\limits_{]0, +\infty[ \hsx \times \hsx ]-\pi, \pi[}  
f(r \cos \theta, r \sin \theta) r \ \td \Lm^2 (r, \theta).
\end{align*}

\chapter{
$\boldsymbol{\S}$\textbf{12.3}.\quad   PROOF}
\setlength\parindent{2em}
\renewcommand{\thepage}{12-\S3-\arabic{page}}

\qquad{\bf 12.3.1:}\quad
{\small\bf SUBLEMMA} \ 
If $T : \R^k \ra \R^n$ \ $(1 \leq k \leq n)$ is Lipschitz continuous and if $E \subset \R^k$ 
is Lebesgue measurable, then $T(E)$ is $\sH^k$-measurable.
\\[-.5cm]

PROOF \ 
It can be assumed that $E$ is bounded.  
Accordingly let $\{K_j\}$ be a sequence of compact sets such that $K_j \subset E$, $K_j \subset K_{j + 1}$, and 

\[
\Lm^k (E) 
\ = \ 
\Lm^k \big( \bigcup\limits_j \ K_j \big)
\quad (\implies \Lm^k \big( E \hsx \backslash \hsx \bigcup\limits_j \ K_j \big) ) 
\ = \ 0.
\]
Since $T$ is continuous, it follows that 

\[
T \big( \bigcup\limits_j \ K_j \big)
\ = \ 
\bigcup\limits_j \ T( K_j)
\]
is Borel (being a countable union of compacta), thus is $\sH^k$-measurable. 
Now write

\[
T(E) 
\ = \ 
T\big( \bigcup\limits_j \ K_j \big) 
\hsx \cup \hsx 
T\big( E \hsx \backslash \hsx \bigcup\limits_j \ K_j \big).
\]
Then
\begin{align*}
\sH^k (T\big( E \hsx \backslash \hsx \bigcup\limits_j \ K_j \big)) \ 
&\leq \ 
\Lip (T)^k \ 
\sH^k \big( E \hsx \backslash \hsx \bigcup\limits_j \ K_j \big)
\\[11pt]
&=\ 
\Lip (T)^k \ \Lm^k \big( E \hsx \backslash \hsx \bigcup\limits_j \ K_j \big) 
\\[11pt]
&=\ 
\Lip (T)^k  \cdot 0
\\[11pt]
&=\ 
0.
\end{align*}
Therefore $T(E)$ is the union of a Borel set and a set of zero $\sH^k$-measure, so $T(E)$ is $\sH^k$-measurable.
\\[-.5cm]

[Note: \ 
There are various easy variations of this theme.]
\\

\qquad{\bf 12.3.2:}\quad
{\small\bf DATA} \ 
$1 \leq k \leq n$, $\Omega \subset \R^k$ a nonempty open set, 
$\Phi : \Omega \ra \R^n$ a one-to-one function of class $C^\prime$, $E \subset \Omega$ a 
Lebesgue measurable set.
\\

\qquad{\bf 12.3.3:}\quad
{\small\bf THEOREM} \ 
$\Phi(E)$ is $\sH^k$-measurable and 
\[
\sH^k (\Phi(E)) 
\ = \ 
\int\limits_E \ \tJ (\Phi) \ \td \Lm^k 
\qquad \text{(area formula)}.
\]
\\[-.75cm]

\qquad{\bf 12.3.4:}\quad
{\small\bf EXAMPLE} \ 
Suppose that $G \subset \Phi(\Omega)$ is Borel, hence $\sH^k$-measurable $-$then
\[
\Phi (\Phi^{-1} (G)) 
\ = \ 
G \ \cap \ \Phi (\Omega) 
\ = \ G,
\]
so
\begin{align*}
\sH^k (G) \ 
&=\ 
\sH^k  (\Phi (\Phi^{-1} (G)) )
\qquad (E = \Phi^{-1} (G) \subset \Omega)
\\[15pt]
&=\ 
\int\limits_{\Phi^{-1} (G)} \ \tJ (\Phi) \ \td \Lm^k .
\end{align*}
\\[-.75cm]

\qquad
\textbf{LEMMA A}
Let $x_0 \in \Omega$ and assume that $\td \hsy \Phi(x_0) \in \Hom (\R^k, \R^n)$ is injective 
$-$then $\forall \ \epsilon > 0$ $(< 1)$ there exists a neighborhood $U \subset \Omega$ of $x_0$ 
such that for all $x^\prime$, $x^{\prime\prime} \in U$,
\begin{align*}
(1 - \varepsilon) \hsx \norm{\td \hsy \Phi(x_0) x^\prime - \td \hsy \Phi(x_0) x^{\prime\prime}} \ 
&\leq \ 
\norm{\Phi(x^\prime) - \Phi(x^{\prime\prime})}
\\[11pt]
&\leq \ 
(1+ \varepsilon) \hsx \norm{\td \hsy \Phi(x_0) x^\prime - \td \hsy \Phi(x_0) x^{\prime\prime}}.
\end{align*}
\\[-1cm]

PROOF \ 
Fix $\varepsilon > 0$ $(< 1)$ and choose $C > 0$:
\[
\norm{\td \hsy \Phi(x_0) x} 
\ \geq \ 
C \hsy \norm{x} 
\qquad (x \in \R^k).
\]
Since $\Phi$ is class $C^\prime$, there exists $\delta > 0$ such that
\[
\norm{x - x_0} \ll \delta
\implies
\norm{\td \hsy \Phi(x) - \td \hsy \Phi(x_0)} 
\ \leq \ 
C \hsy \varepsilon.
\]
So, for $x^\prime$, $x^{\prime\prime} \in U \equiv B(x_0, \delta)^0$,
\begin{align*}
||\Phi(x^\prime) 
&- \Phi(x^{\prime\prime}) - \td \Phi(x_0) (x^\prime - x^{\prime\prime})|| 
\\[15pt]
&=\ 
\bigg|\bigg| \ 
\int\limits_0^1 \ 
\frac{\td}{\td t}
(
\Phi(x^{\prime\prime} + t(x^\prime - x^{\prime\prime})) 
- \td \Phi(x_0) (x^{\prime\prime} + t(x^\prime - x^{\prime\prime}))) \ \td t \ 
\bigg|\bigg|
\\[15pt]
&=\ 
\bigg|\bigg| \ 
\int\limits_0^1 \ 
\left[
\td \Phi(x^{\prime\prime} + t(x^\prime - x^{\prime\prime})) - \td \Phi(x_0)
\right]
\hsx
(x^\prime - x^{\prime\prime}) \ \td t \ 
\bigg|\bigg|
\\[15pt]
&\leq\ 
C \hsy \varepsilon \norm{x^\prime - x^{\prime\prime}}
\\[15pt]
&\leq\ 
\varepsilon \hsy
\norm{\td \Phi(x_0) (x^\prime - x^{\prime\prime})}.
\end{align*}

Therefore
\\[-.25cm]

\qquad \textbullet \quad 
$\norm{\Phi(x^\prime) - \Phi(x^{\prime\prime})}$
\\

\hspace{1.5cm}
$\leq 
\norm{\td \Phi(x_0) (x^\prime - x^{\prime\prime})} 
+ 
\norm{\Phi(x^\prime) - \Phi(x^{\prime\prime}) - \td \Phi(x_0) (x^\prime - x^{\prime\prime})}$
\\

\hspace{1.5cm}
$\leq 
\norm{\td \Phi(x_0) \hsx x^\prime - \td \Phi(x_0) \hsx x^{\prime\prime})}
+
\varepsilon \hsy
\norm{\td \Phi(x_0) \hsx x^\prime - \td \Phi(x_0) \hsx x^{\prime\prime})}$
\\

\hspace{1.5cm}
$\leq 
(1 + \varepsilon) \hsx 
\norm{\td \Phi(x_0) \hsx x^\prime - \td \Phi(x_0) \hsx x^{\prime\prime})}$.
\\


\qquad \textbullet \quad 
$\norm{\Phi(x^\prime) - \Phi(x^{\prime\prime})}$
\\

\hspace{1.5cm}
$\geq 
\norm{\td \Phi(x_0) (x^\prime - x^{\prime\prime})} 
-
\norm{\Phi(x^\prime) - \Phi(x^{\prime\prime}) - \td \Phi(x_0) (x^\prime - x^{\prime\prime})}$
\\

\hspace{1.5cm}
$\geq 
\norm{\td \Phi(x_0) \hsx x^\prime - \td \Phi(x_0) \hsx x^{\prime\prime})}
-
\varepsilon \hsy
\norm{\td \Phi(x_0) \hsx x^\prime - \td \Phi(x_0) \hsx x^{\prime\prime})}$
\\

\hspace{1.5cm}
$\geq 
(1 - \varepsilon) \hsx 
\norm{\td \Phi(x_0) \hsx x^\prime - \td \Phi(x_0) \hsx x^{\prime\prime})}$.
\\

\qquad
\textbf{LEMMA B} \ 
Let $x_0 \in \Omega$ and assume that $\td \Phi (x_0)  \in \Hom(\R^k, \R^n)$ is injective $-$then 
$\forall \ \epsilon > 0$ $(< 1)$ there exists a neighborhood $U \subset \Omega$ of $x_0$ 
such that for each Lebesgue measurable set $E \subset U$, $\Phi(E)$ is $\sH^k$-measurable and

\begin{align*}
(1 - \varepsilon)^{k+1} \ \int\limits_E \ \tJ (\Phi)  \ \td \Lm^k \ 
&\leq \ 
\sH^k \Phi(E))
\\[15pt]
&\leq \ 
(1+ \varepsilon)^{k+1} \ \int\limits_E \ \tJ (\Phi) \ \td \Lm^k.
\end{align*}
PROOF \ 
Since the linear transformation $\td \Phi (x_0) : \R^k \ra \R^n$ is injective, 
\[
\td \Phi (x_0)^{-1} \hsx : \hsx \td \Phi (x_0) \hsx \R^k \ra \R^k.
\]
Given $\varepsilon > 0$ $( < 1)$, choose $\delta > 0$ so small that the conclusion of 
\text{\small LEMMA A} holds, where as there 
$U \equiv B(x_0, \delta)^\circ$ and in addition
\[
(1 + \varepsilon)^{-1} \hsx  \tJ (\Phi) (x) 
\ \leq \ 
\tJ (\Phi) (x_0)
\ \leq \ 
(1 + \varepsilon) \hsx  \tJ (\Phi) (x),
\]
$\Phi$ being of class $C^\prime$ $(\norm{x - x_0} < \delta)$.  
In the relation

\[
\norm{\Phi(x^\prime) - \Phi(x^{\prime\prime})}
\ \leq \  
(1 + \varepsilon)  \norm{\td \hsy \Phi(x_0) x^\prime - \td \hsy \Phi(x_0) x^{\prime\prime}},
\]
take 

\[
\begin{cases}
&x^\prime = \td \Phi (x_0)^{-1} y^\prime\\[4pt]
&x^{\prime\prime} = \td \Phi (x_0)^{-1} y^{\prime\prime}
\end{cases}
\qquad (y^\prime, y^{\prime\prime} \in \td \Phi (x_0) (U)).
\]
Then
\begin{align*}
||
(\Phi \circ \td \Phi (x_0)^{-1}) (y^\prime) 
&
- (\Phi \circ \td \Phi (x_0)^{-1}) (y^{\prime\prime})
|| 
\\[11pt]
&=\ 
\norm{\Phi(x^\prime) - \Phi(x^{\prime\prime})}
\\[11pt]
&\leq \ 
(1 + \varepsilon) \hsx \norm{\td \Phi (x_0)\hsx x^\prime - \td \Phi (x_0) \hsx x^{\prime\prime}}
\\[11pt]
&= \ 
(1 + \varepsilon) \hsx \norm{y^\prime - y^{\prime\prime}}.
\end{align*}
Therefore

\[
\Phi \circ \td \Phi (x_0)^{-1} \hsy : \hsy \td \Phi (x_0) (U) \ra \R^n
\]
is Lipschitz continuous with 

\[\Lip (\Phi \circ \td \Phi (x_0)^{-1}) 
\ \leq \ 
1 + \varepsilon.
\]
Consequently, 
\begin{align*}
\sH^k (\Phi(E)) \ 
&=\ 
\sH^k ((\Phi \circ \td \Phi (x_0)^{-1}) \hsx  (\td \Phi (x_0) (E))))
\\[15pt]
&\leq \ 
(1 + \varepsilon)^k \ 
\sH^k ( \td \Phi (x_0) (E)) 
\hspace{1.25cm} \text{(cf. 2.3.10)}
\\[15pt]
&=\ 
(1 + \varepsilon)^k \  
\tJ (\td \Phi(x_0)) \  \Lm^k (E)
\qquad \text{(cf. 12.1.11)}
\\[15pt]
&=\ 
(1 + \varepsilon)^k \ 
\tJ (\Phi) (x_0) \  \Lm^k (E)
\hspace{1.1cm} \text{(cf. 12.2.2)}
\\[15pt]
&=\ 
(1 + \varepsilon)^k \ 
\tJ (\Phi) (x_0) \ 
\int\limits_E \ 1 \ \td \Lm^k 
\\[15pt]
&=\ 
(1 + \varepsilon)^k  \ 
\int\limits_E \ \tJ (\Phi) (x_0) \  \td \Lm^k 
\\[15pt]
&\leq \ 
(1 + \varepsilon)^{k+1}  \ 
\int\limits_E \ \tJ (\Phi) (x) \  \td \Lm^k 
\qquad 
\ \left(
\norm{x - x_0} < \delta
\right),
\end{align*}
the sought for relation from above.  
To arrive at the estimate from below, in the relation

\[
(1 - \varepsilon) \hsx \norm{\td \Phi(x_0) x^\prime - \td \Phi (x_0) x^{\prime\prime}}
\ \leq \ 
\norm{\Phi(x^\prime) - \Phi(x^{\prime\prime})}
\]
take 

\[
\begin{cases}
&x^\prime =  \Phi^{-1} (y^\prime)\\[4pt]
&x^{\prime\prime} = \Phi^{-1} (y^{\prime\prime}) 
\end{cases}
\qquad (y^\prime, y^{\prime\prime} \in \Phi(U))
\]
to get
\begin{align*}
\norm{(\td \Phi(x_0)  \circ \Phi^{-1}) (y^\prime) - (\td \Phi(x_0)  \circ \Phi^{-1}) (y^{\prime\prime})} \ 
&\leq \ 
(1 - \varepsilon)^{-1} 
\norm{\Phi(x^\prime) - \Phi(x^{\prime\prime})}
\\[11pt]
&=\ 
(1 - \varepsilon)^{-1} 
\norm{y^\prime - y^{\prime\prime}}.
\end{align*}
Therefore

\[
\td \Phi(x_0) \circ \Phi^{-1} \hsx : \hsx \Phi(U) \ra \R^n
\]
is Lipschitz continuous with 

\[
\Lip (\td \Phi(x_0) \circ \Phi^{-1}) 
\ \leq \ 
(1 - \varepsilon)^{-1} .
\]
\\[-1cm]

\noindent
Now manipulate as before: 
\begin{align*}
(1 + \varepsilon)^{-1} \ 
\int\limits_E \ \tJ (\Phi) (x) \ \td \Lm^k \ 
&\leq \ 
\int\limits_E \ \tJ (\td \Phi (x_0)) \ \td \Lm^k 
\\[15pt]
&=\ 
\sH^k (\td \Phi(x_0) (E))
\\[15pt]
&=\ 
\sH^k ((\td \Phi(x_0) \circ \Phi^{-1}) \hsx (\Phi(E)))
\\[15pt]
&\leq \ 
(1 - \varepsilon)^{-k} \ \sH^k (\Phi(E))
\end{align*}

$\implies$

\[
(1 - \varepsilon)^k \hsx (1 + \varepsilon)^{-1}  \ 
\int\limits_E \ \tJ (\Phi) (x) \ \td \Lm^k 
\ \leq \ 
\sH^k (\Phi(E))
\]

$\implies$

\[
(1 - \varepsilon)^{k + 1} \ 
\int\limits_E \ \tJ (\Phi) (x) \ \td \Lm^k 
\ \leq \ 
\sH^k (\Phi(E)) 
\qquad   ((1 + \varepsilon)^{-1} \geq 1 - \varepsilon).
\]
\\[-.75cm]

\qquad{\bf 12.3.5:}\quad
{\small\bf \un{N.B.}}  \ 
$E \subset U$ is Lebesgue measurable and the claim is that $\Phi(E)$ is $\sH^k$-measurable.  
\\[-.5cm]

[To see this, let $T = \Phi \circ \td \Phi (x_0)^{-1}$, thus by construction $T: \td \Phi(x_0) (U) \ra \R^n$ is
Lipschitz continuous.  
And

\[
\Phi(E) 
\ = \ 
T (\td f(x_0) (E)).
\]
But $\td f(x_0)$ is Lipschitz continuous, so $\td f(x_0) (E)$ is $\sH^k$-measurable, thus the same is true of $T (\td f(x_0) (E))$.]
\\

\qquad
\textbf{LEMMA C} \ 
Suppose that $\forall \ x \in \Omega$, $\td \Phi(x)$ is injective $-$then

\[
\sH^k (\Phi(E))
\ = \ 
\int\limits_E \ \tJ (\Phi) (x) \ \td \Lm^k .
\]

PROOF \ 
Fix $\varepsilon > 0$ $(< 1)$ and cover $\Omega$ with countably many $U_i \subset \Omega$ such that for any 
$E \subset U_i$, 
\text{\small LEMMA B} is in force.
Given now an $E \subset \Omega$, define inductively

\[
E_1 
\ = \ 
E \cap U_1 , \hsx \ldots, \hsx 
E_i = 
\left(
E \cap U_i
\right)
\hsx \backslash \hsx \bigcup\limits_{j = 1}^{i - 1} \ B_j.
\]
Then the $E_i$ are pairwise disjoint and 

\[
E 
\ = \ 
\bigcup\limits_{i = 1}^\infty \ E_i.
\]
Proceeding, apply 
\text{\small LEMMA B} to each $E_i$, thus
\begin{align*}
(1 - \varepsilon)^{k + 1} \ \int\limits_{E_i} \ \tJ (\Phi) \ \td \Lm^k \ 
&\leq \ 
\sH^k (\Phi(E_i)) 
\\[15pt]
&\leq \ 
(1 + \varepsilon)^{k + 1} \ \int\limits_{E_i} \ \tJ (\Phi) \ \td \Lm^k 
\end{align*}
or still, upon summing over $i$ and bearing in mind that $\Phi$ is one-to-one,

\begin{align*}
(1 - \varepsilon)^{k + 1} \ \int\limits_{E} \ \tJ (\Phi) \ \td \Lm^k \ 
&\leq \ 
\sH^k (\Phi(E)) 
\\[15pt]
&\leq \ 
(1 + \varepsilon)^{k + 1} \ \int\limits_E \ \tJ (\Phi) \ \td \Lm^k. 
\end{align*}
Finish by sending $\varepsilon$ to 0.
\\

\qquad{\bf 12.3.6:}\quad
{\small\bf \un{N.B.}} \ 
$\forall \ i$, $\Phi(E_i)$ is $\sH^k$-measurable, thus the same is true of $\Phi(E)$.
\\

\qquad
\textbf{LEMMA D} \ 
Let $\Sigma \subset \Omega$ be the set of $x \in \Omega$  with the property that $\td \Phi(x)$ is not injective, 
hence that $\tJ (\Phi) (x) = 0$ $-$then
\[
\sH^k (\Phi(\Sigma))
\ = \ 0.
\]

PROOF \ 
Since the matter is local, it can be assumed that $\Omega$ is bounded and that $\td \Phi$ is bounded in $\Omega$, 
say $\norm{\td \Phi (x)} \leq M$ for all $x \in \Omega$.  
Given $\varepsilon > 0$, consider the function 

\[
\Phi_\varepsilon \hsy : \hsy \Omega \ra \R^n \times \R^k
\]
defined by the rule
\[
\Phi_\varepsilon (x) 
\ = \ 
(\Phi(x), \varepsilon \hsy x)
\qquad (x \in \Omega),
\]
so 
\[
\Phi
\ = \ 
\Pi \circ \Phi_\varepsilon,
\]
where
\[
\Pi \hsy : \hsy \R^n \times \R^k \ra \R^n
\]
is the projection operator given by $\Pi(y,x) = y$, a Lipschitz continuous function with Lipschitz constant 1 (i.e., $\Lip(\Pi) = 1$).  
Since $\forall  \ x \in \Omega$, $\td \Phi_\varepsilon (x)$ is injective, it follows from 
\text{\small LEMMA C} that 
\begin{align*}
\sH^k (\Phi(\Sigma))\ 
&=\ 
\sH^k (\Pi ( \Phi_\varepsilon (\Sigma)))
\\[15pt]
&\leq \ 
(\Lip \hsy \Pi)^k \ \sH^k (\Phi_\varepsilon (\Sigma))
\\[15pt]
&=\ 
\sH^k (\Phi_\varepsilon(\Sigma))
\\[15pt]
&=\ 
\int\limits_\Sigma \ \tJ (\Phi_\varepsilon) \ \td \Lm^k.
\end{align*}
To estimate this integral, use Cauchy-Binet to produce a constant $C > 0$ such that $\forall \ x \in \Omega$, 

\[
((\tJ (\Phi_\varepsilon) (x))^2 
\ \leq \ 
(\tJ (\Phi) (x))^2 + C^2 \hsy \norm{\td \Phi (x)}^2 \hsy \varepsilon^2.
\]
In particular, if $x \in \Sigma$, then 
\begin{align*}
(\tJ (\Phi_\varepsilon) (x))^2 \ 
&\leq \ 
C^2 \hsy \norm{\td \Phi (x)}^2 \hsy \varepsilon^2
\\[11pt]
&\leq \ 
C^2 \hsy M^2 \hsy \varepsilon^2
\end{align*}
or still, 

\[
\tJ (\Phi_\varepsilon) (x)
\ \leq \ 
C \hsy M \hsy \varepsilon.
\]
Therefore
\begin{align*}
\sH^k (\Phi(\Sigma)) \ 
&\leq \ 
\int\limits_\Sigma \ \tJ (\Phi_\varepsilon) \ \td \Lm^k
\\[15pt]
&\leq \ 
C \hsy M \hsy \varepsilon \hsy \Lm^k (\Sigma).
\end{align*}
Now let $\varepsilon \ra 0$ to get

\[
\sH^k (\Phi(\Sigma)) 
\ = \ 0.
\]
\\[-.75cm]

PROOF OF THEOREM \ 
Given a Lebesgue measurable set $E \subset \Omega$, write
\[
E 
\ = \ 
E \hsx \backslash \hsx \Sigma \cup E \cap \Sigma.
\]
Then

\[
\Phi(E) 
\ = \ 
\Phi(E \hsx \backslash \hsx \Sigma)  \cup \Phi(E \cap \Sigma).
\]
Owing to \text{\small LEMMA D}, 
$\Phi(E \cap \Sigma)$ is a set of zero $\sH^k$-measure.
On the other hand, $E \hsx \backslash \hsx \Sigma \subset \Omega \hsx \backslash \hsx \Sigma$ (an open set), 
hence $\Phi(E \hsx \backslash \hsx \Sigma)$ is $\sH^k$-measurable (cf. \text{\small LEMMA C} and 12.3.6).  
Therefore $\Phi (E)$ is $\sH^k$-measurable.  
And finally
\allowdisplaybreaks
\begin{align*}
\sH^k( \Phi (E)) \ 
&=\ 
\sH^k (\Phi(E \hsx \backslash \hsx \Sigma))
\\[15pt]
&=\ 
\int\limits_{E \hsx \backslash \hsx \Sigma} \ \tJ (\Phi) \ \td \Lm^k
\\[15pt]
&=\ 
\int\limits_E \ \tJ (\Phi) \ \td \Lm^k.
\end{align*}
The supposition that $\Phi$ is one-to-one can be dropped.
\\

\qquad{\bf 12.3.7:}\quad
{\small\bf THEOREM (AREA FORMULA)} \ 
If $\Phi \hsy : \hsy \Omega \ra \R^n$ is a function of class $C^\prime$ and if $E \subset \Omega$ is a Lebesgue measuable set, 
then $\Phi(E)$ is $\sH^k$-measurable and 

\[
\int\limits_{\R^n} \ \sH^0 (E \cap \Phi^{-1} (y)) \ \td \sH^k (y)
\ = \ 
\int\limits_E  \  \tJ (\Phi)  \ \td \Lm^k. 
\]
\\[-1cm]

[Note: \ 
If $\Phi$ is one-to-one, then matters reduce to
\[
\sH^k (\Phi(E)) 
\ = \ 
\int\limits_E  \  \tJ (\Phi)  \ \td \Lm^k.]
\]
\\[-1cm]

\qquad{\bf 12.3.8:}\quad
{\small\bf \un{N.B.}} \ 
The arrow
\begin{align*}
y \ra \sH^0 (E \cap \Phi^{-1}(y)) \ 
&=\ 
\#\{x \in E \hsy : \hsy \Phi(x) = y\}
\\[8pt]
&\equiv \ 
 N( \Phi, E, y) 
\end{align*}
is the multiplicity function of $\Phi$ at $y \in \R^n$ w.r.t. $E$ and the assignment 
\[
y \ra  N( \Phi, E, y) 
\]
defines an $\sH^k$-measurable function (cf. 11.2.4.) (recall that $\Phi$ has property (N)).
\\[-.5cm]

[Note: \ 
\[
\Phi(E) 
\ = \ 
\{y \in \R^n \hsy : \hsy N( \Phi, E, y) \neq 0\},
\]
so the integral over $\R^n$ can be replaced by an integral over $\Phi (E)$.]
\\

\qquad{\bf 12.3.9:}\quad
{\small\bf CHANGE OF VARIABLES} \ 
Suppose that  $\Phi \hsy : \hsy \Omega \ra \R^n$ is  class $C^\prime$ and $u : \Omega \ra [0, +\infty]$ is Lebesgue measurable $-$then 
the assignment
\[
y \hsx \lra \hsx \sum\limits_{x \in \Phi^{-1} (y)} \  u(x)
\]
defines an $\sH^k$-measurable function and 

\[
\int\limits_\Omega  \ u \hsy \tJ (\Phi)  \ \td \Lm^k 
\ = \ 
\int\limits_{\R^n} \ 
\bigg(
\sum\limits_{x \in \Phi^{-1} (y)} \  u(x)
\bigg)
\  \td \sH^k (y).
\]

[The proof is canonical, given what we know.  
Thus start with $u = \chisubE$ $(E \subset \Omega$ Lebesgue measurable) and note that 

\[
\sum\limits_{x \in \Phi^{-1} (y)} \ \chisubE(x) 
\ = \ 
\sH^0 (E \cap \Phi^{-1} (y)).
\]
Therefore
\begin{align*}\ 
\int\limits_\Omega  \ \chisubE \hsy \tJ (\Phi)  \ \td \Lm^k \ 
&= \ 
\int\limits_E \ \tJ (\Phi)  \ \td \Lm^k 
\\[15pt]
&= \ 
\int\limits_{\R^n} \ \sH^0 (E \cap \Phi^{-1} (y)) \  \td \sH^k (y)
\qquad \text{(area formula)}
\\[15pt]
&= \ 
\int\limits_{\R^n} \ 
\bigg(
\sum\limits_{x \in \Phi^{-1} (y)} \  \chisubE (x)
\bigg)
\  \td \sH^k (y).
\end{align*}
By linearity, this settles the contention for simple functions, thence \ldots \hsx .]
\\

\qquad{\bf 12.3.10:}\quad
{\small\bf SCHOLIUM} \ 
Suppose that $\Phi \hsy : \hsy \Omega \ra \R^n$ is  class $C^\prime$ and $v : \R^n \ra [0, +\infty]$ is
 $\sH^k$-measurable $-$then for every Lebesgue measurable set $E \subset \Omega$, 
 
\[
\int\limits_E \ (v \circ \Phi) \hsx \tJ (\Phi) \ \td \Lm^k 
\ = \ 
\int\limits_{\R^n} \ v ( y) \hsx N(\Phi, E, y) \ \td \sH^k (y).
\]
\\[-.75cm]

\qquad \un{Step 1:}\ 
Take $E = \Omega$ and $v = \chisubV$ $(V \subset \R^n)$, $V$ open $-$then
\begin{align*}
\int\limits_\Omega \ (\chisubV \circ \Phi) \hsx \tJ (\Phi) \ \td \Lm^k \ 
&=\ 
\int\limits_{\Phi^{-1}(V)} \  \tJ (\Phi) \ \td \Lm^k
\\[15pt]
&=\ 
\int\limits_{\R^n} \ N(\Phi, \Phi^{-1}(V), y) \  \td \sH^k (y)
\qquad \text{(area formula)}.
\end{align*}
And

\[
\int\limits_{\R^n}
\ = \ 
\int\limits_{\Phi(\Phi^{-1}(V))} \ 
\ = \ 
\int\limits_{V \cap \Phi(\Omega)}
\]

$\implies$
\begin{align*}
\int\limits_{\R^n} \ N(\Phi, \Phi^{-1}(V), y) \  \td \sH^k (y) \ 
&=\ 
\int\limits_V \ N(\Phi, \Phi^{-1}(V), y) \  \td \sH^k (y)
\\[15pt]
&=\ 
\int\limits_V \ N(\Phi, \Omega, y) \  \td \sH^k (y)
\\[15pt]
&=\ 
\int\limits_{\R^n} \ \chisubV \hsx N(\Phi, \Omega, y) \  \td \sH^k (y).
\end{align*}
\\[-1cm]

[Note: \ 
\[
\int\limits_{\R^n}
\ = \ 
\int\limits_{V \cap \hsx \Phi(\Omega)}.
\]
\\[-1cm]

Meanwhile
\begin{align*}
\int\limits_{\R^n} \ 
\ \geq\  \ 
\int\limits_V 
\ \geq\  \ \ 
\int\limits_{V \cap \hsx \Phi(\Omega)}
\end{align*}
\\[-.75cm]

$\implies$
\[
\int\limits_{\R^n}
\ = \ 
\int\limits_V \ .\big]
\]
\\[-.75cm]

\qquad \un{Step 2:}\ 
Take $E \subset \Omega$ compact and $v$ a simple function constant on open sets.
\\[-.25cm]

\qquad \un{Step 3:}\ 
Take $E \subset \Omega$ compact and $v$ an arbitrary simple function. 
\\[-.25cm]

\qquad \un{Step 4:}\ 
Take $E \subset \Omega$ compact and $v \geq 0$ an arbitrary measurable function.
\\[-.25cm]

\qquad \un{Step 5:}\ 
Take $E \subset \Omega$ Lebesgue measurable and $v \geq 0$ an arbitrary measurable function.

\chapter{
$\boldsymbol{\S}$\textbf{12.4}.\quad THE DIFFERENTIABLE CASE}
\setlength\parindent{2em}
\renewcommand{\thepage}{12-\S4-\arabic{page}}

\qquad
The central conclusion of the preceding $\S$ is the  fact that $\Phi(E)$ is $\sH^k$-measurable and 
\[
\sH^k (\Phi(E)) 
\ = \ 
\int\limits_E \ 
\tJ (\Phi) \ 
\td \Lm^k.
\]
Here $1 \leq k \leq n$, $\Omega \subset \R^k$ is a nonempty open set, $\Phi : \Omega \ra \R^n$ is a one-to-one function of class $C^\prime$, 
and $E \subset \Omega$ is a Lebesgue measurable set.

It turns out that one can drop the assumption that $\Phi$ is class $C^\prime$, it being enough to suppose that $\Phi$ is merely differentiable 
(as well as one-to-one).
\\

\qquad{\bf 12.4.1:}\quad
{\small\bf WHITNEY APPROXIMATION PRINCIPLE} \ 
There exists a sequence of disjoint closed sets $F_j \subset \Omega$ and a sequence of $C^\prime$ functions 
$\Phi_j : \R^k \ra \R^n$ such that in $F_j$
\[
\Phi 
\ = \ 
\Phi_j 
\quad \text{and} \quad 
\tJ(\Phi) 
\ = \ 
\tJ(\Phi_j).
\]
Moreover
\[
\Lm^k(\Omega \backslash F)
\ = \ 
0,
\]
where $F = \bigcup\limits_j \ F_j$.
\\

\qquad{\bf 12.4.2:}\quad
{\small\bf LEMMA} \ 

\[
\sH^k(\Phi (\Omega \backslash F))
\ = \ 
0.
\]

PROOF \ 
Write

\[
\Omega \backslash F
\ = \ 
\bigcup\limits_{\ell = 1}^\infty \ E_\ell,
\]
$E_\ell$ being the set of all points $x \in \Omega \backslash F$ such that 

\[
\frac{\norm{\Phi(x) - \Phi(x^\prime)}}{\norm{x - x^\prime}} 
\ \leq \
\ell
\]
for all $x^\prime \in \Omega$ with 

\[
0 
\ < \ 
\norm{x - x^\prime}
\ \leq \
\frac{1}{\ell}.
\]

\noindent
Claim: \ The restriction of $\Phi$ to $E_\ell$ is locally Lipschitz.  
For suppose that $x$, $x^\prime$ belong to a compact $K \subset E_\ell$ and $\ds\norm{x - x^\prime} \leq \frac{1}{\ell}$ $-$then 

\[
\norm{\Phi(x) - \Phi(x^\prime)}
\ \leq \
\ell \hsx \norm{x - x^\prime}
\]
by the very definition of $E_\ell$.  
On the other hand, if $\ds\norm{x - x^\prime} > \frac{1}{\ell}$, then

\begin{align*}
\norm{\Phi(x) - \Phi(x^\prime)} \ 
&\leq \ 
2 \max\limits_K \norm{\Phi}
\\[11pt]
&\leq \ 
2 \max\limits_K \norm{\Phi} \cdot 1
\\[11pt]
&\leq \ 
2 \max\limits_K \norm{\Phi} \hsy \ell \hsy \norm{x - x^\prime}.
\end{align*}
Hence the claim.  
Consequently
\begin{align*}
\sH^k (\Phi (K)) \ 
&\leq \ 
\Lip (\restr{\Phi}{K} ; K) \hsy \sH^k (K)
\\[11pt]
&= \ 
\Lip (\restr{\Phi}{K} ; K) \hsy \Lm^k (K),
\end{align*}
where 

\[
\Lip (\restr{\Phi}{K} ; K) 
\ \leq \ 
\ell(1 + 2 \max\limits_K \norm{\Phi}).
\]
But
\[
K \subset E_\ell \subset \Omega \backslash F
\]

\hspace{.5cm}
$\implies$
\[
\Lm^k (K) 
\ \leq \ 
\Lm^k (E_\ell) 
\ \leq \ 
\Lm^k (\Omega \backslash F) 
\ = \ 
0.
\]
Therefore

\[
\sH^k (\Phi(K)) 
\ = \ 
0.
\]
Now let $K \uparrow E_\ell$ invade $E_\ell$ to get

\[
\sH^k (\Phi(E_\ell))
\ = \ 
0,
\]
so in the end
\begin{align*}
\sH^k (\Phi(\Omega \backslash F)) \ 
&\leq \ 
\sum\limits_{\ell = 1}^\infty \ \sH^k (\restr{\Phi}{E_\ell}))
\\[15pt]
&=\ 
0.
\end{align*}
\\[-.75cm]

\qquad{\bf 12.4.3:}\quad
{\small\bf APPLICATION} \ 
\\[-.5cm]

\hspace{3cm}
$
\Omega 
\ = \ 
F \cup \Omega \backslash F
$

\hspace{1.5cm} $\implies$

\hspace{3cm}
$
\Phi(\Omega) 
\ = \ 
\Phi(F) \cup \Phi(\Omega \backslash F)
$

\hspace{1.5cm} $\implies$

\hspace{3cm}
$
\Phi(E) 
\ = \ 
\Phi(E) \cap \Phi(\Omega)
\ = \ 
\Phi(E) \cap \Phi(F) \cup \Phi(E) \cap \Phi(\Omega \backslash F)
$

\hspace{1.5cm} $\implies$

\hspace{3cm}
$
\sH^k (\Phi(E))
\ = \ 
\sH^k (\Phi(E) \cap \Phi(F)) \ + \  \sH^k (\Phi(E) \cap \Phi(\Omega \backslash F)).
$
\\[-.5cm]

And

\vspace{-1cm}
\allowdisplaybreaks
\begin{align*}
\sH^k (\Phi(E) \cap  \Phi(\Omega \backslash F)) \
&\leq \ 
\sH^k (\Phi(\Omega \backslash F))
\\[11pt]
&=\ 
0.
\end{align*}

Now compute: 

\vspace{-1cm}
\allowdisplaybreaks
\begin{align*}
\sH^k (\Phi(E)) \ 
&= \ 
\sH^k (\Phi(E) \cap \Phi(F)) 
\\[15pt]
&=\ 
\sum\limits_j \ 
\sH^k (\Phi(E) \cap \Phi(F_j))
\\[15pt]
&=\ 
\sum\limits_j \ 
\sH^k (\Phi(E) \cap \Phi_j(F_j))
\\[15pt]
&=\ 
\sum\limits_j \ 
\int\limits_{\Phi^{-1} (\Phi(E)) \cap F_j} \ 
\tJ(\Phi_j) 
\ \td \Lm^k
\\[15pt]
&=\ 
\sum\limits_j \ 
\int\limits_{\Phi^{-1} (\Phi(E)) \cap F_j} \ 
\tJ(\Phi) 
\ \td \Lm^k
\\[15pt]
&=\ 
\sum\limits_j \ 
\int\limits_{E \cap F_j} \ 
\tJ(\Phi) 
\ \td \Lm^k
\\[15pt]
&=\ 
\int\limits_{E \cap F} \ 
\tJ(\Phi) 
\ \td \Lm^k
\\[15pt]
&=\ 
\int\limits_{E \hsy \cap \hsy  F  \hsy \cup \hsy  E \hsy  \cap \hsy  \Omega \backslash F} \ 
\tJ(\Phi) 
\ \td \Lm^k
\\[15pt]
&=\ 
\int\limits_E \ \tJ(\Phi) \ \td \Lm^k.
\end{align*}


\chapter{
$\boldsymbol{\S}$\textbf{12.5}.\quad  THE LIPSCHITZ CASE}
\setlength\parindent{2em}
\renewcommand{\thepage}{12-\S5-\arabic{page}}

\qquad
\qquad{\bf 12.5.1:}\quad
{\small\bf DATA} \ 
$1 \leq k \leq n$, $\Phi : \R^k \ra \R^n$ a Lipschitz continuous function, $E \subset \R^k$ a 
Lebesgue measurable set. 
\\

\qquad{\bf 12.5.2:}\quad
{\small\bf RAPPEL} \ 
Owing to Rademacher, $\tJ(\Phi)$ is defined $\Lm^k$-almost everywhere.
\\

\qquad{\bf 12.5.3:}\quad
{\small\bf THEOREM (AREA FORMULA)} \ 
$\Phi(E)$ is $\sH^k$-measurable and 

\[
\int\limits_{\R^n} \ 
\sH^0(E \cap \Phi^{-1} (y)) \ \td \sH^k(y) 
\ = \ 
\int\limits_E \ \tJ(\Phi) \ \td \Lm^k.
\]
\\[-.75cm]

\qquad{\bf 12.5.4:}\quad
{\small\bf \un{N.B.}}  \ 
There is an a priori estimate

\[
\int\limits_{\R^n} \ 
\sH^0(E \cap \Phi^{-1} (y)) \ \td \sH^k(y) 
\ \leq \ 
(\Lip (\Phi))^k \hsx \Lm^k (E).
\]
\\[-.75cm]

\qquad{\bf 12.5.5:}\quad
{\small\bf REMARK} \ 
$\Phi$ has property (N), thus the assignment

\[
y \ra N(\Phi, E, y)
\]
defines an $\sH^k$-measurable function (cf. 11.2.4).
\\

\qquad{\bf 12.5.6:}\quad
{\small\bf LEMMA} \ 
$\forall \ \varepsilon > 0$, there exists a closed set $F_\varepsilon \subset \R^k$ and a $C^\prime$ 
function $\Phi_\varepsilon : \R^k \ra \R^n$ such that in $F_\varepsilon$, 

\[
\Phi 
\ = \ 
\Phi_\varepsilon 
\quad \text{and} \quad 
\tD \Phi = \tD \Phi_\varepsilon.
\]
Moreover
\[
\Lm^k (\R^k \backslash F_\varepsilon) 
\ < \ 
\varepsilon.
\]
\\[-.25cm]

Granted this and bearing in mind the $C^\prime$ version of the area formula, we have
\begin{align*}
\int\limits_E \ 
\tJ (\Phi)
\ \td \Lm^k
&=\ 
\int\limits_{F_\varepsilon \cap E} \ 
\tJ (\Phi)
\ \td \Lm^k 
\hsx + \hsx 
\int\limits_{E \backslash F_\varepsilon} \ 
\tJ (\Phi)
\ \td \Lm^k
\\[15pt]
&=\ 
\int\limits_{F_\varepsilon \cap E} \ 
\tJ (\Phi_\varepsilon)
\ \td \Lm^k
\hsx + \hsx 
\int\limits_{E \backslash F_\varepsilon} \ 
\tJ (\Phi)
\ \td \Lm^k
\\[15pt]
&=\ 
\int\limits_{\R^n} \ 
\sH^0(F_\varepsilon \cap E \cap \Phi_\varepsilon^{-1}(y))
\ \td \sH^k (y)
\hsx + \hsx 
\int\limits_{E \backslash F_\varepsilon} \ 
\tJ (\Phi)
\ \td \Lm^k
\\[15pt]
&=\ 
\int\limits_{\R^n} \ 
\sH^0(F_\varepsilon \cap E \cap \Phi^{-1}(y))
\ \td \sH^k (y)
\hsx + \hsx 
\int\limits_{E \backslash F_\varepsilon} \ 
\tJ (\Phi)
\ \td \Lm^k.
\end{align*}
Now send $\varepsilon$ to 0, noting that $\Lm^k (E \backslash F_\varepsilon) \ra 0$ 
(use monotone convergence).
\\

\qquad{\bf 12.5.7:}\quad
{\small\bf EXAMPLE} \ 
Given a Lipschitz continuous function $f : \R^k \ra \R$, put

\[
\Phi (x) \ = \ (x, f(x)) 
\quad (x \in \R^k).
\]
Then $\Phi : \R^k \ra \R^{k + 1}$ is Lipschitz continuous, one-to-one, and

\[
\sH^k (\Phi (E)) 
\ = \ 
\int\limits_E \ 
\tJ(\Phi) \ \td \Lm^k 
\ = \ 
\int\limits_E \ 
\sqrt{1 + \norm{\tD f}^2} \ 
\td \Lm^k.
\]

\chapter{
$\boldsymbol{\S}$\textbf{12.6}.\quad  THE SOBOLEV CASE}
\setlength\parindent{2em}
\setcounter{chapter}{12}
\renewcommand{\thepage}{12-\S6-\arabic{page}}

\qquad 
Let $\Omega$ be a nonempty open subset of $\R^n$.  
Given a continuous function $f \in W^{1, p} (\Omega)$ $(1 \leq p < +\infty)$ and a Lebesgue measurable set 
$E \subset \Omega$, put

\[
\Gr_f (E) 
\ = \ 
\{(x, f(x) : \hsy x \in E\}
\subsetx 
\R^{n + 1}.
\]
\\[-1cm]

\qquad{\bf 12.6.1:}\quad 
{\small\bf THEOREM} \ 

\[
\sH^n (\Gr_f (E)) 
\ = \ 
\int\limits_E \ \sqrt{1 + \norm{\nabla f}^2} 
\ \td \Lm^n.
\]

Per 7.1.5., write

\[
\Omega
\ = \ 
\left(
\bigcup\limits_{k = 1}^\infty  \ E_k
\right)
\cup Z,
\]
where the $E_k$ are Lebesgue measurable sets such that 
$\restr{f}{E_k}$ is Lipschitz and $Z$ has Lebesgue measure 0.  
Extend $\restr{f}{E_k}$ to a Lipschitz function 
$f_k : \R^n \ra \R$ $-$then $\norm{\nabla f_k} = \norm{\nabla f}$ almost everywere in $E_k$.  
Now apply 12.5.7. to get

\[
\sH^n (\Gr_f (E \cap E_k)) 
\ = \ 
\int\limits_{E \cap E_k} \ \sqrt{1 + \norm{\nabla f}^2} 
\ \td \Lm^n.
\]
Put $\ds E_0 = \bigcup\limits_{k = 1}^\infty \ E_k$ and sum over $k$, hence
\begin{align*}
\sH^n (\Gr_f (E \cap E_0)) \ 
&=\ 
\int\limits_{E \cap E_0} \ \sqrt{1 + \norm{\nabla f}^2} 
\ \td \Lm^n
\\[15pt]
&=\ 
\int\limits_E \ \sqrt{1 + \norm{\nabla f}^2} 
\ \td \Lm^n.
\end{align*}
It remains to pass from 

\[
\sH^n (\Gr_f (E \cap E_0)) 
\quad \text{to} \quad 
\sH^n (\Gr_f (E)) 
\]
and for this, it need only be shown that 

\[
\sH^n (\Gr_f (E \backslash E_0)) 
\ = \ 0.
\]
\\[-1cm]

\qquad{\bf 12.6.2:}\quad
{\small\bf LEMMA} \ 
Let $f \in W^{1, p} (\Omega)$ $-$then $\sH^n (\Gr_f (S)) = 0$ if $S \subset \Omega$ is a set of Lebesgue measure 0.
\\[-.5cm]

[It suffices to make the verification in 
\[
W_\locx^{1, 1} (\Omega) 
\quad (\supset W^{1, p} (\Omega)).]
\]

\chapter{
$\boldsymbol{\S}$\textbf{12.7}.\quad  THE APPROXIMATE CASE}
\setlength\parindent{2em}
\renewcommand{\thepage}{12-\S7-\arabic{page}}

\qquad
Suppose that $\Phi : \Omega \ra \R^n$ is approximately differentiable almost everywhere in $\Omega$ $-$then
using approximate partial derivatives, one can form $\tJ_\ap (\Phi)$.
\\

\qquad{\bf 12.7.1:}\quad
{\small\bf LEMMA} \ 
(cf. 12.4.1) \ 
There exists an increasing sequence $K_1 \subset K_2 \subset \ldots$ of compact subsets of $\Omega \ (\subset \R^k)$ for which

\[
\Lm^k (\Omega \backslash K) 
\ = \ 
0
\quad (K = \bigcup\limits_j \ K_j)
\]
and a sequence of $C^\prime$ functions $\Phi_j : \R^k \ra \R^k$ such that in $K_j$, 
\[
\Phi 
\ = \ 
\Phi_j 
\quad \text{and} \quad
\tJ_\ap (\Phi) 
\ = \ 
\tJ (\Phi_j) .
\]
\\[-1cm]

\qquad{\bf 12.7.2:}\quad
{\small\bf NOTATION} \ 
Given $y \in \R^n$, let $m^j(y)$ be the cardinality of $\Phi^{-1} (y) \cap K_j$.
\\

\qquad{\bf 12.7.3:}\quad
{\small\bf LEMMA} \ 
$m^{j} (-)$ is Borel measurable and $\forall \ y$, 

\[
m^1(y) 
\ \leq \ 
m^2 (y) 
\ \leq \ 
\ldots \hsx .
\]
Put

\[
m(y) 
\ \equiv \ 
\lim\limits_{y \ra \infty} \ m^j(y).
\]
Then $m(-)$ is Borel measurable.
\\

\qquad{\bf 12.7.4:}\quad
{\small\bf THEOREM} \ 

\[
\int\limits_{\R^n} \ m(y) \ \td \sH^k (y)
\ = \ 
\int\limits_\Omega \ \tJ_\ap (\Phi) \ \td \Lm^k.
\]

It suffices to show that 

\[
\int\limits_{\R^n} \ m(y) \ \td \sH^k (y)
\ \geq \ 
\int\limits_\Omega \ \tJ_\ap (\Phi) \ \td \Lm^k
\]
and 

\[
\int\limits_{\R^n} \ m(y) \ \td \sH^k (y)
\ \leq \ 
\int\limits_\Omega \ \tJ_\ap (\Phi) \ \td \Lm^k.
\]
The second point being the easier of the two, note that

\allowdisplaybreaks
\begin{align*}
\int\limits_\Omega \ \tJ_\ap (\Phi) \ \td \Lm^k  \ 
&\geq \ 
\int\limits_{K_j} \ \tJ_\ap(\Phi) \td \Lm^k
\\[15pt]
&=\ 
\int\limits_{K_j} \ \tJ(\Phi_j) \td \Lm^k
\\[15pt]
&=\ 
\int\limits_{\R^n} \ N(\Phi_j, K_j, y)\ \td \sH^k (y)
\\[15pt]
&=\ 
\int\limits_{\R^n} \ \sH^0 (\Phi_j^{-1} (y) \cap K_j) \ \td \sH^k (y)
\\[15pt]
&=\ 
\int\limits_{\R^n} \ m^j(y) \ \td \sH^k (y)
\\[15pt]
&
\underset{(j \ra \infty)}{\longrightarrow} \ 
\int\limits_{\R^n} \ m(y) \ \td \sH^k (y).
\end{align*}
\\[-.5cm]

\qquad{\bf 12.7.5:}\quad
{\small\bf \un{N.B.}}\ 
Under the supposition that 

\[
\int\limits_\Omega \ \tJ_\ap (\Phi) \ \td \Lm^k 
\ < \ 
+\infty,
\]
the ``$m$'' is independent of the choice of data, i.e., if
\[
\begin{cases}
& m_1(-) \longleftrightarrow \{K_j^1\}
\\[4pt]
& m_2(-) \longleftrightarrow \{K_j^2\}
\end{cases}
,
\]
then $m_1 = m_2 \hsy \sH^k$ almost everywhere.
\\[-.5cm]


[Let

\[
m_3 (-) \longleftrightarrow \{K_j^1 \cap K_j^2\}.
\]
Then

\[
m_3(y) \ \leq \ 
\begin{cases}
&m_1(y)
\\[4pt]
&m_2(y)
\end{cases}
\qquad (y \in \R^n).
\]
But

\allowdisplaybreaks
\begin{align*}
\int\limits_\Omega \ \tJ_\ap (\Phi) \ \td \Lm^k \ 
&=\ 
\int\limits_{\R^n} \ m_3(y) \ \td \sH^k (y)
\\[26pt]
&=\ 
\begin{cases}
&\ds\int\limits_{\R^n} \ m_1(y) \ \td \sH^k (y)
\\[26pt]
&\ds\int\limits_{\R^n} \ m_2(y) \ \td \sH^k (y)
\end{cases}
\\[11pt]
\implies \hspace{1cm}
\\[11pt]
&\hspace{0.65cm}
\begin{cases}
&\ds\int\limits_{\R^n} \ (m_1 - m_3) \td \sH^k  \ = \ 0 
\\[26pt]
&\ds\int\limits_{\R^n} \ (m_2 - m_3) \td \sH^k  \ = \ 0 
\end{cases}
\\[15pt]
\implies \hspace{1cm}
\\[7pt]
&\hspace{0.65cm}
\begin{cases}
&m_1 - m_3  \ = \ 0 \\[11pt]
&m_2 - m_3  \ = \ 0 
\end{cases}
\qquad 
\text{$\sH^k$ almost everywhere}
\\[7pt]
\implies \hspace{1cm}
\\[7pt]
&\hspace{1.6cm}
m_1 = m_2 \ \text{$\sH^k$ almost everywhere}.
\end{align*}

\newpage

\centerline{\textbf{\large REFERENCES}}
\setcounter{page}{1}
\setcounter{theoremn}{0}
\renewcommand{\thepage}{References-\arabic{page}}
\vspace{0.75cm}

\[
\text{BOOKS}
\]

\begin{rf}
Adams, Robert, \textit{Sobolev Spaces}, 
Academic Press, New York, 1975.
\end{rf}

\begin{rf}
Ambrosio, Luigi, Fusco, Nicola, and Pallara, Diego, 
\textit{Functions of Bounded Variation and Free Discontinuity Problems}, 
Oxford Science Publications, 2000.
\end{rf}

\begin{rf}
Apostol, Tom., 
\textit{A Modern Approach to Advanced Calculus, Second Edition,}
Pearson, 1974.
\end{rf}

\begin{rf}
Aumann, Georg, 
\textit{Reelle Funktionen}, 
Springer-Verlag, Berlin, 1969.
\end{rf}

\begin{rf}
Blank, J., Exner, P., and Havlicek, M., 
\textit{Hilbert Space Operators in Quantum Physics}, 
American Institute of Physics, New York, 1994.
\end{rf}

\begin{rf}
Bogachev, V. I., 
\textit{Measure Theory, Volume I}, 
Springer-Verlag, Berlin Heidelberg, 2007.
\end{rf}

\begin{rf}
Bogachev, V. I., 
\textit{Measure Theory, Volume II}, 
Springer-Verlag, Berlin Heidelberg, 2007.
\end{rf}

\begin{rf}
Dieudonn\'e, J., 
\textit{Foundations of Modern Analysis}, 
Academic Press, New York, 1960.
\end{rf}

\begin{rf}
Dunford, N., and Schwartz, J., 
\textit{Linear Operators, Part I: General Theory}, 
Interscience Publishers, New York, 1958.
\end{rf}

\begin{rf}
Dunford, N., and Schwartz, J., 
\textit{Linear Operators, Part II: Special Theory}, 
Interscience Publishers, New York, 1963.
\end{rf}

\begin{rf}
Evans, Lawrence C. and Gariepy, Ronald F., 
\textit{Measure Theory and Fine Properties of Functions}, 
Chapman and Hall, 2015.
\end{rf}

\begin{rf}
Federer, Herbert, 
\textit{Geometric Measure Theory}, 
Springer-Verlag, New York, 1969.
\end{rf}

\begin{rf}
Fleming, Wendell, 
\textit{Functions of Several Variables, Second Edition}, 
Springer-Verlag, New York, 1977.
\end{rf}

\begin{rf}
Flett, T. M.,
\textit{Mathematical Analysis}, 
McGraw-Hill, England, 1966.
\end{rf}

\begin{rf}
Folland, Gerald B., 
\textit{Real Analysis: Modern Techniques and Their Applications, Second Edition}, 
Wiley, 2007.
\end{rf}

\begin{rf}
Fonseca, I., and Gangbo, W., 
\textit{Degree Theory in Analysis and Applications}, 
Oxford Science Publications, Oxford, 1995.
\end{rf}

\begin{rf}
Gel'Fand, I. M., and Shilov, G. E., 
\textit{Generalized Functions, Volume 3, Theory of Differential Equations},
Academic Press, New York, 1967.
\end{rf}

\begin{rf}
Giaquinta, Mariano and Modica, Giuseppe, 
\textit{Mathematical Analysis, An Introduction to Functions of Several Variables}, 
Birkh\"auser, 2007.
\end{rf}

\begin{rf}
Giaquinta, Mariano and Modica, Giuseppe, 
\textit{Mathematical Analysis, Foundations and Advanced Techniques for Functions of Several Variables}, 
Birkh\"auser, 2011.
\end{rf}

\begin{rf}
Goffman, Casper,
\textit{Real Functions}, 
Rinehart, New York, 1953.
\end{rf}

\begin{rf}
Goffman, Casper, Nishiura, Toga, and Waterman, Daniel, 
\textit{Homeomorphisms in Analysis}, 
American Mathematical Society, 1997.
\end{rf}

\begin{rf}
Heinonen, Juha, 
\textit{Lectures on Analysis on Metric Spaces}, 
Springer-Verlag, New York, 2001.
\end{rf}

\begin{rf}
Hewitt, Edwin and Stromberg, Karl, 
\textit{Real and Abstract Analysis},
Springer-Verlag, New York, 1969.
\end{rf}

\begin{rf}
Hirsch, Morris,  
\textit{Differential Topology}, 
Springer-Verlag, New York, 1976.
\end{rf}

\begin{rf}
Kashiwara, Masaki and Schapira, Pierre, 
\textit{Categories and Sheaves}, 
Springer-Verlag, 2005.
\end{rf}

\begin{rf}
Leoni, Giovanni, \textit{A First Course in Sobolev Spaces, Second Edition}, 
American Mathematical Society,  2017. 
\end{rf}

\begin{rf}
Lin, Fanghua and Yang, Xiaoping, 
\textit{Geometric Measure Theory}, 
International Press, Somerville, Massachusetts, 2010.
\end{rf}

\begin{rf}
Lindenstrauss, J., Preiss, D., and Tiser, J., 
\textit{Fr\'echet Differentiability of Lipschitz Functions and Porous Sets in Banach Spaces}, 
Princeton University Press, 2012.
\end{rf}

\begin{rf}
Loomis, Lynn H. and Sternberg, Shlomo, 
\textit{Advanced Calculus, Revised Edition}, 
Jones and Bartlett Publishers, Boston, 1990.
\end{rf}

\begin{rf}
Morrey, Jr., Charles B., 
\textit{Multiple Integrals in the Calculus of Variations},  
Springer-Verlag, New York, 1966.
\end{rf}

\begin{rf}
Nachbin, L., 
\textit{The Haar Integral}, 
D. Van Nostrand Co., Princeton, 1965.
\end{rf}

\begin{rf}
Roberts, A. and Vargerg, D., 
\textit{Convex Functions}, 
Academic Press, 1973.
\end{rf}

\begin{rf}
Rogers, C. A., 
\textit{Hausdorff Measures}, 
Cambridge University Press, 1970, Reprinted 1998.
\end{rf}

\begin{rf}
Rudin, Walter, 
\textit{Principles of Mathematical Analysis, Third Edition}, 
McGraw-Hill Book Company, 1976.
\end{rf}

\begin{rf}
Stein, Elias M., 
\textit{Singular Integrals and Differentiability Properties of Functions}, 
Princeton University Press, 1970.
\end{rf}

\begin{rf}
Torchinsky, Alberto, 
\textit{Real-Variable Methods in Harmonic Analysis}, 
Academic Press, New York, 1986.
\end{rf}

\begin{rf}
Treves, F., 
\textit{Topological Vector Spaces, Distributions and Kernels}, 
Academic Press, New York, 1967.
\end{rf}

\begin{rf}
Yeh, J., 
\textit{Real Analysis, Theory of Measure and Integration}, 
World Scientific Publishing Co. Pte. Ltd. Singapore, 2014.
\end{rf}

\begin{rf}
Ziemer, William P., 
\textit{Weakly Differentiable Functions}, 
Springer-Verlag, New York, 1989.
\end{rf}

\newpage


\setcounter{theoremn}{0}

\[
\text{ARTICLES}
\]

\begin{rf}
Avdispahi\'c, M., Concepts of generalized bounded variation and the theory of Fourier series, 
\textit{Internat. J. Math. Sci.} \textbf{9} (1986), 223-244. 
\end{rf}

\begin{rf}
Bruckner, A. M., Density-preserving homeomorphisms and the theorem of Maximoff, 
\textit{Quart. J. Math. Oxford Ser. (2)}  \textbf{21} (1970), 337-347. 
\end{rf}

\begin{rf}
Buczolich, Z., Density points and bi-Lipschitz functions in $\R^m$, 
\textit{Proc. Amer. Math. Soc.} \textbf{116} (1992), 53-59.
\end{rf}

\begin{rf}
Fleming, W. H., Functions whose partial derivatives are measures, 
\textit{Illinois J. Math.} \textbf{4} (1960), 452-478.
\end{rf}

\begin{rf}
Goffman, C., Non-parametric surfaces given by linearly continuous functions, 
\textit{Acta. Math.}  \textbf{103} (1960), 269-291.
\end{rf}

\begin{rf}
Goffman, C., A characterization of linearly continuous functions whose partial derivatives are measures, 
\textit{Acta. Math.} \textbf{117} (1967), 165-190.
\end{rf}

\begin{rf}
Goffman, C. and Neugebauer, C. J., On approximate derivatives, 
\textit{Proc. Amer. Math. Soc.}  \textbf{11} (1960), 962-966.
\end{rf}

\begin{rf}
Goffman, C. and Liu, F-C, Derivative measures, 
\textit{Proc. Amer. Math. Soc.}  \textbf{78} (1980), 218-220.
\end{rf}

\begin{rf}
Martio, O. and Ziemer, W. P.,  Lusin's Condition (N) and Mappings with Nonnegative Jacobians, 
\textit{Mich. Math. J.}  \textbf{39} (1992), 495-508.
\end{rf}

\begin{rf}
Ponomarev, S. P., An example of an $\text{ACTL}^p$ homeomorphism not absolutely continuous in the sense of Banach, 
\textit{Dokl. Akad. Nauk SSSR}  \textbf{201} (1971) 1053-1054 (Russian);  
translation in \textit{Soviet Math. Dokl.} \textbf{12} (1971) 1788-1790.
\end{rf}

\begin{rf}
Reshetnyak, Yu. G., Property N for the space mappings of class $W_{n, \loc}^1$,  
\textit{Sibirsk. Mat. Zh.}  \textbf{28} (1987),  149-153 (Russian); 
translation in \textit{Siberian Math. J.}, \textbf{28} (1987), 810-818.
\end{rf}

\begin{rf}
Hughs, R. E., Functions of BVC type, 
\textit{Proc. Amer. Math. Soc.}  \textbf{12} (1961),  698-701.
\end{rf}

\begin{rf}
Serrin, J., On the differentiability of functions of several variables, 
\textit{Arch. Rational Math. Anal.} \textbf{7} (1961), 359-372.
\end{rf}

\end{document}